\pdfoutput=1
\documentclass[11pt,twoside,reqno]{amsart}
\usepackage{amssymb}
\usepackage{mathtools}
\usepackage{bm}
\usepackage{dsfont}
\usepackage{enumitem}

\usepackage[top=1.5in, bottom=1.5in, left=1.23in, right=1.23in]{geometry}

\newcommand{\C}{\mathbb{C}}
\newcommand{\D}{\mathbb{D}}
\newcommand{\E}{\mathbb{E}}
\newcommand{\F}{\mathcal{F}}
\newcommand{\G}{\mathbb{G}}
\newcommand{\I}{\mathbb{I}}
\newcommand{\K}{\mathbb{K}}
\newcommand{\N}{\mathbb{N}}
\newcommand{\Q}{\mathbb{Q}}
\newcommand{\R}{\mathbb{R}}
\newcommand{\T}{\mathbb{T}}
\newcommand{\Z}{\mathbb{Z}}

\newcommand{\nn}{\mathrm{n}}
\newcommand{\ind}[1]{\mathds{1}_{{#1}}}

\numberwithin{equation}{section}

\newtheorem{condition}[equation]{Condition}
\newtheorem{conjecture}[equation]{Conjecture}
\newtheorem{corollary}[equation]{Corollary}
\newtheorem{lemma}[equation]{Lemma}
\newtheorem{proposition}[equation]{Proposition}
\newtheorem{theorem}[equation]{Theorem}

\theoremstyle{definition}

\newtheorem{definition}[equation]{Definition}
\newtheorem{example}[equation]{Example}
\newtheorem{problem}[equation]{Problem}
\newtheorem{question}[equation]{Question}
\newtheorem{remark}[equation]{Remark}

\DeclareMathOperator{\Log}{Log}

\usepackage{hyperref} 
\hypersetup{colorlinks=true, linktocpage=true,
	linkcolor=blue, citecolor=magenta, urlcolor=blue}

\title[The multilinear circle method and a question of Bergelson]{The multilinear circle method \\ and a question of Bergelson}

\author[D.~Kosz]{Dariusz Kosz}
\address{Faculty of Pure and Applied Mathematics, Wroc{\l}aw University of Science and Tech\-nology, Wybrze{\.z}e Stanis{\l}awa Wyspia{\'n}skiego 27, 50-370 Wroc{\l}aw, Poland}
\email{dariusz.kosz@pwr.edu.pl}

\author[M.~Mirek]{Mariusz Mirek}
\address{Department of Mathematics, Rutgers University, Piscataway, NJ 08854-8019, USA \& Instytut Matematyczny, Uniwersytet Wroc{\l}awski, Plac Grunwaldzki 2/4, 50-384 Wro\-c{\l}aw, Poland}
\email{mariusz.mirek@rutgers.edu}

\author[S.~Peluse]{Sarah Peluse}
\address{Department of Mathematics, Stanford University, 450 Jane Stanford Way, Building 380, Stanford, CA 94305, USA}
\email{speluse@stanford.edu}

\author[R.~Wan]{Renhui Wan}
\address{Ministry of Education Key Laboratory of NSLSCS, School of Mathematical Sciences, Nanjing Normal University, Nanjing 210023, People’s Republic of China}
\email{wrh@njnu.edu.cn}

\author[J.~Wright]{James Wright}
\address{School of Mathematics and Maxwell Institute for Mathematical Science, James Clerk Maxwell Building, The King’s Buildings, Peter Guthrie Tait Road, Edinburgh, EH9 3FD, UK}
\email{J.R.Wright@ed.ac.uk}

\thanks{Dariusz Kosz was partially supported by Basque Government grant BERC 2022-2025, by Spanish State Research Agency grant CEX2021-001142-S, and by National Science Centre of Poland grant SONATA BIS 2022/46/E/ST1/00036. Mariusz Mirek was partially supported by NSF grant DMS-2154712 and by NSF CAREER grant DMS-2236493. Sarah Peluse was partially supported by NSF grant DMS-2516641 and by a Sloan Research Fellowship. James Wright was partially supported by a Leverhulme Research Fellowship RF-2023-709$\backslash$9}

\dedicatory{Dedicated to the memory of John L.~Wright.}

\begin{document}

	\begin{abstract}
		Let $k\in \mathbb Z_+$ and $(X, \mathcal B(X), \mu)$ be a probability space equipped with a family of commuting invertible measure-preserving transformations $T_1,\ldots, T_k \colon X\to X$. Let $P_1,\ldots, P_k\in\mathbb Z[\rm n]$ be polynomials with integer coefficients and distinct degrees. We establish pointwise almost everywhere convergence of the multilinear polynomial ergodic averages
		\[
		\frac{1}{N}\sum_{n=1}^Nf_1\big(T_1^{P_1(n)}x\big)\cdots f_k\big(T_k^{P_k(n)}x\big), \qquad x\in X,
		\]
		as $N\to\infty$ for any functions $f_1, \ldots, f_k\in L^{\infty}(X)$. Besides a couple of results in the bilinear setting $k=2$, and then only in the single transformation case $T_1 = T_2$, this is the first pointwise result for general polynomial multilinear ergodic averages in arbitrary measure-preserving systems. This answers a question of Bergelson from 1996 in the affirmative for any polynomials with distinct degrees, and makes progress on the Furstenberg--Bergelson--Leibman conjecture.
		
		In this paper, we build a versatile \emph{multilinear circle method} by developing the Ionescu--Wainger multiplier theorem for the set of canonical fractions, which gives a positive answer to a question of Ionescu and Wainger from 2005. We also establish multilinear $L^p$-improving bounds and an inverse theorem in higher order Fourier analysis for averages over polynomial corner configurations, which we use to establish a multilinear analogue of Weyl's inequality and its real counterpart, a Sobolev smoothing estimate. 
	\end{abstract}
	
	\date{\today}
	
	\maketitle
	
	\tableofcontents
	
	\section{Introduction}
	\label{sec:1}
	\subsection{A brief history} A fundamental problem in ergodic theory is to understand the convergence, both in norm and pointwise almost everywhere, of multilinear polynomial ergodic averages. This line of inquiry started in the early 1930s with von Neumann's mean ergodic theorem \cite{vN} and Birkhoff's pointwise ergodic theorem \cite{BI}. Significant advances, which we will summarize shortly, have been made in this area of research over the last century.
	
	In 1996, Bergelson \cite[Question~9]{Ber1} formulated the following question.
	\begin{question}[Bergelson, 1996]
		\label{con:B}
		Let $k\in \Z_+$ and $(X, \mathcal B(X), \mu)$ be a probability space endowed with a family of commuting invertible measure-preserving transformations $T_1,\ldots, T_k \colon \\ X\to X$. Let $P_1,\ldots, P_k\in\Z[\rm n]$ be any polynomials with integer coefficients. Is it true that for any functions $f_1, \ldots, f_k\in L^{\infty}(X)$ the multilinear polynomial ergodic averages
		\begin{align} \label{eq:283}
		A_{N; X, T_1,\ldots, T_k}^{P_1,\ldots, P_k}(f_1,\ldots, f_k)(x)
		\coloneqq \E_{n\in [N]}\prod_{i \in [k]} f_i \big (T_i^{P_i(n)}x \big), \qquad x\in X,
		\end{align}
		converge pointwise almost everywhere on $X$ with respect to $\mu$ as $N\to \infty$?
	\end{question}
	
	Here and throughout the paper we use the notation $[N] \coloneqq (0, N]\cap\Z$ for any real number $N\ge1$ and
	$\E_{y\in Y}f(y) \coloneqq \frac{1}{\#Y}\sum_{y\in Y}f(y)$ for any finite set $Y\neq\emptyset$ and any function $f \colon Y\to\C$.
	
	One of the main results of this paper is the following theorem.
	\begin{theorem}
		\label{thm:B}
		The answer to Question~\ref{con:B} is yes for any $k\in\Z_+$ and any polynomials $P_1,\ldots, P_k\in\Z[\rm n]$ with distinct degrees.
	\end{theorem}
	
	In order to understand the origins of Bergelson's question, one has to go back to 1977, when Furstenberg \cite{Fur0} gave an ergodic theoretic proof of Szemer{\'e}di's theorem \cite{Sem1}, which asserts that every subset of the integers with positive upper density must contain arbitrarily long arithmetic progressions. In Furstenberg's approach \cite{Fur0}, multilinear averages of the form \eqref{eq:283} with $P_1(n)=n, \ldots, P_k(n)=kn$ and $T_1=\dots=T_k$ served as a natural tool to detect recurrent points and, consequently, arithmetic progressions in subsets of integers with positive upper density. This gave rise to a new field now known as ergodic Ramsey theory. 
	
	Not long afterwards, Bergelson had the great insight to initiate a challenging program with the goal of establishing a polynomial extension of Szemer{\'e}di's theorem by studying the asymptotic behavior of the corresponding multilinear polynomial ergodic averages. This led to his foundational paper \cite{Ber0} on weakly mixing PET (polynomial ergodic theorem), in which the van der Corput differencing technique was cemented as one of the primary tools in the field and was used to establish $L^2(X)$ norm convergence of multilinear averages \eqref{eq:283} to the product of integrals $\int_X f_1\cdots\int_X f_k$ provided that $T_1=\dots=T_k=T$ is a weakly mixing measure-preserving transformation on $X$ and $P_1,\dots,P_k$ are pairwise essentially distinct, i.e., $P_i-P_j$ is nonconstant whenever $i\neq j$.
	
	Almost two decades after Furstenberg's paper \cite{Fur0}, Bergelson and Leibman \cite{BL1} achieved the main goal of Bergelson's program by establishing (among other results --- see the discussion below Conjecture~\ref{con:FBL} and inequality \eqref{eq:2}) the following far-reaching polynomial extension of the classical Poincar{\'e} recurrence theorem and the multidimensional Szemer{\'e}di theorem of Furstenberg and Katznelson \cite{FurKa}.
	
	\begin{theorem}[The Bergelson--Leibman polynomial Szemer{\'e}di theorem, 1996]
		\label{thm:BLcorrespondence}
		Let $k\in \Z_+$ and $(X, \mathcal B(X), \mu)$ be a probability space
		endowed with a family of commuting invertible measure-preserving
		transformations $T_1,\ldots, T_k \colon \\ X\to X$. Let
		$P_1,\ldots, P_k\in\Z[\rm n]$ be polynomials with vanishing constant
		terms. Then, for any $A\in \mathcal B(X)$ satisfying $\mu(A)>0$, one has
		\begin{align}
		\label{eq:284}
		\liminf_{N\to \infty}\E_{n\in[N]}\mu\left(T_1^{-P_1(n)}(A)\cap\dots\cap T_k^{-P_k(n)}(A)\right)>0.
		\end{align}
	\end{theorem}
	
	Theorem~\ref{thm:BLcorrespondence} is merely a special case of a more general multidimensional polynomial Szemer{\'e}di theorem by Bergelson and Leibman~\cite{BL1}, which is formulated in \eqref{eq:2} below. This theorem sparked interest in understanding the asymptotic behavior of multilinear polynomial ergodic averages \eqref{eq:283} as $N\to \infty$; see Conjecture~\ref{con:FBL} and inequality \eqref{eq:2}. Therefore, Theorem~\ref{thm:B} naturally contributes to Bergelson's program of understanding the asymptotic behavior of multilinear polynomial ergodic averages.
	
	\subsection{Norm convergence} Bergelson's Question~\ref{con:B} was initially about the convergence of \eqref{eq:283} in both the $L^2(X)$ norm and pointwise almost everywhere. Nowadays, the $L^2(X)$ norm convergence of \eqref{eq:283} is fairly well understood due to groundbreaking work of Walsh \cite{W}.  Prior to Walsh's work, there was an extensive body of research towards establishing $L^2(X)$ norm convergence for \eqref{eq:283} in the single transformation case $T_1=\dots=T_k=T$. This includes fundamental work for linear polynomials due to Host and Kra \cite{HK} and, independently, Ziegler \cite{Z1}. We also refer to Leibman's appendix in \cite{Ber2} where the constructions of characteristic factors from the papers by Host and Kra \cite{HK} and Ziegler \cite{Z1} are compared and shown to be equivalent. The case of general polynomials was addressed in work by Leibman \cite{Leibman}, as well as by Frantzikinakis and Kra \cite{FraKra}, and Host and Kra \cite{HK2}. In the single transformation case, in contrast to Walsh's paper \cite{W}, one can identify the limiting function in \eqref{eq:283} thanks to the theory of Host--Kra factors~\cite{HK} and equidistribution on nilmanifolds~\cite{Leibman}. For instance, the limit in the single transformation case for linearly independent polynomials was established in \cite{FraKra} (see also \cite{FraKra1, HK2}), whereas the cases for three polynomials or multiples of one polynomial were addressed in \cite{Fra1}. We also refer to \cite{Fra}, where precise references and additional information on this topic can be found.
	
	The case of arbitrary commuting measure-preserving transformations is much more complicated and for linear polynomials was subsequently studied by Tao~\cite{Tao}, Austin \cite{A2}, and Host \cite{H}. In \cite{CFH}, Chu, Frantzikinakis, and Host established $L^2(X)$ norm convergence for averages of the form~\eqref{eq:283} when the polynomials have distinct degrees. Finally, Walsh \cite{W} established norm convergence of \eqref{eq:283} in the general case, even handling noncommutative transformations $T_1, \ldots, T_k$ generating a nilpotent group. For more on this topic, we refer the interested reader to the discussion below Conjecture~\ref{con:FBL}, as well as the articles \cite{A1, PZK} for alternative proofs and generalizations of Walsh's result and the survey articles \cite{Ber1, Ber2, Fra} that include comprehensive historical background and an extensive literature on the subject of norm convergence and its combinatorial applications.
	
	Although Walsh's result establishes $L^2(X)$ norm convergence for \eqref{eq:283}, the question of identifying the limit for arbitrary transformations and polynomials remains widely unanswered, save the case of a single transformation and linear polynomials~\cite{HK, Z1} or a single transformation and linearly independent polynomials~\cite{FraKra, FraKra1,  HK2, Fra1}.
	
	Identifying the limit for general polynomial ergodic averages is a well-known open problem in ergodic theory, even in the linear case. The $L^2(X)$ limit for \eqref{eq:283} in the case of commuting transformations and linearly independent polynomials has very recently been identified in a breakthrough paper by Frantzikinakis and Kuca \cite{FK1}. They proved the following.
	
	\begin{theorem}[Frantzikinakis--Kuca, 2022]
		\label{thm:FK1}
		Let $k\in \Z_+$ and $(X, \mathcal B(X), \mu)$ be a probability space endowed with a family of commuting invertible measure-preserving transformations $T_1,\ldots, T_k \colon X\to X$. Let $P_1,\ldots, P_k\in\Z[\rm n]$ be linearly independent polynomials. Then the rational Kronecker factor is characteristic for the averages~\eqref{eq:283}. In particular, if $T_1,\ldots, T_k$ are totally ergodic, i.e., $T_i^n$ is ergodic for any $n\in\Z_+$ and $i\in[k]$, then we have
		\begin{align}
		\label{eq:285}
		\lim_{N\to\infty}\bigg\|A_{N; X, T_1,\ldots, T_k}^{P_1,\ldots,P_k} (f_1,\ldots,f_k)-\prod_{i\in[k]}\int_{X}f_i(y)d\mu(y)\bigg\|_{L^2(X)}=0
		\end{align}
		for any $f_1, \ldots, f_k\in L^{\infty}(X)$. 
	\end{theorem}
	
	Combining \eqref{eq:285} with our Theorem~\ref{thm:B} shows that the time averages \eqref{eq:283} also converge pointwise almost everywhere to the product of space averages for totally ergodic transformations. More precisely, we have the following corollary.
	
	\begin{corollary}
		\label{cor:lim}
		Let $k\in \Z_+$ and $(X, \mathcal B(X), \mu)$ be a probability space endowed with a family of commuting invertible measure-preserving transformations $T_1,\ldots, T_k \colon X\to X$ that are totally ergodic. Let $P_1,\ldots, P_k\in\Z[\rm n]$ be polynomials with distinct degrees. Then we have
		\begin{align*}
		\lim_{N\to \infty}A_{N; X, T_1,\ldots, T_k}^{P_1,\ldots, P_k}(f_1,\ldots, f_k)(x)
		=\prod_{i\in[k]}\int_{X}f_i(y)d\mu(y)
		\end{align*}
		for any $f_1, \ldots, f_k\in L^{\infty}(X)$ and for almost every $x\in X$.
	\end{corollary}
	
	As we have seen in the brief overview above, there has been more than four decades of tremendous effort from both the ergodic and combinatorial perspectives to understand the asymptotic nature of multilinear polynomial ergodic averages \eqref{eq:283} in the $L^2(X)$ norm.
	
	\subsection{Pointwise convergence} The state of knowledge is dramatically worse for pointwise almost everywhere convergence of multilinear polynomial ergodic averages. The phenomenon of pointwise convergence, while the most natural mode of convergence, is very subtle and can differ significantly from norm convergence, making the study of pointwise convergence problems quite challenging. To begin to see the differences, let us consider the so-called uniform averages, which are given by \eqref{eq:283} with shifted averages $\E_{n\in[N]\setminus[M]}\coloneqq\frac{1}{N-M}\sum_{n=M+1}^N$ in place of the ordinary averages $\E_{n\in[N]}\coloneqq\frac{1}{N}\sum_{n=1}^N$. These two averages coincide for $M=0$. The asymptotic behavior of uniform averages in the $L^2(X)$ norm is pivotal in combinatorial applications; we refer to~\cite{Ber1, Ber2, BM} for comprehensive motivations and an extensive literature.
	
	By the dominated convergence theorem, the pointwise almost everywhere convergence of the uniform averages $\E_{n\in[N]\setminus[M]}\prod_{i \in [k]} f_i \big (T_i^{P_i(n)}x \big)$ as $N - M \to \infty$ implies their norm convergence for all $f_1, \ldots, f_k\in L^{\infty}(X)$ on a probability space $(X, \mathcal B(X), \mu)$. Interestingly, pointwise almost everywhere convergence for uniform averages can fail, even for linear averages, as shown by del Junco and Rosenblatt~\cite{dJR}. This stands in sharp contrast to $L^2(X)$ norm convergence of uniform averages, which, for instance, follows from Zorin-Kranich's paper \cite{PZK} in a very general setting.
	
	However, returning to the averages in \eqref{eq:283}, as indicated in Bergelson's question (see also Conjecture~\ref{con:FBL} below), it is conjectured that pointwise convergence holds, highlighting a significant distinction between ordinary averages $\E_{n\in[N]}$ and uniform averages $\E_{n\in[N]\setminus[M]}$. Pointwise almost everywhere convergence has only been established for a few special cases of averages $A_{N; X, T_1,\ldots, T_k}^{P_1,\ldots, P_k}(f_1,\ldots, f_k)=\E_{n\in[N]}\prod_{i \in [k]} f_i \big (T_i^{P_i(n)}x \big)$.
	
	\begin{enumerate}[label*={\arabic*}.]
		\item The case $k=1$ with $P_{1}(n)=n$ corresponds to the classical Birkhoff ergodic theorem \cite{BI}.
		
		\item The case $k=1$ with an arbitrary polynomial $P_{1}\in \Z[\mathrm n]$ was an open problem of Bellow \cite{Bel} and Furstenberg \cite{Fur3}, and was solved by Bourgain in a series of fundamental papers \cite{B1, B2, B3} in the late 1980s.
		
		\item Soon afterwards, Bourgain \cite{B4} also proved pointwise almost everywhere convergence in the bilinear setting when $k=2$ with $T_1=T_2$, $P_{1}(n)=an$, and
		$P_{2}(n)=bn$ for any $a, b\in\Z$, providing a positive answer to a question of Furstenberg \cite{Fur2}.
		
		\item Recently, the second author along with Krause and Tao \cite{KMT} established pointwise almost everywhere convergence when $k=2$ with $T_1=T_2$, $P_1(n)=n$, and an arbitrary polynomial $P_2\in \Z[{\rm n}]$ of degree at least two.  Throughout the paper, we will call
		\[
		\E_{n\in[N]}f_1(T^nx)f_2(T^{P_2(n)}x)
		\]
		the \emph{Furstenberg--Weiss averages}, even though Furstenberg and Weiss studied
		$L^2(X)$ norm convergence of such averages in \cite{FurWei} only for quadratic polynomials $P_2$.
		
		\item Pointwise convergence of multilinear polynomial averages was established for some special classes of measure-preserving systems, such as exact endomorphisms and $K$-automorphisms \cite{DL} and nilsystems \cite{Leibman2}. The context of commuting transformations along linear orbits for distal systems was studied in \cite{HSY} and \cite{DS}. 
	\end{enumerate}
	
	Aside from these results, nothing further is known about pointwise almost everywhere convergence for multilinear polynomial averages \eqref{eq:283} in general measure-preserving systems.
	
	Our paper forges new ground on three fronts: we treat averages with an arbitrary degree of multilinearity $k$, with any commuting transformations, and with all polynomials allowed to be nonlinear. All prior work in general measure-preserving systems treated only the linear or bilinear case $k\leq 2$, a single transformation, and at most one nonlinear polynomial. We obtain, for example, the first pointwise almost everywhere convergence result for averages with commuting transformations $S, T \colon X\to X$ corresponding to the ``sqorners'' configuration
	\begin{align}
	\label{eq:305}
	\E_{n\in [N]} f_1 (S^{n}x)f_2 (T^{n^2}x), \qquad x\in X,
	\end{align}
	recently considered by the third author along with Prendiville and Shao~\cite{PPS}. The methods used in our paper are robust enough to suggest that we can handle more general classes of ergodic averages than those considered here. We will elaborate on this later in the introduction, when we describe the \emph{multilinear circle method} that we develop in this article.
	
	\subsection{Statement of the main results} Throughout this paper, the triple $(X, \mathcal B(X), \mu)$ denotes a $\sigma$-finite measure space and $\K$ denotes either $\Z$ or $\R$. Correspondingly, $\K[{\rm n}]$ denotes the space of all formal polynomials $P({\rm n})$ with coefficients and indeterminate ${\rm n}$ in $\K$. Each polynomial $P\in\K[{\rm n}]$ is always identified with a map $\K\ni n\mapsto P(n)\in\K$.
	
	Let $k \in\Z_+$. Given a family ${\mathcal T} = \{T_1,\ldots, T_k\}$ of invertible commuting measure-preserving transformations on $X$, complex-valued measurable functions $f_1,\ldots, f_k$ on $X$, a family of polynomials $\mathcal P = \{P_1,\ldots, P_k\} \subset \Z[\mathrm n]$, and a real number $N \geq 1$, we define, as in \eqref{eq:283}, the corresponding multilinear polynomial ergodic average by
	\begin{align}
	\label{eq:98}
	A_{N; X, {\mathcal T}}^{\mathcal P}(f_1,\ldots, f_k)(x) \coloneqq \E_{n\in [N]}\prod_{i \in [k]} f_i\big(T_i^{P_i(n)}x\big),
	\qquad x\in X,
	\end{align}
	and its truncated version by
	\begin{align}
	\label{eq:102}
	\tilde{A}_{N; X, {\mathcal T}}^{\mathcal P}(f_1,\ldots, f_k)(x) \coloneqq \E_{n\in [N]\setminus[N/2]}\prod_{i \in [k]} f_i\big(T_i^{P_i(n)}x\big), \qquad x\in X.
	\end{align}
	
	We will often abbreviate $A_{N; X, {\mathcal T}}^{{\mathcal P}}$ to $A_{N; X}^{{\mathcal P}}$ and $\tilde{A}_{N; X, {\mathcal T}}^{\mathcal P}$ to $\tilde{A}_{N; X}^{\mathcal P}$ when the transformations are understood. Depending on how explicit we want to be, in some instances, we will write out the averages
	\begin{align*}
	A_{N;X}^{\mathcal P} = A_{N;X}^{P_1,\ldots, P_k}
	\quad {\rm or} \quad 
	A_{N; X, {\mathcal T}}^{\mathcal P} =
	A_{N;X,T_1,\ldots, T_k}^{P_1,\ldots, P_k},\\
	\tilde A_{N;X}^{\mathcal P} =\tilde  A_{N;X}^{P_1,\ldots, P_k}
	\quad {\rm or} \quad 
	\tilde A_{N; X, {\mathcal T}}^{\mathcal P} =
	\tilde  A_{N;X,T_1,\ldots, T_k}^{P_1,\ldots, P_k}.
	\end{align*}
	
	Using the asymptotic notation from Section~\ref{section:2} and the definition of $r$-variational norm ${\bf V}^r$ from \eqref{eq:3}, we now state the main result of this article, i.e., the following quantitative ergodic theorem that implies Theorem~\ref{thm:B} and, consequently, gives an affirmative answer to Bergelson's Question~\ref{con:B} for polynomials with distinct degrees.
	
	\begin{theorem}
		\label{thm:main}
		Let $k\in\Z_+$ and $(X, \mathcal B(X), \mu)$ be a $\sigma$-finite measure space equipped with a family $\mathcal T$ of invertible measure-preserving transformations $T_1,\ldots, T_k \colon X\to X$. Suppose that ${\mathcal P} = \{P_1,\ldots, P_k\} \subset \Z[\mathrm n]$ is a family of polynomials with distinct degrees. Let $f_i \in L^{p_i}(X)$ for $i\in[k]$ and some $1<p_1,\ldots, p_k<\infty$ such that $\frac{1}{p_1}+\cdots+\frac{1}{p_k}=\frac{1}{p}\le 1$, and let $A_{N; X, {\mathcal T}}^{\mathcal P}(f_1,\ldots, f_k)$ be the averages defined in \eqref{eq:98}. Then the following statements hold.
		
		\begin{itemize}
			\item[\rm (i)] \textit{(Mean ergodic theorem)} The functions $A_{N; X, {\mathcal T}}^{\mathcal P}(f_1,\ldots, f_k)$ converge in the $L^p(X)$ norm as $N \to \infty$.
			
			\item[\rm (ii)] \textit{(Pointwise ergodic theorem)} The functions $A_{N; X, {\mathcal T}}^{\mathcal P}(f_1,\ldots, f_k)$ converge pointwise almost everywhere as $N \to \infty$.
			
			\item[\rm (iii)] \textit{(Maximal ergodic theorem)}
			One has
			\begin{align}
			\label{eq:101sup}
			\bigg \|\sup_{N\in\Z_+}\left|A_{N; X, {\mathcal T}}^{\mathcal P}(f_1,\ldots, f_k)\right|\bigg\|_{L^p(X)}\lesssim \prod_{i\in[k]}\|f_i\|_{L^{p_i}(X)}.
			\end{align}
			
			\item[\rm (iv)] \textit{(Variational ergodic theorem)}
			Fix $r>2$ and $\lambda>1$. Then one has
			\begin{align}
			\label{eq:101}
			\left\|{\bf V}^r\left(A_{N; X, {\mathcal T}}^{\mathcal P}(f_1,\ldots, f_k): N\in\D\right) \right\|_{L^p(X)}\lesssim\prod_{i\in[k]}\|f_i\|_{L^{p_i}(X)}
			\end{align}
			if $\mathbb D = \{\lambda_n :n\in\N\} \subset [1,+\infty)$ is $\lambda$-lacunary, i.e., $\inf_{n\in\N}\frac{\lambda_{n+1}}{\lambda_n}\ge\lambda$.
		\end{itemize}
		The same results hold with the truncated averages $\tilde{A}_{N; X, {\mathcal T}}^{\mathcal P}$ in place of $A_{N; X, {\mathcal T}}^{\mathcal P}$.
	\end{theorem}
	
	We now give some remarks about Theorem~\ref{thm:main} and its consequences.
	
	\begin{enumerate}[label*={\arabic*}.]
		\item On a probability space, bounded functions are $L^p$ functions and so the conclusion from part (ii) implies Theorem~\ref{thm:B}.
		
		\item When $k=1$ and $P_1(n)=n$, parts (i)--(iii) follow from works of von Neumann \cite{vN}, Birkhoff \cite{BI}, and Hopf \cite{Hopf}. It is also known that the maximal function is of weak-type $(1,1)$, which can be derived from the corresponding bounds for the Hardy--Littlewood maximal function \cite{bigs} by using the Calder{\'o}n transference principle \cite{C1}. Part (iv) follows from \cite{J+}, and weak-type $(1,1)$ estimates also hold. 
		
		\item When $k=1$ and $P_1\in \Z[{\rm n}]$ is arbitrary, part (i) was established by Furstenberg \cite{Fur2}, parts (ii) and (iii) were established by Bourgain \cite{B1, B2, B3}, and part (iv) follows, for instance, from \cite{MSZ3}. Here, if $\deg P_1\ge 2$, then pointwise almost everywhere convergence may fail for $p_1=1$, as was shown in \cite{BM, LaV1}.
		
		\item When $k=2$, $P_1(n)=an$ and $P_2(n)=bn$ with $a, b\in\Z$, and $T_1=T_2$, part (i) was established by Furstenberg \cite{Fur2} and part (ii) for bounded functions on a probability space was established by Bourgain \cite{B4}. Part (iii) follows from work of Lacey~\cite{LAC}; see also \cite{Demeter}. Finally, part (iv) was established by Do, Oberlin, and Palsson \cite{DOP} for large $r > 2$.
		
		\item When $k=2$, $P_1(n)=n$ and $P_2\in\Z[\rm n]$ is arbitrary with $\deg P_2\ge2$, and $T_1=T_2$, part (i) was originally established by Furstenberg and Weiss \cite{FurWei} in the case $P_2(n)=n^2$, whereas for general polynomials $P_2\in\Z[\rm n]$ with $\deg P_2\ge2$, it follows from the papers of Host and Kra~\cite{HK2} and Frantzikinakis and Kra~\cite{FraKra}. Parts (ii)--(iv) were recently established by the second author with Krause and Tao in \cite{KMT}.
		
		\item The implicit constant in \eqref{eq:101} is allowed to depend on the parameters $p_1,\ldots, p_k,r,\lambda$ and on the polynomials $P_1,\ldots, P_k$.
		
		\item If \eqref{eq:101} is established for some $r > 2$ and, for example, the $2$-lacunary set $\{2^n : n \in \N\}$, then \eqref{eq:101sup} holds for the same $p_1, \dots, p_k$ with an implicit constant that depends only on $p_1,\ldots, p_k,P_1,\ldots, P_k$. The maximal inequality \eqref{eq:101sup} also addresses a question raised by Christ, Durcik, and Roos in \cite[Section~6, Problem~4]{CDR}.
		
		\item If \eqref{eq:101} holds with some $r > 2$ for all $\lambda > 1$ and all $\lambda$-lacunary sets $\mathbb D \subset [1,+\infty)$, then the limit of the averages $A_{N; X, {\mathcal T}}^{\mathcal P}(f_1,\ldots, f_k)(x)$ exists for $\mu$-almost every $x\in X$ as $\mathbb D\ni N\to\infty$. This, in turn, implies pointwise almost everywhere convergence by taking the $2^{1/s}$-lacunary set $\{ 2^{n/s} : n \in \N \}$ for every $s \in \Z_+$; see \cite[Proposition~3.2(i)]{KMT}.
		
		\item Combining \eqref{eq:101sup} with the pointwise almost everywhere convergence of the averages $A_{N; X, {\mathcal T}}^{\mathcal P}(f_1,\ldots, f_k)$ and the dominated convergence theorem yields the norm convergence, which was originally proven by Chu, Frantzikinakis, and Host \cite{CFH}, using different methods.
		
		\item Taking into account the last three items above, it suffices to establish inequality \eqref{eq:101}. We can assume that $N \geq C_0$ for some fixed $C_0 \in \Z_+$.
		
		\item The condition $r>2$ in Theorem~\ref{thm:main} is necessary, as no variational estimate is possible for $r \leq 2$; see \cite[Corollary~12.4]{KMT}. Thus, the range in \eqref{eq:101} is sharp.
		
		\item We are also able to ``break duality'' in Theorem~\ref{thm:main} by handling some ranges of exponents $1<p_1,\ldots, p_k<\infty$ such that $\frac{1}{p_1}+\cdots+\frac{1}{p_k}=\frac{1}{p}> 1$, at the cost of increasing $r$ slightly in \eqref{eq:101}. 
		
		\item Estimates \eqref{eq:101sup} and \eqref{eq:101} imply the analogous estimates with $\tilde{A}_{N; X, {\mathcal T}}^{\mathcal P}$ in place of $A_{N; X, {\mathcal T}}^{\mathcal P}$ and vice versa. This follows by a simple telescoping argument; see \cite[Proposition~3.2(iii)]{KMT}. We will work with $\tilde{A}_{N; X, {\mathcal T}}^{\mathcal P}$ to avoid some technicalities.
		
		\item In the proof of Theorem~\ref{thm:main}, we will use the Calder{\'o}n transference principle \cite{C1}, which reduces the estimates \eqref{eq:101sup} and \eqref{eq:101} in abstract measure-preserving systems $(X, \mathcal{B}(X), \mu)$ to those for the integer shift system, which is $\sigma$-finite; see Example~\ref{ex:1}. This is the reason why we formulate Theorem~\ref{thm:main} for $\sigma$-finite measure spaces. However, from the point of view of applications in combinatorics or elsewhere, where the statistical properties of the ergodic averages \eqref{eq:98} for totally ergodic systems matter, the general setting of $\sigma$-finite measure spaces is not interesting (since the limits are zero when $\mu(X) = \infty$) and only finite measure spaces are important.
		
		\item A continuous analogue of the ``sqorners'' averages \eqref{eq:305} was studied by Christ, Durcik, Kova{\v c}, and Roos in \cite{CDKR} and pointwise convergence for these averages was established. The key estimate is a Sobolev smoothing bound from an earlier paper of Christ, Durcik, and Roos \cite[Theorem~5]{CDR}, which is a special instance of our multilinear Weyl inequality Theorem~\ref{weyl} in the real setting $\K = \R$. Theorem~\ref{weyl} establishes this key smoothing inequality for general $k$-linear polynomial corner averages in the distinct degree case. We also refer to a recent paper by Hsu and Yu-Hsiang Lin~\cite[Proposition~1.2]{HLin}, which provides a short and elementary proof of the smoothing inequality from \cite[Theorem~5]{CDR}. In \cite{CDR}, Christ, Durcik, and Roos use their smoothing inequality to prove maximal estimates and quantitative nonlinear Roth-type theorems for sqorner configurations in ${\mathbb R}^2$. The same applications for distinct degree polynomial corner configurations in ${\mathbb R}^k$ can be derived from our general smoothing inequality Theorem~\ref{weyl}.
	\end{enumerate}
	
	\subsection{The Furstenberg--Bergelson--Leibman conjecture} Both Theorem~\ref{thm:B} and Theorem~\ref{thm:main} yield progress on the Furstenberg--Bergelson--Leibman conjecture, which asserts the following.
	
	\begin{conjecture}[The Furstenberg--Bergelson--Leibman conjecture]
		\label{con:FBL}
		Let $d, k\in\Z_+$ and $(X, \mathcal B(X), \mu)$ be a probability space equipped with a family of invertible measure-preserving transformations $T_1,\ldots, T_d \colon X\to X$ generating a nilpotent group. Consider any polynomials $P_{1, 1},\ldots,P_{i,j},\dots,P_{d,k}\in \Z[\mathrm n]$ with integer coefficients. Then, for any functions $f_1, \ldots, f_k\in L^{\infty}(X)$, the multilinear polynomial ergodic averages
		\begin{align}
		\label{eq:272}
		\E_{n\in[N]}\prod_{j \in [k]}f_j\Big(T_1^{P_{1, j}(n)}\cdots T_d^{P_{d, j}(n)} x\Big), \qquad x\in X,
		\end{align}
		converge pointwise almost everywhere as $N\to \infty$.	
	\end{conjecture}
	
	We give a few remarks about Conjecture~\ref{con:FBL}.
	
	\begin{enumerate}[label*={\arabic*}.]
		\item This conjecture is a very challenging problem in pointwise ergodic theory and modern harmonic analysis. It had been promoted by Bergelson in \cite[Question~9, p.~52]{Ber1} and in \cite[Section~6, p.~838]{Ber2}, and in person by Furstenberg (see Austin's article \cite[p.~6662]{A1}) before it was published by Bergelson and Leibman \cite[Section~5.5, p.~468]{BL2}.
		
		\item Bergelson's Question~\ref{con:B} is a special case of Conjecture~\ref{con:FBL}. Indeed, if $d=k$, $P_{j, j}=P_j\in\Z[\mathrm n]$, $P_{i, j}\equiv0$ whenever $i\neq j$, and the transformations $T_1,\ldots, T_k$ commute, then the averages from \eqref{eq:272} coincide with the averages in \eqref{eq:283}. Hence, Theorem~\ref{thm:B} constitutes progress on Conjecture~\ref{con:FBL} in the commutative case for arbitrary polynomials with distinct degrees.
		
		\item The polynomial Szemer{\'e}di theorem of Bergelson and Leibman \cite{BL1} was proved with the averages \eqref{eq:272} in place of \eqref{eq:284} for invertible commuting measure-preserving transformations $T_1,\ldots, T_d$. To be more precise, if $P_{1, 1},\ldots,P_{i,j},\dots,P_{d,k}\in \Z[\mathrm n]$ are polynomials with vanishing constant terms, then for any $A\in\mathcal B(X)$ satisfying $\mu(A)>0$, they show that
		\begin{align}
		\label{eq:2}
		\liminf_{N\to \infty}\E_{n\in[N]}\mu\bigg(\bigcap_{j\in[k]}T_1^{-P_{1, j}(n)}\cdots T_d^{-P_{d, j}(n)}(A)\bigg)>0.
		\end{align}
		The most general polynomial Szemer{\'e}di theorem for commuting transformations and {\rm IP}-sets to date is due to Bergelson and McCutcheon~\cite{BMcC}. We also refer to \cite[Theorem~6.13]{BM} for a variant of \eqref{eq:2} with the uniform averages $\E_{n\in[N]\setminus[M]}$ in place of $\E_{n\in[N]}$. In \cite{BL2}, Bergelson and Leibman began the study of recurrence theorems in the nilpotent setting, which triggered the formulation of Conjecture~\ref{con:FBL}.
		
		\item The Furstenberg--Bergelson--Leibman conjecture initially asked about convergence for \eqref{eq:272} both in $L^2(X)$ and pointwise almost everywhere. As mentioned previously, Conjecture~\ref{con:FBL} in the context of $L^2(X)$ norm convergence for \eqref{eq:272} was established by Walsh~\cite{W}. Zorin-Kranich~\cite{PZK} also established $L^2(X)$ norm convergence for \eqref{eq:272} with the uniform averages $\E_{n\in[N]\setminus[M]}$ in place of $\E_{n\in[N]}$. 
		
		\item Bergelson and Leibman \cite{BL2} showed that convergence for \eqref{eq:272} may fail if $T_1,\ldots, T_d$ generate a solvable group that is not nilpotent. This suggests that the nilpotent setting is the appropriate one for this conjecture.
		
		\item Conjecture~\ref{con:FBL} in a genuinely nilpotent (i.e., of step at least two) setting is widely open, except one case. Recently, the second author along with Ionescu, Magyar, and Szarek \cite{IMMS} proved this conjecture with $d\in\Z_+$ and $k=1$ for arbitrary polynomials $P_{1,1},\ldots, P_{d,1}\in\Z[{\rm n}]$ and arbitrary invertible measure-preserving transformations $T_1,\ldots, T_{d}$ on a $\sigma$-finite measure space $(X, \mathcal B(X), \mu)$ that generate a nilpotent group of step two; see also \cite{IMSW, MSW0} for some special cases of \cite{IMMS}.
	\end{enumerate}
	
	\subsection{Reduction to the integer shift system}
	In pointwise convergence problems, the most important dynamical system is the integer shift system.
	
	\begin{example}\label{ex:1}
		Consider the $k$-dimensional lattice $(\Z^k, \mathcal B(\Z^k), \mu_{\Z^k})$, where $\mathcal B(\Z^k)$ is the $\sigma$-algebra of all subsets of $\Z^k$ and $\mu_{\Z^k}$ is counting measure on $\Z^k$, equipped with the standard family of shifts $S_1,\ldots, S_k \colon \Z^k\to\Z^k$. More precisely, let $S_i(x) \coloneqq x-e_i$ for all $x\in\Z^k$, where $e_i$ is the $i$-th standard basis vector for each $i\in[k]$. The average $A_{N; X, T_1,\ldots, T_k}^{P_1, \ldots, P_k}$ from \eqref{eq:98} with the family $\{S_1, \dots, S_k\}$ in place of $\{T_1, \dots, T_k\}$ can be rewritten as
		\begin{align}
		\label{eq:99}
		A_{N; \Z^k}^{P_1, \ldots, P_k}(f_1,\ldots, f_k)(x)=\E_{n\in [N]} \prod_{i \in [k]}  f_i(x-P_i(n)e_i),
		\qquad  x\in\Z^k.
		\end{align}
		Its truncated variant \eqref{eq:102} takes the form
		\begin{align}
		\label{eq:103}
		\tilde{A}_{N; \Z^k}^{P_1, \ldots, P_k}(f_1,\ldots, f_k)(x)=
		\E_{n\in [N]\setminus[N/2]} \prod_{i \in [k]}  f_i(x-P_i(n)e_i),
		\qquad  x\in\Z^k.
		\end{align}
		We will often abbreviate $A_{N; \Z^k}^{P_1, \ldots, P_k}$ to
		$A_{N; \Z^k}^{{\mathcal P}}$ and $\tilde A_{N; \Z^k}^{P_1, \ldots, P_k}$ to
		$\tilde A_{N; \Z^k}^{{\mathcal P}}$. 
	\end{example}
	
	In view of the Calder{\'o}n transference principle \cite{C1} (or, more precisely, following the argument from \cite[Proposition~3.2(ii)]{KMT} or \cite[Theorem~1.6]{Kosz}), it will suffice to work with the integer shift system and establish \eqref{eq:101} with $A_{N; \Z^k}^{\mathcal P}$ in place of $A_{N; X, \mathcal T}^{\mathcal P}$. This will allow us to utilize the algebraic structure of $\Z^k$ and employ Fourier methods on $\Z^k$, which are not available in abstract measure preserving systems in general.
	
	The Calder{\'o}n transference principle allows us to transfer the quantitative estimates \eqref{eq:101sup} and \eqref{eq:101} from the integer shift system to corresponding estimates for $A_{N; X, \mathcal T}^{\mathcal P}$ in abstract measure-preserving systems $(X, \mathcal{B}(X), \mu)$. We emphasize that the Calder{\'o}n transference principle \cite{C1} only transfers quantitative bounds that imply pointwise almost everywhere convergence, but does not transfer pointwise almost everywhere convergence itself. In fact, in the integer shift system, pointwise convergence is implied by norm convergence, since the $\ell^\infty(\Z^k)$ norm is dominated by the $\ell^2(\Z^k)$ norm. Hence, we will only be concerned with proving quantitative bounds for $A_{N; \Z^k}^{\mathcal P}$ or $\tilde A_{N; \Z^k}^{\mathcal P}$ in the integer shift system, not pointwise convergence on $\Z^k$. For technical reasons we will only work with the truncated averages $\tilde A_{N; \Z^k}^{\mathcal P}$.
	
	After these reductions, our main result reads as follows.
	
	\begin{theorem}
		\label{thm:main1}
		Let $k\in\Z_+$ and ${\mathcal P} = \{P_1,\ldots, P_k\} \subset \Z[\mathrm n]$ be a family of polynomials with distinct degrees. Let $f_i\in \ell^{p_i}(\Z^k)$ for $i\in[k]$ and some $1<p_1,\ldots, p_k<\infty$ such that $\frac{1}{p_1}+\cdots+\frac{1}{p_k}=\frac{1}{p}\le 1$, and let $\tilde{A}_{N; \Z^k}^{\mathcal P}(f_1,\ldots, f_k)$ be the averages defined in \eqref{eq:103}. Fix $r>2$ and $\lambda>1$. Then one has
		\begin{align}
		\label{eq:104}
		\left\|{\bf V}^r\left(\tilde{A}_{N; \Z^k}^{\mathcal P}(f_1,\ldots, f_k): N\in\D\right)\right\|_{\ell^p(\Z^k)}\lesssim\prod_{i\in[k]}\|f_i\|_{\ell^{p_i}(\Z^k)},
		\end{align}
		if $\mathbb D=\{\lambda_n :n\in\N\} \subset [1,+\infty)$ is $\lambda$-lacunary. The same result holds with the averages $A_{N; \Z^k}^{\mathcal P}(f_1,\ldots, f_k)$, defined in \eqref{eq:99}, in place of $\tilde{A}_{N; \Z^k}^{\mathcal P}(f_1,\ldots, f_k)$.
	\end{theorem}
	
	Arguing as in \cite[Proposition~3.2(ii)--(iii)]{KMT}, one sees that Theorem~\ref{thm:main1} implies Theorem~\ref{thm:main}. Thus, it suffices to prove Theorem~\ref{thm:main1}.
	
	To establish the $r$-variational inequality in \eqref{eq:104}, we will develop a new robust method that we call the \emph{multilinear circle method}. This method can be viewed as a classical Hardy--Littlewood--Ramanujan circle method in $\ell^p(\Z^k)$ spaces. A bilinear variant of the circle method was recently developed by the second author with Krause and Tao \cite{KMT} in the context of pointwise convergence for the Furstenberg--Weiss averages $\E_{n\in[N]}f_1(T^nx)f_2(T^{P(n)}x)$ for any $P\in\Z[{\rm n}]$ with $\deg P\ge2$. The arguments from \cite{KMT} are limited to the Furstenberg--Weiss averages. Here we can handle genuinely multilinear averages involving polynomials with distinct degrees and arbitrary commuting transformations. Although the arguments in this paper are inspired by those in \cite{KMT}, we will present a conceptually different approach to overcome several new difficulties that arise in the more general setting. The key tools that we develop to prove inequality \eqref{eq:104} that make up our multilinear circle method are the following.
	
	\begin{enumerate}[label*={\arabic*}.]
		\item\label{T1} \emph{An Ionescu--Wainger multiplier theorem for the set of canonical fractions}; see Theorem~\ref{thm:IW}. This is a multifrequency multiplier theorem, which we prove in Section~\ref{sec:IW} for the set of \textit{canonical fractions} (see \eqref{eq:4} for their definition), giving a positive answer to a question of Ionescu and Wainger from \cite[Remark~3, p.~361]{IW}. This new tool will enable us to implement the circle method directly in $\ell^p(\Z^k)$ spaces.
		
		\item\label{T2} \emph{An inverse theorem for averages over distinct degree polynomial corner configurations}; see Theorem~\ref{pp}. The inverse theorem is a powerful result in additive combinatorics, offering important structural information that reveals the minor and major arc structure for the multilinear operators $\tilde{A}_{N; \Z^k}^{\mathcal P}(f_1,\ldots, f_k)$. We prove it in Section~\ref{sec:inverse} in both the integer and real settings, and it constitutes a joint generalization of an inverse theorem of the third author~\cite{P2} for averages over distinct degree polynomial progressions and of the third author along with Prendiville and Shao~\cite{PPS} for averages over the two-dimensional ``sqorners''
		configuration $(x_1,x_2),(x_1+n,x_2),(x_1,x_2+n^2)$.
		
		\item\label{T3} \emph{A multilinear $\ell^{p}(\Z^k)$-improving inequality}; see Theorem~\ref{thm:improv}. This is a new inequality in the polynomial corners setting that will allow us to relax the $\ell^{\infty}(\Z^k)$ bounds that arise in the inverse theorem to $\ell^{p}(\Z^k)$ bounds, and to handle genuinely multilinear cases for $k\ge 3$. We prove this inequality in Section~\ref{sec:improving} in both the integer and real settings.
		
		\item\label{T4} \emph{Multilinear Weyl and Sobolev smoothing inequalities}; see Theorem~\ref{weyl}. The multilinear Weyl inequality \eqref{eq:298} is a key tool to control the minor arc contribution in our multilinear circle method. The multilinear Sobolev smoothing inequality will be used to understand the major arc contribution. These inequalities are proved in Section~\ref{sec:weyl}, and are sometimes referred to as \textit{smoothing estimates} in the literature. A bilinear Weyl inequality was recently proved in \cite{KMT} as a consequence of the inverse theorem of the third author \cite{P2}. However, the bounds obtained in \cite{KMT} are logarithmic in scale, which is insufficient for the methods in this paper. Here, we make quantitative improvements to the bilinear estimates and derive a multilinear Weyl inequality with polynomial bounds that are consistent with the bounds in the classical Weyl inequality for exponential sums. This is possible thanks to our Ionescu--Wainger multiplier theorem for the set of canonical fractions, inverse theorem for polynomial corner configurations, and multilinear $\ell^p(\Z^k)$-improving inequality. 
	\end{enumerate}
	
	Finally, in Section~\ref{sec:ergodic}, we will use these tools to develop the multilinear circle method in the context of Theorem~\ref{thm:main1}. An important new feature of our argument is that we do not require $p$-adic methods, which were employed in~\cite{KMT}. This answers a question of Magyar \cite{Mag}, who asked whether the use of $p$-adic and adelic harmonic analysis is necessary in \cite{KMT}.
	
	\subsection{The multilinear circle method} The details of the multilinear circle method in the context of Theorem~\ref{thm:main1} will be presented in Section~\ref{sec:ergodic}. We now briefly describe its key features, starting with a few basic concepts from the classical circle method.
	
	\subsubsection{A first glimpse of the circle method} The use of the classical circle method to attack pointwise convergence problems for linear ergodic averages with polynomial orbits originates in Bourgain's papers \cite{B1, B2, B3}, and can be summarized as follows:
	
	\begin{enumerate}[label*={(\alph*)}]
		\item\label{I1} To control the minor arc contribution, we apply Plancherel's theorem and Weyl's inequality for exponential sums.
		
		\item\label{I2} To control the major arc contribution, we use multifrequency harmonic analysis in the spirit of the Ionescu--Wainger multiplier theorem; see Theorem~\ref{thm:IW} and Theorem~\ref{thm:IWvar}.
	\end{enumerate}
	
	We now illustrate the ideas from \ref{I1} and \ref{I2} in the context of $r$-variational estimates \eqref{eq:104} in the linear case $k=1$ for finitely supported $f\in \ell^2(\Z)$. This will highlight the differences between the classical circle method and the multilinear circle method.
	
	Let $P\in\Z[\rm n]$ be a polynomial of degree $d\ge2$. For the sake of discussion we will work with $A_{N; \Z}^{P}$ instead of its truncation. Using the Fourier transform, we note that
	\begin{align*}
	\mathcal F_{\Z}A_{N; \Z}^{P}(f)(\xi)=m_{N}(\xi)\mathcal F_{\Z}f(\xi), \qquad \xi\in\T,
	\end{align*}
	where the multiplier $m_N$ is the exponential sum
	\begin{align*}
	m_{N}(\xi) \coloneqq \E_{n\in [N]}e(\xi P(n)).
	\end{align*}
	The classical circle method can be used to understand the nature of the multiplier $m_{N}$. This will require the concepts of canonical fractions and their corresponding major arcs. Given $N_1, N_2\in \R_+$, we define the set of \emph{canonical fractions} by
	\begin{align}
	\label{eq:4}
	\mathcal R_{\le N_1} \coloneqq \left\{\frac{a}{q}\in \Q/\Z: q\in[N_1] \text{ and } \gcd(a, q)=1\right\},
	\end{align}
	where $\frac{a}{q}$ represents the equivalence class $\frac{a}{q} + \Z$, and the corresponding set of \emph{major arcs} by
	\begin{align*}
	\mathfrak M_{\le N_2}(\mathcal R_{\le N_1}) \coloneqq \bigcup_{\theta\in \mathcal R_{\le N_1}}[\theta-N_2, \theta+N_2].
	\end{align*}
	The set of \emph{minor arcs} is then defined as the complement of the set of major arcs in $\T$. We formulate Weyl's estimate for the multiplier $m_N(\xi)$ as follows: for every $C\in\R_+$ there exists a small constant $c\in(0, 1)$ such that we have
	\begin{align}
	\label{eq:288}
	\xi \in \T \setminus \mathfrak M_{\le N^{-d}\delta^{-C}}(\mathcal R_{\le \delta^{-C}}) \implies |m_{N}(\xi)|\le c^{-1}(\delta^c+N^{-c}) 
	\end{align}
	for all $N\ge 1$ and $\delta\in (0, 1]$. In fact, \eqref{eq:288} is the classical Weyl sum estimate for normalized exponential sums; see, for example, \cite[Exercise~1.1.21, p.~16]{Tho}.
	
	If $N\ge1$ is sufficiently large in terms of $\delta$, say $N>3\delta^{-3C}$, then the intervals that comprise the set of major arcs are narrow and disjoint in $\T$.  Taking a smooth even cutoff function $\eta \colon \R\to[0, 1]$ such that $\ind{[-1/4, 1/4]}\le \eta\le \ind{(-1/2, 1/2)}$, we define for $N_1, N_2\in\R_+$ the smooth projection operators $\Pi[\le N_1, \le N_2] \colon \ell^2(\Z)\to \ell^2(\Z)$ by setting
	\begin{align}
	\label{eq:290}
	\mathcal F_{\Z}\left(\Pi\left[\le N_1, \le N_2\right]f\right)(\xi)
	\coloneqq \sum_{\theta\in \mathcal R_{\le N_1}}\eta\left(N_2^{-1}(\xi-\theta)\right)\mathcal F_{\Z}f(\xi).
	\end{align}
	These projections will be called the Ionescu--Wainger projections and will allow us to effectively localize to the major arcs. Their boundedness properties will be extensively studied in Section~\ref{sec:IW} in a more general context. By \eqref{eq:290} and Plancherel's theorem, the operator $\Pi[\le \delta^{-C}, \le N^{-d}\delta^{-C}]$ is a contraction on $\ell^2(\Z)$ provided that the arcs in $\mathfrak M_{\le N^{-d}\delta^{-C}}(\mathcal R_{\le \delta^{-C}})$ are disjoint. We also see that $\Pi[\le \delta^{-C}, \le N^{-d}\delta^{-C}]$ is bounded on $\ell^p(\Z)$ for all $p\in[1, \infty]$, since 
	\begin{align}
	\label{eq:292}
	\left\|\Pi \left[\le \delta^{-C}, \le N^{-d}\delta^{-C}\right] f \right\|_{\ell^p(\Z)}
	\lesssim \# \mathcal R_{\le \delta^{-C}} \|f\|_{\ell^p(\Z)} \le \delta^{-2C} \|f\|_{\ell^p(\Z)}.
	\end{align}
	This bound is not very useful, but much better bounds will follow from the Ionescu--Wainger multiplier theorem proved in Section~\ref{sec:IW}; see Theorem~\ref{thm:IW}.
	
	Taking $\delta=4^{1/C} N^{-\varepsilon}$ with $\varepsilon\in(0, 1)$ sufficiently small, say $\varepsilon C<1/4$, and $N$ sufficiently large (so that, in particular, $\delta \leq 1$), we gain a negative power of $N$ in \eqref{eq:288} and this, combined with Plancherel's theorem, yields
	\begin{align}
	\label{eq:289}
	\left\|A_{N; \Z}^{P}\left(f-\Pi\left[\le N^{C \varepsilon}, \le N^{-d+C\varepsilon}\right]f\right)\right\|_{\ell^2(\Z)}
	\lesssim N^{-c\varepsilon}\|f\|_{\ell^2(\Z)},
	\end{align}
	since $\F_\Z(f-\Pi[\le N^{C \varepsilon}, \le N^{-d+C\varepsilon}]f)$ vanishes on the major arcs
	\[
	\mathfrak M_{\le N^{-d+C\varepsilon}/4}(\mathcal R_{\le N^{C\varepsilon}}) \supseteq \mathfrak M_{\le N^{-d}\delta^{-C}}(\mathcal R_{\le \delta^{-C}}).
	\]
	If we use \eqref{eq:289} and the fact that $N\in\mathbb D$, where $\mathbb D$ is lacunary, inequality \eqref{eq:104} for $k=1$ and $p=2$ is reduced to proving  the following inequality
	\begin{align}
	\label{eq:291}
	\left\|{\bf V}^r\left(A_{N; \Z}^{P}\left(\Pi\left[\le N^{C \varepsilon}, \le N^{-d+C\varepsilon}\right]f\right): N\in\D\right)\right\|_{\ell^2(\Z)}
	\lesssim\|f\|_{\ell^{2}(\Z)}.
	\end{align}
	To estimate \eqref{eq:291}, we split $\Pi[\le N^{C \varepsilon}, \le N^{-d+C\varepsilon}]$ dyadically into pieces corresponding to reduced fractions $\theta=\frac{a}{q}$ with $q\simeq 2^l$. Then, for a major arc frequency $\xi\in\T$ such that $\eta(N^{d-C\varepsilon}(\xi-\theta))\neq0$, we approximate $m_N(\xi)$ by a product $G(\theta) \mathfrak m_N(\xi-\theta)$ of arithmetic and continuous symbols given by 
	\begin{align*}
	G(\theta) \coloneqq \E_{n\in [q]}e(\theta P(n)) \quad \text{for} \quad \theta=\frac{a}{q} \text{ with } q\simeq 2^l,
	\end{align*}
	and
	\begin{align*}
	\mathfrak m_N(\xi) \coloneqq \int_0^1e(\xi P(Nt))dt
	\quad \text{for} \quad \xi\in\T,
	\end{align*}
	respectively. After a further factorization, the arithmetic part is summable over $l\in\N$ because $|G(\theta)|\lesssim q^{-c}\lesssim 2^{-cl}$, whereas the continuous part can be handled by Theorem~\ref{thm:IWvar}. This completes the outline of the proof of \eqref{eq:291}.
	
	\subsubsection{Weyl's inequality in $\ell^p(\Z)$ spaces and basic Ionescu--Wainger theory}
	The argument presented above for estimating \eqref{eq:104} when $k=1$ and $f\in\ell^2(\mathbb{Z})$ can be extended to $f\in\ell^p(\Z)$ for all $p\in(1, \infty)$. However, there are various challenges: for example, it is not clear how to make use of Weyl's inequality, as Plancherel's theorem is not available in $\ell^p(\Z)$ when $p \neq 2$. A way to overcome this difficulty is to proceed as follows:
	\begin{enumerate}[label*={(\roman*)}]
		\item One works directly with the Ionescu--Wainger projections from \eqref{eq:290}, as they are bounded on $\ell^p(\Z)$ for $p\in(1, \infty)$ and localized to the major arcs.
		
		\item Instead of working with the exponential sum $m_N$ itself, one works with the corresponding averaging operator $A_{N; \Z}^{P}$.
	\end{enumerate}
	
	The $\ell^p(\Z)$ norms of the Ionescu--Wainger projections $\Pi[\le \delta^{-C}, \le N^{-d}\delta^{-C}]$ have reasonably good growth in terms of the size of the set of canonical fractions $\mathcal R_{\le \delta^{-C}}$ due to our Ionescu--Wainger multiplier theorem. Namely, for any $p\in(1, \infty)$ and $\rho\in(0, 1)$, whenever
	\begin{align}
	\label{eq:296}
	N>C_p\delta^{-C_p}
	\end{align}
	for some large $C_p\in\R_+$, it follows from Theorem~\ref{thm:IW} that
	\begin{align}
	\label{eq:293}
	\left\|\Pi\left[\le \delta^{-C}, \le N^{-d}\delta^{-C}\right]f\right\|_{\ell^p(\Z)}
	\lesssim_{p, \rho} \delta^{-C\rho} \|f\|_{\ell^p(\Z)}.
	\end{align}
	This bound represents a significant quantitative improvement compared to \eqref{eq:292} and will be proved in Section~\ref{sec:IW} as a part of our Ionescu--Wainger multifrequency multiplier theory.
	
	Below are a few comments about this inequality.
	\begin{enumerate}[label*={\arabic*}.]	
		\item In an impactful and influential paper, Ionescu and Wainger \cite{IW} established a deep multiplier theorem for the set of so-called \textit{Ionescu--Wainger fractions} (see \eqref{eq:8} below for their definition), in order to establish $\ell^p(\Z^k)$ bounds for discrete singular integral Radon transforms. The Ionescu--Wainger multiplier theorem quickly became the main tool in the study of discrete analogues in harmonic analysis \cite{M1, MSZ3}.
		
		\item The set of Ionescu--Wainger fractions is defined as
		\begin{align}
		\label{eq:8}
		\tilde{\mathcal R}_{\le N} \coloneqq \left\{ \frac{a}{q} \in\Q/\Z: q\in P_{\le N} \text{ and } \gcd(a, q)=1 \right\},
		\end{align}
		where $P_{\le N}$ is a subtle set of natural numbers with certain prime power factorizations. In particular, the original Ionescu--Wainger theory \cite{IW} implies inequality \eqref{eq:293}, as well as the conclusion of Theorem~\ref{thm:IW}, for the projections defined as in \eqref{eq:290}, with the set of Ionescu--Wainger fractions $\tilde{\mathcal R}_{\le \delta^{-C}}$ in place of the set of canonical fractions $\mathcal R_{\le \delta^{-C}}$, provided that, instead of condition \eqref{eq:296}, we assume that 
		\begin{align}
		\label{eq:297}
		\log N\gtrsim \delta^{-C\gamma}
		\end{align}
		holds for some arbitrarily small $\gamma\in (0, 1)$. 
		
		\item Ionescu--Wainger theory \cite{IW} was originally developed for scalar-valued multipliers. An important aspect of the theory was that the $\ell^p(\Z)$ norm of the Ionescu--Wainger multipliers corresponding to the set of Ionescu--Wainger fractions $\tilde{\mathcal R}_{\le \delta^{-C}}$ was controlled by a multiple of $\log(\delta^{-C}+1)^D$, where $D \coloneqq \lfloor 2\gamma^{-1}\rfloor+1$. The proof in \cite{IW} is based on an intricate inductive argument that takes advantage of superorthogonality phenomena. A slightly different proof with the factor $\log(\delta^{-C}+1)$ in place of $\log(\delta^{-C}+1)^D$ was given in \cite{M1}. The latter proof relied on certain recursive arguments instead of the induction used in \cite{IW}. This approach helped  clarify the role of the underlying square functions and orthogonalities; see also \cite[Section~2]{MSZ3}. Ionescu--Wainger theory, among other topics, was discussed by Pierce \cite{Pierce} in the context of superorthogonality phenomena.  Finally, we refer the reader to the recent paper of Tao \cite{T}, where a uniform bound in place of $\log(\delta^{-C}+1)$ was obtained.
		
		\item The central question about Ionescu--Wainger theory for the set of canonical fractions from \eqref{eq:4} had remained unanswered since it was first formulated in \cite[Remark~3, p.~361]{IW}. In Section~\ref{sec:IW}, the question is answered affirmatively.	
	\end{enumerate}
	
	For any $p\in(1, \infty)$ and $\rho\in(0, 1)$, if $N>C_p\delta^{-C_p}$, then \eqref{eq:293} implies
	\begin{align}
	\label{eq:287}
	\left\|A_{N; \Z}^{P}\left(f-\Pi\left[\le \delta^{-C}, \le N^{-d}\delta^{-C}\right]f\right)\right\|_{\ell^p(\Z)}
	\lesssim_{p, \rho} \delta^{-C\rho} \|f\|_{\ell^p(\Z)}.
	\end{align}
	For $p=2$, by Plancherel's theorem and \eqref{eq:288}, we have the stronger bound
	\begin{align}
	\label{eq:295}
	\left\|A_{N; \Z}^{P}\left(f-\Pi\left[\le \delta^{-C}, \le N^{-d}\delta^{-C}\right]f\right)\right\|_{\ell^2(\Z)}
	\lesssim (\delta^c+N^{-c}) \|f\|_{\ell^2(\Z)}.
	\end{align}
	Thus, for any $p\in(1, \infty)$, by interpolating \eqref{eq:287} and \eqref{eq:295}, we obtain that there exist constants $C'_p \in \Z_+$ and
	$c_p\in(0, 1)$ such that if $N>C'_p\delta^{-C'_p}$, then
	\begin{align}
	\label{eq:294}
	\left\|A_{N; \Z}^{P}\left(f-\Pi\left[\le \delta^{-C}, \le N^{-d}\delta^{-C}\right]f\right)\right\|_{\ell^p(\Z)}
	\lesssim_{p} (\delta^{c_p}+N^{-c_p})\|f\|_{\ell^p(\Z)}.
	\end{align}
	Hence, inequality \eqref{eq:294} can be thought of as a generalization of Weyl's inequality for exponential sums \eqref{eq:288} to $\ell^p(\Z)$ spaces for all $p\in(1, \infty)$.  Here, it is important that the bounds in \eqref{eq:294} are consistent with the bounds in the classical Weyl inequality for exponential sums \eqref{eq:288}. An essential new tool that makes these estimates available is the Ionescu--Wainger multiplier theorem for the set of canonical fractions. Thus, Ionescu--Wainger theory should be understood as a tool that enables us to interpret exponential sum estimates in terms of purely functional analytic language. This interpretation is very useful in the context of the multilinear Weyl inequality stated in \eqref{eq:298} below.
	
	\subsubsection{The multilinear Weyl inequality and the minor arc contribution}
	In the multilinear setting (when $k\ge 2$), Plancherel's theorem and Weyl's inequality \eqref{eq:288} are no longer viable to effectively control the contribution from the minor arcs. As a result, it is not possible to implement the classical circle method as we discussed above.
	
	However, in view of the discussion about the generalization of Weyl's inequality to $\ell^p(\Z)$ spaces, it is natural to formulate the multilinear version using functional analysis. Namely, we expect that if $\mathcal P = \{P_1, \ldots ,P_k\}\subset\Z[\nn]$ is such that the degrees $d_i \coloneqq \deg P_i$ are distinct, and if $1<p_1,\ldots, p_k<\infty$ are such that $\frac{1}{p_1}+\cdots+\frac{1}{p_k}=\frac{1}{p}\le 1$, then for each fixed $C\in\R_+$ the following holds with some small $c\in (0, 1)$ possibly depending on $k,\mathcal P, p_1, \dots, p_k, p, C$. For every $f_i \in \ell^{p_i}(\Z^k)$ for all $i\in[k]$ and for all $N \geq 1$ and $\delta\in(0, 1]$, if $f_j \in \ell^{2}(\Z^k) \cap \ell^{p_j}(\Z^k)$ for some $j\in[k]$ and the $j$-th Fourier transform $\F_{j, \Z^k} f_j$ vanishes on the major arcs ${\mathfrak M}_{\leq N^{-d_j}\delta^{-C}}^j(\mathcal R_{\leq \delta^{-C}})$, then
	\begin{align}
	\label{eq:298}
	\left\| A^{P_1,\ldots, P_k}_{N; \Z^k}(f_1,\ldots, f_k) \right\|_{\ell^p(\Z^k)} 
	\lesssim (\delta^c+N^{-c}) \prod_{i \in [k]} \|f_i\|_{\ell^{p_i}(\Z^k)}.
	\end{align}
	Here, $\F_{j, \Z^k}$ and ${\mathfrak M}_{\leq N^{-d_j}\delta^{-C}}^j(\mathcal R_{\leq \delta^{-C}})$ denote, respectively, the Fourier transform and the set of major arcs defined with respect to the $j$-th variable.  Inequality \eqref{eq:298} is a multilinear variant of inequality \eqref{eq:294} and will be referred to as the \textit{multilinear Weyl inequality} throughout the paper. It will play a key role in facilitating our \textit{multilinear circle method}, namely:
	\[
	\text{multilinear Weyl's inequality}
	\ \leftrightsquigarrow \
	\text{control of the minor arc contribution.}
	\]
	
	The proof of inequality \eqref{eq:298}, as well as its real counterpart, the multilinear Sobolev smoothing inequality, is given in Section~\ref{sec:weyl}; see Theorem~\ref{weyl}.
	
	We now describe the key components of our \textit{multilinear circle method} in the context of establishing Theorem~\ref{thm:main1}. We begin by explaining how the multilinear Weyl inequality from \eqref{eq:298} controls the contribution from the minor arcs. If we take $\delta \simeq N^{-\varepsilon}$ with $\varepsilon\in(0, 1)$ sufficiently small, say $\varepsilon C<1/4$, and $N$ sufficiently large, then we gain a negative power of $N$ in \eqref{eq:298}. Hence, if we split $f_i = f_{i, N}^0+f_{i, N}^1$ for each $i \in [k]$, where
	\begin{align*}
	f_{i, N}^0 & \coloneqq f_i-\Pi\left[\le N^{C \varepsilon}, \le N^{-d_i+C\varepsilon}\right]f_i,\\
	f^1_{i, N} & \coloneqq \Pi\left[\le N^{C \varepsilon}, \le N^{-d_i+C\varepsilon}\right]f_i,
	\end{align*}
	then, in view of \eqref{eq:298} and from the fact that $N\in\mathbb D$, where $\mathbb D$ is lacunary, inequality \eqref{eq:104} is reduced to proving that
	\begin{align}
	\label{eq:299}
	\left\|{\bf V}^r\left(A_{N; \Z^k}^{P_1,\ldots, P_k}\left(f_{1, N}^1,\ldots, f_{k, N}^1\right): N\in\D\right)\right\|_{\ell^p(\Z^k)}\lesssim\prod_{i\in[k]}\|f_i\|_{\ell^{p_i}(\Z^k)},
	\end{align}
	where $f^1_{1, N},\ldots, f^1_{k, N}$ are supported on major arcs. Here, in this informal discussion, we work with the original average \eqref{eq:99} instead of its truncation \eqref{eq:103}. The idea of the proof of inequality \eqref{eq:299} will be discussed momentarily. 
	
	Our multilinear Weyl inequality with polynomial bounds \eqref{eq:298} is a new tool that we will also use in the analysis of the major arc contribution. A bilinear Weyl inequality was recently proved in \cite{KMT} for the Furstenberg--Weiss averages as a consequence of an inverse theorem of the third author \cite{P2}. Namely, under suitable assumptions, for some sufficiently large $\widetilde C \in \Z_+$, we have
	\begin{align}
	\label{eq:300}
	\left\| \E_{n\in[N]}f_1(x-n)f_2(x-P(n))\right\|_{\ell^p(\Z)} 
	\lesssim (\delta^c+ \langle \log N \rangle^{-c \widetilde C}) \prod_{i \in [2]} \|f_i\|_{\ell^{p_i}(\Z)}
	\end{align}
	for any polynomial $P\in\Z[\rm n]$ with degree at least two. However, these bounds are logarithmic in scale, and thus insufficient for the purpose of this paper. The logarithmic factor in the inequality \eqref{eq:300} is a result of applying Ionescu--Wainger theory with the Ionescu--Wainger fractions \eqref{eq:8}. The construction of the Ionescu--Wainger fractions forces condition \eqref{eq:297}, which subsequently leads to the presence of the logarithmic factor in \eqref{eq:300}.  In our approach, by applying Ionescu--Wainger theory with the set of canonical fractions \eqref{eq:4} as developed in Section~\ref{sec:IW}, we can quantitatively improve the bilinear estimates from \eqref{eq:300} and obtain a multilinear Weyl inequality with polynomial bounds in \eqref{eq:298} that match the bounds from the classical Weyl inequality for exponential sums \eqref{eq:288}. It is critical here that Ionescu--Wainger theory with the set of canonical fractions holds under condition \eqref{eq:296}, which enables us to obtain a negative power of $N$ in \eqref{eq:298}.
	
	Apart from Ionescu--Wainger theory for the set of canonical fractions, other essential tools needed to derive \eqref{eq:298} include an inverse theorem for averages over polynomial corner configurations with polynomial bounds, and a multilinear $\ell^p(\Z^k)$-improving inequality; see Theorem~\ref{pp} and Theorem~\ref{thm:improv}, respectively. In order to prove inequality \eqref{eq:298}, we first apply the inverse theorem, which reveals the major arc structure corresponding to the set of canonical fractions. By combining this theorem with the Hahn--Banach separation theorem, we can also reveal the major and minor arc structure for the adjoint operators corresponding to $A^{P_1,\ldots, P_k}_{N; \Z^k}(f_1,\ldots, f_k)$. Next, we use our version of Ionescu--Wainger theory to maintain the distinction between major and minor arcs and to preserve the polynomial bounds in terms of $\delta$ obtained from the inverse theorem throughout this process. Finally, the $\ell^p(\Z^k)$-improving inequality is applied to relax the $\ell^{\infty}(\Z^k)$ bounds arising in the inverse theorem to $\ell^{p}(\Z^k)$ bounds that we need here. This ultimately implies the multilinear Weyl inequality with polynomial bounds \eqref{eq:298}, as desired.
	
	\subsubsection{The inverse theorem for averages over polynomial corner configurations}
	In 2019, the third author developed a new technique in higher order Fourier analysis that she used to prove power-saving bounds in the finite field version of the polynomial Szemer\'edi theorem for linearly independent polynomial progressions~\cite{Pff}, answering a question of Bourgain and Chang~\cite{BC}. This method, now called \textit{degree lowering}, was improved and adapted to the integer setting by the third author and Prendiville~\cite{PP1}, proving the first quantitative bounds for the size of subsets of $[N]$ lacking nontrivial copies of the \textit{nonlinear Roth configuration} $x,x+n,x+n^2$. The third author then extended this result to all distinct degree polynomial progressions~\cite{P2}. The key technical result of~\cite{P2}, from which quantitative bounds in the polynomial Szemer\'edi theorem follow by a standard density increment iteration, is an inverse theorem asserting that if averages of $1$-bounded functions over fixed distinct degree polynomial progressions are large, then the underlying $1$-bounded functions must correlate with functions whose Fourier transforms are supported on major arcs.
	
	\begin{theorem}[Peluse, 2020]\label{thm:P}
		Let $k\in\mathbb{Z}_+$ and $P_1,\dots,P_k\in\mathbb{Z}[\rm n]$ be polynomials with vanishing constant terms and distinct degrees $d_1\coloneqq\deg{P_1}<\dots<d_k\coloneqq\deg{P_k}$. Then there exist large constants $C_1,C_2 \in \Z_+$ such that the following holds. Assume that $\delta\in(0, 1]$ and $N \geq C_{1} \delta^{-C_{1}}$. If $f_0, f_1,\ldots, f_k \in \ell^\infty(\Z)$ are $1$-bounded functions supported on $[ \pm O (N^{d_k})]$ such that
		\begin{align*}
		\bigg|\bigg\langle f_0, \E_{n\in[N]}\prod_{i\in[k]}f_i(x-P_i(n))\bigg\rangle\bigg| \geq \delta N^{d_k},
		\end{align*}
		then, using the projections from \eqref{eq:290}, one has
		\begin{align*}
		\left|\left\langle f_1, \Pi [\le C_{2} \delta^{-C_{2}},
		\le C_{2} \delta^{-C_{2}}N^{-d_1}]f_1\right\rangle\right|
		\ge C_{1}^{-1} \delta^{C_{1}} N^{d_k}.
		\end{align*}
	\end{theorem}
	When $k=2$, one can obtain strong structural information about $f_0$ and $f_2$ from Theorem~\ref{thm:P}; see~\cite{PP2} for details. This structural information was then used by the second author along with Krause and Tao in~\cite{KMT} to prove their bilinear Weyl inequality for Furstenberg--Weiss averages \eqref{eq:300}. When $k>2$, the argument from~\cite{PP2} breaks down, and so structural information cannot be deduced for $f_2,\dots,f_k$; see the discussion in the last section of~\cite{PP2}.
	
	In order to prove our multilinear Weyl inequality, we require an extension of Theorem~\ref{thm:P} to higher dimensions that simultaneously produces strong structural information about each function appearing in the average. Such an inverse theorem is proved in Section~\ref{sec:inverse} in both the integer and real settings; see Theorem~\ref{pp}. This inverse theorem simultaneously generalizes work of the third author~\cite{P2} proving Theorem~\ref{thm:P} and of the third author along with Prendiville and Shao proving an inverse theorem for averages over sqorners
	\[
	(x_1, x_2), (x_1+n, x_2), (x_1, x_2+n^2).
	\]
	
	\subsubsection{The multilinear $\ell^p(\Z^k)$-improving inequality}
	In Section~\ref{sec:improving}, we establish the following multilinear $\ell^p(\Z^k)$-improving inequality. For $k\ge 2$ and $j\in [k]$, there are exponents $1< p < q < 2$ and $1<p_1,\ldots, p_k<\infty$ satisfying $\frac{1}{p_1}+\cdots+\frac{1}{p_k}=\frac{1}{p}$ and $p_j = 2$ such that the following scale-invariant bound
	\begin{align}
	\label{eq:302}
	\| A_{N; \mathbb Z^k}^{\mathcal P}(f_1,\ldots, f_k) \|_{\ell^q(\mathbb Z^k)}
	\lesssim_{\mathcal P} N^{-D(\frac{1}{p}-\frac{1}{q})} 
	\prod_{i \in [k]} \|f_{i}\|_{\ell^{p_i}(\mathbb Z^k)}
	\end{align}
	holds with $D \coloneqq d_1 + \cdots + d_k$.  The bound \eqref{eq:302} plays a key role in relaxing the $\ell^{\infty}(\Z^k)$ bounds required to apply the inverse theorem to $\ell^{p}(\Z^k)$ bounds. It is essential that our multilinear $\ell^p(\Z^k)$-improving inequality is obtained without any $\varepsilon$-loss in the scale $N$. Otherwise, we could not derive the multilinear Weyl inequality with polynomial bounds in \eqref{eq:298}.
	
	Such improving bounds have a long history in the continuous setting (with $\R$ in place of $\Z$) but it is only recently that they were established in the integer setting and then only for linear averages. In \cite{KMT}, the bounds \eqref{eq:302} were needed for particular bilinear averages and were deduced from the linear case of Han--Kova{\v c}--Lacey--Madrid--Yang in \cite{HKLMY} by passing to the dual operator. In the general case considered here, genuinely multilinear $\ell^p(\Z^k)$-improving estimates are needed, and it is no longer possible to reduce matters to the linear case.
	
	To prove \eqref{eq:302} with no $\varepsilon$-loss, we marry Christ's {\it refinements method} \cite{Ch} (which has proved very successful for $L^p$-improving estimates in the continuous setting; see \cite{stovall}, where multilinear $L^p$ improving bounds are established for Radon-like transforms), adapted to the multilinear setting, with the Vinogradov mean value theorem \cite{BDG, WOOLEY}.
	
	\subsubsection{Major arcs contribution: paraproduct-type decomposition}
	The objective now is to understand estimate \eqref{eq:299}. Our aim is to reduce \eqref{eq:299} to estimates for certain model operators, which will be more manageable. We may normalize each $f_i$ so that $\|f_i\|_{\ell^{p_i}(\Z^k)} = 1$. We proceed in two steps:
	\begin{enumerate}[label*={(\roman*)}]
		\item We begin with trimming the size of denominators. This means that each function $f^1_{i, N}=\Pi[\le N^{C \varepsilon}, \le N^{-d_i+C\varepsilon}]f_i$ can be split into finitely many pieces $f^{1, l_i}_{i, N}$, where each piece corresponds to projections onto major arcs centered at fractions in $\mathcal R_{\le N^{C \varepsilon}}$ whose denominators have size of order $2^{l_i}$ for some $l_i\in\N$ such that $2^{l_i}\le N^{C \varepsilon}$. We have seen a similar reduction in the linear case above.
		
		\item Next, we perform a certain paraproduct decomposition of each $f^{1, l_i}_{i, N}$ to match the scales of the bumps in the definition of the Ionescu--Wainger projections. In other words, we replace $\eta(N^{d_i-C\varepsilon}(\xi-\theta))$ in the definition of $f^{1, l_i}_{i, N}$ by $\eta(N^{d_i}(\xi-\theta))$ for $\theta=\frac{a}{q}$ with $q\simeq 2^{l_i}$.  This can be done by a simple telescoping argument, where we further split $f^{1, l_i}_{i, N}$ into finitely many highly oscillatory and nonoscillatory pieces $f^{1, l_i, s_i}_{i, N}$ for some $s_i\in\N$ such that $2^{s_i}\le N^{C \varepsilon}$. If $s_i>0$, then $f^{1, l_i, s_i}_{i, N}$ is highly oscillatory, which means that its $i$-th Fourier transform is supported on an annulus of size $N^{d_i}2^{-s_i}$ centered at $\theta=\frac{a}{q}$. Otherwise, if $s_i=0$, then $f^{1, l_i, 0}_{i, N}$ is nonoscillatory and corresponds to the bumps $\eta(N^{d_i}(\xi-\theta))$ for $\theta=\frac{a}{q}$ with
		$q\simeq 2^{l_i}$.
	\end{enumerate}
	This reduction is standard, and it was also performed in \cite{KMT}. By using this decomposition, inequality \eqref{eq:299} will follow if we can show that, for all $l_1, s_1,\ldots, l_k, s_k\in\N$ and for some $c\in (0, 1)$, the following inequality holds
	\begin{align}
	\label{eq:303}
	\left\|{\bf V}^r\left(A_{N; \Z^k}^{P_1,\ldots, P_k}\left(f^{1, l_1, s_1}_{1, N},\ldots, f^{1, l_k, s_k}_{k, N}\right): N\in\D\right)\right\|_{\ell^p(\Z^k)}\lesssim 2^{-c(l+s)},
	\end{align}
	with $l \coloneqq \max\{l_i: i\in[k]\}$ and $s \coloneqq \max\{s_i: i\in[k]\}$, where each $f_i$ is normalized so that $\|f_i\|_{\ell^{p_i}(\Z^k)} = 1$ and $ f^{1, l_i, s_i}_{i, N} \equiv 0$ if $N^{C\varepsilon} < 2^{\max\{l_i,s_i\}}$. Here our approach deviates significantly from \cite{KMT}. We distinguish two cases:
	\begin{itemize}
		\item[(i)] the \textit{high frequency} case $s \gtrsim l$;
		\item[(ii)] the \textit{low frequency} case $l \gtrsim s$.
	\end{itemize}
	By standard approximations on the major arcs, inequality \eqref{eq:303} can be reduced to showing
	\begin{align}
	\label{eq:304}
	\left\|{\bf V}^r\left(A_{2^u; \Z^k}^{P_1,\ldots, P_k}\left(f^{1, l_1, s_1}_{1, N},\ldots, f^{1, l_k, s_k}_{k, N}\right): N\in\D\right)\right\|_{\ell^p(\Z^k)}\lesssim 2^{-c(l+s)},
	\end{align}
	where $u \coloneqq  100k(s+1)$ in the high frequency case, and $u \coloneqq  100k(l+1)$ in the low frequency case. An important feature of this reduction is that the operator $A_{2^u; \Z^k}^{P_1,\ldots, P_k}$ is independent of the parameter $N\in\mathbb D$.
	
	\subsubsection{The major arc contribution: high frequency case}
	In the high frequency case, the multifrequency problem becomes, in fact, a single frequency problem. By using Taylor expansions, we can separate the nonoscillatory pieces from the highly oscillatory pieces. The nonoscillatory pieces are handled by a seminorm variant of the Ionescu--Wainger theorem; see Theorem~\ref{thm:IWvar}. This theorem is interesting in its own right, as it eliminates the need to consider small and large scales as done previously. The highly oscillatory pieces are controlled by the multilinear Sobolev smoothing inequality, which is a real variant of the multilinear Weyl inequality. This inequality allows us to gain a negative power of $2^{s}$ for the highly oscillatory pieces, which yields the desired decay in \eqref{eq:304} in the high frequency case thanks to $s \gtrsim l$.
	
	This is a novel aspect of our argument and shows that the Sobolev smoothing inequality and the multilinear Weyl inequality, as well as the inverse theorem in both the integer and real settings, should always be proved in tandem.
	
	\subsubsection{The major arc contribution: low frequency case}
	To control the low frequency case, we will use a multiparameter norm interchanging inequality \eqref{eq:84} combined with the multilinear Weyl inequality \eqref{eq:298}. This is an important part of our argument, where the multilinear Weyl inequality is also useful for controlling the contribution from the major arcs. Precisely, it will allow us to gain a negative power of $2^l$, which yields the desired decay in \eqref{eq:304} in the low frequency case thanks to $l \gtrsim s$.
	
	A multiparameter norm interchanging inequality, to a certain extent, linearizes maximal functions in a way that makes it possible to apply, for instance, the multilinear Weyl inequality or similar tools. By a simple interpolation, we can assume that $r=\infty$ in \eqref{eq:304} and focus on maximal functions. Our approach, which relies on the multiparameter norm interchanging inequality, also eliminates the need to consider small and large scales.
	
	Finally, we highlight that our argument heavily depends on the polynomial bounds in the multilinear Weyl inequality \eqref{eq:298}. In the bilinear case treated in \cite{KMT}, the parameter $u$ from \eqref{eq:304} is exponential in $l$, while in our case it is linear. Furthermore, the Weyl inequality from \cite{KMT} was proved with decay $\delta^c+ \langle \log N \rangle^{-c\widetilde C}$, which includes a logarithmic term that requires $u$ to be exponential in $l$. However, thanks to our improvement in the Ionescu--Wainger theorem, our Weyl inequality holds with $\delta^c+N^{-c}$, allowing us to have a linear dependence between $l$ and $u$. This is crucial for handling the low frequency case and arranging the sizes of supports in certain model operators in terms of \( u \) in such a way that their contribution to the norm of these model operators becomes negligible compared to the size of \( cl \) with $c$ arising in \eqref{eq:298}. If $u$ were exponential in $l$, then our method would not work.
	
	\subsection{Final remarks and open problems}
	We have given an outline of the proof of inequality \eqref{eq:104}, discussing the crucial tools and steps of our multilinear circle method. We wrap up the introduction by giving a few concluding remarks and open problems.
	
	\begin{enumerate}[label*={\arabic*}.]
		\item The multiparameter norm interchanging inequality approach allows us to avoid the $p$-adic methods used in \cite{KMT} and significantly simplifies the argument. Moreover, the $p$-adic approach from \cite{KMT} (especially sharp $L^p$-improving inequalities in the $p$-adic setting), even if adjusted to our context, would not be sufficient to handle the major arc estimates in the low frequency setting.
		
		\begin{problem}
			Is a $p$-adic analogue of Theorem~\ref{thm:main1} (or its maximal counterpart, which should be simpler) true?
		\end{problem}
		
		\item Our Ionescu--Wainger theory is not uniform with respect to the size of the set of canonical fractions \eqref{eq:4}. However, based on Tao's recent paper \cite{T}, we know that the Ionescu--Wainger multiplier theorem for the set of Ionescu--Wainger fractions \eqref{eq:8} holds with uniform bounds.
		
		\begin{problem}
			Does our Ionescu--Wainger multiplier theorem hold with bounds independent of the size of the set of canonical fractions?
		\end{problem}
		
		\item The methods outlined in this paper are robust and fairly versatile.  Specifically, the input from harmonic analysis, such as the multilinear $\ell^p(\Z^k)$-improving inequality or the proof of Theorem~\ref{thm:main1}, does not rely on the assumption that the polynomials have distinct degrees. This assumption is only used in the proof of the inverse theorem and is not used anywhere else in the paper. In other words, if we treated Theorem~\ref{pp} as a black box, we could derive all other conclusions from this paper without needing to assume that our polynomials have distinct degrees. With this in mind, we hope that the methods outlined in the paper will be used in other contexts to address similar questions. 	
		
		\begin{problem}
			Is Theorem~\ref{thm:main} for the weighted operators
			\begin{align}
			\label{eq:5}
			\sum_{n\in[\pm N]\setminus\{0\}}\frac{1}{n}\prod_{i \in [k]} f_i\big(T_i^{P_i(n)}x\big)
			\end{align}
			in place of the averages $A_{N; X, {\mathcal T}}^{\mathcal P}$ from \eqref{eq:98} true? The operators from \eqref{eq:5} are multilinear truncated singular integral variants of $A_{N; X, {\mathcal T}}^{\mathcal P}$. It is also interesting to know any $\ell^p(\Z^k)$ bounds for the full multilinear singular integral operators, i.e., the operators corresponding to \eqref{eq:5} with summation over $[\pm N]\setminus\{0\}$ replaced with $\Z\setminus\{0\}$.
		\end{problem}
		
		\item The assumption that our polynomials have distinct degrees is only used in the proof of the inverse theorem, as it is essential to run the degree lowering argument in Section~\ref{sec:inverse}. However, the proof of Theorem~\ref{pp} also demonstrates that the assumption that the polynomials have distinct degrees in the multilinear setting may have a similar impact as curvature does in the linear theory of maximal and singular Radon transforms in the discrete \cite{M1, MSZ3} or Euclidean \cite{bigs} settings. This is apparent when we start with large multilinear averages \eqref{hyp} and apply the PET induction scheme and the degree lowering argument and end up with an expression as in \eqref{eq:276} that corresponds to modulated linear polynomial averages that can be handled by linear methods, which exploit the curvature. This interesting phenomenon requires further investigation. In view of this remark, it is natural to ask the following.
		
		\begin{problem}
			Is it possible to establish an inverse theorem like in Theorem~\ref{pp} for linearly independent polynomials with polynomial bounds in $\delta$? That would be an important step toward proving Theorem~\ref{thm:B} for linearly independent polynomials.
		\end{problem}
		
		\item A well-known open problem in ergodic theory is to establish  pointwise almost everywhere convergence of
		bilinear averages 
		\[
		\E_{n\in [N]}f_1(S^{n}x)f_2(T^{n}x), \qquad x\in X,
		\]
		corresponding to ``corner'' configurations, for two commuting measure-preserving transformations $S, T \colon X\to X$; see \cite[Problem~19]{Fra}. In the single transformation setting, the pointwise almost everywhere convergence of the trilinear averages
		\[
		\E_{n\in [N]}f_1(T^{n}x)f_2(T^{2n}x)f_3(T^{3n}x), \qquad x\in X,
		\]
		corresponding to the arithmetic progressions of length four, is another well-known open problem; see \cite[Problem~11]{Fra}. Our methods break down in these cases, as the minor arc contributions are no longer negligible, and we are therefore unable to say anything in this direction. In the trilinear case, identities such as
		\[
		e(\alpha x^2) \cdot e(-3\alpha (x-n)^2) \cdot e(3\alpha (x-2n)^2) \cdot e(-\alpha(x-3n)^2)=1
		\]
		speak to the presence of a \emph{quadratic} modulation invariance, which would most likely require developing quadratic time--frequency analysis in the spirit of \cite{LAC, DLTT}; this is out of reach at the moment. However, we think that our methods represent a significant step in understanding the Furstenberg--Bergelson--Leibman conjecture.
		
		\begin{problem}
			How useful are the methods presented in this paper in understanding Question~\ref{con:B} for linearly dependent
			polynomials?
		\end{problem}		
	\end{enumerate}
	
	\section*{Acknowledgments} We warmly thank Vitaly Bergelson, Nikos Frantzikinakis, Ben Green, Noah Kravitz, Borys Kuca, Tomasz Z. Szarek, and Terence Tao for helpful comments on an earlier draft of this paper. The first author gratefully acknowledges the support of the Basque Center for Applied Mathematics, where he carried out his postdoctoral research under the supervision of Renato Luc\`{a} and Luz Roncal, as well as the support of the Wroc{\l}aw University of Science and Technology through the Academia Iuvenum membership. The second and third authors gratefully acknowledge the support and hospitality of the Institute for Advanced Study, which helped facilitate this collaboration. Finally, we thank the referees for careful reading of the manuscript and useful remarks that led to the improvement of the presentation.
	
	\section{Notation and useful tools}\label{section:2}
	We now set up notation that will be used throughout the paper. 
	
	\subsection{Basic notation}  The set of positive integers and the set of nonnegative integers will be denoted, respectively, by $\Z_+ \coloneqq \{1, 2, \ldots\}$ and $\N \coloneqq \{0,1,2,\ldots\}$. For $d\in\Z_+$ the sets $\Z^d$, $\R^d$, $\mathbb C^d$ and $\T^d \coloneqq (\R/\Z)^d$ have their standard meaning. We also denote $\R_+ \coloneqq (0, \infty)$.
	
	For $a = (a_1,\ldots, a_d) \in \Z^d$ and $q\in\Z_+$, we denote by $\gcd(a,q)$ the greatest common divisor of $a$ and $q$, i.e., the largest  $n \in\Z_+$ that divides $q$ and all components $a_1, \ldots, a_d$. Clearly, any vector in $\Q^d$ has a unique representation as $\frac{a}{q}$ with $q\in \Z_{+}$, $a \in \Z^d$, and $\gcd(a,q)=1$.
	
	We use $\ind{A}$ to denote the indicator function of a set $A$. For a statement $S$ we write $\ind{S}$ to denote its indicator, equal to $1$ if $S$ is true and $0$ if $S$ is false. In particular, we have  $\ind{A}(x)=\ind{x\in A}$.
	
	Throughout the paper $C>0$ is an absolute constant that may change from occurrence to occurrence. For two nonnegative quantities $A, B$ we write $A \lesssim B$ if there is an absolute constant $C>0$ such that $A\le CB$. We also write $A \simeq B$ when $A \lesssim B\lesssim A$. We will write $\lesssim_{\delta}$ or $\simeq_{\delta}$ to emphasize that the implicit constant depends on $\delta$. For a function $f \colon X\to \C$ and a positive-valued function $g \colon X\to (0, \infty)$, we write $f = O(g)$ if there exists a constant $C>0$ such that $|f(x)| \le C g(x)$ for all $x\in X$. We will also write $f = O_{\delta}(g)$ if the implicit constant depends on $\delta$.
	
	We will use the Japanese bracket notation
	\[
	\langle z \rangle \coloneqq (1 + |z|^2)^{1/2}\simeq 1+|z|, \qquad z\in\C.
	\]
	For any $x\in\R$ we will use the floor and fractional part functions
	\begin{align*}
	\lfloor x \rfloor \coloneqq \max\{ n \in \Z : n \le x \}
	\quad \text{and} \quad
	\{x\} \coloneqq x-\lfloor x \rfloor.
	\end{align*}
	All logarithms in this paper are taken base $2$, and for any $N \geq 1$ we define 
	\begin{equation}\label{log-scale}
	\Log N \coloneqq \lfloor \log N \rfloor,
	\end{equation}
	i.e., the unique natural number such that $2^{\Log N} \leq N < 2^{\Log N+1}$.
	
	For every $N\in\R_+$ we set
	\begin{align}
	\begin{split}\label{intervals}
	[N]_{\Z} & \coloneqq  (0, N] \cap \Z, \phantom{N],} 
	[\pm N]_{\Z} \coloneqq  [-N, N] \cap \Z, \quad \\
	[N]_{\R} & \coloneqq  (0, N], \phantom{N] \cap \Z,}
	[\pm N]_{\R} \coloneqq  [-N, N].
	\end{split}
	\end{align}
	We shall abbreviate $[N]_{\Z}$ to $[N]$ and $[\pm N]_{{\Z}}$ to $[\pm N]$. We will also write
	\begin{align*}
	\N_{\le N} & \coloneqq  [0, N]\cap\N, \quad \phantom{)}
	\N_{\ge N} \coloneqq  [N, \infty)\cap\N, \\
	\N_{< N} & \coloneqq  [0, N)\cap\N, \quad \phantom{]}
	\N_{> N} \coloneqq  (N, \infty)\cap\N.
	\end{align*}
	
	For any $\lambda>1$ we will say that the set
	\begin{align}\label{eq:46}
	\mathbb D=\{\lambda_n : n\in\N\} \subset [1,+\infty)
	\end{align}
	is $\lambda$-lacunary if $\inf_{n\in\N}\frac{\lambda_{n+1}}{\lambda_n}\ge\lambda$.
	
	Finally, $\mathcal C \colon \mathbb C\to \mathbb C$ will be the complex conjugation, i.e., $\mathcal Cz \coloneqq \overline{z}$ for $z\in \mathbb C$.
	
	\subsection{Euclidean spaces} For  $d\in\Z_+$ let $\{e_i\in\R^d:i\in[d]\}$ be the standard basis in $\R^d$. The standard inner product and the corresponding Euclidean norm on $\R^d$ are denoted, respectively, by 
	\begin{align*}
	x\cdot\xi \coloneqq \sum_{i \in [d]} x_i \xi_i 
	\quad \text{and} \quad
	|x| \coloneqq |x|_2 \coloneqq \sqrt{x\cdot x }
	\end{align*}
	for every $x=(x_1,\ldots, x_d)$ and $\xi=(\xi_1, \ldots, \xi_d)\in\R^d$. 
	
	Throughout the paper the $d$-dimensional torus $\T^d$, which unless otherwise stated will be  identified with $[-1/2, 1/2)^d$, is endowed with the periodic norm
	\begin{align*}
	\|\xi\| \coloneqq \Big(\sum_{i \in [d]} \|\xi_i\|^2\Big)^{1/2},
	\qquad
	\xi=(\xi_1,\ldots,\xi_d)\in\T^d,
	\end{align*}
	where $\|\xi_i\| \coloneqq {\rm dist}(\xi_i, \Z)$ for every $\xi_i\in\T$ and $i\in[d]$. However, identifying $\T^d$ with $[-1/2, 1/2)^d$, we see that the norm $\|\,\cdot\,\|$ coincides with the Euclidean norm $|\,\cdot\,|$ restricted to $[-1/2, 1/2)^d$.
	
	Throughout this paper we work with a fixed cutoff function $\eta \colon \R\to[0, 1]$ that is a smooth and even function such that
	\begin{align}\label{eq:126}
	\ind{[-1/4, 1/4]}\le\eta\le \ind{(-1/2, 1/2)}.
	\end{align}
	All constants are permitted to depend on $\eta$.
	
	\subsection{Function spaces}
	All vector spaces in this paper will be defined over the complex numbers $\C$. The triple $(X, \mathcal B(X), \mu)$ is a measure space $X$ with $\sigma$-algebra $\mathcal B(X)$ and $\sigma$-finite measure $\mu$. The space of all $\mu$-measurable complex-valued functions defined on $X$ will be denoted by $L^0(X, \mu)$. The space of all functions in $L^0(X, \mu)$ whose modulus is integrable with $p$-th power is denoted by $L^p(X, \mu)$ for $p\in(0, \infty)$, whereas $L^{\infty}(X, \mu)$ denotes the space of all essentially bounded functions in $L^0(X, \mu)$. These notions can be extended to functions taking values in a finite-dimensional normed vector space $(B, \|\cdot\|_B)$. For instance, we define
	\begin{align*}
	L^{p}(X,\mu;B)
	\coloneqq \left\{F\in L^0(X, \mu;B):\|F\|_{L^{p}(X, \mu;B)} \coloneqq \left\|\|F\|_B\right\|_{L^{p}(X, \mu)}<\infty\right\},
	\end{align*}
	where $L^0(X, \mu;B)$ denotes the space of measurable functions from $X$ to $B$ (up to almost everywhere equivalence). One can extend these notions to infinite-dimensional $B$, at least when $B$ is separable. In this paper, however, we will always be able to work in finite-dimensional settings by appealing to  standard approximation arguments. For ease of notation, we shall abbreviate $L^p(X, \mu)$ to $L^p(X)$ and $L^{p}(X,\mu;B)$ to $L^{p}(X;B)$ for $p\in[0, \infty]$.
	
	In our case we will mainly work with $X=\R^d$ or $X=\T^d$ endowed with Lebesgue measure, or with $X=\Z^d$ endowed with counting measure. If $X$ is endowed with counting measure, then we abbreviate $L^p(X)$ to $\ell^p(X)$ and $L^p(X; B)$ to $\ell^p(X; B)$.
	
	For a continuous linear map $T \colon B_1 \to B_2$ between two normed vector spaces $B_1$ and  $B_2$, we use $\|T\|_{B_1 \to B_2}$ to denote its operator norm.
	
	We say that $f\in L^0(X)$ is $1$-bounded if $f\in L^{\infty}(X)$ and $\|f\|_{L^\infty(X)}\le 1$. We will mainly work with $1$-bounded functions in Sections~\ref{sec:inverse}~and~\ref{sec:weyl}.
	
	For any $F\in L^0(X, \mu)$ and a measurable set $A\subseteq X$ with $0<\mu(A)<\infty$ we will use the convenient notation
	\begin{align}
	\label{eq:30}
	\E_{t \in A}^{\mu}F(t) \coloneqq \frac{1}{\mu(A)}\int_AF(t)d\mu(t).
	\end{align}
	If $X$ is endowed with counting measure $\mu(A) \coloneqq \#A$, then for $\emptyset\neq A\subseteq X$ we will abbreviate $\E_{t \in A}^{\mu}F(t)$ to $\E_{t \in A}F(t)$, which is the average of $F$ over $A$, i.e., $\E_{t \in A}F(t) \coloneqq \frac{1}{\# A}\sum_{t\in A}F(t)$. Sometimes, to simplify notation, we will write $|A|$ instead of $\#A$ for the cardinality of a set $A$; see, for example, Section~\ref{sec:inverse}.
	
	\subsection{Convolutions} Let $(\G, +)$ be a locally compact abelian group (LCA group) equipped with a Haar measure $\lambda_{\G}$. For two functions $f, g\in L^1(\G)$ we define their convolution by 
	\begin{align*}
	f*_{\G}g(x) \coloneqq f*g(x) \coloneqq\int_{\G}f(x-y)g(y)d\lambda_{\G}(y).
	\end{align*}
	For the sake of simplicity, we will usually abbreviate $*_{\G}$ to $*$;  the meaning of this notation will always be clear from the context and will not cause any confusion. Since $\G$ is abelian, we readily see that $f*g=g*f$.
	
	Fix $k\in\Z_+$ and, for each $i\in[k]$, let $(\G_i, +)$ be a locally compact abelian group equipped with a Haar measure $\lambda_{\G_i}$. We consider their product $\G=\G_1\times\cdots\times \G_k$ with the product measure $\lambda_{\G}=\lambda_{\G_1}\times\cdots\times\lambda_{\G_k}$. Given $i \in[k]$, for any $f\in L^1(\G)$ and $g\in L^1(\G_i)$ we may define their $i$-th convolution by 
	\begin{align*}
	f*_{\G_i} g(x) \coloneqq f*_i g(x) \coloneqq \int_{\G_i} f(x_1,\ldots, x_{i-1}, x_i-y, x_{i+1},\ldots, x_k)g(y)d\lambda_{\G_i}(y),  
	\end{align*}
	for $x=(x_1,\ldots, x_k)\in\G$. In other words, $f*_i g(x)$ coincides with the function $f$ except the $i$-th variable where it is defined as a convolution with the function $g$. We will also write $f*_{\G_i} g(x)=f*_{\G, i} g(x)$ if $\G_1=\cdots=\G_k=\G$. We will mainly work with $\G=\R^d$ endowed with Lebesgue measure or with $\G=\Z^d$ endowed with counting measure.	
	
	\subsection{Fourier transform}
	Although we will use Fourier analysis only on $\R^d$, $\T^d$, and $\Z^d$, it will be convenient to set out some abstract harmonic analysis notation to perform this analysis in a unified fashion and to avoid repetition. We shall write $e(z) \coloneqq e^{2\pi {\bm i} z}$ for every $z\in\C$, where ${\bm i}^2=-1$. 
	
	As before, let $(\G, +)$ be an LCA group equipped with a Haar measure $\lambda_{\G}$. It is well known (see, for instance, \cite{Rudin}) that every LCA group $\G$ has a Pontryagin dual $\hat \G = (\hat \G,+)$, an LCA group with a Haar measure $\lambda_{\hat \G}$ and a pairing, i.e., a continuous bihomomorphism $\G \times \hat \G\ni (x,\xi) \mapsto \langle x , \xi\rangle\in\T$, such that the Fourier transform $\F_{\G} \colon L^1(\G) \to C(\hat \G)$ given by
	\[
	\F_{\G} f(\xi) \coloneqq \int_{\G} f(x) e(\langle x , \xi\rangle)d\lambda_{\G}(x), \qquad \xi\in\hat \G,
	\]
	extends to a unitary map from $L^2(\G)$ to $L^2(\hat \G)$. We have Plancherel's identity
	\[
	\|\F_\G f\|_{L^2(\hat \G)}=\|f\|_{L^2(\G)}, \qquad f\in L^2(\G),
	\]
	and the inverse Fourier transform $\F_\G^{-1} \colon L^2(\hat \G) \to L^2(\G)$ is given by the formula
	\[
	\F_\G^{-1} f(x) = \int_{\hat \G} f(\xi) e(-\langle x , \xi\rangle)d\lambda_{\hat \G}(\xi), \qquad x\in\G,
	\]
	whenever $f\in L^1(\hat \G) \cap L^2(\hat \G)$. We will mainly work with concrete pairs $(\G,\hat \G)$ of Pontryagin dual LCA groups:
	\begin{itemize}
		\item [(i)] Consider $\G = \R^d$ with Lebesgue measure $\lambda_{\R^d} = dx$. Then $\hat \G = \R^d$ with Lebesgue measure $\lambda_{\R^d} = d\xi$ is a Pontryagin dual, with pairing $\langle x , \xi\rangle \coloneqq x \cdot \xi$. Moreover, for any $f \in L^1(\R^d)$ we have
		\begin{align*}
		\mathcal F_{\R^d} f(\xi)  \coloneqq  \int_{\R^d} f(x) e(x\cdot\xi) d x, \qquad \xi\in\R^d.
		\end{align*}
		
		\item[(ii)]  Consider $\G = \Z^d$ with counting measure $\lambda_{\Z^d}$. Then $\hat \G = \T^d$ with Lebesgue measure $\lambda_{\T^d} = d\xi$ is a Pontryagin dual, with pairing $\langle x , \xi\rangle \coloneqq x \cdot \xi$. Moreover, for any $f \in \ell^1(\Z^d)$ we have
		\begin{align*}
		\mathcal F_{\Z^d}f(\xi) \coloneqq \sum_{x \in \Z^d} f(x) e(x\cdot\xi), \qquad \xi\in \T^d.
		\end{align*}
		
		\item[(iii)]  Consider LCA groups $\G_1,\ldots, \G_k$ with Pontryagin duals $\hat \G_1,\ldots, \hat \G_k$. Then the product $\G=\G_1 \times\cdots \times \G_k$ with product Haar measure $\lambda_{\G}=\lambda_{\G_1}\times\cdots\times\lambda_{\G_k}$ is an LCA group with Pontryagin dual $\hat \G=\hat \G_1 \times\cdots\times \hat \G_k$ with product Haar measure $\lambda_{\hat \G}=\lambda_{\hat \G_1}\times\cdots\times\lambda_{\hat \G_k}$ and pairing $\langle (x_1,\ldots, x_k) , (\xi_1,\ldots, \xi_k)\rangle  \coloneqq  \langle x_1, \xi_1\rangle +\cdots+ \langle x_k,  \xi_k\rangle$. 
	\end{itemize}
	
	For LCA product groups $\G=\G_1 \times\cdots \times \G_k$ we consider partial Fourier transforms. Given $i\in[k]$, for $x=(x_1,\ldots, x_k)\in\G$ and $\xi_i \in \hat \G_i$ we set 
	\[
	x_{(\xi_i)} \coloneqq (x_1,\ldots, x_{i-1},\xi_i,x_{i+1},\ldots, x_k).
	\] 
	Then, for any $f\in L^1(\G)\cap L^2(\G)$, the $i$-th Fourier transform is given by
	\begin{align*}
	\F_{i, \G} f(x_{(\xi_i)}) \coloneqq
	\int_{\G_i} f(x) e(\langle x_i , \xi_i \rangle)d\lambda_{\G_i}(x_i),
	\end{align*}
	for $x=(x_1,\ldots, x_k)\in\G$ and $\xi_i\in\hat \G_i$. In a similar manner, we define the $i$-th inverse Fourier transform. 
	
	For any bounded measurable function $\mathfrak m \colon \hat \G\to\C$ and any test function $f \colon \G\to\C$ we define the Fourier multiplier operator by 
	\begin{align}
	\label{eq:100}
	T_{\G}[\mathfrak m]f(x) \coloneqq \int_{\hat \G}e(-\langle x , \xi\rangle)\mathfrak m(\xi)\mathcal F_{\G}f(\xi)d\lambda_{\hat \G}(\xi), \qquad  x\in\G.
	\end{align}
	One may think of $f \colon \G\to\C$ as a function on $\G$ that is compactly supported (and smooth if $\G=\R^d$), or as any other function for which \eqref{eq:100} makes sense. We call $\mathfrak m$ a Fourier multiplier on $\hat \G$.
	
	When $\G = {\mathbb R}^d$ and ${\mathfrak m} \colon {\mathbb R}^d \to {\mathbb C}$ is a Fourier multiplier on $\R^d$ supported on $[-1/2,1/2)^d$, we at times abuse notation and identify ${\mathfrak m}$ with its periodization 
	\[
	{\mathfrak m}_{\rm per}(\xi) \coloneqq \sum_{n \in {\mathbb Z}^d} {\mathfrak m}(\xi + n),
	\]
	thereby also viewing ${\mathfrak m}$ as a Fourier multiplier on ${\mathbb T}^d$; see \cite[Section 2]{MSW} for a detailed comparison of $T_{{\mathbb R}^d}[{\mathfrak m}]$ and $T_{{\mathbb Z}^d}[{\mathfrak m}]$.
	
	If $\G=\G_1 \times\cdots \times \G_k$ is a product of LCA groups and $\mathfrak m \colon \hat \G_i \to\C$ is a bounded measurable function acting only on $\hat \G_i$ for some $i \in[k]$, then for any $f\in L^1(\G)\cap L^2(\G)$ it makes sense to consider $i$-th Fourier multiplier operators, which can be defined by
	\begin{align*}
	T_{i, \G}[\mathfrak m]f(x) \coloneqq \int_{\hat \G_i}e(-\langle x_i , \xi_i\rangle)\mathfrak m(\xi_i)\mathcal F_{i, \G}f(x_{(\xi_i)}) d\lambda_{\hat \G_i}(\xi_i), 
	\end{align*}
	for $x=(x_1,\ldots, x_k)\in\G$.
	
	Finally, for any finite set of frequencies $\Sigma\subseteq \hat \G$ we define
	\begin{align}
	\label{eq:106}
	T_{\G}^{\Sigma}[\mathfrak m]f(x) \coloneqq T_{\G}\Big[\sum_{\theta\in\Sigma}\tau_{\theta}\mathfrak m\Big]f(x), \qquad  x\in\G,
	\end{align}
	where $\tau_{\theta}\mathfrak m(\xi) \coloneqq \mathfrak m(\xi-\theta)$ for $\xi\in\hat \G$. If $\G=\G_1 \times\cdots \times \G_k$ is a product of LCA groups, then we similarly define $T_{i, \G}^{\Sigma}[\mathfrak m]f$ for any finite set $\Sigma\subseteq \hat \G_i$ and $i\in[k]$.
	
	\subsection{Variational norms}
	For any $\mathbb I\subseteq \R$, any family $(\mathfrak a_t)_{t\in\mathbb I}=(\mathfrak a_t: t\in\mathbb I)\subseteq \C$, and any exponent $r \in [1, \infty)$, the $r$-variation seminorm is defined to be
	\begin{align}
	\label{eq:45}
	V^{r}( \mathfrak a_t: t\in\mathbb I) \coloneqq 
	\sup_{J\in\Z_+} \sup_{\substack{t_{0}<\dotsb<t_{J}\\ t_{j}\in\mathbb I}}
	\Big(\sum_{j \in [J]}  |\mathfrak a_{t_{j}}-\mathfrak a_{t_{j-1}}|^{r} \Big)^{1/r},
	\end{align}
	where the supremum is taken over all finite increasing sequences in $\mathbb I$, and is set by convention to equal zero if $\# \I \leq 1$. Taking limits as $r \to \infty$ we also adopt the convention
	\begin{align*}
	V^{\infty}( \mathfrak a_t: t\in\mathbb I) \coloneqq \sup_{t,t' \in \I} |\mathfrak a_t - \mathfrak a_{t'}|.
	\end{align*}
	
	The $r$-variation norm for $r \in [1, \infty)$ is defined by
	\begin{equation}
	\label{eq:3}
	{\bf V}^r( \mathfrak a_t: t\in\mathbb I)
	\coloneqq \sup_{t\in\I}|\mathfrak a_t|+
	V^r( \mathfrak a_t: t\in\mathbb I).
	\end{equation}
	Note that the ${\bf V}^r$ norm is nonincreasing in $r$, and comparable to the $\ell^\infty$ norm when $r=\infty$. We also observe the simple triangle inequality
	\begin{equation*}
	{\bf V}^{r}( \mathfrak a_t: t\in\mathbb I)\lesssim {\bf V}^{r}( \mathfrak a_t: t\in\mathbb I_1)+{\bf V}^{r}( \mathfrak a_t: t\in\mathbb I_2)
	\end{equation*}
	whenever $\I = \I_1 \sqcup \I_2$ is an ordered partition of $\I$, which means that $t_1<t_2$ for all $t_1 \in \I_1, t_2 \in \I_2$. If $\mathbb I\subseteq \R$ is countable, then we have
	\begin{equation*}
	{\bf V}^{r}( \mathfrak a_t: t\in\mathbb I)\lesssim \Big(\sum_{t\in\mathbb I}|\mathfrak a_t|^r\Big)^{1/r}.
	\end{equation*}
	One easily establishes the algebra property
	\begin{equation}\label{var-alg}
	{\bf V}^{r}( \mathfrak a_t\mathfrak b_t: t\in\mathbb I)\lesssim
	{\bf V}^{r}( \mathfrak a_t: t\in\mathbb I){\bf V}^{r}( \mathfrak b_t: t\in\mathbb I)
	\end{equation}
	for any scalar sequences $(\mathfrak a_t)_{t \in \I}, (\mathfrak b_t)_{t \in \I}$.
	
	Finally, the $r$-variational spaces will be defined by
	\begin{align*}
	{\bf V}^r(\mathbb I)\coloneqq &\left\{(\mathfrak a_t: t\in\mathbb I)\subseteq\C: {\bf V}^{r}(\mathfrak a_t: t\in\mathbb I) <\infty\right\},
	\end{align*}
	and we will also write $(\mathfrak a_t)_{t\in\mathbb I}\in L^p(X; {\bf V}^r(\mathbb I))$ if and only if
	\[
	\|(\mathfrak a_t)_{t\in\mathbb I}\|_{L^p(X; {\bf V}^r(\mathbb I))}=\|(\mathfrak a_t: t\in\mathbb I)\|_{L^p(X; {\bf V}^r(\mathbb I))}=\left\|{\bf V}^r(\mathfrak a_t: t\in\mathbb I)\right\|_{L^p(X)}<\infty.
	\]
	We will use a similar convention for the spaces ${V}^r(\mathbb I)$.
	
	We also recall the well known Rademacher--Menshov inequality, which asserts that for any  $n_0, m\in\N$ with $n_0< 2^m$ and any sequence of complex numbers $(\mathfrak a_n: n\in\N)$ we have
	\begin{align}
	\label{eq:164}
	V^{2}( \mathfrak a_n: n \in U)
	\leq \sqrt{2}\sum_{i \in \N_{\leq m}} \Big(\sum_{j \in [2^{m-i}]}\Big| \ind{U_{j}^i\subseteq U} \sum_{n \in U_{j}^i} (\mathfrak a_{n+1}-\mathfrak a_{n})\Big|^2\Big)^{1/2},
	\end{align}
	where $U \coloneqq U(n_0,m) \coloneqq [n_0, 2^m) \cap \Z$ and $U_j^i \coloneqq [(j-1)2^i, j2^i) \cap \Z$ for $i, j\in\Z$; see 
	\cite[Lemma~2.5, p.~534]{MSZ2} for the proof. 
	
	\section{Ionescu--Wainger multiplier theorem for canonical fractions} \label{sec:IW} The results of this section will be stated for future reference for arbitrary dimension $d\in\Z_+$, although only $d=1$ will be used throughout this paper.
	
	\subsection{Magyar--Stein--Wainger sampling method}
	An important tool in discrete analogues in harmonic analysis is a quantitative variant of the Shannon sampling theorem from a seminal paper of Magyar, Stein, and Wainger \cite{MSW}. Recalling the definition of $T_{\Z^d}^{\Sigma}[\mathfrak m]$ from \eqref{eq:106}, their sampling result can be subsumed under the following useful proposition.
	\begin{proposition}
		\label{prop:msw}
		For $d\in\Z_+$ there is a constant ${\bf C}_{\rm MSW}(d) \in \R_+$ such that the following holds. Let $p \in [1,\infty]$ and $q\in\Z_+$. If $\mathfrak m \colon \R^d \to L(B_1, B_2)$ is a measurable function supported on $[-1/2,1/2)^d/q$, whose values are bounded linear operators between two finite-dimensional Banach spaces $B_1$ and $B_2$, then
		\[
		\Big\|T_{\Z^d}^{(\frac{1}{q} \Z / \Z)^d}[\mathfrak m]\Big\|_{\ell^{p}(\Z^d;B_1)\to \ell^{p}(\Z^d;B_2)}\le {\bf C}_{\rm MSW}(d) \|T_{\R^d}[\mathfrak m]\|_{L^{p}(\R^d;B_1)\to L^{p}(\R^d;B_2)}.
		\]
	\end{proposition}
	
	We now make a few remarks about Proposition~\ref{prop:msw}.
	
	\begin{enumerate}[label*={\arabic*}.]
		\item The sampling principle developed by Magyar, Stein, and Wainger serves as a ``bridge'' between discrete analogues in harmonic analysis and classical harmonic analysis. It enables the deduction of discrete problems from their continuous counterparts.
		
		\item The proof of Proposition~\ref{prop:msw} can be found in \cite[Corollary~2.1, p.~196]{MSW}. We also refer to \cite{MSZ1} for a generalization of Proposition~\ref{prop:msw} to real interpolation spaces.
		
		\item All Banach spaces that arise in our problems are infinite-dimensional and separable. Despite Proposition~\ref{prop:msw} being formulated for finite-dimensional Banach spaces, it can still be applied in our situation using standard approximation arguments that reduce infinite-dimensional settings to finite-dimensional settings.
		
		\item We emphasize that $B_1$ and $B_2$ are general (finite-dimensional) Banach spaces and thus, in view of the previous remark, Proposition~\ref{prop:msw} includes maximal functions and can also accommodate $r$-variation seminorms as well as many other seminorms.
		
		\item A fundamental question in discrete harmonic analysis is whether it is possible to replace the set $(\frac{1}{q} \Z / \Z)^d$ appearing in Proposition~\ref{prop:msw} with other sets of fractions $\Sigma$, while still maintaining mapping properties of $T_{\Z^d}^{\Sigma}[\mathfrak m]$ that depend in a good way on the size of $\Sigma$, under a suitable support condition for $\mathfrak m$. This was achieved for the first time in a groundbreaking paper of Ionescu and Wainger \cite{IW} for a certain family of reduced fractions $\Sigma$ whose denominators have special factorization properties.
	\end{enumerate}
	
	In this paper, our aim is to prove a variant of the Ionescu--Wainger multiplier theorem for the set of canonical fractions that arises from the classical circle method. This will be critical in our proof of the multilinear Weyl inequality; see Theorem~\ref{weyl}.
	
	\subsection{A refined variant of the Ionescu--Wainger multiplier theorem}
	
	For $d\in\Z_+$ and $N\ge1$ define $1$-periodic sets of the so-called \textit{canonical fractions} by
	\begin{align}
	\label{IWeq:372}
	\mathcal R_{\le N}^d  \coloneqq  \left\{ \frac{a}{q} \in(\Q / \Z)^d:  q \in [N] \text{ and } \gcd(a, q)=1\right\},
	\end{align}
	where $\gcd(a, q) \coloneqq \gcd(a_1, \dots, a_d, q)$ denotes the greatest common divisor of $a_1, \dots, a_d, q$. To be precise, the elements of ${\mathcal R}^d_{\le N}$ are equivalence classes $\theta = \frac{a}{q} + {\Z}^d$ in $\T^d = (\R/\Z)^d$, with a well-defined denominator $q$ such that every representative $\frac{b}{q} \in \theta$ satisfies $\gcd(b,q) = 1$. We often identify $\frac{a}{q} \in \Q^d$ with the equivalence class $\frac{a}{q} + \Z^d$.
	
	Using the definition of $T_{\Z^d}^{\Sigma}[\mathfrak m]$ from \eqref{eq:106} we now formulate the main result of this section. Here, $p_0' \coloneqq p_0 / (p_0-1)$ denotes the conjugate exponent of $p_0$.
	
	\begin{theorem} \label{thm:IW}
		Let $p_0' \leq p \leq p_0$ for some $p_0 \in 2\Z_+$. Then there exists a constant ${\bf C}_{\rm IW}(p_0) \in \R_+$ such that for every $d \in \Z_+$ and every $N \in \N_{\ge 100}$ the following is true. Assume that 
		\[
		0 < \vartheta \leq (2 p_0 N^{p_0})^{-1}
		\] 
		and let $\mathfrak m \colon \R^{d} \to L(H_1,H_2)$ be a measurable function supported on $[-\vartheta, \vartheta]^d$, whose values are bounded linear operators between two separable Hilbert spaces $H_{1}$ and $H_{2}$. Denote $\rho_N \coloneqq \rho_N(p_0) \coloneqq {\bf C}_{\rm IW}(p_0) \frac{\log \log \log N}{\log \log N}$. Then
		\begin{align}
		\label{IWeq:376}
		\Big\|T_{\Z^d}^{\mathcal R_{\le N}^d}[\mathfrak m]\Big\|_{\ell^p(\Z^d;H_1)\to \ell^p(\Z^d;H_2)}
		\lesssim_{d,p_0} 
		N^{\rho_N}
		\|T_{\R^d}[\mathfrak m]\|_{L^{p_0}(\R^d;H_1)\to L^{p_0}(\R^d;H_2)}.
		\end{align}	
	\end{theorem}
	
	A few remarks about Theorem~\ref{thm:IW}.
	
	\begin{enumerate}[label*={\arabic*}.]
		
		\item Theorem~\ref{thm:IW} gives an affirmative answer to a question of Ionescu and Wainger from \cite[Remark~3, below Theorem~1.5, p.~361]{IW}, which asks about \eqref{IWeq:376} in the scalar setting $H_1=H_2=\mathbb C$  with $C_{d,p_0,\rho} N^{\rho}$ for arbitrary $\rho \in \R_+$ in place of $N^{\rho_N}$. Obviously, \eqref{IWeq:376} yields the bound postulated in \cite{IW},  since $\lim_{N\to \infty} \rho_N = 0$. 
		
		\item Although Proposition~\ref{prop:msw} plays an important role in the proof of Theorem~\ref{thm:IW}, the latter result is only applicable to multipliers defined in the Hilbert space setting. However, in the next section, we will show that Theorem~\ref{thm:IW} can be applied to some Banach space settings. This can be done for a broad range of Banach spaces that, in some manner, exhibit similarities to Hilbert spaces.
		
		\item The advantage of using Theorem~\ref{thm:IW} is that square function estimates in the discrete setting can be directly deduced from the corresponding continuous square function estimates, which, as we will soon see, is very useful in many applications.
		
		\item The norm on the right hand side of \eqref{IWeq:376}, unlike the support hypothesis, is scale-invariant. To be more precise, if $\mathfrak m$ is replaced with $\mathfrak m(A \, \cdot \,)$ for any invertible linear transformation $A \colon \R^d\to \R^d$, then the norm in \eqref{IWeq:376} remains unchanged.
		
		\item It would be interesting to know whether inequality \eqref{IWeq:376} holds with an absolute constant independent of $N$. Although in most applications inequality \eqref{IWeq:376} is sufficient, having \eqref{IWeq:376} with a constant independent of $N$ would simplify some arguments.
		
		\item In the proof of Theorem~\ref{thm:IW} we will use variants of two properties introduced in \cite{IW}: a uniqueness property and the so-called $\mathcal O$ property. These two properties efficiently detect strong orthogonalities between the Ionescu--Wainger fractions, which are rationals whose denominators have certain factorization properties. The novelty of our approach lies in using the set of canonical fractions \eqref{IWeq:372} instead of the Ionescu--Wainger fractions \cite{IW}, which are quite complex.
	\end{enumerate}
	
	\subsubsection{\textbf{Proof of Theorem~\ref{thm:IW}, simple reductions}} The case $p_0=2$ holds by Plancherel's theorem. By duality and interpolation we can assume that $p = p_0 = 2r$ for $r \in \N_{\geq 2}$. Replacing $\mathfrak m$ by $\tilde{\mathfrak m} \colon \R^d \to L(H_1\oplus H_2, H_1\oplus H_2)$, where $\tilde{\mathfrak m}(x)(h_1,h_2) \coloneqq (0,\mathfrak m(x)(h_1))$ for $x \in \R^d$, we can take $H_1=H_2=H$. We normalize $\|T_{\R^d}[\mathfrak m]\|_{L^{2r}(\R^d; H) \to L^{2r}(\R^d; H)} = 1$. Finally, we can assume that $N \geq N_0$ for a large integer $N_0 \in \N_{\geq 10^{10}}$, since $\mathcal R^d_{\le N}$ produces at most $N^{d+1}$ operators, each having a norm  bounded by ${\bf C}_{\rm MSW}(d)$ in view of Proposition~\ref{prop:msw}. 
	
	\subsubsection{\textbf{Uniqueness property and property $\mathcal O_N$}}
	We follow \cite{IW} in using the uniqueness property and property $\mathcal O$, two important concepts introduced therein for detecting orthogonality.
	
	\begin{definition}[Uniqueness property]
		\label{up}
		A sequence $(v_1,\ldots, v_n)$ has the \emph{uniqueness property} if there is $i\in[n]$ such that $v_i\neq v_j$ for all $j\in[n]\setminus\{i\}$. In other words, the element $v_i$ occurs in the sequence $(v_1,\ldots, v_n)$ exactly once. This definition will primarily be used for sequences of rational numbers, but it is also applicable to more general sequences.
	\end{definition}
	
	Using Definition~\ref{up} we can formulate an abstract orthogonality principle.
	
	\begin{proposition}[Abstract orthogonality principle]
		\label{UP}
		Given $r \in \Z_+$, there exists a constant ${\bf C}_{\rm UP}(r) \in \R_+$ such that the following is true. Let $(X,\mathcal B(X), \mu)$ be a measure space, let $H$ be a separable Hilbert space, and let $V_{1},\dots,V_{r} \subset \Z_+$ be finite sets, not necessarily disjoint. For every pair $(i,v)$ with $v \in V_{i}$ let $F_{i,v} \in L^{2r}(X;H)$ and suppose that for each $2r$-tuple
		\begin{equation*}
		(v_{1,1}, v_{1,2}, \ldots, v_{r,1}, v_{r,2})
		\in V_{1}^{2} \times \dotsm \times V_{r}^{2}
		\end{equation*}
		with the uniqueness property the following condition holds
		\begin{align}
		\label{IWeq:6}
		\int_{X} \prod_{i \in [r]} \langle F_{i,v_{i,1}}(x), F_{i,v_{i,2}}(x) \rangle_{H} \, d\mu(x) = 0. 
		\end{align}
		Then
		\begin{equation*}
		\int_{X} \prod_{i \in [r]} \Big \| \sum_{v \in V_{i}} F_{i,v}(x) \Big\|_{H}^{2} \, d\mu(x)
		\leq {\bf C}_{\rm UP}(r)
		\int_{X} \prod_{i \in [r]} \Big( \sum_{v \in V_{i}} \| F_{i,v}(x)\|_{H}^{2} \Big) \, d\mu(x).
		\end{equation*} 
	\end{proposition}
	
	\begin{proof}
		We refer to \cite[Corollary~2.24]{MSZ3} for a detailed proof; see also \cite[Lemma~2.2, p.~363]{IW} for a scalar version of this principle.
	\end{proof} 
	
	Our aim will be to split the set of canonical fractions $\mathcal R^d_{\le N}$ into a controlled number of pieces for which the uniqueness property implies condition \eqref{IWeq:6}.
	
	Here our approach deviates from \cite{IW} and \cite{MSZ3}. In the latter two papers an essential role was played by the Ionescu--Wainger rationals whose denominators were factored into two parts, one being a relatively large product of powers of small prime divisors, and the other consisting of a product of a small number of powers of large prime divisors.
	
	In contrast to \cite{IW} and \cite{MSZ3}, instead of working with the Ionescu--Wainger rationals, we shall work directly with the set of canonical fractions $\mathcal R^d_{\le N}$ from \eqref{IWeq:372}. We introduce a new concept of $N$-lifted composites, which will allow us to effectively split $\mathcal R^d_{\le N}$. 
	
	\begin{definition}[$N$-lifted composites]
		\label{N-lifted}
		Let $\mathbb P$ be the set of prime numbers. For $N\in\Z_{+}$ an \emph{$N$-lifted composite} $Q_q$ of an integer $q\in [N]$ is defined by setting 
		\begin{align*}
		Q_q \coloneqq 
		\prod_{p \in \mathbb{P} \, : \, q \in p \Z}
		p^{\lfloor \log_p N \rfloor}
		\end{align*}
		with the convention $Q_1 \coloneqq 1$. In other words, $q$ and $Q_q$ are divisible by exactly the same prime numbers, while $Q_q$ is a product of integers of the form $p^{\lfloor \log_p N \rfloor}$ for some $p \in \mathbb{P}  \cap \N_{\leq N}$. Obviously, $q$ divides $Q_q$ for each $q\in [N]$.
	\end{definition}
	
	We need to adapt the concept of property $\mathcal O$ from \cite{IW} to fit our needs. 
	
	\begin{definition}[Property $\mathcal O_N$]
		\label{propO}
		For $N \in \N_{\geq N_0}$ we shall say that $\Lambda \subset \Z_+$ has \emph{property $\mathcal O_N$} if $\Lambda \subseteq S_{1}\dotsm S_{k} \coloneqq \{s_{1}\dotsm s_{k} : s_1\in S_{1},\ldots, s_k\in S_{k}\}$ for some  $k\in [ \frac{5 \log N}{\log \log N}]$, where $S_1, \dots, S_k$ are disjoint subsets of 
		$
		\{ p^{\lfloor \log_p N \rfloor} : p \in \mathbb P \cap \N_{\leq N} \}.
		$
	\end{definition}
	
	The following lemma explains the restriction $k\in \big[\frac{5 \log N}{\log \log N}\big]$. 
	
	\begin{lemma} \label{L1}
		If $N_0$ is sufficiently large, then for every $N \in \N_{\geq N_0}$ and distinct $p_1, \dots, p_{k_N} \in \mathbb P$ with $k_N \coloneqq \lfloor \frac{5 \log N}{\log \log N} \rfloor$ we have $p_1  \cdots p_{k_N} > N$.
	\end{lemma}
	
	\begin{proof}
		For $n \in \Z_+$, its primorial $n \#$ is the product of all $p \in \mathbb{P}  \cap \N_{\leq n}$. By the prime number theorem, we have $\lim_{n \to \infty} \sqrt[n]{n \#} = e$. Take $N \in \Z_+$ such that $n_N \coloneqq \lfloor \log N \rfloor + 1$ is large. Then $n_N \# > 2^{n_N} > N$ and $\#(\mathbb{P}  \cap \N_{\leq n_N}) \leq \lfloor \frac{5 \log N}{\log \log N} \rfloor = k_N$. The latter implies $p_1 \cdots p_{k_N} \geq n_N \#$, as desired. 
	\end{proof}
	
	\subsubsection{\textbf{Partitioning of the set of canonical fractions}}
	For $N \in \N_{\geq N_0}$ and $\frac{a}{q} \in \mathcal R^d_{\leq N}$ the corresponding composite $Q_q$ has at most $\lfloor \frac{5 \log N}{\log \log N} \rfloor$ factors $p^{\lfloor \log_p N \rfloor}$, while each $p$ is one of the first $\lfloor \frac{2 N}{\log N} \rfloor$ primes. This motivates the following counting lemma originating in \cite{M1}.
	
	\begin{lemma} \label{L2} 
		There exists a constant ${\bf C}_{\rm sur} \in \Z_+$ such that the following holds. For $N \in \N_{\geq N_0}$ let $V$ be a set of size $l \in [ \frac{2 N}{\log N} ]$ and let $k\in [ \frac{5 \log N}{\log \log N}]$ satisfy $k \leq l$. Then, denoting $J \coloneqq \lfloor {\bf C}_{\rm sur} 2^{{\bf C}_{\rm sur} \log N / \log \log N} \rfloor$, we can find $J$ surjective (not necessarily distinct) functions $ g_{1},\dotsc,g_{J} \colon V \to [k] $ such that $\max \{\#g_{j}(E) : j \in [J]\}  = k$ for every $E \subseteq V$ with $\#E = k$.
	\end{lemma}
	
	\begin{proof}
		We refer to \cite[Lemma~2.16]{MSZ3}.
	\end{proof}
	
	Proceeding in two steps, we split $\mathcal R^d_{\le N}$ into at most $N^{\frac{C}{\log\log N}}$ pieces, each satisfying property $\mathcal O_N$.
	
	\paragraph{\bf Step~1}
	For a large integer $N\in\Z_+$, by Lemma~\ref{L1}, we can split $\mathcal R^d_{\le N}$ into disjoint subsets $ \mathcal R^d_{N,0} \coloneqq \mathcal R^d_{\le 1}$ and $\mathcal R^d_{N,k}$ for $k\in [\frac{5 \log N}{\log \log N}]$, consisting of fractions $\frac{a}{q} \in \mathcal R^d_{\le N}$ whose denominators $q$ have exactly $k$ prime divisors. We may assume, without loss of generality, that $k>0$, since the case $k=0$ may be readily handled by Proposition~\ref{prop:msw}. We also assume that $\mathcal R^d_{N,k}$ is nonempty. 
	
	\paragraph{\bf Step~2} Next, we fix $k\in [\frac{5 \log N}{\log \log N}]$ as above and apply Lemma~\ref{L2} with 
	\[
	V \coloneqq \left\{ p^{\lfloor \log_p N \rfloor} : p \in \mathbb P \cap \N_{\leq N} \right\} = \left\{ Q_{p} : p \in \mathbb P \cap \N_{\leq N} \right\},
	\]
	obtaining surjections $g_1,\ldots, g_J \colon V\to [k]$  with $J \coloneqq \lfloor {\bf C}_{\rm sur} 2^{{\bf C}_{\rm sur} \log N / \log \log N} \rfloor$. Given $j\in[J]$, we define $S_i^j \coloneqq g_j[\{i\}]^{-1}$ for $i \in [k]$ and, consequently, we set
	\begin{align}
	\label{IWeq:8}
	\mathcal R^d_{N,k,j} \coloneqq \bigg\{\frac{a}{q}\in \mathcal R^d_{N,k}: Q_q\in S_1^j\cdots S_k^j\setminus\bigcup_{i\in [j-1]}S_1^i\cdots S_k^i\bigg\}, 
	\end{align}
	where if $j=1$, then the union above is understood to be empty. For each $j \in [J]$ the set $\{ Q_q : \frac{a}{q} \in \mathcal R^d_{N,k,j}\} \subset \Z_+$ has property $\mathcal O_N$. Moreover, $\mathcal R^d_{N,k}=\bigcup_{j\in[J]}\mathcal R^d_{N,k,j}$, which immediately follows from Lemma~\ref{L2}.
	
	Our task is to prove \eqref{IWeq:376} with any nonempty $\mathcal R^d_{N,k,j}$ in place of $\mathcal R^d_{\le N}$. To facilitate matters, we introduce some useful notation and terminology.
	
	\subsubsection{\textbf{\texorpdfstring{$\mathbb P$}{TEXT}-irreducibility and factorization}}
	
	We now present a new concept of \emph{$\mathbb P$-irreducible fractions}, which will be crucial in verifying condition \eqref{IWeq:6}. 
	
	\begin{definition}[Quotients, reduced fractions and $\mathbb P$-irreducible fractions]
		Fix $d, q\in \Z_+$. Identifying $\frac{a}{q}$ with $\frac{a}{q} + \Z^d$, we define the following sets.
		\begin{enumerate}[label*={\arabic*}.]
			\setlength\itemsep{0.25em}
			\item The set of \emph{quotients} or \emph{$q$-quotients} is defined by $\mathcal Q(q) \coloneqq \frac{1}{q} \N_{<q}^d$.
			\item The set of \emph{reduced fractions} is defined by $\mathcal R(1) \coloneqq \mathcal Q(1)$ or 
			\[
			{\mathcal R}(q) \coloneqq \left\{ \frac{a}{q} \in \mathcal Q(q) : \gcd(a,q) = 1 \right\} \quad \text{for} \quad q>1.
			\]
			\item The set of \emph{$\mathbb P$-irreducible fractions} is defined by $\mathcal I(1) \coloneqq \mathcal Q(1)$ or
			\[
			{\mathcal I}(q) \coloneqq \left\{ \frac{a}{q} \in \mathcal Q(q) : \frac{q}{\gcd(a, q)} \, \text{has the same prime divisors as } q \right\} \quad \text{for} \quad q>1.
			\]
		\end{enumerate}	
	\end{definition}
	
	It is easy to see that $\mathcal I(q)\supseteq \mathcal R(q)$ and that $\mathcal I(q)$ is a disjoint union of the sets $\mathcal R(\tilde q)$, where $\tilde q$ runs over the set of  all positive divisors of $q$ having exactly the same prime divisors as $q$.
	
	\begin{example}
		For $d=1$ and $q=12$, identifying $\frac{a}{q}$ with $\frac{a}{q} + \Z$, we obtain
		\begin{align*}
		\quad 
		\mathcal R(12)  = \frac{1}{12} \{1,5,7,11\} \quad {\rm and} \quad
		\mathcal I(12)  = \frac{1}{12} \{1,2,5,7,10,11\}.
		\end{align*}
	\end{example}
	
	For $\Lambda \subset \Z_+$ define 
	\[
	\mathcal R(\Lambda) \coloneqq \bigcup_{q \in \Lambda} \mathcal R(q)
	\quad \text{and}
	\quad 
	\mathcal I(\Lambda) \coloneqq \bigcup_{q \in \Lambda} \mathcal I(q).
	\]
	Obviously, we have $\mathcal R(\Lambda)\subseteq \mathcal I(\Lambda)$ and the set of canonical fractions \eqref{IWeq:372} can be written as $\mathcal R_{\le N}^d=\mathcal R([N])$ by using this definition. Moreover, an important feature of $\mathcal R(\Lambda)$ used in \cite{IW} and \cite{MSZ3} is the factorization property asserting that if $\gcd(q, q') = 1$ for all $(q, q') \in \Lambda \times \Lambda'$, then
	\[
	\mathcal R(\Lambda \Lambda') = \mathcal R(\Lambda) \oplus \mathcal R(\Lambda').
	\]
	The direct sum $\oplus$ indicates that each element in the sumset ${\mathcal R}(\Lambda) + {\mathcal R}(\Lambda')$ has a unique representation $\frac{a}{q} + \frac{a'}{q'}$ as an element of $\T^d$. The identity above is a simple consequence of the Chinese remainder theorem. We now strengthen this factorization property to the $\mathbb P$-irreducible fractions $\mathcal I(\Lambda)$.
	
	\begin{definition}[$\mathbb P$-separated sets]
		\label{def:p-sep}
		A set $\Lambda\subseteq\Z_+$ is \emph{$\mathbb P$-separated} if for any distinct elements $q_1, q_2 \in \Lambda$ there is $p \in \mathbb P$ that divides exactly one of them.
	\end{definition}
	
	We now prove the factorization property for $\mathcal I(\Lambda)$.
	\begin{lemma}
		\label{lem:fact}
		Assume that $\Lambda, \Lambda'\subseteq\Z_+$ are $\mathbb P$-separated and $\gcd(q, q') = 1$ for all $(q, q') \in \Lambda \times \Lambda'$. Then the map $\Phi \colon \mathcal I(\Lambda) \times \mathcal I(\Lambda')\to \mathcal I(\Lambda \Lambda')$ given by
		\[
		\Phi\Big(\frac{a}{q}, \frac{a'}{q'}\Big) \coloneqq \frac{a}{q}+ \frac{a'}{q'} 
		\quad \text{for} \quad
		\Big( \frac{a}{q}, \frac{a'}{q'} \Big) \in \mathcal I(\Lambda) \times \mathcal I(\Lambda')
		\]
		is a bijection. In particular, we have $
		\mathcal I(\Lambda \Lambda') = \mathcal I(\Lambda) \oplus \mathcal I(\Lambda').$ 
	\end{lemma}
	
	\begin{proof}
		By the Chinese remainder theorem, it follows that $\Phi$ is surjective, and from Definition~\ref{def:p-sep} we deduce that $\Phi$ is injective. The second part of the lemma now follows easily from the first part.
	\end{proof}
	
	To establish \eqref{IWeq:376} in Theorem~\ref{thm:IW}, it suffices to prove
	\begin{align}\label{IWeq:k-j}
	\Big\|T_{\Z^d}^{{\mathcal R}^d_{N,k,j}}[\mathfrak m] f\Big\|_{\ell^{2r}(\Z^d;H)}
	\lesssim_{d,p_0} N^{{\bf C}_{\rm IW}(p_0) \frac{\log \log \log N}{\log \log N}}
	\|f\|_{\ell^{2r}(\Z^d;H)}
	\end{align}
	since we have normalized $\|T_{\R^d}[\mathfrak m]\|_{L^{2r}(\R^d;H)\to L^{2r}(\R^d;H)} =1$. Indeed summing \eqref{IWeq:k-j} over $k\in [ \frac{5 \log N}{\log \log N}]$ and $j\in [J]$ incurs an acceptable, additional factor
	\[
	\frac{5 \log N}{\log \log N} \cdot {\bf C}_{\rm sur} 2^{{\bf C}_{\rm sur} \log N/\log \log N}.
	\]
	
	We fix $k$ and $j$. The multiplier corresponding to  $T_{\Z^d}^{ \mathcal R^d_{N,k,j}}[\mathfrak m]$ is given by
	\begin{equation*}
	\Delta_{N,k,j}(\xi)
	\coloneqq 
	\sum_{q \in S_1 \cdots S_k}
	\sum_{\frac{a}{q} \in\mathcal I(q)}
	\ind{N,k,j}\Big(\frac{a}{q}\Big) \,
	\mathfrak m\Big(\xi - \frac{a}{q}\Big),
	\end{equation*}
	where, as in \eqref{IWeq:8}, the sets $S_1 \coloneqq S_1^j \coloneqq g_j[\{1\}]^{-1},\ldots, S_k \coloneqq S_k^j \coloneqq g_j[\{k\}]^{-1}$ form a partition of $V$ determined by the surjection $g_j$, and $\ind{N,k,j}(\frac{a}{q}) \coloneqq 1$ if $\frac{a}{q}$ is $\mathbb P$-irreducible and upon reduction belongs to $\mathcal R^d_{N,k,j}$, and $\ind{N,k,j}(\frac{a}{q}) \coloneqq 0$ otherwise. We note that the sets of ${\mathbb P}$-irreducible fractions ${\mathcal I}(q)$ are pairwise disjoint as $q$ varies over the product set $S_1 \cdots S_k$. 
	
	\subsubsection{\textbf{Denominators}} Now we exploit orthogonality between denominators. In the argument below we will use Lemma~\ref{lem:fact} repeatedly. Note that
	\begin{align*}
	\Big\| T_{\Z^d}^{{\mathcal R}^d_{N,k,j}}[\mathfrak m] f \Big\|^{2r}_{\ell^{2r}(\Z^d;H)}
	=
	\sum_{x \in \Z^d} 
	\Big \| 
	\sum_{q \in S_{1} \cdots S_{k}}
	\sum_{u \in\mathcal I(q)} 
	f_{u}(x) \Big\|_{H}^{2r},
	\end{align*}
	where
	$f_u \coloneqq \mathcal F^{-1}_{\Z^d} ( \ind{N,k,j}(u) \, \mathfrak m(\, \cdot \, - u) \mathcal F_{\Z^d} f )$. 
	Let ${\mathcal S}_{\emptyset} \coloneqq \{1\}$ and ${\mathcal S}_L \coloneqq \prod_{l \in L} S_l$ for every nonempty $L \subseteq [k]$. We observe that
	\[
	\sum_{x \in \Z^d} 
	\Big \| 
	\sum_{q \in S_{1} \cdots S_{k}}
	\sum_{u \in\mathcal I(q)} 
	f_{u}(x) \Big\|_{H}^{2r}
	= 
	\sum_{x \in \Z^d} \Big( \sum_{q \in {\mathcal S}_{\emptyset}} \Big\| F_{q,\emptyset}(x) \Big\|_H^2 \Big)^r,
	\]
	where $F_{q, L} \coloneqq \sum_{q' \in {\mathcal S}_{[k] \setminus L}} \sum_{u \in \mathcal I(qq')} f_u$ for every $q \in {\mathcal S}_L$. Now, proceeding as in \cite[Section~2.4]{MSZ3}, we can show for each $l \in [k]$ that
	\begin{align} \label{denom}
	\sum_{x \in \Z^d} \Big(
	\sum_{q \in {\mathcal S}_{[l-1]}} 
	\Big\| F_{q, [l-1]}(x) \Big\|_{H}^{2} \Big)^r
	\leq {\bf C}_{\rm UP}(r)
	\sum_{x \in \Z^d} \Big(
	\sum_{q \in {\mathcal S}_{[l]}} 
	\Big\| F_{q, [l]}(x) \Big\|_{H}^{2} \Big)^r
	\end{align}
	with ${\bf C}_{\rm UP}(r)$ from Proposition~\ref{UP}. After $k$ iterations of \eqref{denom}, we obtain 
	\begin{align*}
	\sum_{x \in \Z^d} \Big( 
	\sum_{q \in {\mathcal S}_{\emptyset}}
	\Big\| F_{q,\emptyset}(x) \Big\|_{H}^{2} \Big)^r
	\leq {\bf C}_{\rm UP}(r)^{k}
	\sum_{x \in \Z^d} \Big( 
	\sum_{q \in {\mathcal S}_{[k]}}
	\Big\| F_{q, [k]}(x) \Big\|_{H}^{2} \Big)^r
	\end{align*}
	or, equivalently,
	\begin{align} \label{denom-II}
	\Big\| T_{\Z^d}^{{\mathcal R}^d_{N,k,j}}[\mathfrak m] f\Big\|_{\ell^{2r}(\Z^{d};H)}^{2r}
	\le {\bf C}_{\rm UP}(r)^{k} \sum_{x \in \Z^d} 
	\Big(\sum_{q \in S_{1} \cdots S_{k}}
	\Big \| \sum_{u \in\mathcal I(q)} 
	f_{u}(x) \Big\|_{H}^{2}\Big)^r.
	\end{align}
	The constant ${\bf C}_{\rm UP}(r)^{k/2r} \leq C_{r} 2^{C_r \log N / \log \log N}$ is acceptable, since $k \le \frac{5 \log N}{\log \log N}$. 
	
	\subsubsection{\textbf{Numerators}}
	Now we exploit orthogonalities between numerators to bound the right-hand side of \eqref{denom-II}. In the argument below, as before, we will use Lemma~\ref{lem:fact} repeatedly. The right-hand side of \eqref{denom-II} is equal to
	\begin{align*}
	\sum_{x \in \Z^d} 
	\Big(\sum_{q \in S_{1} \cdots S_{k}}
	\Big \| \sum_{u \in\mathcal I(q)} 
	f_{u}(x) \Big\|_{H}^{2}\Big)^r
	= \sum_{x \in \Z^d} \Big( 
	\sum_{q \in {\mathcal S}_{[k]}}
	\Big\| F_{q, [k]}(x) \Big\|_{H}^{2} \Big)^r.
	\end{align*}
	Set $F_{q,q',q''} \coloneqq \sum_{u'' \in \mathcal I(q'')} \big \| \sum_{u' \in\mathcal I(q')} \sum_{u \in\mathcal I(q)} f_{u+u'+u''} \big \|_{H}^{2}$ and observe that
	\[
	\sum_{x \in \Z^d} \Big( 
	\sum_{q \in {\mathcal S}_{[k]}}
	\Big\| F_{q,[k]}(x) \Big\|_{H}^{2} \Big)^r
	=
	\sum_{x \in \Z^d} 
	\sum_{q \in {\mathcal S}_{\emptyset}}
	\Big( \sum_{q' \in {\mathcal S}_{[k]}}
	\sum_{q'' \in {\mathcal S}_{\emptyset}}
	F_{q,q',q''}(x) \Big)^r.
	\]
	Denote ${\mathcal F}_{L,L'} \coloneqq \sum_{q \in {\mathcal S}_L} \big( \sum_{q' \in {\mathcal S}_{L' \setminus L}} \sum_{q'' \in {\mathcal S}_{[k] \setminus L'}} F_{q,q',q''} \big)^r$ for $L \subseteq L' \subseteq [k]$. We can rewrite the right-hand side above as $\| {\mathcal F}_{\emptyset, [k]} \|_{\ell^1(\Z^d)}$. Thus, we obtain
	\begin{align}\label{denom-III}
	\sum_{x \in \Z^d} 
	\Big(\sum_{q \in S_{1} \cdots S_{k}}
	\Big \| \sum_{u \in\mathcal I(q)} 
	f_{u}(x) \Big\|_{H}^{2}\Big)^r
	= \| {\mathcal F}_{\emptyset, [k]} \|_{\ell^1(\Z^d)}.
	\end{align}
	Assuming $L \subsetneq L' \subseteq [k]$ and writing $L_+ = L \cup \{l\}$ and $L'_- = L' \setminus \{l\}$ with $l = \min L' \setminus L$, we proceed as in \cite[Section~2.5]{MSZ3} to show that
	\begin{align} \label{numer}
	\| {\mathcal F}_{L, L'} \|_{\ell^1(\Z^d)} \leq
	C_r \big( \| {\mathcal F}_{L_+, L'} \|_{\ell^1(\Z^d)} + \| {\mathcal F}_{L, L'_-}\|_{\ell^1(\Z^d)} \big)
	\end{align}
	with $C_r \simeq_r 1$ by using Proposition~\ref{UP}. After $k$ iterations of \eqref{numer}, we obtain
	\begin{align}
	\label{IWeq:9}
	\| {\mathcal F}_{\emptyset, [k]} \|_{\ell^1(\Z^d)} 
	\leq C_r^k
	\sum_{L \subseteq [k]}
	\| {\mathcal F}_{L,L} \|_{\ell^1(\Z^d)}.
	\end{align}
	The implied constant in \eqref{IWeq:9} is acceptable, since $k\in [\frac{5 \log N}{\log \log N}]$.
	
	\subsubsection{\textbf{Square function estimates}}
	Gathering the estimates and identities \eqref{denom-II}, \eqref{denom-III}, and \eqref{IWeq:9}, we reduce matters to proving
	\begin{align*}
	\sum_{L \subseteq [k]}
	\| {\mathcal F}_{L,L} \|_{\ell^1(\Z^d)} \lesssim_{r,d}  N^{  \frac{2r C \log \log \log N}{\log \log N}}
	\|f\|_{\ell^{2r}(\Z^{d};H)}^{2r}.
	\end{align*}
	Since there are $2^k$ subsets $L \subseteq [k]$, we may treat each summand above separately at the expense of introducing an acceptable factor $C_r 2^{C_r \log N / \log \log N}$.
	
	We fix $L$. Set $\mathcal I' \coloneqq \mathcal I(\mathcal S_{[k] \setminus L}) \cap \mathcal R^d_{\le N}$. Since $\ind{N,k,j}(u+u') = 0$ for every $u \in \mathcal I(\mathcal S_L)$ and $u' \in \mathcal I(\mathcal S_{[k] \setminus L}) \setminus \mathcal R^d_{\le N}$, we are reduced to proving the inequality
	\begin{align}
	\label{IWeq:10}
	\sum_{x \in \Z^d} 
	\sum_{q \in \mathcal S_L}
	\Big( 
	\sum_{u' \in \mathcal I'}
	\Big \| 
	\sum_{u \in\mathcal I(q)} 
	f_{u+u'}(x) \Big\|_{H}^{2} \Big)^r
	\leq C_{r,d}^{2r} N^{\frac{2r C \log \log \log N}{\log \log N}}
	\|f\|_{\ell^{2r}(\Z^{d};H)}^{2r}
	\end{align}
	that will be interpreted as a linear operator norm bound from $\ell^{2r}(\Z^{d};H)$ to $\ell^{2r}(\Z^{d} \times \mathcal  S_{L};\ell^{2}( \mathcal I' ;H))$. The reduction to \eqref{IWeq:10} uses the fact that the sets $\mathcal I(q) \cap \mathcal I(q')$ are empty for any distinct $q, q' \in \mathcal S_{[k]\setminus L}$.
	
	The following lemma is key to establishing \eqref{IWeq:10}. This is the only place where the expression $\frac{\log \log \log N}{\log \log N}$ appears in the exponent. It would be interesting to determine whether $\frac{1}{\log \log N}$ could be used instead.
	
	\begin{lemma} \label{L4}
		There exists a constant $\mathbf C_{\rm div} \in \R_+$ such that for all $q \in \Z_+$ and $N \in \N_{\geq N_0}$, if $q$ has $k \in [\frac{5 \log N}{\log \log N}]$ prime divisors, then $q$ has at most $\lfloor \mathbf C_{\rm div} N^{\mathbf C_{\rm div} \frac{\log \log \log N}{\log \log N}} \rfloor $ positive divisors not greater than $N$. 
	\end{lemma}
	
	\begin{proof}
		Let $\{p_1,\ldots, p_k\}\subset \mathbb P$ be the set of prime divisors of $q$. Our aim will be to count of all possible tuples $(\alpha_1,\ldots, \alpha_k) \in \N^k$ such that $p_1^{\alpha_1}\cdots p_k^{\alpha_k}\le N$. We may assume that $\alpha_1 + \dots + \alpha_k\le \lfloor 2 \log N \rfloor$, since otherwise $p_1^{\alpha_1}\cdots p_k^{\alpha_k} \geq 2^{\lfloor 2 \log N \rfloor + 1} \ge 4^{\log N}> N$. Thus, we have $\max\{ \alpha_1,\ldots, \alpha_k \} \le \lfloor 2 \log N \rfloor$. Note that $\N_{\le 2 \log N}$ can be split into a finite number of intervals of the form
		\[
		I_i \coloneqq \N_{\le 2 \log N}\cap [(i-1) \log \log N, i\log \log N),
		\] 
		each having about $\log \log N$ elements except the last one, which may be shorter. Now to each tuple $(\alpha_1,\ldots, \alpha_k)$ we assign a unique tuple of integers $(i_1, \dots, i_k)$ such that $(\alpha_1,\ldots, \alpha_k)\in I_{i_1}\times\cdots\times I_{i_k}$. If $(\alpha_1,\ldots, \alpha_k)\in I_{i_1}\times\cdots\times I_{i_k}$, then
		\[
		(i_1+ \cdots +i_k - k) \frac{\log \log N}{2} \le \alpha_1+\cdots+\alpha_k\le \lfloor 2 \log N \rfloor,
		\]
		which implies that $i_1+ \cdots +i_k \leq k + \frac{4 \log N}{\log \log N} \leq \lfloor \frac{10 \log N}{\log \log N} \rfloor$. Taking
		\[
		\mathcal A \coloneqq \left\{(\alpha_1,\ldots, \alpha_k)\in\N_{\le 2 \log N}^k: p_1^{\alpha_1}\cdots p_k^{\alpha_k}\le N\right\}
		\]
		we obtain the desired bound 
		\[
		\# \mathcal A\le \binom{\lfloor 10 \log N / \log \log N \rfloor + k}{k} \lfloor 2 \log \log N \rfloor^{k}\le C N^{C \frac{\log \log \log N}{\log \log N}}.
		\]
		Indeed, the number of tuples $(i_1, \dots, i_k)\in\N^k$ with $i_1+ \cdots +i_k\leq \lfloor \frac{10 \log N}{\log \log N} \rfloor$ is at most $\binom{\lfloor 10 \log N / \log \log N \rfloor + k}{k} \leq 2^{15 \log N / \log \log N}$, and for each fixed $(i_1, \dots, i_k)$ the size of $I_{i_1}\times\cdots\times I_{i_k}$ is at most $\lfloor 2 \log \log N \rfloor^{k} \leq C N^{C \frac{\log \log \log N}{\log \log N}}$.
	\end{proof}
	
	We say that $\tilde q \in \Z_+$ is an \emph{admissible divisor} of $q$ if it divides $q$, satisfies $\tilde q \in [N]$, and has exactly the same prime divisors as $q$. By Lemma~\ref{L4} we can split $\mathcal I(q) \cap \mathcal R^d_{\le N}$ into disjoint sets $\mathcal R(q_1), \dots,\mathcal R(q_{m})$ for some $m \in [ \mathbf C_{\rm div} N^{\mathbf C_{\rm div} \frac{\log \log \log N}{\log \log N}} ]$, where $q_1, \dots, q_{m}$ are admissible divisors of $q$.  By the triangle inequality followed by H{\"o}lder's inequality, it suffices to show that
	\begin{align*}
	\sum_{x \in \Z^d} 
	\sum_{q \in {\mathcal S}_L}
	\Big( 
	\sum_{u' \in \mathcal I'}
	\Big \| 
	\sum_{u \in\mathcal R(\tilde q)} 
	\tilde f_{u+u'} (x) \Big\|_{H}^{2} \Big)^r
	\leq C_{r,d}^{2r} N^{2r C  \frac{\log \log \log N}{\log \log N}}
	\|f\|_{\ell^{2r}(\Z^{d};H)}^{2r}
	\end{align*}
	holds uniformly for all possible choices of admissible divisors $\tilde q$ of $q$, where the function $\tilde f_{u+u'}$ is a variant of $f_{u+u'}$ without the indicator function $\ind{N,k,j}$. Indeed, this follows because $\ind{N,k,j}(u+u')$ is constant on each $\mathcal R(\tilde q)$ when $u' \in \mathcal I'$ is fixed, and equal to $0$ if $u \in \mathcal I(q) \setminus \mathcal R^d_{\le N}$.
	
	Let $\phi \colon \R^d \to [0,1]$ be smooth and such that $\ind{[-\frac{6}{5}, \frac{6}{5}]^d} \leq \phi \leq \ind{[-\frac{9}{5}, \frac{9}{5}]^d}$. Let $\psi \colon \R^d \to [0,\infty)$ be smooth, supported on $[-\frac{1}{5}, \frac{1}{5}]^d$, and such that $\| \psi \|_{L^1(\R^d)} = 1$. Define $\zeta \coloneqq \phi*\psi$ so that $\ind{[-1, 1]^d} \leq \zeta \leq \ind{[-2, 2]^d}$. Proceeding as in \cite[Section~2.6]{MSZ3}, for each $q\in {\mathcal S}_L$ and each admissible divisors $\tilde q$ of $q$, we obtain
	\begin{align*}
	\sum_{x \in \Z^d} 
	\Big( 
	\sum_{u' \in \mathcal I'}
	\Big \| 
	\sum_{u \in\mathcal R(\tilde q)} 
	\tilde f_{u+u'}(x) \Big\|_{H}^{2} \Big)^r
	\leq C_{r,d}^{2r} 
	\sum_{x \in \Z^d} 
	\Big( 
	\sum_{u' \in \mathcal I'}
	\Big \| 
	\sum_{u \in\mathcal R(\tilde q)} 
	f^\zeta_{u+u'}(x) \Big\|_{H}^{2} \Big)^r,
	\end{align*}
	where $f^\zeta_{u+u'}$ is a variant of $\tilde f_{u+u'}$ with $\mathfrak m$ replaced by a nice bump multiplier $\zeta$.
	
	\subsubsection{\textbf{Square function estimates for the nice bump multiplier}}
	Finally, for each fixed choice of admissible divisors $\tilde q$ of $q\in {\mathcal S}_L$, and every $t\in[1, \infty]$, we shall show the following inequality
	\begin{align*}
	\Bigg(\sum_{x \in \Z^d} 
	\sum_{q \in \mathcal S_L}
	\Big( 
	\sum_{u' \in \mathcal I' }
	\Big \| 
	\sum_{u \in\mathcal R(\tilde q)} 
	f^\zeta_{u+u'}(x) \Big\|_{H}^{2} \Big)^t \Bigg)^{\frac{1}{2t}}
	\leq C_{d} N^{C \frac{\log \log \log N}{\log \log N}}
	\|f\|_{\ell^{2t}(\Z^{d};H)}
	\end{align*}
	that will be interpreted as a linear operator norm bound from $\ell^{2t}(\Z^{d};H)$ to $\ell^{2t}(\Z^{d} \times \mathcal  S_{L};\ell^{2}( \mathcal I' ;H))$. By interpolation it suffices to consider $t \in \{1, \infty\}$. 
	
	The case $t=1$ holds with an absolute constant independent of $N$ by Plancherel's theorem. Indeed, the functions $\mathcal F_{\Z^d} f^\zeta_{u+u'}$ have disjoint supports, since $\tilde q \le N$, $r\ge2$, and  $\vartheta \leq (4r N^{2r})^{-1}\le (8 N^{4})^{-1}$, and for distinct $q, q'\in {\mathcal S}_L$ the corresponding admissible divisors $\tilde q, \tilde q'$ have distinct sets of prime divisors. 
	
	In the remaining case $t=\infty$, it suffices to show that
	\begin{align}
	\label{IWeq:12}
	\sup_{x\in\Z^{d}} \sup_{q \in {\mathcal S}_{L}} \sum_{u' \in \mathcal I' } 
	\Big\| \sum_{u\in \mathcal Q(\tilde q)} f^\zeta_{u+u'}(x)
	\Big\|_{H}^2
	\leq
	C_{d}^2 \|f\|_{\ell^{\infty}(\Z^{d};H)}^{2}
	\end{align}
	holds uniformly for all $\tilde q \in [N]$ that are (not necessarily admissible) divisors of $q$, and once \eqref{IWeq:12} is established, the desired inequality with $t=\infty$ follows. Indeed, it suffices to apply the inclusion--exclusion formula involving the M\"obius function $\mu$, which asserts that
	\[
	\sum_{u\in \mathcal R(\tilde q) } F(u)
	= \sum_{b \in \mathcal D(\tilde q)} \mu\left(\frac{\tilde q}{b}\right) \sum_{u \in \mathcal Q(b)} F(u),
	\]
	where $\mathcal D(\tilde q)$ is the set of all positive divisors of $\tilde q$, and to use \eqref{IWeq:12} together with Lemma~\ref{L4}, which implies that the number of all positive divisors of a given admissible divisor $\tilde q \in [N]$ is at most $\lfloor \mathbf C_{\rm div} N^{\mathbf C_{\rm div} \frac{\log \log \log N}{\log \log N}} \rfloor$. In order to prove \eqref{IWeq:12}, we proceed as in \cite[Section~2.7]{MSZ3}.
	
	\subsection{A seminorm variant of the Ionescu--Wainger theorem}
	The goal of this subsection is to obtain an $r$-variation variant of Theorem~\ref{thm:IW}; see \eqref{eq:45}. 
	
	\begin{theorem} \label{thm:IWvar} Let
		$p_0' \leq p \leq p_0$ for some $p_0 \in 2\Z_+$ and let $r\in (2, \infty]$. Then there exists a constant ${\bf C}_{\rm IW}(p_0,r) \in \R_+$ such that for every $d \in \Z_+$ and every $N \in \N_{\ge 100}$ the following is true. Assume that $\mathbb D\subseteq [1, \infty)$ is $\lambda$-lacunary for some $\lambda>1$; see \eqref{eq:46}. For each $n\in\mathbb N$, let
		\[
		0 < \vartheta_n \leq \min\{(2 p_0 N^{p_0})^{-1}, \lambda_n^{-1}\}, 
		\]
		and let $\mathfrak m_n \colon \R^{d} \to \mathbb C$ be measurable and supported on $[-\vartheta_n, \vartheta_n]^d$. Moreover, set
		\begin{align}
		\label{eq:54}
		\begin{split}
		{\bf A}_{p, r} & \coloneqq \|(T_{\R^d}[\mathfrak m_n])_{n\in\N}\|_{L^{p}(\R^d)\to L^{p}(\R^d;{\bf V}^r(\N))},\\
		{\bf B}_{p_0} & \coloneqq \sup_{(\omega_n)_{n\in\N}\in\{-1,1\}^{\mathbb N}}
		\Big\|\sum_{n\in\mathbb N} \omega_n T_{\R^d}[\mathfrak m_{n+1}-\mathfrak m_n]\Big\|_{L^{p_0}(\R^d)\to L^{p_0}(\R^d)}.
		\end{split}
		\end{align}
		Denote $\rho'_N \coloneqq \rho'_N(p_0,r) \coloneqq {\bf C}_{\rm IW}(p_0,r) \frac{\log \log \log N}{\log \log N}$. Then
		\begin{align}
		\label{eq:7}
		\Big\|(T_{\Z^d}^{\mathcal R_{\le N}^d}[\mathfrak m_n])_{n\in\N}\Big\|_{\ell^p(\Z^d)\to \ell^p(\Z^d;{\bf V}^r(\N))}
		\lesssim_{\lambda, d,p_0,  r} N^{\rho'_N} ({\bf A}_{p, r}+{\bf B}_{p_0} ).
		\end{align}
	\end{theorem}
	
	Some remarks about Theorem~\ref{thm:IWvar} are in order.
	
	\begin{enumerate}[label*={\arabic*}.]
		\item Theorem~\ref{thm:IWvar} is a seminorm variant of the Ionescu--Wainger theorem stated at the beginning of this section.  We present a fairly general argument, which can be used to deduce similar results for jumps or oscillations (as well as for norms corresponding to real interpolation spaces) in place of $r$-variations. We refer to \cite{MSZ1} for the definitions.
		
		\item The restriction $r\in (2, \infty]$ arises from the fact that the key tool used in our proof is the Rademacher--Menshov inequality from \eqref{eq:164}. It is highly tempting to ask whether Theorem~\ref{thm:IWvar} remains true for arbitrary seminorms in place of $r$-variations like the Magyar--Stein--Wainger sampling principle; see Proposition~\ref{prop:msw}.
		
		\item In our applications we will have to verify that the quantities ${\bf A}_{p, r}$ and ${\bf B}_{p_0}$ are finite. It might be a difficult task in general, but we will be working with multipliers $\mathfrak m_n$ that can be handled using techniques from classical harmonic analysis \cite{bigs}. To prove that ${\bf A}_{p, r}$ is finite we will use a square function argument from \cite{JSW}. This argument will allow us to replace $(T_{\R^d}[\mathfrak m_n])_{n\in\N}$ with a dyadic martingale. The latter will be handled using the $r$-variational inequality of L{\'e}pingle for martingales with $r>2$; see, for instance, \cite{Lep,MSZ1}. We note that the restriction to $r\in (2, \infty]$ is necessary in this case. To verify that ${\bf B}_{p_0}$ is finite we will use classical square function methods in the spirit of Littlewood and Paley; see, for instance, \cite{bigs, RdF, MSZ2}.
	\end{enumerate}
	
	\begin{proof}[Proof of Theorem~\ref{thm:IWvar}]
		Fix a $\lambda$-lacunary set $\mathbb D\subseteq [1, \infty)$ for some $\lambda>1$, and let $\kappa_N$ satisfy $2^{\kappa_N} = 100p_0 (\log\lambda)^{-1}N$. We show for all $f\in \ell^p(\Z^d)$ that
		\begin{align}
		\label{eq:62}
		&\Big\|{\bf V}^r(T_{\Z^d}^{\mathcal R_{\le N}^d}[\mathfrak m_n]f: n\in\N_{<2^{\kappa_N}})\Big\|_{\ell^p(\Z^d)}\lesssim_{\lambda, d,p_0,  r}
		N^{\rho'_N} {\bf B}_{p_0}\|f\|_{\ell^p(\Z^d)},\\
		\label{eq:63}
		&\Big\|{\bf V}^r(T_{\Z^d}^{\mathcal R_{\le N}^d}[\mathfrak m_n]f: n\in\N_{\ge 2^{\kappa_N}})\Big\|_{\ell^p(\Z^d)}\lesssim_{\lambda, d,p_0,  r}
		N^{\rho'_N} {\bf A}_{p, r}\|f\|_{\ell^p(\Z^d)}.
		\end{align}
		Combining \eqref{eq:62} and \eqref{eq:63}, we see that \eqref{eq:7} readily follows. By using the Rademacher--Menshov inequality from \eqref{eq:164} followed by Khinchine's inequality and Theorem~\ref{thm:IW} with ${\bf B}_{p_0}$ from \eqref{eq:54}, we obtain \eqref{eq:62}, since $\kappa_N\simeq \log N$.  Invoking Proposition~\ref{prop:msw} followed by Theorem~\ref{thm:IW}, we obtain \eqref{eq:63} with ${\bf A}_{p, r}$ from \eqref{eq:54}, and the proof is complete.
	\end{proof}
	
	\section{Inverse theorems}
	\label{sec:inverse}
	In this section, we simultaneously prove an inverse theorem for averages over distinct degree corner configurations in both the discrete and continuous settings. This will be essential to establish a multilinear Weyl inequality and a Sobolev smoothing theorem in Section~\ref{sec:weyl}. We begin with fixing necessary notation and terminology.
	
	\subsection{Basic definitions and the statement of the inverse theorem}
	Throughout, $\mathbb K$ is either the set of integers $\Z$ or the set of real numbers $\R$. If $\mathbb K=\Z$, then $\lambda_{\mathbb K}$ is counting measure on $\Z$. If $\mathbb K=\R$, then $\lambda_{\mathbb K}$ is Lebesgue measure on $\R$. Let $\K_+ \coloneqq \{x\in\K: x>0\}$. We write $\lambda_{\mathbb K^k} \coloneqq \lambda_{\mathbb K}^{\otimes k}$ for the product measure on $\K^k$. For a measurable set $E\subseteq \K^k$, we abbreviate $\lambda_{\mathbb K^k}(E)$ to $|E|_{\K^k}$.
	
	Fix $k\in\Z_+$, and consider a polynomial mapping
	\begin{align}
	\label{eq:42}
	{\mathcal P} \coloneqq (P_1,\ldots, P_k) \colon \mathbb K\to \mathbb K^k,
	\end{align}
	where $P_1, \ldots ,P_k\in\mathbb K[\nn]$ are polynomials with distinct degrees such that
	\begin{align}
	\label{eq:40}
	d_1 \coloneqq \deg P_1<\cdots < d_k \coloneqq \deg P_k,
	\end{align}
	and also define
	\begin{align}
	\label{eq:44}
	D \coloneqq D_k \coloneqq d_1+\cdots+d_k.
	\end{align}
	
	For $N\ge1$, recalling the definitions of $[N]_{\mathbb K}$ from \eqref{intervals} and $\E_{t \in [N]_{\mathbb K}}^{\lambda_{\mathbb K}}$ from \eqref{eq:30}, we define multilinear averages for $f_1, \ldots, f_k\in L^0(\K^k, \lambda_{\K^k})$ by
	\begin{align}
	\label{eq:49}
	A_{N; \mathbb K^k}^{\mathcal P}(f_1,\ldots, f_k)(x) \coloneqq \E_{t \in [N]_{\mathbb K}}^{\lambda_{\mathbb K}} \prod_{i\in[k]} f_i(x-P_i(t)e_i),
	\qquad x\in\mathbb K^k,
	\end{align}
	and their truncations by
	\begin{align}
	\label{eq:52}
	\tilde{A}_{N; \mathbb K^k}^{\mathcal P}(f_1,\ldots, f_k)(x) \coloneqq \E_{t \in [N]_{\mathbb K}\setminus [N/2]_{\mathbb K}}^{\lambda_{\mathbb K}} \prod_{i\in[k]} f_i(x-P_i(t)e_i), \qquad x\in\mathbb K^k.
	\end{align}
	For example, in the integer case $\mathbb K=\Z$, we have
	\[
	A_{N; \mathbb Z^k}^{\mathcal P}(f_1,\ldots, f_k)(x) = \lfloor N \rfloor^{-1} \sum_{n\in[N]} \prod_{i\in[k]} f_i(x-P_i(n)e_i),
	\]
	whereas, in the real case $\mathbb K=\R$, we have
	\[
	A_{N; \mathbb R^k}^{\mathcal P}(f_1,\ldots, f_k)(x) = N^{-1} \int_{0}^N \prod_{i\in[k]} f_i(x-P_i(t)e_i)dt.
	\]
	
	Let ${\mathcal C}z \coloneqq \overline z$ for all $z \in \C$. For $f, g\in L^2(\mathbb K^k, \lambda_{\mathbb K^k})$, define the inner product
	\[ 
	\langle f, g \rangle \coloneqq \langle f, g \rangle_{\lambda_{\mathbb K^k}} \coloneqq\int_{\mathbb K^k} f(x) \overline{g(x)} d\lambda_{\mathbb K^k}(x) = \int_{\mathbb K^k} f(x) \mathcal C g(x) d\lambda_{\mathbb K^k}(x).
	\]
	For compactly supported $f_0, f_1,\ldots, f_k \in L^{\infty}(\mathbb K^k, \lambda_{\mathbb K^k})$, note the identities
	\begin{equation}\label{transpose}
	\left\langle A_{N; \mathbb K^k}^{\mathcal P}(f_1,\ldots, f_k), f_0 \right\rangle 
	= \left\langle f_j,  A_{N; \mathbb K^k}^{{\mathcal P}, *j}(f_1,\ldots, f_{j-1}, f_0, f_{j+1},\ldots, f_k) \right\rangle 
	\end{equation}
	and
	\begin{equation}\label{transpose-t}
	\left\langle \tilde{A}_{N; \mathbb K^k}^{\mathcal P}(f_1,\ldots, f_k), f_0 \right\rangle
	= \left\langle f_j, \tilde{A}_{N; \mathbb K^k}^{{\mathcal P}, *j}(f_1,\ldots, f_{j-1}, f_0, f_{j+1},\ldots, f_k) \right\rangle 
	\end{equation}
	for all $j \in [k]$, where the adjoint operators $A_{N; \mathbb K^k}^{{\mathcal P}, *j}$ and $\tilde{A}_{N; \mathbb K^k}^{{\mathcal P}, *j}$ are given by
	\begin{align}
	\label{eq:108a}
	A_{N; \mathbb K^k}^{{\mathcal P}, *j}(g_1,\ldots,g_k)(x) \coloneqq 
	\E_{t \in [N]_{\mathbb K}}^{\lambda_{\mathbb K}}  
	\prod_{i\in[k]}
	g^{*j}_{i}(x- {\textbf{\textit{P}}}^{*j}_i(t))
	\end{align}
	and
	\begin{align}
	\label{eq:108b}
	\tilde{A}_{N; \mathbb K^k}^{{\mathcal P}, *j}(g_1,\ldots,g_k)(x) \coloneqq 
	\E_{t \in [N]_{\mathbb K}\setminus [N/2]_{\mathbb K}}^{\lambda_{\mathbb K}} 
	\prod_{i\in[k]}
	g^{*j}_{i}(x- {\textbf{\textit{P}}}^{*j}_i(t))
	\end{align}
	with $g^{*j}_{i} \coloneqq \mathcal C^{\ind{i \neq j}} g_{i}$ and ${\textbf{\textit{P}}}^{*j}_i(t) \coloneqq \ind{i \neq j} P_{i}(t)e_{i} - P_j(t)e_j$. We often abbreviate $A_{N; \mathbb K^k}^{\mathcal P}$ to $A_{N; \mathbb K^k}$ and $\tilde{A}_{N; \mathbb K^k}^{\mathcal P}$ to $\tilde{A}_{N; \mathbb K^k}$. In some instances, we will write out the averages $A_{N; \mathbb K^k}^{\mathcal P}=A_{N; \mathbb K^k}^{P_1,\ldots,  P_k}$ and $\tilde{A}_{N; \mathbb K^k}^{\mathcal P}=\tilde{A}_{N; \mathbb K^k}^{P_1,\ldots,  P_k}$ depending on how explicit we need to be, and likewise with $A_{N; \mathbb K^k}^{{\mathcal P}, *j}$ and $\tilde{A}_{N; \mathbb K^k}^{{\mathcal P}, *j}$ for $j\in[k]$.
	
	We also fix some notation and terminology from the classical circle method. For a finite set $\Sigma\subset \hat{\K}$ and any $N\in\R_+$, we define the set of major arcs by
	\begin{align}
	\label{eq:78}
	\mathfrak M_{\le N}(\Sigma) \coloneqq \bigcup_{\theta\in \Sigma}[\theta-N, \theta+N].
	\end{align}
	Using \eqref{eq:78}, we define the set of major arcs in the $j$-th component by
	\begin{align}
	\label{eq:79}
	\mathfrak M_{\le N}^j(\Sigma) \coloneqq \K^{j-1}\times \mathfrak M_{\le N}(\Sigma)\times \K^{k-j}.
	\end{align}
	
	Recall the definition of the set of canonical fractions $\mathcal R_{\le N}^d$ from  \eqref{IWeq:372}. For every $N\ge1$, when $d=1$, we shall abbreviate $\mathcal R_{\le N}^d$ to
	\begin{align}
	\label{eq:129}
	\mathcal R_{\le N} \coloneqq \left \{ \frac{a}{q}\in\mathbb Q / \Z: q\in[N] \text{ and } \gcd(a, q)=1 \right \}.
	\end{align}
	We shall also write $\mathcal R_{\le N}^{\K} \coloneqq \mathcal R_{\le N}$ if $\K=\Z$ and $\mathcal R_{\le N}^{\K} \coloneqq\{0\}$ if $\K=\R$.
	
	Let $M_1\ge1$ and $M_2\in \R_+$, and recall the definition of $T_{i, \K^k}^{\Sigma}[\mathfrak m]$ from \eqref{eq:106}. Using this definition, we introduce the Ionescu--Wainger projections by setting
	\begin{align}
	\label{eq:132}
	\Pi_{\K}^i [\le M_1, \le M_2]f(x) \coloneqq T_{i, \K^k}^{\mathcal R_{\le M_1}^{\K}}[\eta_{[\le M_2]}]f(x)
	\quad \text{for} \quad
	i \in[k], \, x\in\K^k,
	\end{align}
	where $\eta \colon \R\to[0, 1]$ is a smooth and even function satisfying \eqref{eq:126}, while
	\begin{align*}
	\eta_{[\le M_2]}(x) \coloneqq \eta(M_2^{-1}x) \quad \text{for} \quad x\in \R.
	\end{align*}
	Note that the projections from \eqref{eq:132} localize the major arcs ${\mathfrak M}_{\leq M_2}^i(\mathcal R_{\leq M_1}^{\K})$. For every $p\in[1, \infty]$ and $f\in L^p(\K^k)$, by the triangle inequality, we have
	\begin{align*}
	\left\|\Pi_{\K}^i [\le M_1, \le M_2]f\right\|_{L^p(\K^k)}\lesssim \#\mathcal R_{\leq M_1}^{\K}\|f\|_{L^p(\K^k)}.
	\end{align*}
	This bound is much worse than the bounds obtained in Section~\ref{sec:IW}, but it holds for all $p\in[1, \infty]$ and will be sufficient for the purpose of this section.
	
	The main result of this section is the following inverse theorem.
	
	\begin{theorem}\label{pp}
		Let $\K$ be either $\Z$ or $\R$. Fix $C_0 \in \Z_+$ and $k\in\Z_+$, and let $\mathcal P$ be a polynomial mapping satisfying \eqref{eq:42}--\eqref{eq:44}. Then there exist large constants $C_1,C_2 \in \Z_+$  depending only on $C_0$ and $\mathcal P$ such that the following holds. Assume that $\delta\in(0, 1]$ and $N \geq C_{1} \delta^{-C_{1}}$, and fix $j \in [k]$. If $f_0, f_1,\ldots, f_k \in L^\infty(\K^k)$ are $1$-bounded functions supported on $I \coloneqq \prod_{i\in[k]}[ \pm C_0 N^{d_i}]_{\K}$ such that
		\begin{equation}\label{hyp}
		\left|\left\langle \tilde{A}_{N;\K^k}^{\mathcal P}(f_1,\ldots, f_k), f_0 \right\rangle\right| \geq \delta N^{D},
		\end{equation}
		then there exists a $1$-bounded function $h_j \in L^\infty(\K^k)$ supported on $ I$ such that
		\begin{align*}
		\left|\left\langle f_j, \Pi_{\K}^j [\le M_1,
		\le M_2]h_j\right\rangle\right|
		\ge C_{1}^{-1} \delta^{C_{1}} N^D,
		\end{align*}
		where $M_1 \coloneqq C_{2} \delta^{-C_{2}}$ and $M_2 \coloneqq C_{2} \delta^{-C_{2}}N^{-d_j}$.
	\end{theorem}
	
	We will prove this theorem in the next few subsections by extending the methods from \cite{P2, PP2, PPS}. In particular, we adapt the strategy from \cite[Subsection~7.1]{PP2} to the multidimensional setting in order to deduce a structural result for each of the $1$-bounded functions $f_1,\dots,f_k$ in Theorem~\ref{pp}. A new difficulty specific to the multidimensional setting arises during the degree lowering part of the argument, coming from the dependence of phases produced by the $U^2$-inverse theorem on variables besides the differencing parameters, which prevents one from simply applying the major arc lemma as in~\cite{P2}. In~\cite{PPS}, the issue could be sidestepped by a simple application of the pigeonhole principle, but this trick is limited to the case of configurations of the form $(x_1,x_2),(x_1+P_1(n),x_2),(x_1,x_2+P_2(n))$ with $P_1$ linear. To overcome the new difficulty in the general distinct degree case, we perform a more elaborate pigeonholing argument and introduce another application of PET induction combined with concatenation inside the main degree lowering argument.
	
	\subsection{Gowers norms and their basic properties}
	We continue fixing necessary notation. We also gather important tools that will be used later.
	
	\subsubsection{\textbf{Uniform and Fej{\'e}r measures}} Let $k\in\Z_+$. For any measurable $I\subseteq \K^k$ such that $0<|I|_{\K^k}<\infty$, we define the uniform and Fej{\'e}r measures by
	\begin{align*}
	d\lambda_{\K^k, I}(x) \coloneqq |I|^{-1}_{\K^k} \cdot \ind{I}(x)d\lambda_{\K^k}(x)
	\quad \text{and} \quad
	d\sigma_{\K^k, I}(x) \coloneqq \kappa_{\K^k, I}(x)d\lambda_{\K^k}(x),
	\end{align*}
	respectively, where $\kappa_{\K^k, I}(x)$ is the Fej{\'e}r kernel defined by
	\begin{align}
	\label{eq:56}
	\kappa_{\K^k, I}(x) \coloneqq 
	|I|^{-2}_{\K^k} \cdot \ind{I} \ast \ind{-I}(x)
	= |I|^{-2}_{\K^k} \, |\{y\in I: x+y\in I\}|_{\K^k}.
	\end{align}
	We observe that $\lambda_{\K^k, I}(\K^k)=1$ and $\sigma_{\K^k, I}(\K^k)=1$. We will mainly work with $k=1$, and we shall write $\K$ in place of $\K^1$ in the above formulae. If $\K=\R$, $N\in \R_+$, and $I=[N]_{\R}$, then, by \eqref{eq:56}, one can easily check that the Fej{\'e}r kernel $\kappa_{\R, I}(x)$ equals $N^{-1}\big(1-N^{-1}|x|\big)$ for $|x|\le N$, and $0$ otherwise. 
	
	\subsubsection{\textbf{Multiplicative discrete derivatives}}
	Let $k\in\Z_+$. For every $f\in L^0(\K^k)$ and $x, h, h'\in\K^k$, we define the multiplicative discrete derivatives by 
	\begin{align*}
	\Delta_h f(x) \coloneqq f(x){\overline{f(x+h)}}
	\quad \text{and} \quad
	\Delta_{(h, h')}' f(x) \coloneqq {\overline{f(x+h)}} f(x+h').
	\end{align*}
	For $s\in \Z_+$ and vectors ${\bm h} = (h_1, \ldots, h_s) \in (\K^k)^s$ and ${\bm h'} = (h_1', \ldots, h_s') \in (\K^k)^s$, we also define the iterated multiplicative discrete derivatives by 
	\begin{align*}
	\Delta_{\bm h}f(x)& \coloneqq \Delta_{h_1, \ldots, h_s} f(x) \coloneqq \Delta_{h_1}(\cdots(\Delta_{h_s} f(x))\cdots),\\
	\Delta_{(\bm h, \bm h')}'f(x)& \coloneqq \Delta_{(h_1, h_1'), \ldots, (h_s, h_s')}' f(x) \coloneqq \Delta_{(h_1, h_1')}'(\cdots (\Delta_{(h_s, h_s')}' f(x))\cdots).
	\end{align*}
	
	\begin{remark}
		\label{remark:1}
		The above formulae may be rewritten more concisely as
		\begin{align}
		\label{product}
		\begin{split}
		\Delta_{\bm h} f(x) &=\prod_{\omega \in \{0,1\}^s} {\mathcal C}^{|\omega|} f(x + \omega \circ {\bm h}),\\
		\Delta_{(\bm h, \bm h')}'f(x)&=\prod_{\omega \in \{0,1\}^s} {\mathcal C}^{|\omega|} f(x + \omega \circ {\bm h}+({\bm 1}-\omega) \circ {\bm h}'),
		\end{split}
		\end{align}
		where we denote ${\bm 1} \coloneqq (1,\dots, 1) \in \{0,1\}^s$ and define $\omega \circ {\bm h} \coloneqq \sum_{i\in [s]} \omega_i h_i$ for every $\omega = (\omega_1, \ldots, \omega_s) \in \{0,1\}^s$ and ${\bm h} = (h_1, \ldots, h_s) \in (\K^k)^s$. If $k=1$, then $\omega \circ {\bm h}$ coincides with the usual inner product $\omega \cdot {\bm h}$ on $\K^s$.
	\end{remark}
	
	\subsubsection{\textbf{The Gowers box and uniformity norms}}
	We use multiplicative discrete derivatives to define the (localized) Gowers box and uniformity norms.
	
	\begin{definition}[Gowers box and uniformity norms]
		\label{def:4}
		For $k,s \in \Z_+$ and $n\in[k]$, let $V\equiv\K^n$ be an $n$-dimensional linear subspace (when $\K=\R$) or sublattice (when $\K=\Z$) of $\K^k$ endowed with the natural measure $\lambda_{\K^n}$. Let $I$ be a  subset of $\K^k$ with finite and positive measure with respect to $\lambda_{\K^k}$ and $H_1,\ldots, H_s$ be subsets of $V$ with finite and positive measure with respect to $\lambda_{\K^n}$. Setting ${\bm H} \coloneqq \prod_{i \in [s]} H_i$, the normalized Gowers box norm of $f\in L^\infty(\K^k)$ supported on $I$ with respect to the sets of translates $H_1,\ldots, H_s$ is given by
		\begin{align}
		\label{eq:57}
		\| f\|_{\square^s_{H_1,\ldots, H_s}(I)} \coloneqq
		\Big( |I|^{-1}_{\K^k} \int_{\K^k}
		\E_{{\bm h}, {\bm h}' \in {\bm H}}^{\lambda_{\K^n}^{\otimes 2s}} \Delta_{({\bm h}, {\bm h}')}' f(x)
		d\lambda_{\K^k}(x) \Big)^{2^{-s}}.
		\end{align} 
		After the change of variables $x \mapsto x - h_1'- \cdots - h_s'$, we can rewrite \eqref{eq:57} as
		\begin{align*}
		\| f\|_{\square^s_{H_1,\ldots, H_s}(I)} \coloneqq
		\Big( |I|^{-1}_{\K^k} \int_{\K^k}
		\E_{{\bm h}, {\bm h}' \in {\bm H}}^{\lambda_{\K^n}^{\otimes 2s}} \Delta_{{\bm h} - {\bm h}'} f(x)
		d\lambda_{\K^k}(x) \Big)^{2^{-s}}.
		\end{align*} 
		For a subset $J \subset V$ with finite and positive measure with respect to $\lambda_{\K^n}$, the Gowers $U^s$-norm of $f\in L^\infty(\K^k)$ supported on $I$ with respect to $J$ is given by
		\begin{align}
		\label{eq:59}
		\|f\|_{U^s_J(I)} \coloneqq \| f\|_{\square^s_{J,\ldots, J}(I)}.
		\end{align}
	\end{definition}
	
	We make some remarks about Definition~\ref{def:4}.
	\begin{enumerate}[label*={\arabic*}.]
		
		\item Although it is not immediately obvious, the quantity on the right-hand side of \eqref{eq:57} is nonnegative and the box norms satisfy the triangle inequality. Proofs of these standard facts and all other statements in this subsection can be found in~\cite[Appendix~B]{GT-primes}; see also~\cite{DKK} for generalizations. 
		
		\item After the change of variables $x\mapsto x-h_1'-\cdots-h_s'$ in \eqref{eq:57}, we can rewrite the $2^s$-th power of the local Gowers box norm as
		\begin{align}
		\label{eq:165}
		\| f\|_{\square^s_{H_1,\ldots, H_s}(I)}^{2^s}
		&= |I|_{\K^k}^{-1} \ \int_{\K^k} \E_{h_1\in V}^{\sigma_{\K^n, H_1}}\cdots \E_{h_s\in V}^{\sigma_{\K^n, H_s}}\Delta_{h_1,\dots,h_s} f(x) \, d\lambda_{\K^k}(x).
		\end{align}
		In fact, our arguments will more often feature \eqref{eq:165} than \eqref{eq:57}.
		
		\item For every $m\in[s]$, from \eqref{product} and \eqref{eq:165} we immediately see that
		\begin{align}
		\label{eq:166}
		\|f\|_{\square_{H_1, \ldots, H_{s}}^{s}(I)}^{2^{s}}
		=\E_{h_{m+1}\in V}^{\sigma_{\K^n, H_{m+1}}}\cdots \E_{h_s\in V}^{\sigma_{\K^n, H_{s}}}
		\|\Delta_{h_{m+1},\ldots, h_s} f\|_{\square_{H_1, \ldots, H_m}^m(I)}^{2^{m}}.
		\end{align}
		We similarly relate the Gowers $U^{s+1}$-norm to the Gowers $U^s$-norm.
		
		\item Finally, we emphasize that we will not usually need the full generality of the multidimensional definitions of \eqref{eq:57} or \eqref{eq:59}. We will mainly work with the case $k\in\Z_+$ and $n=1$. To be more specific, we will use \eqref{eq:57} or \eqref{eq:59} with $V=\K e_j$ for some $j\in[k]$, where $\K e_j \coloneqq \{te_j: t\in \K\}$ and $\{e_i : i\in[k]\}$ is the standard basis of $\R^k$.
	\end{enumerate}
	
	We can similarly define the Gowers box inner product.
	
	\begin{definition}[Gowers box inner products]
		\label{def:GCS}
		For $k,s \in \Z_+$ and $n\in[k]$, let $V,I,H_1,\dots,H_s$ be as in Definition~\ref{def:4}. For every $\omega\in\{0, 1\}^s$ let $f_{\omega}\in L^\infty(\K^k)$ be supported on a set $I_\omega$ with $|I_\omega|_{\K^k} = |I|_{\K^k}$. The Gowers box inner product of $(f_{\omega})_{\omega\in\{0, 1\}^s}$ with respect to $H_1,\ldots, H_s$ is given by
		\begin{align*}
		\left\langle (f_{\omega})_{\omega\in\{0, 1\}^s}\right\rangle_{H_1,\ldots, H_s}^I
		\coloneqq |I|^{-1}_{\K^k} \int_{\K^k}
		\E_{{\bm h}, {\bm h}'\in {\bm H}}^{\lambda_{\K^n}^{\otimes 2s}}\prod_{\omega \in \{0,1\}^s} {\mathcal C}^{|\omega|} f_{\omega}(x + [{\bm h}, {\bm h'}]_\omega)
		d\lambda_{\K^k}(x), 
		\end{align*}
		where ${\bm H} \coloneqq \prod_{i \in [s]} H_i$ and $[{\bm h}, {\bm h'}]_\omega \coloneqq \omega \circ {\bm h}+({\bm 1}-\omega) \circ {\bm h}'$ for all ${\bm h}, {\bm h}'\in {\bm H}$ and $\omega \in \{0,1\}^s$. After the change of variables $x\mapsto x-h_1'-\cdots-h_s'$, we can rewrite this as
		\[
		|I|^{-1}_{\K^k} \int_{\K^k} \E_{{\bm h}\in V^s}^{\sigma_{\K^n, H_1}\otimes\cdots\otimes\sigma_{\K^n, H_s}} \prod_{\omega \in \{0,1\}^s} {\mathcal C}^{|\omega|} f_{\omega}(x + \omega \circ {\bm h}) d\lambda_{\K^k}(x).
		\]
	\end{definition}
	
	Iterating the Cauchy--Schwarz inequality yields the following inequality, which is a local version of \cite[Equation~11.6]{TV}.
	
	\begin{proposition} [The Gowers--Cauchy--Schwarz inequality]
		\label{prop:GCS}
		Let $V,I,H_1,\dots,H_s,{\bm H}$ and $ (f_{\omega})_{\omega\in\{0, 1\}^s}, (I_{\omega})_{\omega\in\{0, 1\}^s}$ be as in Definition~\ref{def:GCS}. Then we have
		\begin{align}
		\label{GCS-in}
		\left|\left\langle (f_{\omega})_{\omega\in\{0, 1\}^s}\right\rangle_{H_1,\ldots, H_s}^I\right|\le 
		\prod_{\omega \in \{0,1\}^s} \|f_{\omega}\|_{\square_{H_1, \ldots, H_{s}}^{s}(I_\omega)}.
		\end{align}
	\end{proposition}
	
	Finally, we will require a general version of the $U^2$-inverse theorem. For $k\in\Z_+$ and the Pontryagin dual $\hat \K^k$ (i.e., $\R^k = \hat \R^k$ or $\T^k = \hat \Z^k$), let $\lambda_{\hat \K^k}$ be Lebesgue measure on $\hat \K^k$. A proof of the following standard result can be found in~\cite[Lemma~5.1]{KPMW} in the case $\K=\R$, and a proof of a quantitatively weaker version can be found in~\cite[Lemma~2.4]{P2} in the case $\K=\Z$.
	
	\begin{proposition}[$U^2$-inverse theorem]
		\label{U2}
		Let $I, H\subset \K$ be intervals with $0 < |I|_\K, |H|_\K < \infty$. If $f\in L^\infty(\K)$ is 1-bounded and supported on $I$, then
		\begin{align*}
		\|f\|_{U^2_{H}(I)}^4 \le |H|_{\K}^{-2}
		\|\mathcal F_{\K}f\|_{L^{\infty}(\hat \K, \lambda_{\hat \K})}^2.
		\end{align*}
	\end{proposition}
	
	\subsection{The main technical reduction}
	The purpose of this subsection is to formulate a slightly stronger variant of Theorem~\ref{pp}, and then to show that the conclusion of the theorem for $j=1$ implies the same conclusion for general $j\in[k]$. To do this, we will need the notion of admissible polynomials.
	
	\begin{definition}
		\label{def:5}
		Fix $0 < \delta \leq 1$, $N \in [1,\infty)$, and $d \in \Z_+$. Let $Q({\rm t})=\sum_{i=0}^dc_{i}(Q){\rm t}^i\in \K[{\rm t}]$ be a polynomial of degree $d$ with the leading coefficient $\ell(Q) \coloneqq c_d(Q)$. We say that $Q$ is \emph{$(d, \delta, N)$-admissible with tolerance $A \in [1, \infty)$} if its coefficients satisfy 
		\begin{align}
		\label{eq:212}
		A^{-1} \delta^A \leq |\ell(Q)| \leq A \delta^{-A}
		\quad \text{and} \quad  
		|c_i(Q)| \le A \delta^{-A} N^{d-i} \quad \text{for} \quad i \in \N_{< d}.
		\end{align}
	\end{definition}
	When the value of $A$ is understood, we just say that polynomials $Q$ that satisfy \eqref{eq:212} are $(d,\delta, N)$-admissible. Now, we formulate a quantitatively uniform extension of Theorem~\ref{pp} 
	for $(d,\delta, N)$-admissible polynomials.
	
	\begin{theorem}\label{ppg}
		Fix $A \in[1, \infty)$,
		$C_0 \in \Z_+$, and $k\in\Z_+$, and let $\mathcal P$ be a polynomial
		mapping satisfying \eqref{eq:42}--\eqref{eq:44}. Then there
		exist large constants $C_1,C_2 \in \Z_+$  depending on $A,C_0, d_k$ such
		that the following holds. Assume that $\delta\in(0, 1]$ and
		$N \geq C_1 \delta^{-C_1}$, and fix $j\in [k]$. If each $P_i$ is
		$(d_i,\delta, N)$-admissible with tolerance $A$,  
		and $f_0, f_1,\ldots, f_k \in L^\infty(\K^k)$ are $1$-bounded functions
		supported on $I \coloneqq \prod_{i\in[k]}[ \pm C_0 N^{d_i}]_{\K}$ such that
		\begin{equation}\label{hypg}
		\left|\left\langle A_{N;\K^k}^{\mathcal P}(f_1,\ldots, f_k), f_0 \right\rangle\right| \geq \delta N^{D},
		\end{equation}
		then there exists a $1$-bounded function $h_j \in L^\infty(\K^k)$ supported on 
		$I$ such that
		\begin{align*}
		\left|\left\langle f_j, \Pi_{\K}^j [\le M_1,
		\le M_2]h_j\right\rangle\right|
		\ge C_1^{-1} \delta^{C_1} N^D,
		\end{align*}
		where $M_1 \coloneqq C_2 \delta^{-C_2}$ and
		$M_2 \coloneqq C_2 \delta^{-C_2}N^{-d_j}$.
	\end{theorem}
	
	A few comments about this theorem are in order.
	\begin{enumerate}[label*={\arabic*}.]
		\item A fixed polynomial $Q \in\K[{\rm t}]$ of degree $d \in \Z_+$ is $(d, \delta, N)$-admissible with tolerance $\max\{|\ell(Q)|, |\ell(Q)|^{-1}\}$ for all $0 < \delta \le 1$ and $N \geq C_1 \delta^{-C_1}$ if $C_1$ is sufficiently large in terms of the coefficients of $Q$. Thus, the entries $P_i$ of ${\mathcal P}$ from \eqref{eq:42} are $(d_i, \delta, N)$-admissible with tolerance $O_{\mathcal P}(1)$. 
		\item In view of the previous remark, Theorem~\ref{ppg} implies Theorem~\ref{pp}. Indeed, condition \eqref{hyp} for the truncated averages $\tilde A_{N;\K^k}^{\mathcal P}$ implies that
		\[
		\left|\left\langle A_{N;\K^k}^{\mathcal P}(f_1,\ldots, f_k), f_0 \right\rangle\right| \gtrsim \delta N^{D}
		\quad \text{or} \quad
		\left|\left\langle A_{N/2;\K^k}^{\mathcal P}(f_1,\ldots, f_k), f_0 \right\rangle\right| \gtrsim \delta N^{D}
		\]
		holds, which, in turn, gives \eqref{hypg} with comparable $N$ and $\delta$. Theorem~\ref{ppg} then yields the desired conclusion of Theorem~\ref{pp}.
		\item The subsections that follow are devoted to proving Theorem~\ref{ppg} in the special case $j=1$. The aim of this subsection is to show that the general case $j \in [k]$ can be derived inductively from this special case.
	\end{enumerate}
	
	Next, we introduce some convenient notation and important terminology that will facilitate working in high dimensions.
	
	\begin{definition} \label{def:nn}
		For $l, m, n\in\Z_+$ and an $l$-dimensional vector $x \in \K^l$, let 
		\begin{align*}
		x_{\le l}& \coloneqq x_{<l+1} \coloneqq (x_{1},\ldots, x_{l})\in \K^l,\\
		x_{>m} & \coloneqq x_{\ge m+1} \coloneqq  (x_{m+1},\ldots, x_{m+l})\in \K^l, \\
		x_{(m,n]} &\coloneqq (x_{m+1},\ldots, x_{n})\in \K^{l} \quad \text{if} \quad n-m=l,\\
		x_{(m,n]^c} &\coloneqq (x_1,\ldots, x_{m},x_{n+1},\ldots, x_{n+l-m})\in \K^{l}\quad \text{if} \quad m< \min\{l,n+1\}.
		\end{align*}
		We use this notation whenever we need to emphasize that the entries of the underlying vector $x \in \K^l$, whose dimension $l$ is clear from the context, are indexed differently from the standard convention $x=(x_1,\ldots,x_l)$. 
		In particular, for $m \in [l+1]$, we define $x_{\{m\}^c}$ by setting
		\[
		x_{\{m\}^c} \coloneqq x_{(m-1,m]^c} \coloneqq (x_1,\ldots, x_{m-1}, x_{m+1},\ldots, x_{l+1})\in \K^{l},
		\]
		with obvious modifications of indexing parameters when $m \in \{1,l+1\}$.
	\end{definition}
	
	\subsubsection{\textbf{Proof of Theorem~\ref{ppg} assuming Theorem~\ref{ppg} for $j=1$}}
	
	We now demonstrate that it is enough to prove Theorem~\ref{ppg} for $j=1$, as the remaining cases $j\in[k]\setminus\{1\}$ can be derived from this case. The statement of Theorem~\ref{ppg} for $j=1$ is as follows.
	
	\begin{theorem}\label{ppp}
		Fix $A \in[1, \infty)$, $C_0 \in \Z_+$, and $k\in\Z_+$, and let $\mathcal P$ be a polynomial mapping satisfying \eqref{eq:42}--\eqref{eq:44}. Then there exist large constants $C_{1,1},C_{1,2} \in \Z_+$ depending on $A,C_0,d_k$ such that the following holds. Assume that $\delta\in(0, 1]$ and $N \geq C_{1,1} \delta^{-C_{1,1}}$. If each $P_i$ is $(d_i,\delta, N)$-admissible with tolerance $A$, and $f_0, f_1,\ldots, f_k \in L^\infty(\K^k)$ are $1$-bounded functions supported on $I \coloneqq \prod_{i\in[k]}[ \pm C_0 N^{d_i}]_{\K}$ such that
		\begin{equation}
		\label{eq:186}
		\left|\left\langle A_{N;\K^k}^{\mathcal P}(f_1,\ldots, f_k), f_0 \right\rangle\right| \geq \delta N^{D},
		\end{equation}
		then we have
		\begin{align}
		\label{eq:197}
		\left|\left\langle f_1, \Pi_{\K}^1 [\le M_1,
		\le M_2]f_1\right\rangle\right|
		\ge C_{1,1}^{-1} \delta^{C_{1,1}} N^D,
		\end{align}
		where $M_1 \coloneqq C_{1,2} \delta^{-C_{1,2}}$ and $M_2 \coloneqq C_{1,2} \delta^{-C_{1,2}}N^{-d_1}$.
	\end{theorem}
	
	\begin{remark}
		\label{rem:23}
		Theorem~\ref{ppp} remains true if we assume that $f_0, f_1,\ldots, f_k \in L^\infty(\mathbb K^k)$ are $1$-bounded functions supported on  $I \coloneqq \prod_{i\in[k]}[ \pm C_0 \delta^{-C_0} N^{d_i}]_{\mathbb K}$.    
	\end{remark}
	
	Assuming momentarily that Theorem~\ref{ppp} holds, we derive Theorem~\ref{ppg}. The proof of Theorem~\ref{ppp} will be presented in the next two subsections.
	
	\begin{proof}[Proof of Theorem~\ref{ppg}]
		Assume that Theorem~\ref{ppp} holds for all $k\in\Z_+$. Fix $k \in \Z_+$ and $j\in[k]$ as in Theorem~\ref{ppg}. By Theorem~\ref{ppp}, we can assume that $k > 1$ and $j > 1$. Set $D_i \coloneqq d_1+\cdots+d_i$ for all $i\in[k]$.   Furthermore, fix $1$-bounded functions $f_0, f_1,\ldots, f_k \in L^\infty(\K^k)$ supported on $I$. We will proceed in two steps, assuming that $N\ge C_* \delta^{-C_*}$ for some large $C_*\in\Z_+$.
		
		\paragraph{\bf Step~1}
		Applying \eqref{hypg}, \eqref{transpose}, the Cauchy--Schwarz inequality, and then \eqref{transpose} again, we obtain
		\begin{align*}
		\delta N^D \le \left|\left\langle f_1, A_{N; \K^k}^{{\mathcal P}, *1}(f_0, f_2,\ldots, f_k) \right\rangle\right|
		\lesssim  N^{D/2}\left|\left\langle f_0, A_{N; \K^k}^{\mathcal P}(F_1, f_2,\ldots, f_k)\right\rangle\right|^{1/2}
		\end{align*}
		with $F_1 \coloneqq A_{N; \K^k}^{{\mathcal P}, *1}(f_0, f_2,\ldots, f_k)$, given explicitly by the formula
		\[
		F_1(x) \coloneqq \E_{y \in [N]_{\mathbb K}}^{\lambda_{\mathbb K}}  f_0(x+P_1 (y)e_1) \prod_{i\in[k]\setminus\{1\}} \mathcal C f_{i}(x-P_{i}(y)e_{i}+P_1(y)e_1)
		\]
		for all $x \in \K^k$. Now, by Theorem~\ref{ppp}, there exist large $C_{1,1}, C_{1,2} \in\Z_+$ such that, whenever $N \geq C_{1,1} \delta^{-C_{1,1}}$, we have
		\begin{align*}
		\left|\left\langle F_1,\Pi_{\K}^1[\le M_1, \le M_2] F_1\right\rangle\right|\ge C_{1,1}^{-1}\delta^{C_{1,1}}N^D,
		\end{align*}
		where $M_1 \coloneqq C_{1,2}\delta^{-C_{1,2}}$ and $M_2 \coloneqq C_{1,2}\delta^{-C_{1,2}}N^{-d_1}$, which is equivalent to
		\begin{align}
		\label{eq:14}
		\left|\left\langle f_0, A_{N; \K^k}^{\mathcal P}\left(\Pi_{\K}^{1}[\le M_1, \le M_2] F_1, f_2,\ldots, f_k\right)\right\rangle\right|
		\ge C_{1,1}^{-1}\delta^{C_{1,1}}N^D.
		\end{align}
		
		\paragraph{\bf Step~2}
		By the pigeonhole principle, at the expense of worsening the lower bound in \eqref{eq:14} by multiplying it by a large power of $\delta$, we can consider the averages over shorter intervals in \eqref{eq:14}. Then, by the mean value theorem, one can replace the function $\Pi_{\K}^{1}[\le M_1, \le M_2] F_1$ with a constant function. It follows that there exists a measurable set $X_1\subseteq [C_0N^{d_1}]_{\K}$ with $\lambda_{\K}(X_1)\gtrsim \delta^{O(1)}N^{d_{1}}$ such that, for every $x_1\in X_1$, we have
		\begin{align*}
		\left|\left\langle f_{0, x_1}^1, A_{N_1, \K^{k-1}}^{\mathcal Q_{>1}}(f_{2, x_1},\ldots, f_{k, x_1})\right\rangle\right|
		\gtrsim\delta^{O(1)}N_1^{D_k-D_1},
		\end{align*}
		where $\delta^{O(1)}N\lesssim N_1\lesssim \delta^{-O(1)}N$ and $\mathcal Q_{>1} \coloneqq \{Q_2,\ldots, Q_k\}$ is a family of polynomials $Q_i$ which are $(d_i, \delta, N)$-admissible with tolerance $A_i\lesssim_{A,C_{1,1},d_k}1$, with $f_{0, x_1}^1(x_{> 1}) \coloneqq f_{0}^1(x_{1}, x_{> 1})$ and $f_{i, x_{1}}(x_{> 1}) \coloneqq f_{i}(x_{1}, x_{> 1})$ for all $i\in [k]\setminus\{1\}$ and $x_{> 1}\in \K^{k-1}$. If $j=2$, then we apply Theorem~\ref{ppp} for each $x_1\in X_1$ (or, more precisely, we apply the result given in  Remark~\ref{rem:23}, since the functions $f_{0, x_1}^1, f_{2, x_1},\ldots, f_{k, x_1}$ are supported on $\prod_{i\in[k]}[\pm C_{2,0} \delta^{-C_{2,0}} N_1^{d_i}]_{\mathbb K}$ for some $C_{2,0} \in \Z_+$) to conclude that there exist constants $C_{2,1}, C_{2,2}\in\Z_+$ depending on $A,C_0,d_k$ such that, whenever $N\ge C_{2,1}\delta^{-C_{2,1}}$, we have
		\begin{align*}
		\left|\left\langle f_{2, x_1},\Pi_{\K}^2[\le M_1, \le M_2] f_{2, x_1}\right\rangle\right|\ge C_{2,1}^{-1}\delta^{C_{2,1}}N^{D_k-D_1},
		\end{align*}
		where $M_1 \coloneqq C_{2,2}\delta^{-C_{2,2}}$ and $M_2 \coloneqq C_{2,2}\delta^{-C_{2,2}}N^{-d_2}$. Integrating the above inequality over $x_1\in X_1$ and using duality, we obtain
		\begin{align*}
		\left|\left\langle f_{2},\Pi_{\K}^2[\le M_1, \le M_2] h_2\right\rangle\right|\gtrsim \delta^{O(1)}N^{D_k},
		\end{align*}
		where $h_2(x) \coloneqq f_{2}(x)g_2(x_1)$ for some $1$-bounded $g_2\in L^{\infty}(\K)$. This shows that Theorem~\ref{ppg} holds for $j=2$, and iterating this argument, we conclude that it holds for all $j\in[k]$ with constants $C_{j, 1}, C_{j, 2}\in\Z_+$ in place of $C_1$ and $C_2$, respectively. Finally, the implied constants $C_{j, 1}, C_{j, 2}$ for each $j\in[k]$ can be taken to be uniform in $j\in[k]$. This completes the proof of Theorem~\ref{ppg}.
	\end{proof}
	
	\subsection{Gowers norm control}
	In this subsection, we will begin the preparation to prove Theorem~\ref{ppp}. The first step involves the PET induction technique, which is a method going back to the work of Bergelson and Leibman~\cite{BL1}, with a coefficient tracking scheme as in~\cite{P2}. The primary tool of PET induction is the following variant of van der Corput's inequality.
	
	\begin{lemma}
		\label{lem:vdC}
		Fix $\mathfrak g \in L^1(\K)$ and an interval $J \subset \K$. Then
		\[
		\left| \E_{y \in J}^{\lambda_\K} \mathfrak g(y) \right|^2 \leq 2 \E_{y \in J}^{\lambda_\K} \E_{h \in \K}^{\sigma_{\K,[H]_\K}} \ind{J-h}(y) \Delta_h\mathfrak g(y)
		\]
		for every $H \in \R_+$ such that $1 \leq H \le |J|_{\K}$.
	\end{lemma}
	
	\begin{proof}
		We first insert an extra average in $h\in[H]_{\K}$ and then apply the Cauchy--Schwarz inequality.
	\end{proof}
	An important consequence of Lemma~\ref{lem:vdC} is the following simple lemma.
	
	\begin{lemma} \label{lem:pet} Fix an interval $J\subset\K$ and $H \in \R_+$ such that $1 \leq H\le |J|_{\K}$. If $\mathfrak g_1\in L^\infty(\K^k)$ and $\mathfrak g_2 \in L^\infty(\K^{k+1})$ are $1$-bounded functions such that
		\begin{align}
		\label{eq:65'}
		\|\mathfrak g_1\|_{L^1(\K^k)}\le \alpha
		\quad \text{and} \quad
		\sup_{y\in\K}\|\mathfrak g_2(\, \cdot \,, y)\|_{L^1(\K^k)}\le \alpha
		\end{align}
		for some $\alpha \in \R_+$, then one can find $\eta \in \C$ such that $|\eta| \leq 2H|J|_{\K}^{-1}$ and
		\[
		\Big|\frac{1}{\alpha}\int_{\K^{k}} \mathfrak g_1(x) \E_{y \in J}^{\lambda_\K} \mathfrak g_2(x, y) d\lambda_{\K^k}(x)\Big|^2 
		\]
		does not exceed the following nonnegative expression
		\[
		\frac{2}{\alpha} \int_{\K^{k}} \E_{y \in J}^{\lambda_\K} \E_{h \in \K}^{\sigma_{\K,[H]_\K}} \mathfrak g_2(x, y)\overline{\mathfrak g_2(x, y+h)} d\lambda_{\K^k}(x) + \eta.
		\]
	\end{lemma}
	\begin{proof}
		We first apply the Cauchy--Schwarz inequality and Lemma~\ref{lem:vdC}, and then observe that $\ind{J-h}$ can be replaced by $\ind{J}$ in the resulting bound at the cost of an error term $\eta \in \C$ satisfying $|\eta| \leq 2H|J|_{\K}^{-1}$.
	\end{proof}
	Fix $L \in \Z_+$ and a family $\{{\textbf{\textit{Q}}}_l : l \in [L] \}$ of polynomial vectors ${\textbf{\textit{Q}}}_l \coloneqq (Q_{l,1}, \dots, Q_{l,k})$. We say that ${\textbf{\textit{Q}}}_l$ is \emph{constant} if each component $Q_{l,i}$ is constant. 
	\begin{proposition}
		\label{prop:pet}
		Fix $N \in [1, \infty)$, $k \in \Z_+$, $L\in\Z_+$, $l_0 \in [L]$, $\alpha \in \R_+$, and a family $\{{\textbf{\textit{Q}}}_l : l \in [L] \}$ of polynomial vectors that are not constant. If $\mathfrak f_0, \mathfrak f_1,\ldots, \mathfrak f_L \in L^0(\K^k)$ are $1$-bounded, with $L^{1}(\K^k)$ norms bounded by $\alpha$, and
		\begin{align*}
		\Big|\frac{1}{\alpha} \int_{\K^k} \mathfrak f_0(x) \E_{y \in [N]_\K}^{\lambda_{\K}}  \prod_{l \in [L]} \mathfrak f_{l} (x- {\textbf{\textit{Q}}}_{l}(y) ) d\lambda_{\K^k}(x) \Big|\ge \varepsilon
		\end{align*} 
		for some $\varepsilon \in (0,1]$, then for every
		$\varepsilon_0 \in (0,\varepsilon^2/10]$ with
		$\varepsilon_0 \varepsilon^2 N \geq 10^3$, we have
		\begin{align}
		\label{PET-2}
		\Big| \frac{1}{\alpha} \int_{\K^k} \E_{y \in [N]_\K}^{\lambda_{\K}}  \prod_{l \in [L]} \mathfrak f_{l} (x- {\textbf{\textit{Q}}}_{l}(y) ) \overline{ \mathfrak f_{l} (x- {\textbf{\textit{Q}}}_{l}(y+h) ) } d\lambda_{\K^k}(x) \Big| \geq \frac{\varepsilon^2}{10}
		\end{align}
		for all $h$ belonging to a measurable set
		$\mathcal H \subseteq [\pm \varepsilon_0 N]_\K \setminus [\pm \varepsilon_0 \varepsilon^2 N/10^2]_\K$
		such that
		$|\mathcal H|_\K \geq \varepsilon_0 \varepsilon^2 N / 10^2$. Thus, for all $h \in \mathcal H$, substituting $x\mapsto x+{\textbf{\textit{Q}}}_{l_0}(y)$ yields
		\begin{align*}
		\Big|\frac{1}{\alpha} \int_{\K^k} \mathfrak f'_0(x)
		\E_{y \in [N]_\K}^{\lambda_{\K}} 
		\prod_{l' \in [L']} \mathfrak f'_{l'} (x- {\textbf{\textit{Q}}}'_{l'}(y) ) d\lambda_{\K^k}(x) \Big|\ge \frac{\varepsilon^2}{10}
		\end{align*} 
		for some integer $L' \leq 2L-1$, where $\mathfrak f'_0(x)$ is the product of all factors independent of $y$, the new family $\{{\textbf{\textit{Q}}}'_{l'} : l' \in [L'] \}$ of polynomial vectors consists exactly of those ${\textbf{\textit{Q}}}_{l} - {\textbf{\textit{Q}}}_{l_0}$ and ${\textbf{\textit{Q}}}_{l}( \,\cdot\, +h) - {\textbf{\textit{Q}}}_{l_0}$ that are not constant, $\mathfrak{f}_{l'}'$ is either $\mathfrak f_l$ when ${\textbf{\textit{Q}}}_{l'}'={\textbf{\textit{Q}}}_{l}-{\textbf{\textit{Q}}}_{l_0}$ or $\overline{\mathfrak f_{l}}$ when ${\textbf{\textit{Q}}}_{l'}'={\textbf{\textit{Q}}}_{l}( \,\cdot\, +h)-{\textbf{\textit{Q}}}_{l_0}$, and, consequently, $\mathfrak f'_0, \mathfrak f'_1, \dots, \mathfrak f'_{L'}$ are $1$-bounded, with $L^1(\K^k)$ norms bounded by $\alpha$.
	\end{proposition}
	
	\begin{proof}
		It suffices to apply Lemma~\ref{lem:pet} with
		$\mathfrak g_1(x)=\mathfrak f_0(x)$ and
		$\mathfrak g_2(x, y)=\prod_{l\in[L]}\mathfrak f_{l}(x-{\textbf{\textit{Q}}}_l(y))$,
		noting that $\mathfrak g_1, \mathfrak g_2$ are $1$-bounded and satisfy \eqref{eq:65'}.
	\end{proof}
	
	Repeated applications of Proposition~\ref{prop:pet} according to the PET induction scheme allow us to control the left-hand side of~\eqref{eq:186} by Gowers box norms.
	\begin{lemma}
		\label{lem:initialPET}
		Under the assumptions of Theorem~\ref{ppp}, we have
		\begin{align}
		\label{eq:32''}
		\mathbb{E}_{\bm{h}\in\K^T}^{\mu_{\K^T,\bm{H}_T}}\|f_k\|_{\square^s_{c_1(\bm{h})[\pm H']_{\K}e_k,\dots,c_s(\bm{h})[\pm H']_{\K}e_k}(I)}^{2^s} \geq C^{-1} \delta^{C}
		\end{align}
		for some integers $s,T,C \geq 2$ bounded in terms of $A,C_0,d_k$, some nonzero multilinear polynomials $c_1(\bm{h}),\dots,c_s(\bm{h})$ of degree $d_k-1$ with coefficients in $\K$ that are bounded in terms of $A$ and $d_k$, some $H' \in \K$ such that $H' \geq C^{-1} \delta^{C} N$, and $\mu_{\K^T,\bm{H}_T} \coloneqq \lambda_{\K,\mathcal{H}_1}\otimes\dots\otimes\lambda_{\K,\mathcal{H}_T}$, where each $\mathcal{H}_i$ is an interval of the form $[\pm H]_{\K}$ for some $H \in \K$ such that $H \geq C^{-1} \delta^{C} N$.
	\end{lemma}
	In fact, Kravitz, Kuca, and Leng~\cite{KKL1} have recently proven this exact Gowers box norm control result in the case $\K=\Z$ when $P_1,\dots,P_k$ have bounded coefficients (as a special case of their Proposition~4.8). Their argument goes through with no change when $P_1,\dots,P_{k}$ are merely admissible, and thus Lemma~\ref{lem:initialPET} holds for this more general class of polynomials we consider here, in both cases $\K=\Z$ and $\K=\R$. We omit the details.
	
	When $\K=\mathbb{Z}$, averages of Gowers box norms such as those appearing on the left-hand side of~\eqref{eq:32''} can be bounded in terms of $U^s$-norms using \textit{concatenation theorems}. The first concatenation theorem was proven by Tao and Ziegler in~\cite{TZ} for the purpose of obtaining asymptotics for polynomial progressions in the primes. Although it was applicable in very broad generality, their concatenation theorem was purely qualitative, as the proof of Tao and Ziegler produced no explicit bounds. Peluse and Prendiville~\cite{PP1} (for the averages of Gowers box norms arising from the nonlinear Roth configuration) and Peluse~\cite{P2} (for more general averages of Gowers box norms) proved the first quantitative concatenation theorems, which had polynomial bounds. Applying an even more general quantitative concatenation theorem of Kravitz, Kuca, and Leng~\cite[Theorem~1.10]{KKL1} (see also \cite{Ku}) with polynomial bounds proves the following key theorem in the case $\K=\mathbb{Z}$.
	
	\begin{theorem}
		\label{thm:PET-K}
		Fix $A \in[1, \infty)$, $C_0 \in \Z_+$, and $k\in\Z_+$, and let $\mathcal P$ be a polynomial mapping satisfying \eqref{eq:42}--\eqref{eq:44}. Then there exists a large constant $C_* \in \Z_+$ depending only on $A, C_0, d_k$ such that the following holds. Assume that $\delta\in(0, 1]$ and $N \geq C_{*} \delta^{-C_{*}}$. If each $P_i$ is $(d_i,\delta, N)$-admissible with tolerance $A$, and $f_0, f_1,\ldots, f_k \in L^\infty(\K^k)$ are $1$-bounded functions supported on $I \coloneqq \prod_{i\in[k]}[ \pm C_0 \delta^{-C_0} N^{d_i}]_{\K}$, then we have the following implication
		\begin{align}
		\label{eq:32'}
		\left|\left\langle A_{N;\K^k}^{\mathcal P}(f_1,\ldots, f_k), f_0 \right\rangle\right| \geq \delta N^{D}
		\implies 
		\|f_k\|_{U_{[ \pm H_*]_{\K} e_k}^s(I)} \geq C_*^{-1} \delta^{C_*},
		\end{align}
		with $H_* \coloneqq 7C_0 \delta^{-C_0} N^{d_k}$ and an integer $s \geq 2$ bounded in terms of $A, C_0, d_k$.
	\end{theorem}
	
	To deduce Theorem~\ref{thm:PET-K} from Lemma~\ref{lem:initialPET} in the case $\K=\mathbb{R}$, one can simply use the pigeonhole principle to find $\bm{h}$ for which $c_1(\bm{h}),\dots,c_s(\bm{h})\gtrsim_{A,d_k}\delta^{O_{A,d_k}(1)}N^{d_k-1}$ and 
	\begin{equation*}
	\|f_k\|_{\square^s_{c_1(\bm{h})[\pm H']_{\K}e_k,\dots,c_s(\bm{h})[\pm H']_{\K}e_k}(I)}^{2^s} \gtrsim_{A,C_0,d_k} \delta^{O_{A,d_k}(1)},
	\end{equation*}
	apply suitable changes of variables to replace $c_i(\bm{h})[\pm H']_{\mathbb{K}}$ by $[\pm H_i]_\K$, where $H_i \simeq_{A,C_0,d_k} \delta^{O_{A,d_k}(1)}N^{d_k}$ for each $i \in [s]$, and then use the Gowers--Cauchy--Schwarz inequality \eqref{GCS-in} and the support of $f_k$ to replace each $[\pm H_i]_\K$ by two intervals $[\pm H_*]_{\mathbb{K}}$ of uniform length, thereby also replacing $s$ by $s+1$. The $U^s$-norm appearing on the left-hand side of \eqref{eq:32'} can be interpreted as an average of one-dimensional $U^s$-norms with respect to the last variable:
	\begin{align*}
	\|f_k\|_{U_{[ \pm H_*]_{\K} e_k}^s(I)}^{2^s}=
	|I_{<k}|_{\K^{k-1}}^{-1} \int_{\K^{k-1}} \|f_k(x_{<k}, \, \cdot \, )\|_{U_{[ \pm H_*]_\K}^s(I_k)}^{2^s} d \lambda_{\K^{k-1}}(x_{<k}).
	\end{align*}
	
	\subsection{Degree lowering}
	
	In this subsection, we carry out the degree lowering argument, following the broad outline of \cite[Sections~7~and~8]{P2} while overcoming new difficulties that arise in the multidimensional setting.
	
	As before, $\K$ is either $\Z$ or $\R$ with the corresponding dual groups $\hat{\Z}=\T$ (identified with $[-1/2, 1/2)$) and $\hat{\R}=\R$. The dual group $\hat{\K}$ is endowed with a canonical norm: $\|\xi\|_{\hat{\K}} \coloneqq {\rm dist}(\xi, \Z)$ when $\K=\Z$ and $\|\xi\|_{\hat{\K}} \coloneqq|\xi|$ when $\K=\R$.
	
	\subsubsection{\textbf{Degree lowering preparation}}\label{sec:prep}
	First, we collect the tools used in the degree lowering argument. The following result is the inverse form of Weyl's inequality when $\K=\Z$, and of van der Corput's lemma when $\K=\R$.
	
	\begin{proposition}
		\label{lem:Weyl1}
		Fix $P\in\K[{\rm t}]$ of degree $d\in\Z_+$, of the form $P(t)=\xi_dt^d+\dots+\xi_1t+\xi_0$. Then there exists a large constant $C\in\Z_+$ depending only on $d$ such that, for every $\varepsilon \in(0, 1]$ and $N\ge C\varepsilon^{-C}$, the following holds. If
		\[
		\left|\E_{t\in[N]_\K}^{\lambda_{\K}}e(P(t))\right|\geq \varepsilon,
		\]
		then there exists a positive integer $q\in[C\varepsilon^{-C}]$ such that
		\[
		\|q\xi_i\|_{\hat{\K}}\le C\varepsilon^{-C}N^{-i}
		\]
		for all $i\in[d]$. If $\K=\R$, then we take $q=1$.
	\end{proposition}
	
	The proof for $\K=\Z$ can be found in \cite[Proposition~4.3]{GT12}, whereas the case $\K=\R$ is the classical van der Corput lemma for oscillatory integrals with polynomial phases; see \cite[Proposition~2, p.~332]{bigs} or \cite[Theorem~1.1, p.~13]{ACK}.
	
	The following simple lemma will also be needed.
	\begin{lemma}\label{lem:pigeonhole}
		Fix $\alpha\in \R$ and $\gamma \in \R_+$. If 
		$
		|\alpha-u|\leq\gamma
		$
		for some $u \in \Q$, then, for any $M \in[1, \infty)$, one can find $m \in [\pm M]$ and
		$\theta\in[-1,1]$ such that
		\[
		\alpha=u+m\frac{\gamma}{M}+\theta\frac{\gamma}{M}.
		\]
	\end{lemma}
	\begin{proof}
		We take $m \coloneqq \lfloor M( \alpha-u)\gamma^{-1}\rfloor$ and $\theta \coloneqq M(\alpha-u)\gamma^{-1}-m$.
	\end{proof}
	
	To state the next lemma, we will need some notation from~\cite{P2}. For any $s\in\Z_+$ and $X \subseteq\K^{2s}$, let $\square_{s}(X)$ denote the set of $3s$-tuples
	\[
	\left(k_{1}^{(1)},\dots,k_{s}^{(1)},k_{1}^{(2)},\dots,k_{s}^{(2)},k_{1}^{(3)},\dots,k_{s}^{(3)}\right)\in\K^{3s}
	\]
	such that
	\[
	\left(k_{1}^{(1)},\dots,k_{s}^{(1)},k_{1}^{(\omega_1+2)},\dots,k_{s}^{(\omega_s+2)}\right)\in X
	\]
	for all $\omega\in\{0,1\}^s$. The following lemma is analogous to \cite[Lemma~7.4]{P2}, and is proved in the same way.
	
	\begin{lemma}[Dual--difference interchange]\label{lem:interchange}
		Fix $C_0,s, d \in \Z_+$. Then there exists a constant $C \in \Z_+$ depending only on $C_0, s, d$ such that the following holds. Assume that $\delta \in (0,1]$ and $N \geq C \delta^{-C}$. Let $C_0^{-1}\delta^{C_0}N^{d}\le H\le C_0\delta^{-C_0}N^{d}$ and let $\mathfrak H \subseteq [H]_{\K}^{2s}$ be a measurable set such that $\lambda_{\K}^{\otimes 2s}(\mathfrak H)\ge C_0^{-1}\delta^{C_0}H^{2s}$. Consider $F(x) \coloneqq \E_{t\in[N]_{\K}}^{\lambda_{\K}}F_t(x)$, where the map 
		\[
		\K \times [N]_\K \ni (x,t) \mapsto F_t(x) \in \C
		\]
		is measurable and, for each $t\in[N]_{\K}$, the map $\K\ni x \mapsto F_t(x)$ is $1$-bounded and supported on $J \coloneqq [\pm C_0 N^{d}]_\K$. If
		\begin{align*}
		\E_{({\bm h},{\bm h}')\in \mathfrak H}^{\lambda_{\K}^{\otimes 2s}}
		\left| |J|^{-1}_{\K} \int_{\K}
		[\Delta'_{(h_i,h_i')_{i\in[s]}}
		F(x)] e(\phi({\bm h},{\bm h}')x)d\lambda_{\K}(x)\right|^2\geq C_0^{-1}\delta^{C_0}
		\end{align*}
		for some measurable function $\phi \colon \mathfrak H \to \hat{\K}$, then
		\begin{align*}
		\E_{{\bm k}\in[H]_{\K}^{3s}}^{\lambda_{\K}^{\otimes 3s}}\ind{\square_s(\mathfrak H)}(\bm k)
		\left| |J|^{-1}_{\K} \int_{\K}G_{{\bm k}}(x) e(\psi({\bm k})x)d\lambda_{\K}(x)\right|^2 \geq C^{-1}\delta^C,
		\end{align*}
		where, for all ${\bm k} = (k_{1}^{(1)},\dots,k_{s}^{(1)},k_{1}^{(2)},\dots,k_{s}^{(2)},k_{1}^{(3)},\dots,k_{s}^{(3)}) \in \square_s(\mathfrak H)$, we have
		\begin{align*}
		G_{{\bm k}}(x) \coloneqq 
		\E_{t\in[N]_{\K}}^{\lambda_{\K}}
		\Delta'_{(k^{(2)}_{i},k^{(3)}_i)_{i\in[s]}}F_t(x)
		\end{align*}
		and $\psi \colon \square_s(\mathfrak H) \to {\hat \K}$ is given by
		\begin{equation*}
		\psi({\bm k}) \coloneqq 
		\sum_{\omega\in\{0,1\}^s}(-1)^{|\omega|}\phi(k_{1}^{(1)},\dots,k_{s}^{(1)},k_{1}^{(\omega_1+2)},\dots,k_{s}^{(\omega_s+2)}).
		\end{equation*}
	\end{lemma}
	
	We will also need the following variant of \cite[Lemma~7.5]{P2}.
	\begin{lemma}\label{lem:GCScor}
		Fix $s, C_0 \in \Z_+$. Then there exists a constant $C \in\Z_+$ depending only on $s, C_0$ such that the following holds. Let $M, H \in[1,\infty)$, $\varepsilon \in (0,1]$, and, for each $i\in[s]$, let $\phi_i \colon \K^{2s}\to\hat{\K}$ be a measurable function independent of the $(s+i)$-th variable. If $C_0^{-1} \varepsilon^{C_0} M \leq H \leq C_0 \varepsilon^{-C_0} M$ and $f\in L^{\infty}(\K)$ is a $1$-bounded function supported on $J \coloneqq [\pm M]_{\K}$ and such that
		\[
		\E_{{\bm h},{\bm h}'\in [H]_{\K}^s}^{\lambda_{\K^s}^{\otimes 2}}
		\bigg||J|^{-1}_{\K} \int_{\K}\Delta'_{({\bm h},{\bm h}')}f(x)e\Big(\sum_{i\in[s]}\phi_i({\bm h},{\bm h}')x\Big)d\lambda_{\K}(x) \bigg|^2
		\geq \varepsilon,
		\]
		then we have
		$\|f\|_{U^{s+1}_{[H]_{\K}}(J)}\ge C^{-1} \varepsilon^{C}$.
	\end{lemma}
	
	The next lemma has no analogue in~\cite{P2}, and addresses the added complexity of executing a pigeonholing argument in the multidimensional setting.
	\begin{lemma}\label{lem:randomphi}
		Fix $s \in \Z_+$. There exists a constant $C \in\Z_+$ depending only on $s$ such that the following holds. Assume that $\delta \in (0,1]$ and $H \geq C \delta^{-C}$, and set $B \coloneqq [H]^{2s}_\K$. Then one can find a measurable function $\phi \colon B \to \hat \K$ and an exceptional set $K_\phi \subseteq \square_s(B)$ with $|K_\phi|_{\K^{3s}} \le 2C \delta^{-2} H^{-1} |\square_s(B)|_{\K^{3s}}$ such that if ${\bm k} \in \square_s(B) \setminus K_\phi$, then
		\begin{equation*}
		\|q\psi({\bm k})\|_{\hat \K} > \delta^{-1}H^{-1}
		\end{equation*}
		holds for all $q \in [\delta^{-1}]$, where $\psi({\bm k})$ is defined as in Lemma~\ref{lem:interchange}.
	\end{lemma}
	\begin{proof}
		We will consider the cases $\K = \Z$ and $\K = \R$ separately.
		
		\paragraph{{\bf Case} $\K = \Z$}
		Define a random function $\phi \colon B \to \T$ by choosing each value $\phi({\bm h},{\bm h}')$ independently and uniformly at random from $\T$. Set
		\begin{equation*}
		\mathfrak{M} \coloneqq \left\{\xi \in\T :\|q \xi\|_\T \leq\delta^{-1}H^{-1} \text{ for some positive integer }q\leq\delta^{-1}\right\}.
		\end{equation*}
		Note that $\lambda_{\T}(\mathfrak{M}) \le C\delta^{-2}H^{-1}$. Moreover, we have $\square_{s}(B)=[H]^{3s}$, and for any $\varepsilon_1,\dots,\varepsilon_{2^s} \in\{-1,1\}$ and distinct pairs $({\bm h}_1,{\bm h}_1'),\dots,({\bm h}_{2^s},{\bm h}'_{2^s}) \in B$, the function $\sum_{i\in [2^s]} \varepsilon_i\phi({\bm h}_i,{\bm h}_i') $ is also uniformly distributed on $\T$. We split $\square_s(B) = K_1 \cup K_2$, where
		\[
        K_1 \coloneqq \left\{{\bm k}\in[H]^{3s}: k_i^{(2)}\neq k_i^{(3)} \text{ for all
		} i\in[s]\right\}
        \]
		and $K_2 \coloneqq [H]^{3s}\setminus K_1$. Then $\E_{\phi}[\ind{\mathfrak{M}} (\psi({\bm k}))]=\lambda_{\T}(\mathfrak M)$ for  ${\bm k}\in K_1$, where $\E_{\phi}$ is the expected value over the possible choices of $\phi$. Consequently, we have
		\begin{align*}
		\E_{\phi} \Bigg[\sum_{{\bm k}\in [H]^{3s}} \ind{\mathfrak{M}}(\psi({\bm k})) \Bigg]
		\le   \sum_{{\bm k}\in K_1} \lambda_{\T}(\mathfrak M)
		+ \sum_{{\bm k}\in K_2} 1
		\leq C \delta^{-2}H^{-1}\cdot H^{3s}  + CH^{3s-1}
		\end{align*}
		by $|K_1|_{\Z^{3s}} \le H^{3s}$ and $|K_2|_{\Z^{3s}} \leq C H^{3s-1}$. Thus, for some $\phi \colon B \to \T$, we have
		\begin{align*}
		\sum_{{\bm k}\in [H]^{3s}} \ind{\mathfrak{M}}(\psi({\bm k}))\le 2C\delta^{-2} H^{3s-1},
		\end{align*}
		and this $\phi$ has the desired properties.
		
		\paragraph{{\bf Case} $\K = \R$} Fix a large integer $n \in \Z_+$ to be determined later, and split $[H]_{\R}^s = \bigcup_{J\in \mathcal J}J$, where $\mathcal J \coloneqq \{J_i\subseteq [H]_{\R}^s: i\in[2^{ns}]\}$ is a family of $2^{ns}$ dyadic cubes, each having volume $2^{-ns}H^s$. The cubes in $\mathcal J$ are essentially disjoint, i.e., their interiors are disjoint. Define $\phi({\bm h},{\bm h}') \coloneqq \sum_{i \in [2^{ns}]} 2^{i+1} \delta^{-1} H^{-1} \ind{J_i}({\bm h}')$. Let $\mathcal J_0$ consist of all pairs $(J,J') \in \mathcal J^2$ of cubes with essentially disjoint sides. We split $\square_s(B) = K_1 \cup K_2$, where 
		\begin{align*}
		K_1 \coloneqq \left\{{\bm k}\in[H]_{\R}^{3s}:  {\bm k}_{(s,2s]} \in J, \,  {\bm k}_{(2s,3s]} \in J' \text{ for some } (J,J') \in  \mathcal J_0 \right\}
		\end{align*}
		and $K_2 \coloneqq [H]_{\R}^{3s}\setminus K_1$. Note that $|\psi({\bm k})| > \delta^{-1} H^{-1}$ for almost all ${\bm k}\in K_1$. If ${\bm k}\in K_2$, then ${\bm k}_{(s,2s]} \in J$ and ${\bm k}_{(2s,3s]} \in J'$ for some $J, J'\in\mathcal J$, which have at least one side in common. This gives $|K_2|_{\R^{3s}} \leq s 2^{-n}H^{3s}\le H^{3s-1}$ provided that $n$ is sufficiently large in terms of $H$ and $s$. Hence,
		\begin{align*}
		\int_{[H]_{\R}^{3s}} \ind{[\pm \delta^{-1} H^{-1}]_{\R}}(\psi({\bm k})) d\lambda_{\R^{3s}}({\bm k})
		=  
		\int_{K_2} \ind{[\pm \delta^{-1} H^{-1}]_{\R}}(\psi({\bm k})) d\lambda_{\R^{3s}}({\bm k})
		\leq H^{3s-1}
		\end{align*}
		and $\phi$ has the desired properties.
	\end{proof}
	
	\subsubsection{\textbf{Conditional degree lowering argument}}
	\label{sec:cdeglow}
	We will now carry out a conditional variant of the degree lowering argument following \cite[Section~8]{P2}. 
	
	We will use the convention that $\K^0 \coloneqq \{0\}$ denotes the trivial vector space endowed with the Dirac delta measure $\lambda_{\K^0}$ at zero. As in Definition~\ref{def:nn}, for any polynomial mapping ${\mathcal P}=(P_1,\ldots, P_k) \colon \K\to \K^k$, we will write ${\mathcal P}_{\le l} \coloneqq (P_1,\ldots, P_l) \colon \K\to \K^l$ and ${\mathcal P}_{>l} \coloneqq (P_{l+1},\ldots, P_k) \colon \K\to \K^{k-l}$ so that ${\mathcal P}=({\mathcal P}_{\le l}, {\mathcal P}_{>l})$ for any $l\in[k-1]$. Moreover, let ${\mathcal P}_{>0} \coloneqq {\mathcal P}_{\le k} \coloneqq {\mathcal P}$. We analogously define ${\mathcal P}_{< l} \coloneqq {\mathcal P}_{\le l-1}$, ${\mathcal P}_{\ge l} \coloneqq {\mathcal P}_{>l-1}$, and ${\mathcal P}_{(m, l]}, {\mathcal P}_{(m, l]^c}$ for $m, l\in[k]$ with $m\le l$.
	
	Our main result here is the following.
	
	\begin{lemma}[Conditional degree lowering lemma]
		\label{prop:cdeglow}
		Fix an integer $s\geq 3$, $A\in[1, \infty)$, $C_0\in \Z_+$, $k, l \in\Z_+ \setminus \{1\}$, $m\in\Z_+$ with $ m \leq l \leq k$, and $d_1, \dots, d_k \in \Z_+$ with $d_1<\cdots<d_k$. Then there exists a large constant $C \in \Z_+$ depending only on $s,A,C_0,d_k$ such that the following holds. Assume that $\delta \in (0,1]$ and $N \geq C \delta^{-C}$. Let ${\mathcal P}=(P_1,\ldots, P_k)$ be a polynomial mapping as in \eqref{eq:42}, where each $P_i\in\K[{\rm t}]$ has degree $d_i$ and is $(d_i,\delta, N)$-admissible with tolerance $A$. Assume that $C_0^{-1} \delta^{C_0} N^{d_m} \leq H \leq C_0 \delta^{-C_0} N^{d_m}$. Let $f_0, f_1,\ldots, f_{m-1} \in L^\infty(\K^m)$ be $1$-bounded measurable functions supported on $I_{\le m} \coloneqq \prod_{i\in[m]}[ \pm C_0 \delta^{-C_0} N^{d_i}]_{\K}$. Given a measurable mapping $\xi \coloneqq \xi_{(m, l]} \coloneqq (\xi_{m+1},\dots,\xi_{ l })\colon \K^{m}\to \hat\K^{l-m}$ and a~frequency vector $\eta \coloneqq \eta_{(l, k]} \coloneqq (\eta_{ l +1},\dots,\eta_k)\in \hat\K^{k-l}$, define
		\begin{equation*}
		F_m(x)\coloneqq
		\E_{t\in[N]_\K}^{\lambda_\K}
		f_0(x+{\textbf{\textit{P}}}_{m}(t))
		\Big( \prod_{i=1}^{m-1} f_{i}(x-{\textbf{\textit{P}}}_{i}(t)+{\textbf{\textit{P}}}_{m}(t)) \Big) 
		\cdot (\xi, \eta)(x,t),
		\end{equation*}
		where ${\textbf{\textit{P}}}_i(t) \coloneqq P_i(t)e_i$ for $i\in[m]$ and 
		\[
		(\xi, \eta)(x,t) \coloneqq
		e \left(\xi(x + {\textbf{\textit{P}}}_{m}(t))\cdot {\mathcal P}_{(m, l]}(t) +\eta\cdot {\mathcal P}_{>l}(t)\right)
		\]
		for $x\in\K^{m}$ and $t\in[N]_\K$. We have the following degree lowering implication
		\begin{equation*}
		\|F_m\|_{U^s_{[H]_\K e_{ m}}(I_{\le m})}\geq C_0^{-1} \delta^{C_0}
		\implies
		\|F_m\|_{U^{s-1}_{[H]_\K e_{ m}}(I_{\le m})} \geq C^{-1} \delta^C
		\end{equation*}
		provided that Condition~\ref{con:major} formulated below is satisfied.
		
		\begin{condition}[$(m, l)$-major arc condition]
			\label{con:major}
			For $A, k, l, m, d_1, \dots, d_k, \mathcal P$ as above and $c_0 \in \Z_+$, there exists a large constant $c \in \Z_+$ depending only on $A, c_0, d_k$ such that the following holds. Assume that $\delta \in (0,1]$ and $N \geq c \delta^{-c}$, and that each $P_i\in\K[{\rm t}]$ is $(d_i,\delta, N)$-admissible with tolerance $A$. Let $g_0, g_1,\ldots, g_{m-1} \in L^\infty(\K^{m-1})$ be $1$-bounded measurable functions supported on $J_{<m}$, where $J_{<1} \coloneqq \{0\} = \K^0$ and $J_{< m} \coloneqq \prod_{i\in[m-1]}[ \pm c_0 \delta^{-c_0} N^{d_i}]_{\K}$ if $m \geq 2$. Given a measurable mapping $\zeta \coloneqq \zeta_{[m, l]} \coloneqq (\zeta_{m},\dots,\zeta_{l}) \colon \K^{m-1}\to\hat\K^{l-m+1}$ and a frequency vector $\theta \coloneqq \theta_{(l, k]} \coloneqq (\theta_{l+1},\dots,\theta_k)\in\hat\K^{k-l}$, define $G_m \in \C$ by
			\begin{align*} 
			\int_{\K^{m-1}} 
			g_0(y)
			\E_{t\in[N]_\K}^{\lambda_\K}
			\Big( \prod_{i=1}^{m-1} g_{i}(y - {\textbf{\textit{P}}}_{i}(t))
			\Big) e(\zeta(y)\cdot {\mathcal P}_{[m, l]}(t)+\theta\cdot {\mathcal P}_{>l}(t)) d\lambda_{\K^{m-1}}(y)
			\end{align*}
			with ${\textbf{\textit{P}}}_{i}(t) \coloneqq P_i(y)e_i$. Set $D_n \coloneqq d_1 + \dots + d_{n}$ for $n\in[k]$ and $D_0 \coloneqq 0$. If
			\begin{align}
			\label{eq:258}
			|G_m|
			\geq c_0^{-1} \delta^{c_0} N^{D_{m-1}},
			\end{align}
			then we can find an integer $q \in [c \delta^{-c}]$ and a measurable set $Y \subseteq J_{<m}$ such that $|Y|_{\K^{m-1}} \geq c^{-1} \delta^{c} N^{D_{m-1}}$ and, for every $y \in Y$, we have
			\begin{align}
			\label{eq:262}
			\sum_{j=m}^lN^{d_{j}} \|q \zeta_{j}(y) \|_{\hat \K}
			+
			\sum_{j=l+1}^kN^{d_{j}} \|q \theta_{j}\|_{\hat \K} 
			\leq c \delta^{-c}.
			\end{align}
			If $\K = \R$, then we take $q = 1$. 
		\end{condition}
	\end{lemma}
	
	Note that, in the special case $m=1$, condition \eqref{eq:258} can be rewritten as
	\[
	\left|g_0(0)\E_{t\in[N]_\K}^{\lambda_\K}
	e \left(\zeta(0)\cdot {\mathcal P}_{[m, l]}(t) + \theta\cdot {\mathcal P}_{>l}(t) \right)\right|\ge c_0^{-1} \delta^{c_0}.
	\]  
	This condition implies \eqref{eq:262} by Lemma~\ref{lem:Weyl1}, since the polynomials $P_1,\ldots, P_k$ have distinct degrees. Thus, the major arc condition for $m=1$ holds. In order to verify Condition~\ref{con:major} for $m\ge 2$, we will adapt an induction argument from \cite{P2} using the conditional degree lowering argument from Lemma~\ref{prop:cdeglow}.
	
	\begin{proof}[Proof of Lemma~\ref{prop:cdeglow}] The proof is broadly similar to the analogous argument in~\cite{P2}, with the main difference being the need for a more elaborate pigeonholing argument to deal with the dependence of phases on multilinear fixed variables. We will proceed in a few steps and, if necessary, we will distinguish between the cases $m\ge 2$ and $m=1$ in each step. In order to ensure consistency in our arguments and notation, when $m=1$, we will identify $\K^1$ with $\K^0\times \K$ by setting $\K\ni x\mapsto (0, x)\in\K^1$. We notice that the spaces of measurable functions $L^0(\K)$ and $L^0(\K^1)$ can also be identified by declaring $L^0(\K)\ni f\mapsto \ind{\{0\}}\otimes f\in L^0(\K^1)$. From now on, when $m=1$, we will make use of these identifications and assume that all functions $f\in L^0(\K)$ are defined on the product space $\{0\}\times \K$ with the product measure $\lambda_{\K^0}\otimes \lambda_{\K}$. We will also abbreviate $F_m$ and $G_m$ to $F$ and $G$, respectively.
		
		\paragraph{\bf Step~1} 
		By~\eqref{eq:57}, we obtain
		\[
		\E_{{\bm h},{\bm h}'\in[H]_\K^{s-2}}^{\lambda_\K^{\otimes 2(s-2)}}
		\|\Delta'_{({\bm h}\circ e_m,{\bm h}'\circ e_m)} F \|^4_{U^2_{[H]_\K e_m}(I_{\le m})} \geq (C_0^{-1} \delta^{C_0})^{2^s},
		\]
		where ${\bm h}\circ e_m \coloneqq  (h_1e_m,\ldots, h_{s-2}e_m)\in(\K^{m})^{s-2}$ for ${\bm h}=(h_1,\ldots, h_{s-2})$.
		
		For any $x_{<m}=(x_1,\dots,x_{m-1})\in \K^{m-1}$, where $x_{<1} \coloneqq 0$, we set
		\[
		F_{x_{<m}}(x_m)\coloneqq F(x_{<m},x_m) \coloneqq F(x_1,\dots,x_{m-1},x_m), \qquad x_m \in \K.
		\]
		We observe that, for some constant $C_1 \in \Z_+$, the function $F_{x_{<m}}$ is supported on the set $J_m \coloneqq [\pm C_1 \delta^{-C_1} N^{d_m}]_\K$ for any $x_{<m}\in J_{<m}$, and $F_{x_{<m}}\equiv 0$ for any $x_{<m} \notin J_{<m}$, where $J_{<m}$ is as in Condition~\ref{con:major} with $C_1$ in place of $c_0$.
		
		By the popularity principle and Proposition~\ref{U2}, for some constant $C_2 \in \Z_+$, there exists a subset $X_{<m} \subseteq J_{<m}$ with $|X_{<m}|_{\K^{m-1}} \geq C_2^{-1} \delta^{C_2} |J_{<m}|_{\K^{m-1}}$ such that if $x_{<m} \in X_{<m}$, then 
		\[
		(x_m,t) \mapsto f_0(x+{\textbf{\textit{P}}}_{m}(t))
		\Big( \prod_{i=1}^{m-1} f_{i}(x-{\textbf{\textit{P}}}_{i}(t)+{\textbf{\textit{P}}}_{m}(t)) \Big) 
		\cdot (\xi, \eta)(x,t),
		\]
		where $x = (x_{<m}, x_m)$, is a measurable function from $\K \times [N]_\K$ to $\C$, and there is a measurable set $H_{x_{<m}} \subseteq [H]_\K^{2(s-2)}$ with
		\[
		\left|H_{x_{<m}}\right|_{\K^{2(s-2)}} \geq C_2^{-1} \delta^{C_2} \Big|[H]_\K^{2(s-2)} \Big|_{\K^{2(s-2)}}
		\]
		such that, for every $({\bm h},{\bm h}')\in H_{x_{<m}}$, we have
		\begin{align}
		\label{eq:256}
		|{\mathcal F}(x_{<m}, {\bm h},{\bm h}'; \xi')|^2
		\geq C_2^{-1} \delta^{C_2} N^{2d_m}
		\end{align}
		for some $\xi' \in\hat{\K}'$, where $\hat{\K}' \subseteq \hat{\K}$ is a finite set independent of $x_{<m} \in X_{<m}$, and
		\begin{align*}
		{\mathcal F}(x_{<m}, {\bm h},{\bm h}'; \xi') \coloneqq \int_{\K} \Delta'_{({\bm h}\circ e_m,{\bm h}'\circ e_m)}
		F_{x_{<m}}(x_m) e(\xi' x_m) d\lambda_\K(x_m).
		\end{align*}
		We can take $\hat{\K}'$ finite, since $\hat{\K}\ni \xi' \mapsto {\mathcal F}(x_{<m}, {\bm h},{\bm h}'; \xi') \in \C$ is continuous.
		
		Moreover, by \eqref{eq:256}, for every $x_{<m} \in X_{<m}$ and $({\bm h},{\bm h}')\in H_{x_{<m}}$, we have
		\begin{align}
		\label{eq:257}
		|{\mathcal F}(x_{<m}, {\bm h},{\bm h}'; \phi_{x_{<m}}({\bm h},{\bm h}'))|^2
		\geq C_2^{-1} \delta^{C_2} N^{2d_m}
		\end{align}
		for some measurable function
		\[
		J_{<m} \times [H]_\K^{s-2} \times [H]_\K^{s-2} \ni (x_{<m}, \bm h, \bm h') \mapsto \phi_{x_{<m}}({\bm h},{\bm h}') \in \hat \K
		\]
		such that $\phi_{x_{<m}} \colon [H]^{2(s-2)}_\K \to \hat\K$ is measurable for every $x_{<m}\in J_{<m}$. 
		
		\paragraph{\bf Step~2} For every ${\bm k}=(k_1^{(1)},\ldots, k_{s-2}^{(1)}, k_1^{(2)},\ldots, k_{s-2}^{(2)}, k_1^{(3)},\ldots, k_{s-2}^{(3)})\in \K^{3(s-2)}$, we define the corresponding projections ${\bm k}^{(\iota)} \coloneqq (k_1^{(\iota)},\dots,k_{s-2}^{(\iota)})$ for $\iota\in\{2, 3\}$ and difference operators $\Delta'_{{\bm k}; e_m}\coloneqq \Delta'_{({\bm k}^{(2)}\circ e_m,{\bm k}^{(3)}\circ e_m)}$.
		
		Fix $x_{<m} \in X_{<m}$. By Lemma~\ref{lem:interchange} applied to \eqref{eq:257}, we can estimate
		\[
		\E_{{\bm k}\in[H]_{\K}^{3(s-2)}}^{\lambda_\K^{\otimes 3(s-2)}}
		\ind{\square_{s-2}(H_{x_{<m}})}({\bm k})
		\left|\int_\K G_{x_{<m},{\bm k}}(x_m)e(\psi_{x_{<m}}({\bm k})x_m) d\lambda_\K(x_m) \right|^2
		\]
		from below by a constant multiple of $\delta^{O(1)} N^{2d_m}$, where
		\[
		\psi_{x_{<m}}({\bm k}) \coloneqq \sum_{\omega\in\{0,1\}^{s-2}}(-1)^{|\omega|} \phi_{x_{<m}}(k_1^{(1)},\dots,k_{s-2}^{(1)},k_1^{(\omega_1+2)},\dots,k_{s-2}^{(\omega_{s-2}+2)})
		\]
		and, denoting $x \coloneqq (x_{<m}, x_m)\in\K^m$ for $x_m \in \K$, we define
		$G_{x_{<m},{\bm k}}(x_m)$ by 
		\[
		\E_{t\in[N]_{\K}}^{\lambda_{\K}}
		\Delta'_{{\bm k}; e_m}f_0(x+{\textbf{\textit{P}}}_{m}(t))
		\Big( \prod_{i=1}^{m-1} \Delta'_{{\bm k}; e_m}f_{i}(x-{\textbf{\textit{P}}}_{i}(t)+{\textbf{\textit{P}}}_{m}(t)) 
		\Big) \Delta'_{{\bm k}; e_m}(\xi, {\eta})(x,t).
		\] 
		
		By the popularity principle, for some constant $C_3 \in \Z_+$, if $x_{<m} \in X_{<m}$, then the set $K_{x_{<m}}$, defined as the set of all ${\bm k}\in\square_{s-2}(H_{x_{<m}})$ such that
		\begin{align*}
		\Big| \int_\K G_{x_{<m},{\bm k}}(x_m)e(\psi_{x_{<m}}({\bm k})x_m) d\lambda_\K(x_m) \Big|
		\ge C_3^{-1} \delta^{C_3} N^{d_m},   
		\end{align*}
		is measurable and satisfies $|K_{x_{<m}}|_{\K^{3(s-2)}} \ge C_3^{-1} \delta^{C_3} |[H]_{\K}^{3(s-2)}|_{\K^{3(s-2)}}$.
		
		\paragraph{\bf Step~3} Consider the following measurable set
		\begin{align*}
		Z \coloneqq \left\{(x_{<m}, {\bm k})\in J_{<m}\times \square_{s-2}([H]_\K^{2(s-2)})
		:  x_{<m}\in X_{<m} \text{ and } {\bm k}\in K_{x_{<m}}\right\}
		\end{align*}
		and its sections $Z_{\bm k} \coloneqq \{x_{<m}\in J_{<m} : (x_{<m}, {\bm k})\in Z\}$. By Fubini's theorem, for some constant $C_4 \in \Z_+$, there is a measurable set $K\subseteq \square_{s-2}([H]_\K^{2(s-2)})$ with $|K|_{\K^{3(s-2)}} \ge C_4^{-1} \delta^{C_4} |\square_{s-2}([H]_\K^{2(s-2)})|_{\K^{3(s-2)}}$ such that if ${\bm k}\in K$, then
		\begin{equation}
		\label{eq:C4-k}
		|Z_{{\bm k}}|_{\K^{m-1}} \geq C_4^{-1} \delta^{C_4} |J_{<m}|_{\K^{m-1}}
		\end{equation}
		and $Z_{\bm k}$ has the property that ${\bm k}\in K$ and $x_{<m}\in Z_{{\bm k}}$ together imply ${\bm k}\in K_{x_{<m}}$.
		
		\paragraph{\bf Step~4} Let $C_5 \in \Z_+$ be a large constant to be specified later. By Lemma~\ref{lem:randomphi}, we can find a measurable function $\phi \colon [H]_\K^{2(s-2)}\to\hat\K$ and an exceptional set $E_{\phi} \subseteq \square_{s-2}([H]_\K^{2(s-2)})$ with
		\begin{equation}
		\label{K-except}
		\left|E_{\phi}\right|_{\K^{3(s-2)}} \leq C(s) C_5^2 \delta^{-2C_5} H^{- 1} \left|\square_{s-2}([H]_\K^{2(s-2)})\right|_{\K^{3(s-2)}},
		\end{equation}
		with a constant $C(s)$ depending only on $s$, such that if ${\bm k}\in\square_{s-2}([H]_{\K}^{2(s-2)}) \setminus E_{\phi}$, then $\|q\psi({\bm k})\|_{\hat\K} > C_5 \delta^{-C_5} H^{-1}$ for all $q \in [C_5 \delta^{-C_5}]$ and $\psi({\bm k})$ as in Lemma~\ref{lem:interchange}. In the case $m\ge2$ we adjust the definition of the function $\psi_{x_{<m}}$ by setting ${\psi}^*_{{x_{<m}}}(\bm k) \coloneqq \psi_{{x_{<m}}}(\bm k)$ if $(x_{<m}, {\bm k}) \in Z$ and ${\psi}^*_{{x_{<m}}}(\bm k) \coloneqq \psi(\bm k)$ otherwise. Let $K^*$ be the set of all ${\bm k} \in K$ such that the map $x_{<m} \mapsto {\psi}^*_{{x_{<m}}}(\bm k)$ is measurable. We have $|K^*|_{\K^{3(s-2)}} = |K|_{\K^{3(s-2)}}$. From now on, we will not distinguish between $\psi_{x_{<m}}$ and ${\psi}^*_{x_{<m}}$, and $\psi_{x_{<m}}$ will refer to the adjusted version. Similarly, $K$ will refer to $K^*$. For all ${\bm k}\in K$, by the definition of $K_{x_{<m}}$ and \eqref{eq:C4-k}, we have
		\begin{align*}
		\int_{\K^{m-1}} \Big| \int_{\K} G_{{x_{<m}},{\bm k}}(x_m)e(\psi_{{x_{<m}}}({\bm k})x_m) d\lambda_{\K}(x_m) \Big| d\lambda_{\K^{m-1}}({x_{<m}}) \gtrsim \delta^{O(1)}N^{D_m}
		\end{align*}
		with some implicit constants that are chosen independently of $C_5$.
		
		\paragraph{\bf Step~5} Next, we will show that
		Condition~\ref{con:major} can, essentially, be applied to the inequality above. By duality, we have a lower bound of the same order for
		\begin{align}
		\label{eq:260}
		\int_{\K^{m-1}} \int_{\K} g_0'(x_{<m})G_{{x_{<m}},{\bm k}}(x_m)e(\psi_{{x_{<m}}}({\bm k})x_m) d\lambda_{\K}(x_m)  d\lambda_{\K^{m-1}}({x_{<m}})
		\end{align}
		with some $1$-bounded $g_0'\in L^{\infty}(\K^{m-1})$. Since $\Delta'_{{\bm k}; e_m}(\xi, {\eta})(x,t)$ has the form
		\[
		e\Big(\sum_{\omega\in\{0, 1\}^{s-2}}
		(-1)^{|\omega|}
		\xi\big(x+{\textbf{\textit{P}}}_m(t)+(\omega \cdot {\bm k}^{(2)}) e_m+(({\bf 1}-\omega) \cdot {\bm k}^{(3)})e_m\big)\cdot {\mathcal P}_{(m, l]}(t)\Big),
		\]
		after the change of variables $x_m\mapsto x_m-P_{m}(t)$ we can rewrite \eqref{eq:260} as
		\begin{align*}
		\int_{\K^{m}}
		f_0'(x)
		\E_{t\in[N]_{\K}}^{\lambda_{\K}}
		\Big( \prod_{i=1}^{m-1} \Delta'_{{\bm k}; e_m}f_{i}(x-{\textbf{\textit{P}}}_{i}(t)) 
		\Big) e(\zeta_{[m, l]}'(x)\cdot {\mathcal P}_{[m, l]}(t))
		d\lambda_{\K^{m}}(x),
		\end{align*}
		where $f_0'(x) \coloneqq g_0'(x_{<m})\Delta'_{{\bm k}; e_m}f_0(x)e(\psi_{x_{<m}}({\bm k})x_m)$, and the measurable mapping $\zeta_{[m, l]}' \coloneqq (\zeta_m',\ldots, \zeta_l') \colon \K^m \to \hat{\K}^{l-m}$ is defined by setting $\zeta_m'(x) \coloneqq -\psi_{x_{<m}}({\bm k})$ and 
		\[
		\zeta_{i}'(x) \coloneqq \sum_{\omega\in\{0, 1\}^{s-2}}
		(-1)^{|\omega|}\xi_i\big(x+(\omega \cdot {\bm k}^{(2)}) e_m+(({\bf 1}-\omega) \cdot {\bm k}^{(3)})e_m\big) \quad \text{for} \quad i\in[l]\setminus[m].
		\]
		
		By the pigeonhole principle, we can find $x_m\in I_m$ for which the integral
		\begin{align*}
		\int_{\K^{m-1}}
		g_0(y)
		\E_{t\in[N]_{\K}}^{\lambda_{\K}}
		\Big( \prod_{i=1}^{m-1} g_{i}(y-{\textbf{\textit{P}}}_{i}(t)) 
		\Big) e(\zeta_{[m, l]}(y)\cdot {\mathcal P}_{[m, l]}(t))
		d\lambda_{\K^{m-1}}(y)
		\end{align*}
		admits a lower bound of order $\delta^{O(1)}N^{D_{m-1}}$, and the functions above are the corresponding sections of the previous ones, i.e., we have $g_0(y) \coloneqq f_0'(y, x_m)$ and $g_i(y) \coloneqq \Delta'_{{\bm k}; e_m}f_{i}(y, x_m)$ for $i\in[m-1]$, and $\zeta_i(y) \coloneqq \zeta_i'(y, x_m)$ for $i\in[l]\setminus[m-1]$ so that $\zeta_{[m, l]} \coloneqq (\zeta_m,\ldots, \zeta_l) \colon \K^{m-1} \to \hat{\K}^{l-m+1}$.
		
		Now, the $(m, l)$-major arc condition can be applied to the last integral, yielding that the following holds for some large $c \in \Z_+$ independent of $C_5$. For each ${\bm k}\in K$, we can find a positive integer $q_{{\bm k}} \in [c \delta^{-c}]$ and a measurable set $Y_{<m, {\bm k}}\subseteq J_{<m}$ such that $|Y_{<m, {\bm k}}|_{\K^{m-1}} \geq c^{-1} \delta^{c} |J_{<m}|_{\K^{m-1}}$ and
		\begin{equation*}
		\|q_{{\bm k}} \psi_{x_{<m}}({\bm k})\|_{\hat \K} \le c \delta^{-c} N^{-d_m},
		\qquad x_{<m} \in Y_{<m, {\bm k}}.
		\end{equation*}
		
		Then, by the pigeonhole principle, we find a  constant $C_6\in\Z_+$ (depending on $c$ and $ C_4$, but not on $C_5$), a measurable set $K'\subseteq K$ with $|K'|_{\K^{3(s-2)}} \geq C_6^{-1} \delta^{C_6} |\square_{s-2}([H]_\K^{2(s-2)})|_{\K^{3(s-2)}}$, and $q \in [C_6 \delta^{-C_6}]$ such that
		\begin{align*}
		\|q \psi_{x_{<m}}({\bm k})\|_{\hat \K} \le C_6 \delta^{-C_6} N^{-d_m}, \qquad {\bm k}\in K', \, x_{<m} \in Y_{<m, {\bm k}}.
		\end{align*}
		We can assume, without loss of generality, that
		\begin{align*}
		Y_{<m, {\bm k}}\coloneqq   \left\{x_{<m}\in J_{<m} : \|q \psi_{x_{<m}}({\bm k})\|_{\hat \K} \le C_6 \delta^{-C_6} N^{-d_m}\right\},
		\end{align*}
		and $|Y_{<m, {\bm k}}|_{\K^{m-1}} \geq C_6^{-1} \delta^{C_6} |J_{<m}|_{\K^{m-1}}$. If $\K=\R$, then we simply take $K'=K$, which is possible because $q_{\bm k}=1$ for all ${\bm k}\in K$.
		
		\paragraph{\bf Step~6} We now specify $C_5 > 2C_0C_6$ to conclude that, for every ${\bm k}\in K'$, if $x_{<m} \in J_{<m}\setminus Z_{{\bm k}}$ and $x_{<m} \in Y_{<m, {\bm k}}$, then, necessarily, ${\bm k} \in E_\phi$. By using \eqref{K-except}, together with the assumption that $N \geq C \delta^{-C}$ for a sufficiently large constant $C\in\Z_+$, and by Fubini's theorem, we obtain that, for some constant $C_7\in\Z_+$, there is a measurable set $X_{<m}'\subseteq X_{<m}$ with
		\[
		|X_{<m}'|_{\K^{m-1}} \ge C_7^{-1}\delta^{C_7} |J_{<m}|_{\K^{m-1}}
		\]
		such that if $x_{<m}\in X_{<m}'$, then there is a measurable set
		$U_{x_{<m}} \subseteq K'$ with
		\[
		|U_{x_{<m}}|_{\K^{3(s-2)}} \ge C_7^{-1}\delta^{C_7} \left|\square_{s-2}([H]_\K^{2(s-2)})\right|_{\K^{3(s-2)}}
		\]
		and with the property that $x_{<m}\in X_{<m}'$ and ${\bm k}\in U_{x_{<m}}$ together imply $x_{<m}\in Y_{<m, {\bm k}}\cap Z_{\bm k}$ and ${\bm k}\in K_{x_{<m}}$.
		
		\paragraph{\bf Step~7} Let $M_0 \coloneqq C_6^{-1} C_8 \delta^{C_6-C_8}$
		for a large $C_8 \in \Z_+$ to be specified later. By
		Lemma~\ref{lem:pigeonhole} and the pigeonhole principle, there are
		measurable sets $X_{<m}''\subseteq X_{<m}'$ with $|X_{<m}''|_{\K^{m-1}} \gtrsim \delta^{O(1)} |J_{<m}|_{\K^{m-1}}$ and $V_{x_{<m}}\subseteq U_{x_{<m}}$ with
		\[
		|V_{x_{<m}}|_{\K^{3(s-2)}} \gtrsim \delta^{O(1)} \left|\square_{s-2}([H]_\K^{2(s-2)})\right|_{\K^{3(s-2)}}
		\]
		such that, for some $a, b\in\Z$
		satisfying $|a|+|b| \lesssim \delta^{-O(1)}$, if $x_{<m}\in X_{<m}''$, then
		\[
		\psi_{x_{<m}}({\bm k})=\frac{a}{q}+\frac{b}{M_0 N^{d_m}}
		+\frac{\theta_{x_{<m}}({\bm k})}{M_0 N^{d_m}}
		\quad \text{with} \quad |\theta_{x_{<m}}({\bm k})|\leq 1
		\]
		for all
		${\bm k}\in V_{x_{<m}}$. Here, it is crucial that $\psi_{x_{<m}}({\bm k})$ is the original function from Step~2, since $x_{<m}\in X_{<m}''$ and ${\bm k}\in V_{x_{<m}}$ together imply $(x_{<m}, {\bm k})\in Z$.
		
		Set $\Omega_1 \coloneqq \{ \omega\in\{0,1\}^{s-2} : \omega_1 = 1 \}$. For $x_{<m}\in X_{<m}''$, define $\eta_{x_{<m}}^{(1)}({\bm k})$ by
		\begin{equation*}
		-\sum_{\omega\in \Omega_1}
		(-1)^{|\omega|}\phi_{x_{<m}}\left(k_1^{(1)},\dots,k_{s-2}^{(1)},k_1^{(\omega_1+2)},\dots,k_{s-2}^{(\omega_{s-2}+2)}\right)
		+\frac{a}{q}+\frac{b}{M_0N^{d_m}}
		\end{equation*}
		so that $\eta_{x_{<m}}^{(1)}$ does not depend on $k_1^{(2)}$. Similarly, for each $i\in [s-2]\setminus\{1\}$, set 
		$\Omega_i \coloneqq \{ \omega\in\{0,1\}^{s-2} : \omega_1 = \dots = \omega_{i-1} = 0, \, \omega_i = 1 \}$
		and define $\eta^{(i)}_{x_{<m}}({\bm k})$ by
		\begin{equation*}
		-\sum_{\omega\in \Omega_i} (-1)^{|\omega|}\phi_{x_{<m}}\left(k_1^{(1)},\dots,k_{s-2}^{(1)},k_1^{(\omega_1+2)},\dots,k_{s-2}^{(\omega_{s-2}+2)}\right)
		\end{equation*}
		so that $\eta_{x_{<m}}^{(i)}$ does not depend on $k_i^{(2)}$. Note that if ${\bm k}\in V_{x_{<m}}$, then
		\begin{equation*}
		\bigg| \phi_{x_{<m}}(k_1^{(1)},\dots,k_{s-2}^{(1)},k_1^{(2)},\dots,k_{s-2}^{(2)})
		- \sum_{i\in[s-2]}\eta_{x_{<m}}^{(i)}({\bm k}) \bigg|
		= \bigg| \frac{\theta_{x_{<m}}({\bm k})}{M_0N^{d_m}} \bigg|
		\leq \frac{1}{M_0N^{d_m}}.
		\end{equation*}
		
		\paragraph{\bf Step~8} By the pigeonhole principle, if $x_{<m}\in X_{<m}''$, then there exists a parameter ${\bm h}''_{x_{<m}}\in[H]^{s-2}_\K$ for which the corresponding set
		\begin{equation*}
		H'_{x_{<m}} \coloneqq \{({\bm h},{\bm h}')\in
		H_{x_{<m}}:({\bm h},{\bm h}', {\bm h}''_{x_{<m}})\in V_{x_{<m}}\}
		\end{equation*}
		satisfies $|H_{x_{<m}}'|_{\K^{2(s-2)}}
		\gtrsim \delta^{O(1)} |[H]^{2(s-2)}_\K|_{\K^{2(s-2)}}$. If $C_8\in\Z_+$ is large, then 
		\begin{equation*}
		\E_{({\bm h},{\bm h}')\in H_{x_{<m}}'}^{\lambda_{\K}^{\otimes 2(s-2)}}
		\Big| \int_\K \Delta'_{(\bm h\circ e_m,{\bm h}'\circ e_m)}
		F_{x_{<m}}(x_m)\prod_{i=1}^{s-2}
		e(\eta^{(i)}_{x_{<m}}({\bm h},{\bm h}', {\bm h}''_{x_{<m}})x_m) d\lambda_\K(x_m) \Big|^2
		\end{equation*}
		admits a lower bound of order $\delta^{O(1)} N^{2d_m}$ thanks to \eqref{eq:257}.
		
		Applying Lemma~\ref{lem:GCScor} for each fixed $x_{<m}\in X_{<m}''$ and integrating over all $x_{<m}\in X_{<m}''$, we conclude, by positivity, that $\|F\|_{U^{s-1}_{[H]_\K e_m}(I_{\le m})} \gtrsim \delta^{O(1)}$.
	\end{proof}
	
	\subsubsection{\textbf{The major arc and degree lowering lemmas}}\label{sec:majarc}
	Our goal is now to prove an unconditional version of the degree lowering lemma. To do this, we will need to verify that Condition~\ref{con:major} always holds.
	
	\begin{lemma}[$(m, l)$-major arc lemma]\label{lem:fixedvar}
		Fix $A\in[1, \infty)$, $C_0 \in \Z_+$, $k, l \in\Z_+ \setminus \{1\}$, $m\in \Z_+$ with $m \leq l \leq k$, and $d_1, \dots, d_k \in \Z_+$ with $d_1<\cdots<d_k$. Then there exists a large constant $C \in \Z_+$ depending only on $A,C_0,d_k$ such that the following holds. Assume that $\delta \in (0,1]$ and $N \geq C \delta^{-C}$. Let ${\mathcal P}=(P_1,\ldots, P_k)$ be a~polynomial mapping as in \eqref{eq:42}, where each $P_i\in\K[{\rm t}]$ has degree $d_i$ and is $(d_i,\delta, N)$-admissible with tolerance $A$. Let $g_0, g_1,\ldots, g_{m-1} \in L^\infty(\K^{m-1})$ be $1$-bounded measurable functions supported on $I_{<m}$, where $I_{<1} \coloneqq \{0\} = \K^0$ and $I_{<m} \coloneqq \prod_{i\in[m-1]}[ \pm C_0 \delta^{-C_0} N^{d_i}]_{\K}$ if $m \geq 2$. Given a measurable mapping $\zeta \coloneqq \zeta_{[m, l]} \coloneqq (\zeta_{m},\dots,\zeta_{l}) \colon \K^{m-1}\to\hat\K^{l-m+1}$ and a frequency vector $\theta \coloneqq \theta_{(l, k]} \coloneqq (\theta_{l+1},\dots,\theta_k)\in\hat\K^{k-l}$, define $G_m \in \C$ by
		\begin{align*} 
		\int_{\K^{m-1}} 
		g_0(y)
		\E_{t\in[N]_\K}^{\lambda_\K}
		\Big( \prod_{i=1}^{m-1} g_{i}(y - {\textbf{\textit{P}}}_{i}(t))
		\Big) e(\zeta(y)\cdot {\mathcal P}_{[m, l]}(t)+\theta\cdot {\mathcal P}_{>l}(t)) d\lambda_{\K^{m-1}}(y)
		\end{align*}
		with ${\textbf{\textit{P}}}_{i}(t) \coloneqq P_i(y)e_i$. Set $D_n \coloneqq d_1 + \dots + d_{n}$ for $n\in[k]$ and $D_0 \coloneqq 0$. If
		\begin{align}
		\label{eq:263}
		|G_m| \geq C_0^{-1} \delta^{C_0} N^{D_{m-1}},
		\end{align}
		then we can find an integer $q \in [C\delta^{-C}]$ and a measurable set $Y \subseteq I_{<m}$ such that $|Y|_{\K^{m-1}} \geq C^{-1} \delta^{C} N^{D_{m-1}}$ and, for every $y \in Y$, we have
		\begin{align}
		\label{eq:264}
		\sum_{j=m}^lN^{d_{j}} \|q \zeta_{j}(y) \|_{\hat \K}
		+
		\sum_{j=l+1}^kN^{d_{j}} \|q \theta_{j}\|_{\hat \K} 
		\leq C \delta^{-C}.
		\end{align}
		If $\K = \R$, then we take $q = 1$. 
	\end{lemma}
	
	Observe that inequality \eqref{eq:263} is exactly inequality \eqref{eq:258} with $c_0 = C_0$. Thus, Lemma~\ref{lem:fixedvar} verifies the $(m, l)$-major arc condition, thereby establishing the degree lowering lemma unconditionally; see Lemma~\ref{lem:deglower} below.
	
	\begin{proof}[Proof of Lemma~\ref{lem:fixedvar}]
		In the discussion following Lemma~\ref{prop:cdeglow}, we have already observed that Lemma~\ref{lem:fixedvar} holds for $m=1$ and any $k, l\in\Z_+\setminus\{1\}$ such that $m\le l\le k$. In other words, the $(1, l)$-major arc lemma is true, establishing the base case of an induction argument. Let $k, l\in\Z_+\setminus\{1\}$ be fixed and let $m\in\Z_+$ be such that $m+1\le l$. Assuming that the $(m, l)$-major arc lemma is true, we proceed inductively to show that the $(m+1, l)$-major arc lemma is also true. Let $g_0, g_1,\ldots, g_{m} \in L^\infty(\K^{m})$ be $1$-bounded measurable functions supported on $I_{\le m} \coloneqq \prod_{i\in[m]}[ \pm C_0 \delta^{-C_0} N^{d_i}]_{\K}$. Let $\zeta_{[m+1, l]} \coloneqq (\zeta_{m+1},\dots,\zeta_{l}) \colon \K^{m}\to\hat\K^{l-m}$ be a measurable mapping and $\theta \coloneqq \theta_{(l, k]} \coloneqq (\theta_{l+1},\dots,\theta_k)\in\hat\K^{k-l}$. Suppose that
		\begin{align}
		\label{eq:263'}
		|G_{m+1}|
		\geq C_0^{-1}\delta^{C_0} N^{D_{m}}.
		\end{align}
		We now prove that \eqref{eq:263'} implies condition \eqref{eq:264} with $m+1$ in place of $m$.
		
		\paragraph{\bf Step~1}
		By making the change of variables $y\mapsto y+{\textbf{\textit{P}}}_{m}(t)$ in $G_{m+1}$, applying the Cauchy--Schwarz inequality to double the $t$ variable, and undoing the change of variables, we obtain that 
		\begin{align*}
		\left|
		\int_{\K^{m}} 
		g_0(y)  \E_{t\in[N]_\K}^{\lambda_\K}
		\Big(
		\prod_{i=1}^{m-1} g_{i}(y-{\textbf{\textit{P}}}_{i}(t))
		\Big) 
		\mathcal C F_{m}(y-{\textbf{\textit{P}}}_{m}(t))
		\cdot (\xi, \theta)(y,t)
		d\lambda_{\K^{m}}(y) \right|
		\end{align*}
		admits a lower bound of order $\delta^{O(1)} N^{D_{m}}$, where $F_m(y)$ is given by
		\[
		\E_{t\in[N]_\K}^{\lambda_\K} g_0(y+{\textbf{\textit{P}}}_m(t)) 
		\Big(
		\prod_{i=1}^{m-1} g_{i}(y + {\textbf{\textit{P}}}_m(t) -{\textbf{\textit{P}}}_{i}(t)) 
		\Big) 
		\cdot (\xi, \theta)(y + {\textbf{\textit{P}}}_m(t),t)
		\]
		for $y \in \K^m$, and we also define
		\begin{align*}
		(\xi, \theta)(y,t) \coloneqq
		e\left(\zeta_{(m, l]}(y)\cdot {\mathcal P}_{(m, l]}(t) + \theta_{(l, k]}\cdot {\mathcal P}_{(l, k]}(t) \right)
		\quad \text{for} \quad t\in[N]_{\K}, \, y\in\K^{m}.
		\end{align*}
		
		\paragraph{\bf Step~2}
		We claim that the estimate from Step~1 implies that
		\begin{align}
		\label{eq:261}
		\|\mathcal C F_m\|_{U^{s}_{[H]_\K e_m}(I_{\le m})} \gtrsim \delta^{O(1)}
		\end{align}
		for some $\delta^{O(1)} N^{d_{m}} \lesssim H \lesssim \delta^{-O(1)} N^{d_{m}}$ and an integer $s \geq 2$ bounded in terms of $A,C_0,d_k$. For this purpose, we need to eliminate the factor $(\zeta, \theta)(y,t)$ by applying van der Corput differencing $d_k+1$ times. More precisely, we apply Proposition~\ref{prop:pet} with the conclusion~\eqref{PET-2} for the expression from Step~1 repeatedly $d_k+1$ times and obtain a lower bound of the same order for
		\begin{align*}
		\mathbb{E}_{\bm{h}\in {\bm H}}^{\lambda_{\K}^{\otimes d_k+1}} \Big| \int_{\K^{m}} \E_{t\in[N]_\K}^{\lambda_\K}
		\prod_{i=1}^{m} \prod_{\omega\in\{0,1\}^{d_k+1}}
		\mathcal{C}^{|\omega|} 
		g_{i}'(x - {\textbf{\textit{P}}}_{i}(t + \omega \cdot {\bm h}))
		d\lambda_{\K^{m}}(x) \Big|
		\end{align*}
		with $g_m' \coloneqq \mathcal C F_m$ and $g_i' \coloneqq g_i$ for $i\in[m-1]$, and with some suitably chosen ${\bm H} \subseteq [\pm N]_{\K}^{d_k+1}$ such that $|{\bm H}|_{\K^{d_k+1}} \gtrsim \delta^{O(1)} |[\pm N]_{\K}^{d_k+1}|_{\K^{d_k+1}}$. We observe that, for each $\bm{h} \in \bm{H}$, the family
		\[
		\{P_{i}(\, \cdot \, + \omega \cdot {\bm h}) : i \in [m], \, \omega\in\{0,1\}^{d_k+1}\}
		\]
		consists of $(d_{m}, \delta, N)$-admissible with a new tolerance constant $A'\lesssim_{A,C_0,d_k}1$. 
		
		Now, we can apply Theorem~\ref{thm:PET-K}, or more precisely, repeat the arguments used to derive it. Namely, \eqref{eq:261} follows by first applying \cite[Proposition~4.8]{KKL1} (which, again, also holds for admissible polynomials and when $\mathbb{K}=\mathbb{R}$) for each $\bm{h} \in \bm{H}$ and then, after averaging over $\bm{h}$, either applying \cite[Theorem~1.10]{KKL1} when $\mathbb{K}=\mathbb{Z}$, or performing a change of variables and applying the Gowers--Cauchy--Schwarz inequality \eqref{GCS-in} when $\mathbb{K}=\mathbb{R}$.
		
		\paragraph{\bf Step~3}
		By the induction hypothesis, the $(m,l)$-major arc condition holds. Thus, using \eqref{eq:261} and Lemma~\ref{prop:cdeglow} repeatedly $s-2$ times, we obtain that
		\[
		\|\mathcal C F_m \|_{U^2_{[H]_\K e_m}(I_{\le m})}^4 = |I_{<m}|^{-1}_{\K^{m-1}} \int_{\K^{m-1}} \|\mathcal C F_m(x_{<m},\, \cdot \,)\|_{U^2_{[H]_{\K}}(I_m)}^4 
		d\lambda_{\K^{m-1}}(x_{<m})
		\]
		admits a lower bound of order $\delta^{O(1)}$. By the popularity principle and Proposition~\ref{U2}, we find a
		measurable function $\zeta_m \colon \K^{m-1}\to \hat{\K}$ such that
		\begin{align}
		\label{eq:267}
		\int_{\K^{m-1}}
		\Big|\int_{\K}F_m(x)e(\zeta_m(x_{<m})x_m)d\lambda_{\K}(x_m)\Big|d\lambda_{\K^{m-1}}(x_{<m})\gtrsim \delta^{O(1)}N^{D_m}.
		\end{align}
		Dualizing the left-hand side of \eqref{eq:267}, expanding $F_m$, and making the change of variables $x\mapsto x-{\textbf{\textit{P}}}_{m}(t)$, we obtain that \eqref{eq:267} can be rewritten as
		\begin{align*}
		\int_{\K^{m}}
		f_0'(x) \E_{t\in[N]_\K}^{\lambda_\K}
		\Big(
		\prod_{i=1}^{m-1}f'_{i}(x-{\textbf{\textit{P}}}_{i}(t))
		\Big) \cdot (\xi, \theta)(x,t)
		d\lambda_{\K^{m}}(x) \gtrsim \delta^{O(1)}N^{D_m},
		\end{align*}
		where $f_0'(x) \coloneqq g_0(x)g_0''(x_{<m})e(\zeta_{m}(x_{<m})x_m)$ for some $1$-bounded function $g_0''\in L^{\infty}(\K^{m-1})$ such that the integrand above is nonnegative, $f_i' \coloneqq g_i$ for each $i \in [m-1]$, and $\xi \coloneqq \xi_{[m, l]} \coloneqq (\xi_{m},\ldots, \xi_{l}) \colon \K^{m}\to \hat{\K}^{l-m+1}$ with $\xi_{m}(x) \coloneqq -\zeta_m(x_{<m})$ and $\xi_i(x) \coloneqq \zeta_i(x)$ for $i\in[l]\setminus[m]$ and $x\in\K^m$.
		
		By the popularity principle applied to the last integral, there exists a measurable set $Y_m\subseteq I_m$ such that $|Y_m|_{\K}\gtrsim \delta^{O(1)}|I_m|_{\K}$ and, setting
		\begin{align*}
		h_i(x) & \coloneqq f'_i(x, x_m) \phantom{\xi'_{i}} \ \text{for} \quad
		x\in\K^{m-1}, \, i\in\N_{<m}, \\
		\zeta_i'(x) & \coloneqq \xi_{i}(x, x_m) \phantom{f'_{i}} \ \text{for} \quad
		x\in\K^{m-1}, \, i\in\N_{\le l}\setminus\N_{<m},
		\end{align*}
		for every $x_m\in Y_m$, we obtain that the integral
		\begin{align*}
		\int_{\K^{m-1}}
		h_0(x) \E_{t\in[N]_\K}^{\lambda_\K}
		\Big(
		\prod_{i=1}^{m-1}h_{i}(x-{\textbf{\textit{P}}}_{i}(t))
		\Big) \cdot (\zeta', \theta)(x,t)
		d\lambda_{\K^{m-1}}(x)
		\end{align*}
		admits a lower bound of order $\delta^{O(1)}N^{D_{m-1}}$. By the induction hypothesis, this estimate implies that, for every $x_m\in Y_m$, we can find a positive integer $q_{x_m} \lesssim \delta^{-O(1)}$ and a measurable set $Y_{x_m} \subseteq I_{<m}$ satisfying $|Y_{x_m}|_{\K^{m-1}} \gtrsim \delta^{O(1)} N^{D_{m-1}}$ such that, for every $y \in Y_{x_m}$, we have
		\begin{align}
		\label{eq:268}
		\sum_{j=m}^lN^{d_{j}} \|q_{x_m} \zeta'_{j}(y) \|_{\hat \K}
		+
		\sum_{j=l+1}^kN^{d_{j}} \|q_{x_m} \theta_{j}\|_{\hat \K} 
		\lesssim \delta^{-O(1)}.
		\end{align}
		
		By the pigeonhole principle, we can find a positive integer $q \lesssim \delta^{-O(1)}$ and a measurable set $Y_m'\subseteq Y_m$ satisfying $|Y_m'|_{\K}\gtrsim \delta^{O(1)}|I_m|_{\K}$ for which the following holds: for every $x_m\in Y_m'$ there exists a measurable set $Y_{x_m} \subseteq I_{<m}$ satisfying $|Y_{x_m}|_{\K^{m-1}} \gtrsim \delta^{O(1)} N^{D_{m-1}}$ such that for every $y \in Y_{x_m}$ inequality \eqref{eq:268} holds with $q$ in place of $q_{x_m}$. If $\K=\R$, then we simply take $Y_m'=Y_m$, which is possible, since $q_{x_m}=1$ for all $x_m\in Y_m$. Moreover, we can assume, without loss of generality, that
		\begin{align*}
		Y_{x_m}\coloneqq \Big \{y\in I_{<m} : \sum_{j=m+1}^lN^{d_{j}} \|q \zeta_{j}(y) \|_{\hat \K}
		+
		\sum_{j=l+1}^kN^{d_{j}} \|q\theta_{j}\|_{\hat \K} 
		\le C_* \delta^{-C_*} \Big\}
		\end{align*}
		for some $1\le q\lesssim \delta^{-O(1)}$ and $C_* \in \Z_+$, and for all $x_m\in Y_{m}'$. Then the set $Y_{x_m}$ is measurable and  \eqref{eq:264} holds with $m+1$ in place of $m$ by taking
		\[
		Y \coloneqq \{(y, x_m)\in I_{<m}\times I_m: x_m\in Y_m' \text{ and } y\in Y_{x_m}\},
		\]
		which is measurable and satisfies $|Y|_{\K^m}\gtrsim \delta^{O(1)}N^{D_m}$ by Fubini's theorem.
	\end{proof}
	
	Having proved the $(m, l)$-major arc lemma, we can now establish the unconditional degree lowering lemma. For $l\in [k]$, a frequency vector $\xi \coloneqq \xi_{>l} \coloneqq (\xi_{l+1},\dots,\xi_k)\in \hat{\K}^{k-l}$, and functions $f_1, \ldots, f_l\in L^\infty(\K^k, \lambda_{\K^k})$ we define a modulated variant of the multilinear average from \eqref{eq:49} by setting
	\begin{align}
	\label{eq:245}
	A_{N, \xi; \mathbb K^k}^{{\mathcal P}_{\le l}, {\mathcal P}_{>l}}(f_1,\ldots, f_l)(x)
	\coloneqq \E_{t \in [N]_{\mathbb K}}^{\lambda_{\mathbb K}} 
	\Big( \prod_{i\in[l]} f_i(x- {\textbf{\textit{P}}}_i(t)) \Big) e(\xi\cdot {\mathcal P}_{>l}(t))
	\end{align}
	for $x \in \K^k$, with ${\textbf{\textit{P}}}_i(t) \coloneqq P_i(t)e_i$. From \eqref{eq:245} we have $A_{N, \xi; \mathbb K^k}^{{\mathcal P}_{\le k}, {\mathcal P}_{>k}}=A_{N;\mathbb K^k}^{{\mathcal P}}$ for any $\xi\in\K^{k}$. The adjoint operators $A_{N, \xi; \mathbb K^k}^{{\mathcal P}_{\le l}, {\mathcal P}_{>l}, *j}$ for $j\in[l]$ are given by
	\begin{align}
	\label{eq:246}
	A_{N, \xi; \mathbb K^k}^{{\mathcal P}_{\le l}, {\mathcal P}_{>l}, *j}(g_1,\ldots,g_l)(x) \coloneqq 
	\E_{t \in [N]_{\mathbb K}}^{\lambda_{\mathbb K}}  
	\Big( \prod_{i\in[l]}
	g^{*j}_{i}(x- {\textbf{\textit{P}}}^{*j}_i(t)) \Big) \overline{e(\xi\cdot {\mathcal P}_{>l}(t))}
	\end{align}
	for $x \in \K^k$, with $g^{*j}_{i} \coloneqq \mathcal C^{\ind{i \neq j}} g_{i}$ and ${\textbf{\textit{P}}}^{*j}_i(t) \coloneqq \ind{i \neq j} P_{i}(t)e_{i} - P_j(t)e_j$; see \eqref{transpose} and \eqref{eq:108a} for analogous formulae for $A_{N;\mathbb K^k}^{{\mathcal P}, *j}$.
	
	Our concluding result here is as follows.
	
	\begin{lemma}[Degree lowering lemma]\label{lem:deglower}
		Fix an integer $s\geq 3$, $A\in[1, \infty)$, $C_0 \in \Z_+$, $k, l \in\Z_+ \setminus \{1\}$ with $l \leq k$, and $d_1, \dots, d_k \in \Z_+$ with $d_1<\cdots<d_k$. Then there exists a large constant $C \in \Z_+$ depending only on $s,A,C_0,d_k$ such that the following holds.  Assume that $\delta \in (0,1]$ and $N \geq C \delta^{-C}$. Let ${\mathcal P}=(P_1,\ldots, P_k)$ be a polynomial mapping as in \eqref{eq:42}, where each $P_i\in\K[{\rm t}]$ has degree $d_i$ and is $(d_i,\delta, N)$-admissible with tolerance $A$. Assume that $C_0^{-1} \delta^{C_0} N^{d_l} \leq H \leq C_0 \delta^{-C_0} N^{d_l}$. Let $f_0, f_1,\dots,f_{l-1} \in L^{\infty}(\K^k)$ be $1$-bounded measurable functions supported on $I \coloneqq \prod_{i \in [k]} [\pm C_0 \delta^{-C_0} N^{d_i}]$. For a frequency vector $\eta \coloneqq (\eta_{l+1},\dots,\eta_k)\in \hat\K^{k-l}$, define
		\begin{align*}
		F^{\eta}(x) \coloneqq A_{N, \eta; \mathbb K^k}^{{\mathcal P}_{\le l}, {\mathcal P}_{>l}, *l}(f_1,\ldots,f_{l-1}, f_0)(x) \quad \text{for} \quad x\in\K^k
		\end{align*}
		with $A_{N, \eta; \mathbb K^k}^{{\mathcal P}_{\le l}, {\mathcal P}_{>l}, *l}$ from
		\eqref{eq:246}. Then the following implication holds:
		\begin{equation*}
		\|F^{\eta}\|_{U^s_{[H]_\K e_l}(I)}\geq C_0^{-1} \delta^{C_0}
		\implies
		\|F^{\eta}\|_{U^{s-1}_{[H]_\K e_l}(I)} \geq C^{-1} \delta^C.
		\end{equation*}
	\end{lemma}
	
	\begin{proof}
		We first apply the popularity principle to reduce the matter from $\K^k$ to
		$\K^l$, and then appeal to the $(l,l)$-major arc lemma and Lemma~\ref{prop:cdeglow}.
	\end{proof}
	
	\subsubsection{\textbf{All together: proof of Theorem~\ref{ppp}}}
	Finally, we can complete the proof of Theorem~\ref{ppp}. We begin with the following intermediate result.
	
	\begin{lemma}\label{lem:replace}
		Fix $A\in[1, \infty)$, $C_0 \in \Z_+$, $k\in\Z_+ \setminus \{1\}$, and $d_1, \dots, d_k \in \Z_+$ with $d_1<\cdots<d_k$. Then there exists a large constant $C \in \Z_+$ depending only on $A,C_0,d_k$ such that the following holds. Assume that $\delta \in (0,1]$ and $N \geq C \delta^{-C}$. Let ${\mathcal P}=(P_1,\ldots, P_k)$ be a polynomial mapping as in \eqref{eq:42}, where each $P_i\in\K[{\rm t}]$ has degree $d_i$ and is $(d_i,\delta, N)$-admissible with tolerance $A$. Let $f_0, f_1,\ldots, f_{k} \in L^\infty(\K^{k})$ be $1$-bounded measurable functions supported on $I \coloneqq \prod_{i\in[k]}[ \pm C_0 \delta^{-C_0} N^{d_i}]_{\K}$. If condition \eqref{eq:186} from Theorem~\ref{ppp} holds, then 
		\begin{align}
		\label{eq:271}
		\left|\left\langle A_{N, \zeta_{>1}; \mathbb K^{k}}^{{\mathcal P}_{\le 1}, {\mathcal P}_{>1}}(f_1), f_0' \right\rangle\right|\geq C^{-1} \delta^C N^{D_{k}}
		\end{align}
		for some $1$-bounded measurable function $f_0' \in L^\infty(\K^{k})$ satisfying $|f_0'|=|f_0|$ and some frequency vector $\zeta_{>1} \coloneqq (\zeta_{2},\dots,\zeta_k)\in\hat\K^{k-1}$.
	\end{lemma}
	
	\begin{proof}
		Set $D_l \coloneqq d_1 + \dots + d_{ l}$ for all $l \in [k]$. Given $l\in[k]$, we will prove that there exists a frequency vector $\zeta_{>l} \coloneqq (\zeta_{l +1},\dots,\zeta_k)\in\hat\K^{k-l}$ such that
		\begin{align}
		\label{eq:1}
		\left|\left\langle A_{N, \zeta_{>l}; \mathbb K^k}^{{\mathcal P}_{\le l}, {\mathcal P}_{>l}}(f_1,\ldots,f_l), f_{0, l} \right\rangle\right|\gtrsim \delta^{O(1)} N^{D_{k}},
		\end{align}
		for some $f_{0, l}\in L^{\infty}(\K^k)$ such that $|f_{0, l}|=|f_0|$. We prove \eqref{eq:1} using backward induction on $l\in[k]$. The base case $l=k$, with $\zeta_{>k}$ and $\mathcal P_{>k}$ being empty vectors, follows from \eqref{eq:186} because $A_{N, \zeta_{>k}; \mathbb K^k}^{{\mathcal P}_{\le k}, {\mathcal P}_{>k}}=A_{N;\mathbb K^k}^{{\mathcal P}}$. We now assume that \eqref{eq:1} holds for some $l\in[k]\setminus\{1\}$ and verify \eqref{eq:1} for $l-1$ in place of $l$. Once the final case $l=1$ is verified, Lemma~\ref{lem:replace} follows. 
		
		\paragraph{\bf Step~1} Define $F_l(x) \coloneqq A_{N, \zeta_{>l}; \mathbb K^k}^{{\mathcal P}_{\le l}, {\mathcal P}_{>l}, *l}(f_1,\ldots,f_{l-1}, f_{0, l})(x)$. By duality and the Cauchy--Schwarz inequality, we may replace $f_l$ with $F_l$ in \eqref{eq:1} and obtain
		\begin{align*}
		\left|\left\langle A_{N, \zeta_{>l}; \mathbb K^k}^{{\mathcal P}_{\le l}, {\mathcal P}_{>l}}(f_1,\ldots,f_{l-1}, F_l), f_{0, l} \right\rangle\right|\gtrsim \delta^{O(1)} N^{D_{ k }}.
		\end{align*}
		As in Step~2 of the proof of Lemma~\ref{lem:fixedvar}, we obtain $\|F_l \|_{U^s_{[H]_\K e_l}(I)}\gtrsim \delta^{O(1)}$ for some suitably bounded integer $s \geq 2$ and $\delta^{O(1)} N^{d_l} \lesssim H \lesssim \delta^{-O(1)} N^{d_l}$. Then 
		\[
		\|F_l \|_{U^2_{[H]_\K e_l}(I)}^4 = |I_{\{l\}^c}|^{-1}_{\K^{k-1}} \int_{\K^{k-1}} \|F_l(x_{\{l\}^c}, \, \cdot \,)\|_{U^2_{[H]_{\K}}(I_l)}^4 d\lambda_{\K^{k-1}}(x_{\{l\}^c}) \gtrsim \delta^{O(1)}
		\]
		by Lemma~\ref{lem:deglower}, where we define $x_{\{l\}^c} \coloneqq (x_1,\ldots, x_{l-1}, x_{l+1},\ldots, x_k)\in \K^{k-1}$ and identify
		$(x_{\{l\}^c},x_l)$ with $(x_1,\dots,x_k)\in\K^k$, and similarly define $I_{\{l\}^c} \coloneqq \prod_{i\in[k]\setminus\{l\}}[ \pm C_0 \delta^{-C_0} N^{d_i}]_{\K}$ and $I_l \coloneqq [ \pm C_0 \delta^{-C_0} N^{d_l}]_{\K}$. We will also write $x_{\{l\}^c} = (x_{<l}, x_{>l})$ and $I_{\{l\}^c} = I_{<l} \times I_{>l}$ for $x_{<l}, x_{>l}, I_{<l}, I_{>l}$ defined analogously.
		
		By the popularity principle and Proposition~\ref{U2}, we find a measurable function $\xi \colon \K^{k-1}\to \hat{\K}$ and a measurable set $X_{>l}\subseteq I_{>l}$ satisfying $|X_{>l}|_{\K^{k-l}}\gtrsim \delta^{O(1)}|I_{>l}|_{\K^{k-l}}$ such that, for every $x_{>l}\in X_{>l}$, there exists a measurable set $X_{<l}^{x_{>l}}\subseteq I_{<l}$ satisfying $|X_{<l}^{x_{>l}}|_{\K^{l-1}}\gtrsim \delta^{O(1)}|I_{<l}|_{\K^{l-1}}$ with the following property: whenever $x_{>l}\in X_{>l}$ and $x_{<l}\in X_{<l}^{x_{>l}}$, we have
		\begin{align}
		\label{eq:274}
		\Big |\int_{\K}F_l(x_{<l},x_l, x_{>l})e(\xi(x_{<l}, x_{>l})x_l)d\lambda_{\K}(x_l)\Big|
		\gtrsim \delta^{O(1)}N^{d_l}.
		\end{align}
		
		\paragraph{\bf Step~2} For some large $C_1 \in \Z_+$ to be specified later, we can find a frequency $\xi_0 \in \hat \K$ such that $ \|q\xi_0\|_{\hat \K} > C_1 \delta^{-C_1} N^{-d_l}$ for all $q \in [C_1 \delta^{-C_1}]$. Indeed, set
		\begin{equation*}
		\mathfrak{M} \coloneqq \left\{\zeta \in \hat \K :\|q \zeta\|_{\hat \K} \leq C_1 \delta^{-C_1} N^{-d_l} \text{ for some positive integer }q\leq C_1 \delta^{-C_1}\right\}.
		\end{equation*}
		The measure of $\mathfrak{M}$ does not exceed $10 C_1^2\delta^{-2C_1}N^{-d_l}$, which is less than $1/2$ provided that $N \geq C \delta^{-C}$ and $C$ is sufficiently large in terms of $C_1$. Let $\xi_0$ be any element of $\hat \K \setminus \mathfrak M$. We define $\xi_l \colon \K^{k-1} \to \hat \K$ by $\xi_ l (x_{<l}, x_{>l}) \coloneqq - \xi (x_{<l}, x_{>l})$ if $x_{>l}\in X_{>l}$ and $x_{<l}\in X_{<l}^{x_{>l}}$, and by $\xi_l (x_{<l}, x_{>l}) \coloneqq - \xi_0$ otherwise. Then, with this measurable
		function $\xi_l$, we obtain for every $x_{>l}\in X_{>l}$ that
		\begin{align*}
		\int_{\K^{l-1}}\left|\int_{\K}F_l(x_{<l},x_l, x_{>l})e(-\xi_l(x_{<l}, x_{>l})x_l)d\lambda_{\K}(x_l)\right|d\lambda_{\K^{l-1}}(x_{<l})
		\gtrsim \delta^{O(1)}N^{D_{l}}.
		\end{align*}
		
		\paragraph{\bf Step~3} For each $x_{>l}\in X_{>l}$, arguing as in Step~3 of the proof of Lemma~\ref{lem:fixedvar}, we find a function $g_0 \in L^\infty(\K^{k})$ with $|g_0|=|f_{0, l}|=|f_0|$ and, by the pigeonhole principle, an element $x_l\in I_l$ such that the integral
		\begin{align*}
		\int_{\K^{l-1}}
		g_0(x)\E_{t\in[N]_{\K}}^{\lambda_{\K}}
		\Big( 
		\prod_{i=1}^{l-1}f_{i}(x-{\textbf{\textit{P}}}_{i}(t))
		\Big) e((\xi_l(x_{<l}, x_{>l}), \zeta_{>l})\cdot {\mathcal P}_{\ge l}(t))
		d\lambda_{\K^{l-1}}(x_{<l})
		\end{align*}
		admits a lower bound of order $\delta^{O(1)}N^{D_{l-1}}$.
		
		Applying the $(l, l)$-major arc lemma and the pigeonhole principle, we find a positive integer $q\lesssim \delta^{-O(1)}$ (if $\K=\R$, then we take $q=1$) and a measurable set $Y^{x_{>l}}_{<l} \subseteq I_{<l}$ such that $|Y^{x_{>l}}_{<l}|_{\K^{l-1}}\gtrsim \delta^{O(1)}|I_{<l}|_{\K^{l-1}}$ and, if $x_{<l}\in Y^{x_{>l}}_{<l}$, then $\|q\xi_l(x_{<l}, x_{>l})\|_{\hat{\K}}\le C_2\delta^{-C_2}N^{-d_l}$ for some absolute constant $C_2\in\Z_+$ independent of $C_1$.  If $C_1>C_2$, then
		\[
		|Y^{x_{>l}}_{<l} \cap X_{<l}^{x_{>l}}|_{\K^{l-1}}=|Y^{x_{>l}}_{<l}|_{\K^{l-1}}\gtrsim \delta^{O(1)}|I_{<l}|_{\K^{l-1}}
		\]
		by the definition of $\xi_l$ and our choice of $\xi_0\in \hat \K \setminus \mathfrak M$. Let $M_0 \coloneqq C_2^{-1} C_3 \delta^{C_2-C_3}$ for a large constant $C_3 \in \Z_+$ to be specified later. We use Lemma~\ref{lem:pigeonhole}, followed by the pigeonhole principle, to find a frequency $\zeta_{l}\in\hat{\K}$ and measurable sets $Z_{>l} \subseteq X_{>l}$ satisfying $|Z_{>l}|_{\K^{l-1}}\gtrsim \delta^{O(1)}|I_{> l}|_{\K^{l-1}}$ and $Z^{x_{>l}}_{<l} \subseteq Y^{x_{>l}}_{<l}$ satisfying $|Z^{x_{>l}}_{<l}|_{\K^{l-1}}\gtrsim \delta^{O(1)}|I_{<l}|_{\K^{l-1}}$ for every $x_{>l}\in Z_{>l}$ and such that 
		\begin{align}
		\label{eq:275}
		|\xi_l(x_{<l}, x_{>l})-\zeta_l|\le M_0^{-1} N^{-d_l}.
		\end{align}
		Such a  choice  of $\zeta_l\in\hat{\K}$ is possible, since the function $\xi_l\colon \K^{k-1}\to \hat{\K}$ is measurable and we have at most $O(\delta^{-O(1)})$ possibilities to choose $\zeta_l\in\hat{\K}$, which is clear from Lemma~\ref{lem:pigeonhole}. 
		Now inequality \eqref{eq:274} implies that
		\begin{align*}
		\Big|\int_{\K}
		g_0'(x)\E_{t\in[N]_{\K}}^{\lambda_{\K}}
		\Big( 
		\prod_{i=1}^{l-1}f_{i}(x-{\textbf{\textit{P}}}_{i}(t))
		\Big) e((\xi_l(x_{<l},x_{>l}),\zeta_{>l})\cdot {\mathcal P}_{\ge l}(t))
		d\lambda_{\K}(x_{l})\Big|
		\end{align*}
		admits a lower bound of order $\delta^{O(1)}N^{d_l}$, where $g_0'\in L^{\infty}(\K^k)$ is a $1$-bounded function satisfying $|g'_0|=|f_0|$. If $C_3$ is sufficiently large, then, using \eqref{eq:275}, we can replace $\xi_l(x_{<l},x_{>l})$ with $\zeta_l$ in the integral above and obtain \eqref{eq:1} with $l-1$ in place of $l$ and with a $1$-bounded function $f_{0, l-1}\in L^{\infty}(\K^k)$ such that $|f_{0, l-1}|=|f_{0, l}|=|f_0|$. This completes the proof of the lemma.
	\end{proof}
	
	Now, we can finally prove Theorem~\ref{ppp}. 
	
	\begin{proof}[Proof of Theorem~\ref{ppp}]
		Recall that $D \coloneqq d_1+\dots+d_k$. For $k > 1$, we apply Lemma~\ref{lem:replace} to find a frequency vector $\zeta_{>1} \in {\hat \K}^{k-1}$ such that \eqref{eq:271} holds. By the popularity principle, we find a measurable set $X\subseteq I_{>1}$ satisfying $|X|_{\K^{k-1}}\gtrsim \delta^{O(1)}|I_{>1}|_{\K^{k-1}}$ and such that for every $x_{>1}\in X$ we have
		\begin{align}
		\label{eq:276}
		\Big| \int_{\K} g_{x_{>1}}(x)
		\E_{t\in[N]_\K}^{\lambda_\K}
		f'(x-P_1(t))
		e(\zeta_{>1} \cdot {\mathcal P}_{>1}(t))
		d\lambda_{\K}(x) \Bigr| \gtrsim\delta^{O(1)}N^{d_1},
		\end{align}
		whenever $N\gtrsim \delta^{-O(1)}$, where $f'(x) \coloneqq f_1(x, x_{>1})$ and $g_{x_{>1}} \in L^{\infty}(\K)$ satisfies $|g_{x_{>1}}| = |f_0|$. For $k=1$, we have $e(\zeta_{>1} \cdot {\mathcal P}_{>1}(t)) \coloneqq 1$ above and so \eqref{eq:276} holds. We will proceed in a few steps, now fixing $x_{>1} \in X$ and assuming that $N\ge C_* \delta^{-C_*}$ for some large $C_*\in\Z_+$.
		
		\paragraph{\bf Step~1}
		Applying the Cauchy--Schwarz inequality to \eqref{eq:276}, we obtain
		\begin{align}
		\label{eq:277}
		\|F_{x_{>1}}\|_{L^2(\K)} \gtrsim\delta^{O(1)}N^{d_1/2},
		\end{align}
		with $F_{x_{>1}}(x) \coloneqq \E_{t\in[N]_{\K}}^{\lambda_\K} f'(x-P_1(t)) e(\zeta_{>1} \cdot {\mathcal P}_{>1}(t))$, since $\|g_{x_{>1}}\|_{L^2(\K)}\lesssim N^{d_1/2}$. By \eqref{eq:277} and Plancherel's theorem, we obtain
		\begin{align}
		\label{eq:278}
		\delta^{O(1)}N^{d_1}\lesssim \|F_{x_{>1}}\|_{L^2(\K)}^2 =
		\|\mathcal F_{\K}(F_{x_{>1}})\|_{L^2(\hat{\K})}^2,
		\end{align}
		where, upon making the change of variables $x\mapsto x+P_1(t)$, we have
		\begin{align*}
		\mathcal F_{\K}(F_{x_{>1}})(\xi)=\mathcal F_{1, \K^k}f_1(\xi, x_{>1})
		\E_{t\in[N]_{\K}}^{\lambda_\K} e(\xi P_1(t)+\zeta_{>1} \cdot {\mathcal P}_{>1}(t)).
		\end{align*}
		
		\paragraph{\bf Step~2} Let $M\in\Z_+$ be a large constant to be specified later. By Proposition~\ref{lem:Weyl1}, there exists a large constant $C_1 \in \Z_+$ depending only on $A$ and $d_k$ such that if $N\ge C_1M^{C_1}\delta^{-C_1M}$ and
		\begin{align}
		\label{eq:280}
		\left| \E_{t\in[N]_{\K}}^{\lambda_\K} e(\xi P_1(t)+\zeta_{>1} \cdot {\mathcal P}_{>1}(t))\right| \ge M^{-1}\delta^M,
		\end{align}
		then there exists an integer $q\in[C_1 M^{C_1}\delta^{-C_1 M}]$ (if $\K=\R$, then we take $q=1$) such that
		\[
		\|q\xi\|_{\hat{\K}} \le C_1 M^{C_1}\delta^{-C_1 M}N^{-d_1}.
		\]
		Here, it is crucial that the polynomials $P_1,\ldots, P_k$ have distinct degrees. 
		
		Set $C' \coloneqq 4C_1 M^{C_1}$. Fix $M_1 \coloneqq C' \delta^{-C'}$, $M_2 \coloneqq C' \delta^{-C'} N^{-d_1}$, and a smooth and even function $\eta \colon \R\to[0, 1]$ satisfying \eqref{eq:126}, and set $\eta_{[\le M_2]}(\xi) \coloneqq \eta(M_2^{-1}\xi)$ for $\xi\in \R$. In order to distinguish between major and minor arcs, let us define 
		\[
		\Xi_{\K}(\xi) \coloneqq \sum_{\theta \in \mathcal R_{\le M_1}^{\K}}\eta_{[\le M_2]}(\xi-\theta)^{1/2}
		\]
		with $\mathcal R_{\le M_1}^{\Z} \coloneqq \mathcal R_{\le M_1}$ and $\mathcal R_{\le M_1}^{\R} \coloneqq \{0\}$, where $\mathcal R_{\le M_1}$ is given in \eqref{eq:129}. If $\K=\Z$, then $\Xi_{\K}^2$ is the Fourier transform of the Ionescu--Wainger projection from \eqref{eq:132}, since $\Xi_{\K}(\xi)^2 = \sum_{\theta \in \mathcal R_{\le M_1}^{\K}}\eta_{[\le M_2]}(\xi-\theta)$ by the disjointness of the supports of $\eta_{[\le M_2]}(\xi-\theta)^{1/2}$ when $\theta \in \mathcal R_{\le M_1}^{\K}$. To ensure disjointness, we used the condition $N\ge C_* \delta^{-C_*}$, assuming that $C_*\in\Z_+$ is large with respect to $C'$.
		
		\paragraph{\bf Step~3} Using $\Xi_{\K}$, we can write $\mathcal F_{\K}(F_{x_{>1}})(\xi)=G_1(\xi, x_{>1})+G_2(\xi, x_{>1})$, where
		\begin{align*}
		G_1(\xi, x_{>1})& \coloneqq \mathcal F_{1, \K^k}f_1(\xi, x_{>1}) \cdot
		\E_{t\in[N]_{\K}}^{\lambda_\K} e(\xi P_1(t)+\zeta_{>1} \cdot {\mathcal P}_{>1}(t)) \cdot (1-\Xi_{\K}(\xi)),\\
		G_2(\xi, x_{>1})& \coloneqq \mathcal F_{1, \K^k}f_1(\xi, x_{>1}) \cdot
		\E_{t\in[N]_{\K}}^{\lambda_\K} e(\xi P_1(t)+\zeta_{>1} \cdot {\mathcal P}_{>1}(t)) \cdot \Xi_{\K}(\xi).
		\end{align*}
		By \eqref{eq:280}, we have $|G_1(\xi, x_{>1})|\le M^{-1}\delta^{M}|\mathcal F_{1, \K^k}f_1(\xi, x_{>1})|$. Thus, we have
		\begin{align*}
		\|G_1(\, \cdot \,, x_{>1})\|_{L^2(\hat{\K})}^2
		\lesssim M^{-1}\delta^{M} N^{d_1}
		\end{align*}
		by Plancherel's theorem, since $f_1(\, \cdot \,, x_{>1})$ is $1$-bounded and supported on a set of size controlled by $N^{d_1}$. If $M$ is large, then the above bound and \eqref{eq:278} yield
		\begin{align*}
		\|\mathcal F_{1, \K^k}f_1(\, \cdot \,, x_{>1})  \cdot \Xi_{\K}(\, \cdot \,)\|_{L^2(\hat{\K})}^2\gtrsim\delta^{O(1)}N^{d_1}
		\end{align*}
		by $|G_2(\xi, x_{>1})|\le |\mathcal F_{1, \K^k}f_1(\xi, x_{>1}) \cdot \Xi_{\K}(\xi)|$. Integrating over $X$, we obtain that
		\begin{align*}
		\left|\int_{\K^{k-1}}\int_{{\hat \K}} \mathcal F_{1, \K^k}(f_1)(\xi, x_{>1}) \overline{\mathcal F_{1, \K^k}(f_1)(\xi, x_{>1})\Xi_{\K}^2(\xi)} d\lambda_{{\hat \K}}(\xi) d\lambda_{\K^{k-1}}(x_{>1})\right|
		\end{align*}
		admits a lower bound of order $\delta^{O(1)}N^{D}$. By Plancherel's theorem, we have
		\begin{align*}
		\left|\left\langle f_1, \Pi_{\K}^1 [\le M_1,
		\le M_2]f_1\right\rangle\right|
		\gtrsim \delta^{O(1)} N^D.
		\end{align*}
		This establishes \eqref{eq:197}
		and completes the proof of Theorem~\ref{ppp}. 
	\end{proof}
	
	\section{Multilinear \texorpdfstring{$L^p$}{TEXT}-improving inequalities}
	\label{sec:improving}
	In this section, we continue to use the notation introduced in Section~\ref{sec:inverse}. Our objective is to establish multilinear $L^p$-improving inequalities on $\Z^k$ and $\R^k$ for the multilinear averages
	\[
	A_{N; \mathbb K^k}^{\mathcal P}(f_1,\ldots, f_k)(x)
	\coloneqq \E_{t \in [N]_{\mathbb K}}^{\lambda_{\mathbb K}} \prod_{i\in[k]} f_i(x-P_i(t)e_i),
	\qquad x\in\mathbb K^k,
	\]
	defined as in \eqref{eq:49} for any polynomial mapping ${\mathcal P} = (P_1, \dots, P_k) \colon \mathbb K \to \mathbb K^k$ with $\deg(P_i) \in \Z_+$ for $i \in [k]$. In contrast to the previous section, the assumption that the degrees of ${\mathcal P}$ are distinct is not needed for the results in this section.
	
	We will be interested in the $L^p$-improving bounds
	\begin{equation}\label{main-bound}
	\left\| A_{N; \mathbb K^k}^{\mathcal P}(f_1,\ldots, f_k) \right\|_{L^q(\mathbb K^k)} \lesssim_{\mathcal P} N^{-D(\frac{1}{p}-\frac{1}{q})} 
	\prod_{i \in [k]} \|f_{i}\|_{L^{p_i}(\mathbb K^k)}
	\end{equation}
	when $p<q$ with $\frac1p = \frac{1}{p_1} + \cdots + \frac{1}{p_k}$, and $D = d_1 + \cdots + d_k$ with $d_i = \deg(P_i)$ for $i \in [k]$. These bounds will apply equally well to the truncated averages \eqref{eq:52}.
	
	Our main result here is the following multilinear $L^p$-improving inequality.
	
	\begin{theorem}\label{thm:improv} Fix $k\in \N_{\geq 2}$ and $j \in [k]$. Then we have \eqref{main-bound} for some $1< p < q < 2$ and $1<p_1,\ldots, p_k<\infty$ such that $\frac{1}{p_1}+\cdots+\frac{1}{p_k}=\frac{1}{p}$ and $p_j = 2$.
	\end{theorem}
	
	Recall the multilinear adjoint operators $A_{N; \mathbb K^k}^{{\mathcal P}, *j}$ defined by
	\[
	A_{N; \mathbb K^k}^{{\mathcal P}, *j}(f_1,\ldots,f_k)(x) \coloneqq 
	\E_{t \in [N]_{\mathbb K}}^{\lambda_{\mathbb K}}  
	\prod_{i\in[k]}
	\mathcal C^{\ind{i \neq j}} f_{i}\left(x- \ind{i \neq j} P_{i}(t)e_{i} + P_j(t)e_j\right).
	\]
	
	\begin{corollary}\label{cor:improvZ}
		Fix $k \in \Z_+$ and $j\in [k]$. Then there exist exponents $1<r<2$ and $1<r_1, \ldots, r_k<\infty$ such that
		$\frac{1}{r_1}+\cdots+\frac{1}{r_k}=\frac{1}{r}$ and
		\begin{equation}\label{adjoint-bound}
		\left\|A_{N; \mathbb K^k}^{{\mathcal P}, *j} (f_1,\ldots, f_k) \right\|_{L^2(\mathbb K^k)} \lesssim_{\mathcal P} N^{-D(\frac{1}{r}-\frac{1}{2})} \prod_{i \in [k]} \|f_i\|_{L^{r_i}(\mathbb K^k)}.
		\end{equation}
	\end{corollary}

	\begin{proof}[Proof of Corollary~\ref{cor:improvZ}] 
		When $k=1$, we can use duality to replace the adjoint operator by $A_{N; \mathbb K}^{\mathcal P}$. If $d_1 = 1$, then \eqref{adjoint-bound} holds by Young's convolution inequality for $r_1 = \frac{4}{3}$, say. If $\mathbb K = \Z$ and $d_1 \geq 2$, then \eqref{adjoint-bound} holds by a result of Han--Kova{\v c}--Lacey--Madrid--Yang in \cite{HKLMY}; see also \cite[Proposition~6.21]{KMT}. We refer to \cite{Ch} and \cite{Gr} for details in the remaining case $\mathbb K = \R$ and $d_1 \geq 2$.
		
		When $k > 1$, let $p_1,\dots,p_k,p,q$ be the exponents from Theorem~\ref{thm:improv} with $p_j=2$. By duality, H\"older's inequality, and \eqref{main-bound}, we obtain
		\begin{align*}
		\left\|A_{N; \mathbb K^k}^{{\mathcal P}, *j} (f_1,\ldots,f_k) \right\|_{L^2(\mathbb K^k)} 
		\lesssim_{\mathcal P} N^{-D(\frac{1}{p} - \frac{1}{q})} \|f_j\|_{L^{q'}(\mathbb K^k)} \prod_{i \in [k] \setminus \{j\}} 
		\|f_i\|_{L^{p_i}(\mathbb K^k)},
		\end{align*}
		where $q' \coloneqq q / (q-1)$ is the conjugate exponent of $q$. Taking $r_j = q'$,
		$r_i = p_i$ for $i \not= j$, and $r$ such that
		$ \frac{1}{r_1} + \dots + \frac{1}{r_k} = \frac{1}{r}$, we obtain
		$\frac{1}{p} - \frac{1}{q} = \frac{1}{r} - \frac{1}{2}$
		so that \eqref{adjoint-bound} holds. From $1 < p < q < 2$ and $\frac{1}{p} - \frac{1}{q} = \frac{1}{r} - \frac{1}{2}$, we see that $1 < r < 2$.
	\end{proof}
	
	\subsection{Vinogradov mean value theorem} The proof of Theorem~\ref{thm:improv} in the case $\mathbb K = {\mathbb Z}$ uses sharp bounds in the Vinogradov mean value theorem \cite{BDG}. For fixed $s,d \in \Z_+$, consider the following system of equations
	\[
	\begin{cases}
	\arraycolsep=1.5pt
	\begin{array}{ccccc}
	(n_1 + \cdots + n_s) & - & (m_1 + \cdots + m_s) & = & 0, \\
	(n_1^2 + \cdots + n_s^2) & - & (m_1^2 + \cdots + m_s^2) & = & 0, \\
	& & & \vdots & \\
	(n_1^d + \cdots + n_s^d) & - & (m_1^d + \cdots + m_s^d) & = & 0.
	\end{array}
	\end{cases}
	\]
	Let $J_{s,d}(N)$ be the number of its integer solutions $(n_1,\ldots, n_s, m_1, \ldots, m_s) \in [N]^{2s}$ for each fixed $N \in \Z_+$.  The Vinogradov
	mean value theorem states that
	\[
	J_{s,d}(N) \lesssim_{\epsilon} N^{\varepsilon} ( N^{s} + N^{2s-d(d+1)/2} )
	\]
	holds for all $\varepsilon > 0$ and, moreover, the $N^{\varepsilon}$ term may be removed in the range $s > d(d+1)/2$. This longstanding conjecture has recently been resolved by Wooley \cite{WOOLEY1} in the cubic case and by Bourgain, Demeter, and Guth in \cite{BDG} in the general case. In  \cite{WOOLEY}, Wooley used his efficient congruencing method to give  an alternative proof of this conjecture. 
	
	In the 1930's, Vinogradov verified the conjecture in the range $s \gtrsim d^2 \log d$. By the usual integral representation of $J_{s,d}(N)$ one has $J_{s,d,\xi}(N) \le J_{s,d}(N)$ for each $\xi \in \Z^d$, where $J_{s,d,\xi}(N)$ refers to the inhomogeneous system with the vector $(\xi_1,\dots,\xi_d)$ instead of $(0,\dots,0)$ on the right-hand side.
	
	Theorem~\ref{thm:improv} is a consequence of the following result.
	
	\begin{proposition}
		\label{restricted-improv}
		Fix $k \in \N_{\geq 2}$ and $s_i \in \N_{\geq D^{*}}$ for all $i\in [k]$, where $D^{*} \coloneqq \sum_{i \in [k]} \frac{d_i(d_i +1)}{2}$, and set $S \coloneqq \sum_{i \in [k]} s_i$. Then, for any distinct $j, j' \in [k]$, we have the following restricted weak-type bound
		\begin{equation}\label{restricted-bound-q}
		\left\| A_{N; \mathbb K^k}^{\mathcal P}(f_1, \ldots, f_k) \right\|_{L^{q,\infty}(\mathbb K^k)} 
		\lesssim_{\mathcal P}  N^{-D(\frac{1}{p} - \frac{1}{q})}
		\prod_{i \in [k]} 
		\|f_i\|_{L^{p_i,1}(\mathbb K^k)}  
		\end{equation}
		with $p_j = 2$, $p_{j'} = \frac{2S}{s_{j'}+1}$, and $p_i = \frac{2S}{s_i}$ for $i \notin \{j, j'\}$, where $1 < p < q < 2$ are chosen such that $\frac{1}{p_1} + \dots + \frac{1}{p_k} = \frac{1}{p}$ and $\frac{1}{q} = \frac{2S-s_j}{2S}$.
	\end{proposition}
	
	In order to deduce Theorem~\ref{thm:improv} from Proposition~\ref{restricted-improv}, we use a multilinear Marcinkiewicz interpolation argument in the spirit of \cite[Appendix]{MuscSch}. Precisely, we interpolate \eqref{restricted-bound-q} with $k+1$ trivial bounds of the form \eqref{main-bound} with $p=q$ and no decay in $N$, which hold for the tuples $(p_1, \dots, p_k)$ belonging to the set
	\[
	\big \{ (1,\infty,\dots,\infty), (\infty,1,\infty,\dots,\infty), \dots, (\infty, \dots, \infty, 1), (\infty, \dots, \infty) \big\}.
	\]
	Then Theorem~\ref{thm:improv} holds for some tuple $(p_1, \dots, p_k)$ and $1 < p < q < 2$ such that $p_j = 2$ and the vector $(1/p_1, \dots, 1/p_k, 1-1/q)$ lies in the interior of the convex hull of the $k+2$ analogous vectors determined by the corresponding exponents used in the interpolation.
	
	\begin{remark}
		If $\mathbb K = \Z$, then \eqref{restricted-bound-q} holds with $s_i \in \N_{> d_i(d_i +1)/2}$ for $i \in [k]$.
	\end{remark}
	
	\subsection{Lifting procedure}
	To deduce \eqref{restricted-bound-q} from Proposition~\ref{restricted-improv}, we may assume,
	without loss of generality, that each $P_i$ has no constant term. Furthermore,
	the bound for $A^{\mathcal P}_{N;\mathbb K^k}$ will follow from the corresponding
	bound for the {\it lifted} multilinear averaging operator
	\[
	A_{N; \mathbb K^k}^{{\mathcal P}, \, \rm lift}(f_1,\ldots, f_k)(x)
	\coloneqq \E_{u \in [N]_{\mathbb K}}^{\lambda_{\mathbb K}} \prod_{i\in[k]} f_i(x-\Gamma_i(u)),
	\qquad x\in\mathbb K^D,
	\]
	where now $f_i \colon \mathbb K^D \to {\mathbb C}$ for $i \in [k]$. Here $\Gamma_i$ is defined by
	\begin{equation}\label{gamma}
	\Gamma_i(u) \coloneqq (\underbrace{0, \dots, 0}_{D_{<i} \text{ times}},u,u^2,\ldots, u^{d_i},\underbrace{0, \dots, 0}_{D_{>i} \text{ times}})\in\K^D,
	\end{equation}
	where $D_{<i} \coloneqq d_1 + \dots + d_{i-1}$ and $D_{>i} \coloneqq d_{i+1} + \dots + d_k$ for $i\in[k]$.
	
	Below we extend an argument from \cite[Section~5]{HKLMY}.
	
	\begin{proposition}\label{lifting}
		Assume that $P_1(0) = \cdots = P_k(0)=0$ and recall that $D^{*} = \sum_{i \in [k]} \frac{d_i(d_i+1)}{2}$. Then \eqref{main-bound} follows from the estimate
		\begin{equation}\label{lifted-bound}
		\left\|A_{N; \mathbb K^k}^{{\mathcal P}, \, \rm lift}(f_1, \ldots, f_k)\right\|_{L^q(\mathbb K^{D})} \lesssim_{\mathcal P}  N^{-D^{*}(\frac{1}{p} - \frac{1}{q})} \prod_{i \in [k]} \|f_i\|_{L^{p_i}(\mathbb K^{D})},
		\end{equation}
		and an analogous transference of restricted weak-type bounds also holds.
	\end{proposition}
	
	\begin{proof} 
		For each $x \in \mathbb K^D$ we use the notation
		\[
		x = (x_{(1)}, \dots, x_{(k)}), \quad \text{where} \quad x_{(i)} \coloneqq (x_{D_{<i}+1}, \dots, x_{D_{<i+1}})\in \K^{d_i},
		\]
		or, alternatively, $x = (x',x'')$ referring to
		$x'' \coloneqq (x_{D_{<2}}, \dots, x_{D_{<k+1}})\in \K^{k}$ and
		\[
		x' = (x'_{(1)}, \dots, x'_{(k)}), \quad \text{where} \quad x'_{(i)} \coloneqq (x_{D_{<i}+1}, \dots, x_{D_{<i+1}-1})\in \K^{d_i-1}.
		\]
		
		\paragraph{\bf Step~1} We \emph{lift} functions $g \colon \mathbb K^k \to {\mathbb C}$ to functions $f \colon \mathbb K^D \to {\mathbb C}$. The definition differs slightly whether $\mathbb K = {\mathbb Z}$ or $\mathbb K = {\mathbb R}$ due to the lack of dilation structure on ${\mathbb Z}$. Let $P_i(u) = a^i_{1} u + \dots + a^i_{d_i} u^{d_i}$ with $a^i_{d_i} \neq 0$. Given $g \colon \Z^k \to {\mathbb C}$ and $r_i \in R_i \coloneqq\N_{<a_{d_i}^i}$, we define the \emph{lift} $f = f_{g,r_1,\ldots,r_k} \colon {\mathbb Z}^{D} \to {\mathbb C}$ of $g$ by
		\[
		f(x) \coloneqq g(Q_1(x_{(1)}) + r_1, \ldots, Q_k(x_{(k)}) + r_k) \ind{E_1}(x_{(1)}') \cdots \ind{E_k}(x_{(k)}'),
		\]
		where $E_i \coloneqq [\pm N]\times \cdots \times [\pm N^{d_i -1}]$ and
		\[
		Q_i(x_{(i)}) \coloneqq Q_i(x'_{(i)},x_{D_{<i+1}}) \coloneqq a^i_1 x_{D_{<i}+1} + \cdots + a^i_{d_i} x_{D_{<i+1}}.
		\]
		In the real case, given $g \colon \R^k \to {\mathbb C}$, we define the \emph{lift} $f = f_{g} \colon {\mathbb R}^{D} \to {\mathbb C}$ of $g$ by
		\[
		f(x) \coloneqq g(Q_1(x_{(1)}), \ldots, Q_k(x_{(k)})) \ind{E_1}(x_{(1)}') \cdots \ind{E_k}(x_{(k)}'),
		\]
		where
		$E_i \coloneqq [\pm N]_{\mathbb R}\times \cdots \times [\pm N^{d_i -1}]_{\mathbb R}$
		and $Q_i$ is as before.
		
		\paragraph{\bf Step~2} Fix $p_0 \in [1, \infty)$. Then
		\begin{align}
		\label{f-g}
		\|f\|_{L^{p_0}({\K}^{D})}^{p_0} \lesssim N^{D'} \|g\|_{L^{p_0}({\K}^k)}^{p_0}   
		\end{align}
		holds with $D' \coloneqq \sum_{i \in [k]} \frac{(d_i -1)d_i}{2}$. Moreover, if $x\in\K^D$ is such that $x'_{(i)} \in E_i'$ for all $i \in [k]$, where $E_i' \coloneqq [N]_{\mathbb K}\times \cdots \times [N^{d_i -1}]_{\mathbb K}$, then we have the equality
		\begin{equation*}
		A^{{\mathcal P}, \, \rm lift}_{N; \mathbb K^k}(f_1,\ldots, f_k)(x) 
		= A_{N; \mathbb K^k}^{\mathcal P} (g_1, \ldots, g_k)
		(Q_1(x_{(1)}) + r_1, \ldots, Q_k(x_{(k)}) + r_k),
		\end{equation*}
		with $r_1=\cdots=r_k=0$ if $\K=\R$.
		
		\paragraph{\bf Step~3} Suppose that  \eqref{lifted-bound} holds. By the above equality, we have
		\[
		N^{D'} \left\|A_{N; \mathbb K^k}^{\mathcal P} (g_1, \ldots, g_k)\right\|_{L^q(\mathbb K^k)}^q \lesssim_{\mathcal P} \left\|A_{N; \mathbb K^k}^{{\mathcal P}, \, \rm lift}(f_1,\ldots, f_k) \right\|_{L^q(\mathbb K^{D})}^q .
		\]
		Combining this with \eqref{lifted-bound} and \eqref{f-g}, we obtain
		\[
		N^{\frac{D'}{q}} \left\|A_{N; \mathbb K^k}^{\mathcal P} (g_1, \ldots, g_k)\right\|_{L^q(\mathbb K^k)} \lesssim_{\mathcal P} N^{-D^*(\frac{1}{p}-\frac{1}{q})} N^{\frac{D'}{p}} \prod_{i \in [k]} \|g_i\|_{L^{p_i}(\mathbb K^k)},
		\]
		which in view of $D^*-D' = \sum_{i \in [k]} d_i = D$ gives
		\eqref{main-bound}, as desired.
		
		The restricted weak-type estimates are handled analogously.
	\end{proof}
	
	\subsection{Refinements}
	By Proposition~\ref{lifting} it remains to prove
	\begin{equation}\label{lifted-restricted-bound}
	\left\|A_{N; \mathbb K^k}^{{\mathcal P}, \, \rm lift}(f_1, \ldots, f_k) \right\|_{L^{q,\infty}(\mathbb K^{D})} \lesssim_{\mathcal P} N^{-D^{*}(\frac{1}{p}-\frac{1}{q})}
	\prod_{i \in [k]} \|f_i\|_{L^{p_i,1}(\mathbb K^{D})} 
	\end{equation}
	for $p_1, \dots, p_k, p, q$ as in \eqref{restricted-bound-q}. We shall use the \emph{refinement method} developed by Christ in \cite{Ch}. From now on, we abbreviate $A_{N; \mathbb K^k}^{{\mathcal P}, \, \rm lift}$ to $A_{N; \mathbb K}$. We denote the corresponding adjoint operators by $A_{N; \mathbb K}^{*i}$.
	
	The bound \eqref{lifted-restricted-bound} equivalently means that 
	\[
	\langle \ind{E_0}, A_{N;\K}(\ind{E_1}, \ldots, \ind{E_k}) \rangle \lesssim N^{-D^*(\frac{1}{p} - \frac{1}{q})} |E_0|_{\K^D}^{\frac{1}{q'}} |E_1|_{\K^D}^{\frac{1}{p_1}}\cdots |E_k|_{\K^D}^{\frac{1}{p_k}}
	\]
	holds for all measurable subsets $E_0, E_1, \dots, E_k \subseteq \mathbb K^{D}$. Denoting the left-hand side by $K$, we can rewrite the above inequality as
	\begin{equation}\label{restricted-rewritten}
	K^{2S} \lesssim N^{-D^*} |E_0|_{\K^D}^{s_j} |E_j|_{\K^D}^S |E_{j'}|_{\K^D}^{s_{j'}+1} 
	\prod_{i \in [k] \setminus \{j,j'\}} |E_i|_{\K^D}^{s_i},
	\end{equation}
	since $\frac{1}{p} = \frac{2S+1-s_j}{2S}$, $\frac{1}{q} = \frac{2S-s_j}{2S}$, $p_j = 2$, $p_{j'} = \frac{2S}{s_{j'}+1}$, and $p_i = \frac{2S}{s_i}$ for $i \notin \{j, j'\}$. 
	
	Let $\alpha_0, \alpha_1, \ldots, \alpha_k \in [0,1]$ be parameters such that
	\[
	K = \alpha_0 |E_0|_{\K^D} = \alpha_1 |E_1|_{\K^D} = \cdots = \alpha_k |E_k|_{\K^D}.
	\]
	Note that, for each $i \in [k]$, we have
	\[
	\left\langle \ind{E_0}, A_{N; \mathbb K}(\ind{E_1}, \ldots, \ind{E_k}) \right\rangle = \left\langle \ind{E_{i}}, A_{N; \mathbb K}^{*i}(\ind{E_1}, \ldots, \ind{E_{i-1}}, \ind{E_0}, \ind{E_{i+1}}, \ldots, \ind{E_k}) \right\rangle
	\]
	where, referring to the mapping \eqref{gamma}, we define
	\[
	A_{N; \mathbb K}^{*i}(g_1,\ldots,g_k)(x) \coloneqq \E_{u \in [N]_{\mathbb K}}^{\lambda_{\mathbb K}} \prod_{i'\in[k]}
	\mathcal C^{\ind{i' \neq i}} g_{i'}\left(x - \ind{i' \neq i} \Gamma_{i'}(u) + \Gamma_{i}(u)\right).
	\]
	We observe that the average value of $A_{N; \mathbb K}(\ind{E_1}, \ldots, \ind{E_k})$ on
	$E_0$ is equal to
	\[
	|E_0|_{\K^D}^{-1} \langle A_{N; \mathbb K}(\ind{E_1}, \ldots, \ind{E_k}), \ind{E_0} \rangle = \alpha_0,
	\]
	and the average value of
	$A_{N; \mathbb K}^{*i}(\ind{E_1}, \ldots, \ind{E_{i-1}}, \ind{E_0}, \ind{E_{i+1}}, \ldots, \ind{E_k})$
	on $E_i$ is $\alpha_{i}$. Let $E_i^0 \coloneqq E_i$. We recursively define the \emph{refinements} of $E_i$ of level $r \in \Z_+$ by
	\begin{align*}
	\arraycolsep=1.5pt
	\begin{array}{ccc}
	E_0^{r} 
	& \coloneqq 
	& \Big \{x \in E_0^{r-1} : A_{N; \mathbb K} 
	(\ind{E_1^{r-1}}, \ind{E_2^{r-1}}, \ldots, \ind{E_{k-1}^{r-1}},
	\ind{E_k^{r-1}})(x) 
	\geq c_{r,0} \alpha_{0} \Big \}, \\
	E_1^{r} 
	& \coloneqq 
	& \Big \{x \in E_1^{r-1} : A_{N; \mathbb K}^{*1}
	( \hspace{0.175cm} \ind{E_0^{r\hspace{0.175cm}}}, \ind{E_2^{r-1}}, \ldots, \ind{E_{k-1}^{r-1}}, 
	\ind{E_k^{r-1}})(x) 
	\geq c_{r,1} \alpha_{1} \Big \}, \\
	& \vdots & \\
	E_k^{r} 
	& \coloneqq 
	& \Big \{x\in E_k^{r-1} : A_{N; \mathbb K}^{*k} 
	(\hspace{0.175cm}\ind{E_1^{r\hspace{0.175cm}}}, \hspace{0.175cm}\ind{E_2^{r\hspace{0.175cm}}},
	\ldots, \ind{E_{k-1}^{r}},
	\hspace{0.175cm} \ind{E_0^{r\hspace{0.175cm}}})(x) 
	\geq c_{r,k} \alpha_k \Big \},
	\end{array}
	\end{align*}    
	where 
	$(c_{r,0}, c_{r,1}, \dots, c_{r,k}) \coloneqq (\frac{1}{2^{(r-1)(k+1)+1}}, \frac{1}{2^{(r-1)(k+1)+2}}, \dots, \frac{1}{2^{r(k+1)}})$.
	
	The next result shows that passing to refinements is essentially lossless. 
	
	\begin{lemma}\label{no-loss} For all $r \in \Z_+$ and $i\in [k]$, we have
		\begin{equation} \label{no-loss-1} \left\langle \ind{E_0^{r}}, A_{N; \mathbb K}(\ind{E_1^{r-1}}, \ldots, \ind{E_k^{r-1}}) \right\rangle \ge c_{r,0} K
		\end{equation}
		and
		\begin{equation} \label{no-loss-2} \left\langle \ind{E_{i}^{r}}, A_{N; \mathbb K}^{*i}(\ind{E_1^{r}}, \ldots, \ind{E_{i-1}^{r}}, \ind{E_0^{r}}, \ind{E_{i+1}^{r-1}}, \dots, \ind{E_k^{r-1}}) \right\rangle \ge c_{r,i} K.
		\end{equation}
	\end{lemma}
	
	\begin{proof} 
		Let $r=1$. We estimate the left-hand side of \eqref{no-loss-1} by
		\[
		\Big\langle \ind{E_0}, A_{N; \mathbb K}(\ind{E_1}, \ldots, \ind{E_k})\Big\rangle - \left\langle \ind{E_0 \setminus E_0^{1}}, A_{N; \mathbb K}(\ind{E_1}, \ldots, \ind{E_k})\right\rangle \geq K - K/2 = K/2,
		\]
		splitting $\ind{E_0^{1}} = \ind{E_0} - \ind{E_0 \setminus E_0^{1}}$ and then using the definition of $E_0^1$. Similarly, for $i=1$, we rewrite the
		left-hand side of \eqref{no-loss-2} as
		\[
		\left\langle \ind{E_1}, A_{N; \mathbb K}^{*1}(\ind{E_0^1}, \ind{E_2},\ldots, \ind{E_k})\right\rangle
		- \left\langle \ind{E_1\setminus E^{1}_1}, A_{N; \mathbb K}^{*1}(\ind{E_0^1}, 
		\ind{E_2}, \ldots, \ind{E_k})\right\rangle
		\]
		which by duality, the previous step, and the definition of $E_1^1$ is estimated by
		\[
		\left\langle \ind{E_0^1}, A_{N; \mathbb K}(\ind{E_1}, \ldots, \ind{E_k})\right\rangle
		- \left\langle \ind{E_1\setminus E^{1}_1}, A_{N; \mathbb K}^{*1}(\ind{E_0^1}, 
		\ind{E_2}, \ldots, \ind{E_k})\right\rangle
		\ge K/2 - K/4.
		\]
		Analogously, for every $i \in [k]$, we have
		\[
		\left\langle \ind{E_i^{1}}, A_{N; \mathbb K}^{*i}(\ind{E_1^1}, \ldots, \ind{E_{i-1}^1},
		\ind{E_0^1}, \ind{E_{i+1}}, \ldots, \ind{E_k})\right\rangle \ge K/2^{i+1}.
		\]
		By induction, the same argument yields \eqref{no-loss-1} and  \eqref{no-loss-2} for every $r \in \Z_+$. 
	\end{proof}
	
	\subsection{Parameter towers and flows}
	From now on, by symmetry, we may assume that $j=1$ and $j'=k$. To establish \eqref{restricted-rewritten} we may assume that $K>0$ so that by Lemma~\ref{no-loss} every refinement $E_i^{r}$ is nonempty. Using $\alpha_i$ as before we can rewrite \eqref{restricted-rewritten} as
	\begin{equation}\label{restricted-rewritten-again}
	\alpha_0^{s_1} \alpha_1^{S}
	\alpha_{2}^{s_{2}}
	\cdots \alpha_{k}^{s_{k}}
	\lesssim N^{-D^*} |E_k|_{\K^D}.
	\end{equation}
	
	We fix some element in one of the refinements, say $z \in E_1^{r}$, and obtain
	\[
	A_{N; \mathbb K}^{*1}(\ind{E_0^{r}}, \ind{E_2^{r-1}}, \ldots, \ind{E_k^{r-1}})
	(z)
	\gtrsim_{k,r} \alpha_1.
	\]
	Hence, we obtain that 
	\[
	|\{u \in [N]_{\mathbb K} : z + \Gamma_1(u) \in E_0^r 
	\text{ and } 
	z + \Gamma_1(u) - \Gamma_i(u) \in E_i^{r-1} 
	\text{ for all } i \in [k] \setminus \{1\} 
	\}|_{\K}
	\]
	admits a lower bound of order $\alpha_1 N$
	and, in particular, the parameter set 
	\[
	I^{z} \coloneqq \{u \in [N]_{\mathbb K} : z + \Gamma_1(u) \in  E_0^{r} \}
	\]
	satisfies $|I^{z}|_{\K} \gtrsim_{k,r} \alpha_1 N$. Next, for each $u \in I^z$, we have $z + \Gamma_1(u) \in  E_0^{r}$. Thus,
	\[
	|\{ v \in [N]_{\mathbb K} : z + \Gamma_1(u) - \Gamma_i(v) \in E_i^{r-1} \text{ for all } i \in [k] \}|_{\K} \gtrsim_{k,r} \alpha_0 N
	\]
	and, in particular, the parameter set
	\[
	I^{z,u} \coloneqq \{v \in [N]_{\mathbb K} : z + \Gamma_1(u) - \Gamma_1(v) \in  E_1^{r-1} \}
	\]
	satisfies $|I^{z,u}|_{\K} \gtrsim_{k,r} \alpha_0 N$. This begins a construction of a parameter tower 
	\[
	{\mathbf P}_{(1,0,1)} \coloneqq \{(u,v) \in [N]_{\mathbb K}^2 : u \in I^{z}, v \in I^{z,u} \}
	\]
	such that $|{\mathbf P}_{(1,0,1)}|_{\K^2} \gtrsim_{k,r} \alpha_1 \alpha_0 N^2$ and a flow
	\[
	\Phi_{(1,0,1)}(u,v) \coloneqq z + \Gamma_1(u) - \Gamma_1(v) = z 
	+ (u - v, \ldots, u^{d_1} - v^{d_1}, 0, \ldots, 0) 
	\]
	such that $\Phi_{(1,0,1)}({\mathbf P}_{(1,0,1)}) \subseteq E_1^{r-1}$. We express this as $E_1^{r} \to E_0^{r} \to E_1^{r-1}$.
	
	When we flow from one set $E_i$ into another $E_{i'}$ via the above process, there is no ambiguity which level of refinement we are using:
	\begin{itemize}
		\item[(i)] if $i < i'$, then necessarily we drop a level, so $E_i^r \to E_{i'}^{r-1}$;
		\item[(ii)] if $i' < i$, then
		necessarily we stay at the same level, so $E_i^r \to E_{i'}^r$.
	\end{itemize}
	
	Our flows will always consist of a uniformly controlled finite number of iterations; hence, there is no need to keep track of the level of refinement as long as we start with $z \in E_i^{r_0}$ for a fixed level $r_0 \lesssim 1$ which is sufficiently large to ensure that the last level is positive.
	
	From now on, to simplify notation, we suppress the refinement level, dropping the superscript $r$. For the flow above, for example, we write $E_1 \to E_0 \to E_1$ and $\Phi_{(1,0,1)}({\mathbf P_{(1,0,1)}}) \subseteq E_1$. Therefore, 
	\[
	|E_1|_{\K^D} \ge |\Phi_{(1,0,1)}({\mathbf P_{(1,0,1)}})|_{\K^D}
	= \int_{\Phi_{(1,0,1)}({\mathbf P_{(1,0,1)}})} d\lambda_{\mathbb K}^{\otimes d_1}(\xi_{(1)})
	\]
	and the idea is to make the change of variables $\xi_{(1)} = \Phi_{(1,0,1)}(u,v)$.
	
	\begin{example}
		For this illustration, let us drop the subscript $(1,0,1)$ and suppose for simplicity that $d_1 =D= 2$. When $\mathbb K = \R$, we can use the usual change of variables formula involving the Jacobian $J_{\Phi}$. Namely, we have
		\begin{equation}\label{E1-I}
		|E_1|_{\R^2} \ge \iint_{{\mathbf P}} |J_{\Phi}(t_1,t_2)| dt_1 dt_2
		= \int_{I^{z}} \int_{I^{z, t_1}} 2 |t_1 - t_2| dt_2 dt_1,
		\end{equation}
		which reduces matters to obtaining an appropriate lower bound for a weighted measure of the parameter set ${\mathbf P}$. The latter can be effectively achieved using sublevel set bounds for $J_{\Phi}$. When $\mathbb K = \Z$, we have
		\[
		|E_1|_{\Z^2} \ge \sum_{\xi_{(1)} \in \Phi({\mathbf P})} 1 =
		\sum_{\xi_{(1)} \in \Phi({\mathbf P})} J_{\xi_{(1)}}^{-1} 
		\sum_{(n,m) \in \mathbf P} \ind{\Phi(n,m) = \xi_{(1)}},
		\]
		where $J_{\xi_{(1)}} \coloneqq |\{ (n,m) \in {\mathbf P} : \Phi(n,m)  = \xi_{(1)} \}|_{\Z^2}$.
		Interchanging the sums gives 
		\begin{equation}\label{E1-II}
		|E_1|_{\Z^2} \ge |\Phi({\mathbf P})|  = \sum_{(n,m) \in {\mathbf P}} J_{\Phi(n,m)}^{-1},
		\end{equation}
		reducing matters to obtaining a uniform lower bound for the ``Jacobian'' $J_{\xi_{(1)}}$, where $J_{\xi_{(1)}}$ counts the number of integer solutions $(n,m) \in [N]^2$ to the system of Diophantine equations
		\[
		\begin{cases}
		z_1 + n\phantom{^2} - m\phantom{^2} = \xi_1, \\
		z_2 + n^2 - m^2 = \xi_2.
		\end{cases}
		\] 
	\end{example}
	
	Returning to the general case, in both settings the ultimate flow will be a concatenation of $k$ subflows
	$({\overline{E_1 \ E_i}})$ or  $({\overline{E_0 \ E_i}})$, where each $({\overline{E_{i'} \ E_i}})$ is itself a flow consisting of some number of $E_{i'}$ and $E_{i}$ alternating.
	When $\mathbb K = {\mathbb Z}$,
	the ultimate flow has the form
	\begin{equation}\label{k-blocks-integer}
	({\overline{E_1 \ E_0}}) \to 
	({\overline{E_1 \ E_2}}) \to 
	({\overline{E_1 \ E_{3}}}) \to 
	\cdots \to ({\overline{E_1 \ E_k}}).
	\end{equation}
	When $\mathbb K = {\mathbb R}$, the ultimate flow is slightly different and has the form
	\begin{equation}\label{k-blocks-real}
	({\overline{E_0 \ E_1}}) \to 
	({\overline{E_0 \ E_2}}) \to 
	\cdots \to ({\overline{E_0 \ E_k}}).
	\end{equation}
	For example, the flow $E_1 \to E_0 \to E_1 \to E_2 \to E_1$ produces a parameter tower 
	\[
	{\mathbf P}_{(1,0,1,2,1)} \coloneqq \{(u,v,u',v') \in [N]_{\mathbb K}^4 : u \in I^{z}, v \in I^{z,u}, u' \in I^{z,u,v}, v' \in I^{z,u,v,u'}  \}
	\]
	with $|{\mathbf P}_{(1,0,1,2,1)}|_{\K^4} \gtrsim_{k,r} \alpha_1 \alpha_2 \alpha_1 \alpha_0 N^4$ and $\Phi_{(1,0,1,2,1)}({\mathbf P}_{(1,0,1,2,1)}) \subseteq E_1$, where
	\[
	\Phi_{(1,0,1,2,1)}(u,v,u',v') \coloneqq z + \Gamma_1(u) - \Gamma_1(v) + \Gamma_1(u') - \Gamma_2(u') - \Gamma_1(v') + \Gamma_2(v').
	\]
	
	\subsection{The integer case} We now prove Proposition
	\ref{restricted-improv} in the integer setting.
	
	\begin{proof}[Proof of Proposition~\ref{restricted-improv} for $\K=\Z$]
		Suppose that the parameters $i$ from \eqref{k-blocks-integer} determine a sequence $\mathcal I$ of the following form
		\[ 
		(\underbrace{1,0,\dots,1,0}_{2s_1 \text{ elements}},  
		\underbrace{1,2,\dots,1,2}_{2s_2 \text{ elements}},
		\underbrace{1,3,\dots,1,3}_{2s_3 \text{ elements}},
		\dots,
		\underbrace{1,k-1,\dots,1,k-1}_{2s_{k-1} \text{ elements}},
		\underbrace{k,1,\dots,1,k}_{2s_k + 1 \text{ elements}}) 
		\]
		with a slightly different final block, where we used the transition $E_{k-1} \to E_k$.
		
		\paragraph{\bf Step~1} Unpacking all this gives a parameter tower
		\[ {\mathbf P}_{\mathcal I} \coloneqq \{ (\textbf{\textit{n}},\textbf{\textit{m}}) \in [N]^{2S} : \textbf{\textit{n}}_1 \in I^{z}, \textbf{\textit{m}}_1 \in I^{z, \textbf{\textit{n}}_1}, 
		\ldots, \textbf{\textit{m}}_S \in I^{z,\textbf{\textit{n}}_1,\textbf{\textit{m}}_1, \dots, \textbf{\textit{n}}_{S}} \},
		\]
		where $r_0 \coloneqq 10S$ and $z \in E_1^{r_0} \subseteq {\mathbb Z}^D$ is fixed. Thus, we have
		\[
		|\mathbf P_{\mathcal I}|_{\Z^D} \gtrsim_{\mathcal P} \alpha_0^{s_1} \alpha_1^S \alpha_2^{s_2} \cdots \alpha_k^{s_k} N^{2S}.
		\]
		The corresponding flow map
		$\Phi_{\mathcal I}(\textbf{\textit{n}},\textbf{\textit{m}})$ can be
		expressed as
		\[
		z - \big( \xi^1_{\textbf{\textit{n}},\textbf{\textit{m}}} - \Phi_{\mathcal I}^{(1)}(\textbf{\textit{n}},\textbf{\textit{m}}), \Phi_{\mathcal I}^{(2)}(\textbf{\textit{n}},\textbf{\textit{m}}), \dots, \Phi_{\mathcal I}^{(k-1)}(\textbf{\textit{n}},\textbf{\textit{m}}), \xi^k_{\textbf{\textit{n}},\textbf{\textit{m}}} - \Phi_{\mathcal I}^{(k)}(\textbf{\textit{n}},\textbf{\textit{m}}) \big),
		\]
		where, writing $\textbf{\textit{n}}=(\textbf{\textit{n}}_{(1)}, \dots, \textbf{\textit{n}}_{(k)}), \textbf{\textit{m}}=(\textbf{\textit{m}}_{(1)}, \dots, \textbf{\textit{m}}_{(k)}) \in [N]^{s_1} \times \dots \times [N]^{s_k}$,
		we let $\Phi_{\mathcal I}^{(i)}(\textbf{\textit{n}},\textbf{\textit{m}}) \coloneqq \Phi_{i}(\textbf{\textit{n}}_{(i)},\textbf{\textit{m}}_{(i)})$ for all $i \in [k]$, with $\Phi_i(u,v)$ defined by
		\[
		\big( (u_1 + \cdots + u_{s_i}) - (v_1 + \cdots + v_{s_i}), \ldots, (u_1^{d_i} + \cdots + u_{s_i}^{d_i}) - (v_1^{d_i} + \cdots + v_{s_i}^{d_i}) \big)
		\]
		for $u,v \in [N]^{s_i}$, and $\xi^1_{\textbf{\textit{n}},\textbf{\textit{m}}} \in \Z^{d_1}$ is independent of $\textbf{\textit{n}}_{(1)}, \textbf{\textit{m}}_{(1)}$, while $\xi^k_{\textbf{\textit{n}},\textbf{\textit{m}}} \in \Z^{d_k}$ depends solely on $\textbf{\textit{n}}_{(k-1)}, \textbf{\textit{m}}_{(k-1)}$ if $k \geq 3$ or simply vanishes if $k=2$.
		
		\paragraph{\bf Step~2} We drop the subscript $\mathcal I$. Since
		$\Phi({\mathbf P}) \subseteq E_k$, we obtain
		\[
		|E_k|_{\Z^D} \ge \sum_{(\textbf{\textit{n}},\textbf{\textit{m}}) \in \mathbf P} J^{-1}_{\Phi(\textbf{\textit{n}},\textbf{\textit{m}})}, 
		\]
		proceeding as in \eqref{E1-II}. Next, we will find a uniform upper bound for $J_{\xi}$, where $\xi \in \Z^{D}$. We will count all integer solutions $(\textbf{\textit{n}},\textbf{\textit{m}}) \in [N]^{2S}$ to the system 
		\[
		\begin{cases}
		\arraycolsep=1.5pt
		\begin{array}{ccl}
		\xi^1_{\textbf{\textit{n}},\textbf{\textit{m}}} - \Phi_1(\textbf{\textit{n}}_{(1)},\textbf{\textit{m}}_{(1)}) & = & z_{(1)} - \xi_{(1)}, \\
		\hfill \Phi_2(\textbf{\textit{n}}_{(2)},\textbf{\textit{m}}_{(2)}) & = & z_{(2)} - \xi_{(2)}, \\
		& \vdots & \\
		\hfill \Phi_{k-1}(\textbf{\textit{n}}_{(k-1)},\textbf{\textit{m}}_{(k-1)}) & = & z_{(k-1)} - \xi_{(k-1)}, \\
		\xi^k_{\textbf{\textit{n}},\textbf{\textit{m}}} - \Phi_k(\textbf{\textit{n}}_{(k)},\textbf{\textit{m}}_{(k)}) & = & z_{(k)} - \xi_{(k)}, 
		\end{array}
		\end{cases}
		\]
		identifying $z, \xi \in {\mathbb Z}^D$ with
		$(z_{(1)}, \ldots, z_{(k)}), \, \xi = (\xi_{(1)}, \ldots, \xi_{(k)}) \in {\mathbb Z}^{d_1} \times \cdots \times {\mathbb Z}^{d_k}$.
		
		\paragraph{\bf Step~3}  Using $s_i > \frac{d_i(d_i+1)}{2}$, the definition of $\Phi_i$, and
		the Vinogradov mean value theorem, we conclude that the subsystem consisting of
		all but the first and last equations has at most
		\[
		J_{s_2, d_2}(N) \cdots J_{s_{k-1}, d_{k-1}}(N) \lesssim N^{2(s_2 + \dots + s_{k-1}) - (d_2(d_2+1) + \dots + d_{k-1}(d_{k-1}+1))/2} 
		\]
		solutions. Once these solutions have been counted, the variables
		$\textbf{\textit{n}}_{(k-1)},\textbf{\textit{m}}_{(k-1)}$ are
		determined and we can bound the number of solutions to the last
		equation by
		$J_{s_k, d_k}(N) \lesssim N^{2 s_k - d_k(d_k+1)/2}$. Once all
		variables $\textbf{\textit{n}}_{(i)},\textbf{\textit{m}}_{(i)}$
		with $i > 1$ have been determined, we can bound the number of solutions
		to the first equation by
		$J_{s_1, d_1}(N) \lesssim N^{2s_1 - d_1(d_1+1)/2}$.  Altogether, we
		have $ J_{\xi} \lesssim N^{2S - D^{*}}$. Thus,
		\[
		|E_k|_{\Z^D} \gtrsim N^{-2S + D^{*}}
		|\mathbf P|_{\Z^D} \gtrsim N^{-2S + D^{*}} \alpha_0^{s_1} \alpha_1^S  \alpha_2^{s_2} \cdots \alpha_k^{s_k} N^{2S}
		\]
		and \eqref{restricted-rewritten-again} holds. This completes the proof of Proposition~\ref{restricted-improv}
		for $\mathbb K = {\mathbb Z}$.
	\end{proof}
	
	\subsection{The real case} We now prove Proposition
	\ref{restricted-improv} in the real setting.
	
	\begin{proof}[Proof of Proposition~\ref{restricted-improv} for $\K=\R$]
		In this case, we will produce exactly $D$ equations in order to follow the Jacobian approach \eqref{E1-I}. In particular, the flow \eqref{k-blocks-real} now depends on the parities of the degrees $d_i$ of the polynomials $P_i$. Precisely, each block 
		$({\overline{E_0 \ E_i}})$ in \eqref{k-blocks-real} ends with $E_i$ and is of length $d_i$, except the last one which ends with $E_k$ but is of length $d_k + 1$. Therefore, the parameters $i$ from \eqref{k-blocks-real} determine a sequence $\mathcal I$ of the following form
		\[
		(\omega_{1}^1,\ldots, \omega_{d_1}^1,  
		\omega_{1}^2,\ldots, \omega_{d_2}^2,
		\dots,
		\omega_{1}^{k-1},\ldots, \omega_{d_{k-1}}^{k-1},
		\omega_{1}^{k},\ldots, \omega_{d_{k}+1}^{k}),
		\]
		where
		\begin{align*}
		(\omega_{1}^{l},\ldots, \omega_{d_{l}}^{l})& \coloneqq
		\begin{cases}
		\ (0,l,\ldots, 0, l) & \text{ if } d_l \text{ is even},\\
		(l, 0, \ldots, l, 0, l) & \text{ if } d_l \text{ is odd},
		\end{cases}
		\quad \text{ for } \quad l\in[k-1],\\
		(\omega_{1}^{k},\ldots, \omega_{d_{k}+1}^{k})& \coloneqq
		\begin{cases}
		(k,0,\ldots, k, 0, k) & \text{ if } d_k \text{ is even},\\
		\ (0, k, \ldots,  0, k) & \text{ if } d_k \text{ is odd}.
		\end{cases}
		\end{align*}
		For example, if all $d_i$ are even, then $\mathcal I$ takes the form
		\[ 
		(\underbrace{0,1,\dots,0,1}_{d_1 \text{ elements}},  
		\underbrace{0,2,\dots,0,2}_{d_2 \text{ elements}},
		\dots,
		\underbrace{0,k-1,\dots,0,k-1}_{d_{k-1} \text{ elements}},
		\underbrace{k,0,\dots,k,0,k}_{d_k + 1 \text{ elements}}), 
		\]
		whereas if all $d_i$ are odd, then $\mathcal I$ takes the form
		\[ 
		(\underbrace{1,0,\dots,1,0,1}_{d_1 \text{ elements}},  
		\underbrace{2,0,\dots,2,0,2}_{d_2 \text{ elements}},
		\dots,
		\underbrace{k-1,0,\dots,k-1,0,k-1}_{d_{k-1} \text{ elements}},
		\underbrace{0,k,\dots,0,k}_{d_k + 1 \text{ elements}}).
		\]
		
		\paragraph{\bf Step~1} Set $r_0 \coloneqq 10D$. Unpacking all this gives a parameter tower
		\[
		{\mathbf P}_{\mathcal I} \coloneqq \{ \textbf{\textit{t}} \in [N]_{\R}^{D} : 
		\textbf{\textit{t}}_1 \in I^{z}, \textbf{\textit{t}}_2 \in I^{z, \textbf{\textit{t}}_1}, \dots,  \textbf{\textit{t}}_D \in I^{z, \textbf{\textit{t}}_1, \textbf{\textit{t}}_2, \dots, \textbf{\textit{t}}_{D-1}} \},
		\]
		for some fixed $z \in E_0^{r_0}$ if $d_1$ is even or $z \in E_1^{r_0}$ if $d_1$ is odd, and the flow map
		\[
		\Phi_{\mathcal I}(\textbf{\textit{t}}) \coloneqq
		z - 
		\big(
		\xi^1_{\textbf{\textit{t}}} + \Phi_1(\textbf{\textit{t}}_{(1)}), 
		\xi^2_{\textbf{\textit{t}}} + \Phi_2(\textbf{\textit{t}}_{(2)}),
		\dots,
		\xi^k_{\textbf{\textit{t}}}
		+ \Phi_k(\textbf{\textit{t}}_{(k)}) \big),
		\]
		where, writing
		$\textbf{\textit{t}}=(\textbf{\textit{t}}_{(1)}, \dots, \textbf{\textit{t}}_{(k)}), \xi_{\textbf{\textit{t}}} = (\xi^1_{\textbf{\textit{t}}}, \dots, \xi^k_{\textbf{\textit{t}}}) \in [N]_{\R}^{d_1} \times \dots \times [N]_{\R}^{d_k}$
		as before, for each $i \in [k]$ such that $d_i$ is even we have 
		\[
		\Phi_i(u_1, \dots, u_{d_i}) \coloneqq \big( u_1 - u_2 + \dots + u_{d_i-1} - u_{d_i}, \ldots, u_1^{d_i} - u_2^{d_i} + \dots + u_{d_i-1}^{d_i} - u_{d_i}^{d_i} \big),
		\]
		whereas for each $i \in [k]$ such that $d_i$ is odd we have 
		\[
		\Phi_i(u_1, \dots, u_{d_i}) \coloneqq - \big( u_1 - u_2 + \dots - u_{d_i-1} + u_{d_i}, \ldots, u_1^{d_i} - u_2^{d_i} + \dots - u_{d_i-1}^{d_i} + u_{d_i}^{d_i} \big),
		\]
		and $\xi^1_{\textbf{\textit{t}}}$ vanishes, whereas if $i \neq 1$, then
		$\textbf{\textit{t}} \mapsto \xi^i_{\textbf{\textit{t}}}$ is a
		polynomial that may depend on the coordinates
		$\textbf{\textit{t}}_{(i')}$ with $i' < i$ but is independent of
		$\textbf{\textit{t}}_{(i')}$ with $i' \geq i$.
		
		\paragraph{\bf Step~2} We drop the subscript $\mathcal I$. Since $\Phi({\mathbf P}) \subseteq E_k$, we obtain
		\begin{equation}\label{Ek-II}
		|E_k|_{\R^D} \geq |\Phi({\mathbf P})|_{\R^D} \gtrsim \int_{\mathbf P} |J_{\Phi}({\textbf{t}}) | d{\textbf{t}},
		\end{equation}
		and the Jacobian matrix $\Phi'$, consisting of the first order partial derivatives of the polynomial map $\Phi$, takes the following form
		\[
		\arraycolsep=1.4pt\def\arraystretch{1.2}
		\left(\begin{array}{c|c|c|c|c}
		\mathbf W_1 & \mathbf 0 & \mathbf 0 & \, \cdots \, & \mathbf 0  \\
		\hline
		\mathbf M_{2,1} & \mathbf W_2 & \mathbf 0 & \, \cdots \, & \mathbf 0 \\
		\hline
		\, \mathbf M_{3,1} \, & \, \mathbf M_{3,2} \, & \mathbf W_{3} & \, \cdots \, & \mathbf 0 \\
		\hline
		\vdots & \vdots & \vdots & \, \ddots \, & \vdots \\
		\hline
		\mathbf M_{k,1} & \mathbf M_{k,2} & \, \mathbf M_{k,3} \, & \, \cdots \, & \, \mathbf W_k \,  
		\end{array}\right)
		\]
		and for each block $\mathbf W_i$, after removing signs from every second row, we have
		\[
		|\det \mathbf W_i(\textbf{\textit{t}}) | = 
		\left| \, \det \left(
		\arraycolsep=2pt\def\arraystretch{1.25}
		\begin{array}{ccccc}
		1 & 2 \textbf{\textit{t}}_{D_{<i}+1} & 3 \textbf{\textit{t}}_{D_{<i}+1}^2 & \, \cdots \, & d_i \textbf{\textit{t}}_{D_{<i}+1}^{d_i-1}  \\
		1 & 2 \textbf{\textit{t}}_{D_{<i}+2} & 3 \textbf{\textit{t}}_{D_{<i}+2}^2 & \, \cdots \, & d_i \textbf{\textit{t}}_{D_{<i}+2}^{d_i-1} \\
		1 & 2 \textbf{\textit{t}}_{D_{<i}+3}  & 3 \textbf{\textit{t}}_{D_{<i}+3}^2 & \, \cdots \, & d_i \textbf{\textit{t}}_{D_{<i}+3}^{d_i-1}  \\
		\vdots & \vdots & \vdots & \, \ddots \, & \vdots \\
		1 & 2 \textbf{\textit{t}}_{D_{<i+1}} & 3 \textbf{\textit{t}}_{D_{<i+1}}^2 & \, \cdots \, & d_i \textbf{\textit{t}}_{D_{<i+1}}^{d_i-1}   
		\end{array} \right) \right| 
		\]
		and therefore, using properties of the Vandermonde matrix, we conclude that
		\[
		|J_{\Phi}({\textbf{t}})| = |\det \Phi'({\textbf{t}}) | = \prod_{i\in[k]} \Bigg( d_i!
		\prod_{\substack{j_i, j_i' \in [D_{<i+1}] \setminus [D_{<i}]\\  j_i< j_i'}} |\textbf{\textit{t}}_{j_i} - \textbf{\textit{t}}_{j_i'}| \Bigg).
		\]
		
		\paragraph{\bf Step~3} Consider the innermost integral in \eqref{Ek-II}. For fixed $\textbf{\textit{t}}_{1}, \dots, \textbf{\textit{t}}_{D-1}$ we integrate $|J_{\Phi}({\textbf{t}})|$ over $\textbf{\textit{t}}_{D} \in I^{z,\textbf{\textit{t}}_{1}, \dots, \textbf{\textit{t}}_{D-1}} \subseteq [0,N]$, where $|I^{z,\textbf{\textit{t}}_{1}, \dots, \textbf{\textit{t}}_{D-1}}|_{\R} \gtrsim_{k,r_0} \alpha_k N$. Let $I \subseteq I^{z,\textbf{\textit{t}}_{1}, \dots, \textbf{\textit{t}}_{D-1}}$ be such that $|I|_{\R} \gtrsim_{\mathcal P} \alpha_k N$ and $|\textbf{\textit{t}}_{j} - \textbf{\textit{t}}_{D}| \gtrsim_{\mathcal P} \alpha_k N$ for all $j \in [D-1] \setminus [D-d_k]$ and $\textbf{\textit{t}}_{D} \in I$. We pull out all factors independent of $\textbf{\textit{t}}_{D}$ and estimate the remaining part as follows
		\[
		\int_{I^{z,\textbf{\textit{t}}_{1}, \dots, \textbf{\textit{t}}_{D-1}}} 
		|\textbf{\textit{t}}_{D-d_k+1} - \textbf{\textit{t}}_{D}| \cdots 
		|\textbf{\textit{t}}_{D-1} - \textbf{\textit{t}}_{D}|
		d\textbf{\textit{t}}_{D}
		\gtrsim_{\mathcal P} |I|_{\R} (\alpha_k N)^{d_k-1} 
		\gtrsim_{\mathcal P} (\alpha_k N)^{d_k} 
		\]
		with some constant that is uniform in $\textbf{\textit{t}}_{1}, \dots, \textbf{\textit{t}}_{D-1}$. Continuing this way, we obtain the uniform lower bound $(\alpha_0 N)^{d_k-1}$ for the analogous integral over $\textbf{\textit{t}}_{D-1}$ and after $d_k$ such steps, using $\alpha_i \in (0,1]$, we end up with the integrand
		\[
		\prod_{i\in[k-1]} \Bigg( d_i!
		\prod_{\substack{j_i, j_i' \in [D_{<i+1}] \setminus [D_{<i}]\\ j_i< j_i'}} |\textbf{\textit{t}}_{j_i} - \textbf{\textit{t}}_{j_i'}| \Bigg)(\alpha_0 \alpha_1 \cdots \alpha_k N)^{\frac{d_k(d_k+1)}{2}}.
		\]
		Repeating this for other blocks leads to
		$|E_k|_{\R^D} \gtrsim_{\mathcal P} (\alpha_0 \alpha_1 \cdots \alpha_k N)^{D^*}$ and, since $s_i \ge D^{*}$ for every $i \in [k]$, this implies $|E_k|_{\R^D} \gtrsim_{\mathcal P} \alpha_0^{s_1} \alpha_1^S \alpha_2^{s_2} \cdots \alpha_k^{s_k} N^{D^{*}}$ so that \eqref{restricted-rewritten-again} holds. This completes the proof of Proposition~\ref{restricted-improv} for $\mathbb K = \R$.
	\end{proof}
	
	\section{Multilinear Weyl inequality and Sobolev smoothing theorem}
	\label{sec:weyl}
	
	In this section, we will formulate and prove a multilinear Weyl inequality, which asserts that the averaging operator $\tilde A^{P_1(\nn),\ldots, P_k(\nn)}_{N;\Z^k}(f_1,\ldots, f_k)$ is negligible when the $j$-th Fourier transform of $f_j$, for at least one $j\in[k]$, vanishes on appropriate major arcs. We will also need a multilinear Sobolev smoothing inequality, which can be thought of as a continuous variant of the multilinear Weyl inequality for the averaging operators $\tilde A^{P_1({\rm t}),\ldots, P_k({\rm t})}_{N;\R^k}(f_1,\ldots, f_k)$.
	
	\subsection{Multilinear smoothing inequalities}
	Using notation from Section~\ref{sec:inverse}, in particular \eqref{transpose}--\eqref{eq:129}, we formulate the multilinear Weyl inequality and the multilinear Sobolev smoothing inequality as the following unified theorem.
	
	\begin{theorem}\label{weyl}
		Let $\K$ be either $\Z$ or $\R$. Fix $C_1, C_2\in\R_+$ and $k\in\Z_+$, and let $\mathcal P = (P_1, \dots, P_k)$ be a polynomial mapping satisfying \eqref{eq:42}--\eqref{eq:44}. Let $1<p_1,\ldots, p_k<\infty$ be exponents such that $\frac{1}{p_1}+\cdots+\frac{1}{p_k}=\frac{1}{p}\le 1$. Then there exists a small constant $c \in (0,1)$, depending only on $k,\mathcal P, p_1, \dots, p_k, p, C_1, C_2$, such that the following holds. Assume that $\delta\in(0, 1]$ and $N \geq 1$. Let $f_i \in L^{p_i}(\K^k)$ for all $i\in[k]$. If $f_j \in L^{2}(\K^k)$ for some $j\in[k]$ and the $j$-th Fourier transform $\F_{j, \K^k} f_j$ vanishes on the major arcs ${\mathfrak M}_{\leq N^{-d_j}\delta^{-C_2}}^j(\mathcal R_{\leq \delta^{-C_1}}^{\K})$, then
		\begin{align}
		\label{eq:43}
		\left\| \tilde A^{P_1,\ldots, P_k}_{N; \K^k}(f_1,\ldots, f_k) \right\|_{L^p(\K^k)} 
		\leq c^{-1} (\delta^c+N^{-c}) \prod_{i \in [k]} \|f_i\|_{L^{p_i}(\K^k)}.
		\end{align}
		The same conclusion holds for $A^{P_1,\ldots, P_k}_{N; \K^k}$ in place of $\tilde A^{P_1,\ldots, P_k}_{N; \K^k}$ in \eqref{eq:43}.
	\end{theorem}
	If $\K=\Z$, then \eqref{eq:43} is the multilinear Weyl inequality. If $\K=\R$, then \eqref{eq:43} is the multilinear Sobolev smoothing inequality.  Theorem~\ref{weyl} plays an essential role in our arguments and will be used repeatedly. Here, in contrast to \cite{KMT}, we need this theorem for $\K=\Z$ and $\K=\R$. In \cite{KMT}, a bilinear variant of \eqref{eq:43} for the Furstenberg--Weiss averages and $\K=\Z$ was established with $\delta^c + \langle \Log N \rangle^{-c\widetilde C}$ in place of $\delta^c+N^{-c}$; see \eqref{eq:300}. Although in \cite{KMT} the decay $\delta^c + \langle \Log N \rangle^{-c\widetilde C}$ was sufficient, in our situation we need $\delta^c+N^{-c}$ in \eqref{eq:43}. This bound can be derived using Theorem~\ref{thm:IW}. We abbreviate $\tilde A^{P_1,\ldots, P_k}_{N; \K^k}$ to $\tilde A_{N; \K^k}$.
	
	\subsection{Structural theorem for adjoints}
	We now prove a structural theorem for the adjoint multilinear operators $\tilde{A}_{N; \K^k}^{*j}(f_1,\ldots, f_k)$ from \eqref{eq:108b}. The key tool will be the following variant of the Hahn--Banach theorem. 
	
	\begin{lemma}\label{hahn} Let $(X, \mathcal B(X), \mu)$ be a $\sigma$-finite measure space. Given $A,B > 0$, $G\in L^2(X)$, and $\Phi\subseteq L^2(X)$, suppose that the following inverse theorem holds: if $f \in L^2(X)$ satisfies $\|f\|_{L^\infty(X)}\le 1$ and $|\langle f, G \rangle| > A$, then $|\langle f, \phi \rangle| > B$ for some $\phi \in \Phi$. Then $G$ lies in the closed convex hull of the set
		\begin{equation*}
		V \coloneqq \big\{ \lambda\phi\in L^2(X): \phi \in \Phi \text{ and } |\lambda| \leq A/B \big\}
		\cup \big\{ h \in L^2(X): \|h\|_{L^1(X)} \leq A \big\}.
		\end{equation*}
	\end{lemma}
	
	\begin{proof}
		We refer to \cite[Lemma~6.9]{KMT}.
	\end{proof}
	
	Lemma~\ref{hahn} combined with the inverse theorem from Section~\ref{sec:inverse} will reveal a major and minor arc structure for the multilinear operators $\tilde{A}_{N;\K^k}^{*j}(f_1,\ldots, f_k)$. In the case $\K=\R$, it still makes sense to refer to the major and minor arcs, even though we have exactly one major arc around $0$. The following theorem will play an essential role in the proof of Theorem~\ref{weyl}.
	
	\begin{theorem}\label{struct}
		Let $\K$ be either $\Z$ or $\R$. Fix $C_0' \in \Z_+$ and $k\in\Z_+$, and let $\mathcal P = (P_1,\ldots, P_k)$ be a polynomial mapping satisfying \eqref{eq:42}--\eqref{eq:44}. Then there exists a large constant $C' \in\Z_+$, depending only on $\mathcal P$ and $C_0'$, such that the following holds. Assume that $\delta_0 \in (0,1]$ and $N \geq C' \delta_0^{-C'}$, and fix $j \in [k]$. Then for all $1$-bounded functions $f_1,\ldots,f_k\in L^{\infty}(\K^k)$ supported on $\prod_{i\in[k]}[ \pm N_i]_{\K}$, where $N_i \coloneqq C_0' N^{d_i}$ for each $i\in[k]$, one can decompose
		\begin{equation}\label{decomp}
		\tilde{A}_{N; \K^k}^{*j}(f_1,\ldots,  f_k)
		= F + E,
		\end{equation} 
		where $F \in L^2(\K^k)$ and $\mathcal F_{j, \K^k}F$ is supported on the set
		\begin{align}
		\label{eq:108}
		\prod_{i\in[j-1]}[\pm C' N_i]_{\K}\times
		\mathfrak M_{\le C'\delta_0^{-C'}N^{-d_j}}\Big(\mathcal R_{\le C'\delta_0^{-C'}}^{\K}\Big)
		\times \prod_{i\in[k]\setminus[j]}[\pm C' N_i]_{\K}.
		\end{align}
		Furthermore, we have
		\begin{equation}\label{faq-bound}
		\| F\|_{L^\infty(\K^k)} \le C' \delta_0^{-C'}
		\quad\text{and}\quad
		\|F\|_{L^1(\K^k)} \le C' \delta_0^{-C'}N^{D},
		\end{equation}
		where $D \coloneqq d_1+\cdots+d_k$, while the error term satisfies
		\begin{equation}\label{e-bound}
		\|E\|_{L^1(\K^k)} \le C_1' \delta_0  N^D.
		\end{equation}
		The same conclusion holds for $A_{N; \K^k}^{*j}$ in place of $\tilde{A}_{N; \K^k}^{*j}$ in \eqref{decomp}.
	\end{theorem}
	
	In the case $\K=\Z$, by \eqref{eq:108}, the $j$-th Fourier transform of $F$ in the sum on the right-hand side of \eqref{decomp} is supported on the major arcs corresponding to the canonical fractions $\mathcal R_{\le N_1}$ with $N_1 = C' \delta_0^{-C'}$. This part can be thought of as the major arcs part of $\tilde{A}_{N; \Z^k}^{*j}(f_1,\ldots, f_k)$, whereas the error function $E$ from \eqref{decomp} is negligible thanks to \eqref{e-bound}. From Theorem~\ref{minorstruct} below it will become clear that $E$ can be thought of as the minor arcs part of $\tilde{A}_{N; \Z^k}^{*j}(f_1,\ldots, f_k)$.
	
	\begin{proof}[Proof of Theorem~\ref{struct}]
		The proof is a combination of Lemma~\ref{hahn} and Theorem~\ref{pp}; see \cite[Corollary~6.10]{KMT} or \cite[Corollary~7.13]{KPMW} for details.
	\end{proof}
	
	\subsection{Proof of multilinear Weyl inequality}
	Here we will use the Ionescu--Wainger projections $\Pi^j_{\K} [\leq M_1, \leq M_2]$ from \eqref{eq:132} with $M_1 \coloneqq \delta^{-C}$ and $M_2 \coloneqq \delta^{-C}N^{-d_j}$, where $\delta\in (0, 1)$ is small and $C \in \R_+$ is large depending on $C_1, C_2$ from Theorem~\ref{weyl}. We make the following important remark.
	
	\begin{remark}
		\label{remark:2}
		Given a large integer ${p_0} \in 2\Z$, there exists a large constant $C_{p_0} \ge 1$, possibly depending on $C$ but independent of $N$ and $\delta$, such that if
		\begin{align}
		\label{eq:133}
		N \ge C_{p_0}  \delta^{-C_{p_0} },
		\end{align}
		then for each fixed $\varepsilon \in (0, 1)$ and each $f\in L^{q} (\K^k)$ with $q \in [p_0',p_0]$, we have
		\begin{align}
		\label{eq:148}
		\left\|\Pi_{\K}^j\left[\le \delta^{-C}, \le N^{-d_j}\delta^{-C}\right]f\right\|_{L^{q} (\K^k)}\lesssim \delta^{-\varepsilon C}\|f\|_{L^{q} (\K^k)}.
		\end{align}
		
		\begin{itemize}
			\item[(i)] If $\K=\R$, then $\mathcal R_{\le M_1}^{\R}=\{0\}$ and so \eqref{eq:148} holds for all $q\in [1,\infty]$ by Minkowski's convolution inequality with $\|\eta\|_{L^1(\R)}$ in place of $\delta^{-\varepsilon C}$.
			\item[(ii)] If $\K=\Z$, then the situation is more difficult, but  condition \eqref{eq:133} allows us to use  Theorem~\ref{thm:IW} to derive inequality  \eqref{eq:148}.
		\end{itemize}
		Moreover, for intervals $I,J \subset \R$ of finite length and $f$ as before, we have
		\begin{align}
		\label{eq:151}
		\left\|\Pi_{\K}^j\left[\le \delta^{-C}, \le N^{-d_j}\delta^{-C}\right]f\right\|_{L^{q} (J_j)}
		\lesssim  \delta^{-\varepsilon C} \bigg\langle \frac{\mathrm{dist}(I,J)}{N^{d_j}\delta^{C} } \bigg\rangle^{-10} \|f\|_{L^{q} (I_j)},
		\end{align}
		where $I_j \coloneqq \K^{j-1}\times (I\cap \K)\times\K^{k-j}$ and $J_j \coloneqq \K^{j-1}\times (J\cap\K)\times\K^{k-j}$.
		\begin{itemize}
			\item[(iii)] If $\K=\R$, then \eqref{eq:151} follows, since $|\F_{\R}^{-1}\eta_{[\le M_2]}(x)|\lesssim_{\eta}
			M_2 \langle M_2 x\rangle^{-20}$.
			
			\item[(iv)]  If $\K=\Z$, then
			\eqref{eq:151} follows by repeating the proof of \cite[Lemma~5.17, p.~1045]{KMT}, with Theorem~\ref{thm:IW} in place of \cite[Theorem~5.7, p.~1041]{KMT}.
		\end{itemize}
	\end{remark}
	
	The next theorem combined with the Ionescu--Wainger projections \eqref{eq:132} will explain why the error term $E$ from \eqref{decomp} satisfying \eqref{e-bound} can be thought of as the minor arcs piece of the adjoint operator $\tilde{A}_{N; \K^k}^{*j}(f_1,\ldots, f_k)$. In what follows, we choose a very large integer $p_0\in 2\Z$, depending on $k, \mathcal P, p_1, \dots, p_k, p, C$ but independent of $N, \delta$, such that all interpolation arguments using \eqref{eq:148} or \eqref{eq:151} from Remark~\ref{remark:2} will refer only to $q \in [p_0',p_0]$.          
	
	\begin{theorem}
		\label{minorstruct}
		Let $\K$ be either $\Z$ or $\R$. Fix $C\in\R_+$ and $k\in\Z_+$, and let $\mathcal P = (P_1,\ldots, P_k)$ be a polynomial mapping satisfying \eqref{eq:42}--\eqref{eq:44}. Let $1<p_1,\ldots, p_k<\infty$ be exponents such that $\frac{1}{p_1}+\cdots+\frac{1}{p_k}=\frac{1}{p}< 1$. Then there exists a small constant $c_0 \in (0,1)$, depending only on $k, \mathcal P, p_1, \dots, p_k, p, C$, such that the following holds. Assume that $\delta \in (0,1]$ and $N \geq 1$. Let $f_i \in L^{p_i}(\K^k)$ for all $i\in[k]$. If \eqref{eq:133}
		holds, then for each $j \in [k]$ we have
		\begin{align}
		\label{eq:109}
		\Big\| \big(\mathds{I} - \bm \Pi_{\K}^j\big) \tilde{A}_{N; \K^k}^{*j}(f_1,\ldots, f_k)\Big\|_{L^p(\K^k)}
		\le c_0^{-1} \delta^{c_0}\prod_{i\in[k]}\|f_i\|_{L^{p_i}(\K^k)},
		\end{align}
		where $\mathds{I}$ and $\bm \Pi_{\K}^j$ are the identity operator and $\Pi_{\K}^j[\le \delta^{-C}, \le N^{-d_j}\delta^{-C}]$, respectively. The same conclusion holds for $A_{N; \K^k}^{*j}$ in place of $\tilde{A}_{N; \K^k}^{*j}$ in \eqref{eq:109}.
	\end{theorem}
	
	\begin{proof}
		The proof is fairly involved. We will adapt the arguments from \cite{KMT}. Due to subtle differences between our proof and the proof of \cite[Theorem~5.12(i)]{KMT} we present the key details. We assume that $\delta$ is small; otherwise \eqref{eq:109} follows from \eqref{eq:148} and H{\"o}lder's inequality. We allow explicit and implicit constants to depend on $k, \mathcal P, p_1, \dots, p_k, p, C$ but not on $N, \delta$. 
		
		\paragraph{\bf Step~1} Fix $C_0'' \in \Z_+$. We show, with a small constant $c_1\in (0, 1)$, the estimate
		\begin{align}
		\label{eq:139}
		\Big\| \big(\mathds{I} - 
		\bm \Pi_{\K}^j \big) \tilde{A}_{N; \K^k}^{*j}(f_1,\ldots, f_k) \Big\|_{L^2(\K^k)}
		\lesssim \delta^{c_1}N^{D/2} \prod_{i\in[k]}\|f_i\|_{L^\infty(\K^k)}
		\end{align}
		for any $1$-bounded $f_1,\ldots, f_k\in L^\infty(\K^k)$ supported on $\prod_{i\in[k]}[\pm C_0'' N^{d_i}]_{\K}$.
		
		Consider $C'$ from Theorem~\ref{struct} applied with $C_0' = C_0''$, and choose $\delta_0 = \delta^{\gamma}$ for some $\gamma\in(0, 1)$ such that $\gamma C' < C / 2$, where $\delta,C$ are as in \eqref{eq:139}. Then
		\begin{align}
		\label{eq:110}
		\tilde{A}_{N; \K^k}^{*j}(f_1,\ldots, f_k)=F + E,
		\end{align}
		where $F, E$ are as in Theorem~\ref{struct} and satisfy \eqref{eq:108}--\eqref{e-bound} with $\delta_0, C'$ as above. If $\delta$ is small, then $\gamma C' < C / 2$ implies $4 C' \delta^{-\gamma C'} < \delta^{-C}$ so that by the support properties of $\mathcal F_{j, \K^k}F$, see \eqref{eq:108}, we have $(\mathds{I} - \bm \Pi_{\K}^j) F\equiv 0$. Using \eqref{eq:110} and the latter properties of $F$ we conclude that the proof of \eqref{eq:139} can be reduced to showing the following estimate
		\begin{align}
		\label{eq:141}
		\|(\mathds{I} - \bm \Pi_{\K}^j)
		E\|_{L^2(\K^k)}
		\lesssim \delta^{c_1}N^{D/2}.
		\end{align}
		Note that
		$\|\tilde{A}_{N; \K^k}^{*j}(f_1,\ldots, f_k)\|_{L^\infty(\K^k)} \lesssim 1$. By \eqref{decomp}, \eqref{faq-bound}, and \eqref{e-bound} we have
		\begin{align}
		\label{eq:142}
		\|E\|_{L^\infty(\K^k)}\lesssim \delta^{-C'\gamma}
		\quad \text{and} \quad
		\|E\|_{L^1(\K^k)}\lesssim \delta^{\gamma} N^D.
		\end{align}
		
		Now if $p_0\in2\Z$ is very large, so that $p_0' \in (1,2)$ is sufficiently
		close to $1$, then by interpolating the bounds in \eqref{eq:142}, we
		obtain 
		\begin{align}
		\label{eq:143}
		\|E\|_{L^{p_0'}(\K^k)}\lesssim \delta^{\gamma/2} N^{D/p_0'}.
		\end{align}
		By our choice of $\delta$ in \eqref{eq:133} and \eqref{eq:148} with some fixed $\varepsilon\in(0, 1)$, we have
		\begin{align*}
		\left\| \big(\mathds{I} - \bm \Pi_{\K}^j \big)
		\tilde{A}_{N; \K^k}^{*j}(f_1,\ldots, f_k)\right\|_{L^{p_0}(\K^k)}
		\lesssim \delta^{-\varepsilon C}N^{D/p_0},
		\end{align*}
		since $\tilde{A}_{N; \K^k}^{*j}(f_1,\ldots, f_k)$ vanishes outside $\prod_{i\in[k]}[\pm O(N_i)]_{\K}$. Thus, we have
		\begin{align}
		\label{eq:144}
		\|
		(\mathds{I} - \bm \Pi_{\K}^j)
		E\|_{L^{p_0}(\K^k)}\lesssim \delta^{-\varepsilon C}N^{D/p_0}.
		\end{align}
		Using \eqref{eq:148} again, this time together with \eqref{eq:143}, we
		obtain
		\begin{align}
		\label{eq:145}
		\|(\mathds{I} - \bm \Pi_{\K}^j)
		E\|_{L^{p_0'}(\K^k)}
		\lesssim \delta^{-\varepsilon C}\delta^{\gamma/2} N^{D/p_0'}.
		\end{align}
		Interpolating \eqref{eq:144} with \eqref{eq:145} and adjusting $\varepsilon\in(0, 1)$, we establish \eqref{eq:141} with $c_1 = \frac{\gamma}{5}$. Consequently, \eqref{eq:139} follows. 
		
		\paragraph{\bf Step~2} We now relax the $L^\infty(\K^k)$ control in \eqref{eq:139} to the control on $L^p(\K^k)$ spaces. More precisely, assuming \eqref{eq:133}, for fixed $1<p_1,\ldots, p_k<\infty$ such that $\frac{1}{p_1}+\cdots+\frac{1}{p_k}=\frac{1}{2}$ we show, with a small constant $c_2\in(0, 1)$, the estimate
		\begin{align}
		\label{eq:146}
		\left\| \big(\mathds{I} - \bm \Pi_{\K}^j \big)
		\tilde{A}_{N; \K^k}^{*j}(f_1,\ldots, f_k)\right\|_{L^2(\K^k)}
		\lesssim \delta^{c_2}\prod_{i\in[k]}\|f_i\|_{L^{p_i}(\K^k)}
		\end{align}
		for any $f_i\in L^{p_i}(\K^k)$ supported on $\prod_{i\in[k]}[\pm C_0'' N^{d_i}]_{\K}$ for all $i\in[k]$.
		
		The main tool in proving \eqref{eq:146} is the multilinear $L^p(\K^k)$-improving inequality proved in Section~\ref{sec:improving}, which establishes the estimate
		\begin{align*}
		\left\|\tilde{A}_{N; \K^k}^{*j}(g_1,\ldots,g_k) \right\|_{L^2(\K^k)}
		\lesssim N^{-D(\frac{1}{r}-\frac{1}{2})} \prod_{i\in[k]} \|g_i\|_{L^{r_i}(\K^k)}
		\end{align*}
		and therefore
		\begin{equation}
		\label{eq:147}
		\left\| \big(\mathds{I} - \bm \Pi_{\K}^j \big)
		\tilde{A}_{N; \K^k}^{*j}(g_1,\ldots,g_k) \right\|_{L^2(\K^k)}
		\lesssim N^{-D(\frac{1}{r}-\frac{1}{2})} \prod_{i\in[k]} \|g_i\|_{L^{r_i}(\K^k)}
		\end{equation}
		holds for some $1<r<2$ and $1<r_1,\ldots, r_k<\infty$ such that $\frac{1}{r_1}+\cdots+\frac{1}{r_k}=\frac{1}{r}$; see Corollary~\ref{cor:improvZ}. Here, it is crucial that $1 < r < 2$. The idea of utilizing $L^p(\K^k)$-improving inequalities in this context was initiated in \cite{KMT}. Although the bilinear case $k=2$ of \eqref{eq:147}, used in \cite{KMT}, can be deduced from the linear $L^p(\K^k)$-improving inequality \cite{HKLMY}, the general case $k\ge3$ requires a genuine multilinear $L^p(\K^k)$-improving inequality established in Section~\ref{sec:improving}.
		
		We interpolate the bounds in \eqref{eq:139} and \eqref{eq:147} taking $\theta \coloneqq \frac{r}{2}<1$ and $q_i$ such that $\frac{1}{q_i} = \frac{\theta}{r_i}+\frac{1-\theta}{\infty}$ for each $i\in[k]$ so that $\frac{1}{q_1}+\cdots+\frac{1}{q_k}=\frac{\theta}{r}=\frac{1}{2}$. Since $-D(\frac{1}{r}-\frac{1}{2})\theta + \frac{D}{2}(1-\theta) = 0$,
		multilinear interpolation yields
		\begin{align*}
		\left\| \big(\mathds{I} - \bm \Pi_{\K}^j\big)
		\tilde{A}_{N; \K^k}^{*j}(f_1,\ldots, f_k)\right\|_{L^2(\K^k)}
		\lesssim
		\delta^{c_1(1-\theta)} 
		\prod_{i\in [k]} \|f_i\|_{L^{q_i}(\K^k)}.
		\end{align*}
		By Plancherel's theorem and H{\"o}lder's inequality, for all $1\le p_1,\ldots, p_k\le\infty$ such that $\frac{1}{p_1}+\cdots+\frac{1}{p_k}=\frac{1}{2}$ we also have
		\begin{align*}
		\left\| \big(\mathds{I} - \bm \Pi_{\K}^j\big)
		\tilde{A}_{N;\K^k}^{*j}(f_1,\ldots, f_k)\right\|_{L^2(\K^k)}
		\lesssim \prod_{i\in[k]}\|f_i\|_{L^{p_i}(\K^k)}.
		\end{align*}
		
		Interpolating the last two estimates yields \eqref{eq:146}, as desired. To remove the support condition in inequality  \eqref{eq:146} we use the off-diagonal decay estimate \eqref{eq:151} and proceed much the same way as in \cite[Corollary~6.24]{KMT}.
	\end{proof}
	
	Finally, we establish the multilinear Weyl inequality and the Sobolev smoothing inequality stated at the beginning of this section. The key tool will be Theorem~\ref{minorstruct}, which will give us the desired decay in $\delta$ and $N$.
	
	\begin{proof}[Proof of Theorem~\ref{weyl}]
		We can assume that $\delta^{-1}$ and $N$ are large; otherwise \eqref{eq:43} holds by H{\"o}lder's inequality. We consider separately $p>1$ and $p=1$.
		
		\paragraph{\bf Step~1} Let $p > 1$. Let
		$\F_{j, \K^k} f_j$ vanish on the major arcs
		${\mathfrak M}_{\leq N^{-d_j}\delta^{-C_2}}^j(\mathcal R_{\leq \delta^{-C_1}}^{\K})$
		for some $C_1, C_2\in\R_+$. Since $N^{-d_j}\delta^{-C} \leq N^{-d_j}\delta^{-C_2}$ and $\delta^{-C} \leq \delta^{-C_1}$ for $C \coloneqq \min\{C_1, C_2\}$, we see that $\F_{j, \K^k} f_j$ vanishes on ${\mathfrak M}_{\leq N^{-d_j}\delta^{-C}}^j(\mathcal R_{\leq \delta^{-C}}^{\K})$.
		
		\paragraph{\bf Case~1} Consider $N \ge C_{p_0} \delta^{-C_{p_0}}$ as in \eqref{eq:133}. By duality it suffices to prove
		\[
		| \langle \tilde A_{N; \K^k}(f_1,\ldots, f_k), h \rangle| \lesssim \delta^c \|h\|_{L^{p'}(\K^k)} \prod_{i\in[k]}\|f_i\|_{L^{p_i}(\K^k)}
		\]
		for $h \in L^{p'}(\K^k)$ and $f_i\in L^{p_i}(\K^k)$ for all $i\in [k]$, with $f_j \in L^{p_j}(\K^k)\cap L^2(\K^k)$ such that $\F_{j, \K^k} f_j$ vanishes on ${\mathfrak M}_{\leq N^{-d_j}\delta^{-C}}^j(\mathcal R_{\leq \delta^{-C}}^{\K})$. As before, we abbreviate $\Pi_{\K}^j[\le \delta^{-C}, \le N^{-d_j}\delta^{-C}]$ to $\bm \Pi_{\K}^j$. By \eqref{transpose} and properties of $\F_{j, \K^k} f_j$ we have
		\begin{align*}
		\left| \left\langle \tilde A_{N; \K^k}(f_1,\ldots, f_k), h \right\rangle\right|
		= \left|\left\langle \big(\mathds{I}-\bm \Pi_{\K}^j \big)
		\tilde{A}_{N; \K^k}^{*j}(f_1,\ldots, f_{j-1}, h,f_{j+1},\ldots, f_k), f_j \right\rangle\right|
		\end{align*}
		and the desired bound holds with $c=c_0$ by H\"older's inequality and \eqref{eq:109}.
		
		\paragraph{\bf Case~2} Consider
		$N \le C_{p_0} \delta^{-C_{p_0} }$. The proof follows from the previous case, with $\delta_0$ such that $\delta_0^{-1} = (N/C_{p_0})^{1/C_{p_0}} \le \delta^{-1}$, giving Theorem~\ref{weyl} for $p>1$. 
		
		\paragraph{\bf Step~2} Let $p=1$. In this case linear averages (with $k=1$) do not arise so we have $k\ge 2$. We assume that $\F_{j, \K^k} f_j$ vanishes on ${\mathfrak M}_{\leq N^{-d_j}\delta^{-C}}^j(\mathcal R_{\leq \delta^{-C}}^{\K})$ for some suitable $C \in\R_+$, as in Step~1. The key is to show that 
		\begin{align} \label{q<1} 
		\|{\tilde A}_{N;\K^k}(f_1,\ldots, f_k)\|_{L^q(\K^k)} \lesssim \prod_{i\in [k]} \|f_i\|_{L^{q_i}(\K^k)} 
		\end{align}
		for some choice of $1<q_1, \ldots, q_k < \infty$ with $q_j =2$ and $\frac{1}{q_1} + \cdots + \frac{1}{q_k} = \frac{1}{q}$ for $q<1$ close to $1$. Of course, if $q\ge 1$, then \eqref{q<1} holds trivially by H\"older's inequality but we need this for $L^q(\K^k)$ with $q<1$, which is not a Banach space. Interestingly, inequality  \eqref{q<1} follows by an application of the $L^p$-improving bound \eqref{adjoint-bound} from Corollary~\ref{cor:improvZ}. Using the multilinear interpolation result for $q<1$ in \cite[Theorem~1]{GrTao}, we can interpolate the bound \eqref{q<1} with the bounds \eqref{eq:43} in Theorem~\ref{weyl} for $p>1$ established in Step~1 to conclude
		\[
		\|{\tilde A}_{N;\K^k}(f_1, \ldots, f_k) \|_{L^1(\K^k)} \lesssim (\delta^c + N^{-c}) \prod_{i\in [k]} \|f\|_{L^{p_i}(\K^k)}
		\]
		for some $c>0$ and for all $1<p_1, \ldots, p_k < \infty$ such that $\frac{1}{p_1} + \cdots + \frac{1}{p_k} = 1$, as desired. This completes the proof of Theorem~\ref{weyl}.
	\end{proof}
	
	\section{Multilinear circle method: proof of pointwise ergodic theorem}
	\label{sec:ergodic}
	In this section, we develop a multilinear circle method in the context of quantitative multilinear pointwise ergodic theorems. Although our method is robust enough to handle oscillation, or even full $r$-variational estimates, we only illustrate its strength in the context of the long $r$-variational multilinear pointwise ergodic bound \eqref{eq:104}. This will be sufficient to deduce the Furstenberg--Bergelson--Leibman conjecture for commuting transformations along polynomials with distinct degrees; see Theorem~\ref{thm:main1} as well as Theorem~\ref{thm:main}.
	
	The bilinear circle method was developed for the first time in \cite{KMT}. Here we shall present a novel different approach relying on the following tools:
	
	\begin{enumerate}[label*={\arabic*}.]
		\item We use the Ionescu--Wainger projections from \eqref{eq:132} to make a distinction between minor and major arcs. These projections now correspond to the set of canonical fractions \eqref{eq:129}, which is a significant new input.
		
		\item An essential ingredient in obtaining the minor arc estimates is the multilinear Weyl inequality \eqref{eq:43} with a polynomial decay of the form $\delta^c+N^{-c}$; see Theorem~\ref{weyl}. This result represents a significant improvement over the logarithmic decay $\delta^c+\langle \log N\rangle^{-c\widetilde C}$ for the bilinear Furstenberg--Weiss averages from \cite{KMT}. The key tools we use to derive \eqref{eq:43} with a polynomial decay, which were not available in \cite{KMT}, are our Ionescu--Wainger multiplier theorem for the set of canonical fractions and multilinear $L^p$-improving inequality; see Theorem~\ref{thm:IW} and Theorem~\ref{thm:improv}, respectively.
		
		\item Another important ingredient is the Sobolev smoothing inequality, which is a continuous variant of the multilinear Weyl inequality; see Theorem~\ref{weyl}. The multilinear Sobolev smoothing inequality is utilized to obtain the major arc estimates in the high frequency case. In the bilinear situation the major arc estimates involving high frequency cases can be handled by appealing to the bilinear Weyl inequality \eqref{eq:43} and using a simple integration by parts argument. However, the general multilinear case necessitates the full strength of the Sobolev smoothing inequality as the integration by parts trick from \cite{KMT} is limited to the bilinear averages. This is a new phenomenon which was not apparent previously.
		
		\item The most challenging part is obtaining estimates for major arcs in the low frequency regime. A novel tool that we propose for controlling maximal functions is a powerful multiparameter norm interchanging inequality \eqref{eq:84}. It is important to note that polynomial decay in the multilinear Weyl inequality \eqref{eq:43} is crucial for dealing with this case; see Remark~\ref{rem:5} below. A one-parameter norm interchanging inequality originates in \cite{KMT} as a consequence of the Rademacher--Menshov inequality \eqref{eq:164}. Our multiparameter norm interchanging inequality allows us to replace the $p$-adic methods used in \cite{KMT} and significantly simplifies the argument. Moreover, we emphasize that the $p$-adic approach (especially sharp $L^p$-improving inequalities in the $p$-adic setting), even if adjusted to our context, would not be sufficient to handle the major arc estimates in the low frequency regime.
	\end{enumerate}
	
	Given $1<p_1,\ldots, p_k<\infty$ such that $\frac{1}{p_1}+\cdots+\frac{1}{p_k}=\frac{1}{p}\le 1$, we will prove
	\begin{align}
	\label{eq:127}
	\|{\bf V}^r(\tilde A_N(f_1,\ldots, f_k) :N\in\D) \|_{\ell^p(\Z^k)}
	\lesssim 1
	\end{align}
	for all finitely supported $f_1\in \ell^{p_1}(\Z^k),\ldots, f_k\in \ell^{p_k}(\Z^k)$, which are normalized so that $\|f_1\|_{\ell^{p_1}(\Z^k)}=\cdots=\|f_k\|_{\ell^{p_k}(\Z^k)}=1$. Then, by scaling and density arguments combined with the monotone convergence theorem, we readily obtain inequality \eqref{eq:104} from Theorem~\ref{thm:main1} for all $f_1\in \ell^{p_1}(\Z^k),\ldots, f_k\in \ell^{p_k}(\Z^k)$. The proof of \eqref{eq:127} is long and intricate. We begin with making a distinction between minor and major arcs and for this purpose we will use the Ionescu--Wainger projections from \eqref{eq:132} with dyadic parameters.
	
	\subsection{Dyadic Ionescu--Wainger projections}
	Let $\eta \colon \R\to[0, 1]$ be a smooth and even function satisfying \eqref{eq:126}. For any $n, \xi\in\R$, we set
	\begin{align*}
	\eta_{\le n}(\xi)\coloneqq \eta_{[\le 2^n]}(\xi) = \eta(2^{-n}\xi).
	\end{align*}
	For any $l\in\N$, using definition of canonical fractions \eqref{eq:129}, we set
	\begin{align*}
	\Sigma_{\leq l}  \coloneqq  
	\mathcal R_{\le 2^l}
	\quad \text{and} \quad
	\Sigma_l  \coloneqq  \Sigma_{\leq l} \setminus \Sigma_{\leq l-1}
	\end{align*}
	with $\Sigma_{\leq-1}$ denoting the empty set. Then we have
	\begin{align}
	\label{eq:373}
	\# \Sigma_{\leq l} \le 2^{2l}.
	\end{align}
	Similarly, for any $i\in[k]$, $l\in\N$, and $m\in\Z$, using definitions \eqref{eq:78} and \eqref{eq:79}, we introduce the dyadic ``major arc'' by
	\begin{align*}
	{\mathcal M}_{\leq l, \leq m}  \coloneqq
	\mathfrak M_{\le 2^m}(\Sigma_{\le l})
	\quad \text{and} \quad
	{\mathcal M}_{\leq l, \leq m}^i \coloneqq
	\mathfrak M_{\le 2^m}^i(\Sigma_{\le l}).
	\end{align*}
	We note that ${\mathcal M}_{\le l, \leq m} \subseteq {\mathcal M}_{\le l', \leq m'}$ whenever both $l \leq l'$ and $m \leq m'$ hold. Moreover, if
	$m \leq - 2l - 2$, then the arcs
	$[\theta - 2^m, \theta + 2^m]$ that comprise
	${\mathcal M}_{\leq l, \leq m}$ are pairwise disjoint. We also define
	\begin{align*}
	{\mathcal M}_{l, \leq m}  \coloneqq  {\mathcal M}_{\leq l, \leq m} \setminus {\mathcal M}_{\leq l-1, \leq m}
	\quad \text{and} \quad
	{\mathcal M}_{l,m} \coloneqq  {\mathcal M}_{l,\leq m} \setminus {\mathcal M}_{l,\leq m-1}.
	\end{align*}
	
	Referring to \eqref{eq:132}, we define the dyadic Ionescu--Wainger projections $\Pi_{\leq l, \leq m}^i \colon \ell^2(\Z^k) \to \ell^2(\Z^k)$ by setting 
	\begin{align}
	\label{eq:20}
	\Pi_{\leq l, \leq m}^i f(x) \coloneqq \Pi^i_\Z [\le 2^l, \le 2^m]f(x) = T_{i, \Z^k}^{\Sigma_{\le l}}[\eta_{\le m}]f(x),\qquad  x\in\Z^k,
	\end{align}
	for each $f \in \ell^2(\Z^k)$. In particular, $\Pi_{\leq l, \leq m}^i$ is a self-adjoint (and real symmetric) operator, which is the identity operator except in the $i$-th variable.
	
	\begin{remark}
		The following properties are clear from \eqref{eq:20}.
		\label{rem:2}
		\begin{enumerate}[label*={(\roman*)}]
			\item\label{IW4} The function $\mathcal F_{i, \Z^k} (\Pi_{\leq l, \leq m}^i f)$ is supported on the set ${\mathcal M}_{\leq l, \leq m}^i$, and if $\mathcal F_{i, \Z^k} f$ vanishes on that set, then $\Pi_{\leq l, \leq m}^i f \equiv 0$.
			\item\label{IW5} If $m \leq - 2l - 2$ and $\mathcal F_{i, \Z^k} f$ is supported on ${\mathcal M}_{\leq l, \leq m-2}^i$, then $\Pi_{\leq l, \leq m}^i f = f$.
			\item\label{IW6} If $m \leq - 2l - 2$, then $\Pi_{\leq l, \leq m}^i$ is a contraction on $\ell^2(\Z^k)$.
		\end{enumerate}
	\end{remark}
	
	A consequence of Theorem~\ref{thm:IW} is the following important bound.
	
	\begin{proposition}
		\label{mult} For every $k\in\Z_+$ and $p\in(1, \infty)$ there exists a constant $C_p\in\R_+$ such that for every $l\in\N$, $m \in \Z$, and $i\in[k]$, if
		\begin{equation}\label{mult-hyp}
		m \leq - 6 \max\{p, p'\}(l+1),
		\end{equation}
		then 
		\begin{align}
		\label{eq:116}
		\| \Pi_{\leq l, \leq m}^i f \|_{\ell^p(\Z^k)} \le C_{p} \left(2^{C_p \frac{l \log \log l}{\log l}}\ind{l\ge 10}+\ind{l<10}\right)\|f\|_{\ell^p(\Z^k)}.
		\end{align}
	\end{proposition}
	\begin{proof}
		The proof of \eqref{eq:116} is a direct consequence of Theorem~\ref{thm:IW}.
	\end{proof}
	
	\subsection{Minor arc estimates} As in Section~\ref{sec:inverse}, we fix a polynomial mapping
	\begin{align*}
	\mathcal P \coloneqq (P_1,\ldots, P_k) \colon \R\to \R^k,
	\end{align*}
	with $P_1,\ldots, P_k\in\Z[{\rm n}]$ such that $d_1 \coloneqq \deg P_1<\cdots < d_k \coloneqq \deg P_k$ and $D \coloneqq d_1+\cdots+d_k$; see \eqref{eq:42}--\eqref{eq:44}. We also fix $1<p_1,\ldots, p_k<\infty$ such that $\frac{1}{p_1}+\cdots+\frac{1}{p_k}=\frac{1}{p}\le 1$ and choose $p_0\in 2\Z_+$ satisfying
	\begin{align}
	\label{eq:33}
	1<p_0'<p_1,\ldots, p_k<p_0<\infty.
	\end{align}
	The parameter $p_0$ from \eqref{eq:33} will serve as an upper bound for $\max\{p,p'\}$ from condition \eqref{mult-hyp}  to ensure that Proposition~\ref{mult} can be applied. One can think that the parameter $p_0$ from \eqref{eq:33} is a sufficiently large even integer, which, if necessary, may be further adjusted depending on the interpolation arguments that will be used throughout this section. 
	
	Given a small fixed parameter 
	\begin{align*}
	0<\alpha<(10^6 kd_kp_0)^{-1}
	\end{align*}
	and a large constant $C_0\in \Z_+$ depending on $p_0, \mathcal P, \alpha$, we define the quantities 
	\begin{align*}
	l_{(N)} \coloneqq \Log N^\alpha \quad \text{and} \quad L_{(N)} \coloneqq \Log N-l_{(N)}
	\end{align*}
	for any $N\ge C_0$, using the log-scale notation \eqref{log-scale}. We take $C_0 \geq 2^{5/\alpha}$ so that $l_{(N)} \in \N_{\geq 5}$ for $N\ge C_0$. For each $i\in [k]$ and $\varepsilon_i \in \{0,1\}$, define
	\begin{align*}
	\Pi_{\leq l_{(N)}, \leq -d_i L_{(N)}}^{i, \varepsilon_i} f_i \coloneqq 
	\begin{cases}
	\Pi_{\leq l_{(N)}, \leq -d_i L_{(N)}}^{i} f_i& \text{ if } \varepsilon_i=1,\\
	f_i -\Pi_{\leq l_{(N)}, \leq -d_i L_{(N)}}^{i} f_i& \text{ if } \varepsilon_i=0.
	\end{cases}
	\end{align*}
	
	Observe that there exists a small absolute constant $c\in(0, 1)$ such that
	\begin{align}
	\label{eq:11}
	\begin{gathered}
	\| \tilde A_N(\Pi_{\leq l_{(N)}, \leq -d_1L_{(N)}}^{1, \varepsilon_1}f_1,\ldots, \Pi_{\leq l_{(N)}, \leq -d_kL_{(N)}}^{k, \varepsilon_k}f_k) \|_{\ell^p(\Z^k)} 
	\lesssim N^{-c} 
	\end{gathered}
	\end{align}
	whenever $(\varepsilon_1,\ldots, \varepsilon_k)\neq(1,\ldots, 1)$.
	Indeed, assuming $\varepsilon_j = 0$ for some $j \in [k]$ we apply Theorem~\ref{weyl} with $C_1=1$, $C_2=d_j$, and $\delta=2^{5-l_{(N)}} \in (0,1]$, and then \eqref{eq:11} follows from Proposition~\ref{mult}. Writing $\tilde A_N(f_1,\ldots, f_k)$ as 
	\begin{gather*}
	\sum_{(\varepsilon_1,\ldots, \varepsilon_k)\in\{0, 1\}^k}
	\tilde A_N\left(\Pi_{\leq l_{(N)}, \leq -d_1L_{(N)}}^{1, \varepsilon_1}f_1,\ldots, \Pi_{\leq l_{(N)}, \leq -d_kL_{(N)}}^{k, \varepsilon_k}f_k\right)
	\end{gather*}
	we see that, in view of \eqref{eq:11}, the proof of \eqref{eq:127} is reduced to estimating
	\begin{align}
	\label{eq:12}
	\left\|{\bf V}^r\left(\tilde A_N\left(\Pi_{\leq l_{(N)}, \leq -d_1L_{(N)}}^{1}f_1,\ldots, \Pi_{\leq l_{(N)}, \leq -d_kL_{(N)}}^{k}f_k\right):N\in\D\right) \right\|_{\ell^p(\Z^k)} 
	\end{align}
	for all lacunary sequences $\mathbb D \subset [C_0, \infty)$ and normalized functions $f_1, \dots, f_k$, where the functions $\Pi_{\leq l_{(N)}, \leq -d_1L_{(N)}}^{1}f_1,\ldots, \Pi_{\leq l_{(N)}, \leq -d_kL_{(N)}}^{k}f_k$ are all restricted to major arcs.
	
	\subsection{Major arc estimates} 
	\label{sec:4}
	Our aim is to utilize certain model operators, which will be more manageable. We begin with trimming the size of denominators. To this end, we rewrite the average from \eqref{eq:12} as
	\begin{align}
	\label{eq:maj}
	\sum_{l_1\in\N_{\le l_{(N)}}} \dots \sum_{l_k\in\N_{\le l_{(N)}}}
	\tilde A_N\left(\Pi_{l_1, \leq -d_1L_{(N)}}^{1}f_1,\ldots, \Pi_{l_k, \leq -d_kL_{(N)}}^{k}f_k\right),
	\end{align}
	where $\Pi_{l_i, \leq -d_i L_{(N)}}^{i} \coloneqq T_{i, \Z^k}^{\Sigma_{l_i}}[\eta_{\le -d_i L_{(N)}}]$ for every $i\in[k]$. Next, we match the scales of the bumps $\eta_{\leq -d_i \Log N  + d_i l_{(N)}}$ appearing in \eqref{eq:maj} by splitting them into nonoscillatory and highly oscillatory pieces
	\[
	\eta_{\leq -d_i\Log N}+\sum_{s_i \in [l_{(N)}]}\big(\eta_{\le -d_i (\Log N  - s_i )}-\eta_{\le -d_i (\Log N  - s_i+1)}\big).
	\]
	For $i \in [k]$ we define 
	\begin{align}
	\label{eq:34}
	\eta_{N}^{i, s_i} \coloneqq 
	\begin{cases}
	\eta_{\leq -d_i(\Log N  - s_i)}-\eta_{\leq -d_i(\Log N  - s_i+1)} & \text{ if } s_i>0,\\
	\eta_{\leq -d_i\Log N} & \text{ if } s_i=0,\\
	\end{cases}
	\end{align}
	and set
	\begin{align}
	\label{eq:31}
	\tilde{\Pi}_{l_i,s_i}^{i,N} \coloneqq T_{i, \Z^k}^{\Sigma_{l_i}}[\eta_{N}^{i, s_i}].
	\end{align}
	Using \eqref{eq:34}, \eqref{eq:31}, and the previous identities, we write \eqref{eq:maj} as
	\begin{gather*}
	\sum_{l_1, s_1\in\N_{\le l_{(N)}}} \dots \sum_{l_k, s_k\in\N_{\le l_{(N)}}}
	\tilde A_N\left(\tilde{\Pi}_{l_1,s_1}^{1,N}f_1,\ldots, \tilde{\Pi}_{l_k,s_k}^{k,N}f_k\right).
	\end{gather*}
	Set $\D_{l, s} \coloneqq \{N\in\D: N \geq C_0 \text{ and } l_{(N)}\ge \max\{l, s\}\}$, where
	\begin{align}
	\label{eq:38}
	l \coloneqq \max\{l_i:i\in[k]\} \quad \text{and} \quad s \coloneqq \max\{s_i:i\in[k]\}.
	\end{align}
	It suffices to show, for any $s_1, l_1,\ldots, s_k, l_k\in\N$, that
	\begin{align}
	\label{eq:15}
	\begin{gathered}
	\left\|{\bf V}^r\left(\tilde A_N\left(\tilde{\Pi}_{l_1,s_1}^{1,N}f_1,\ldots, \tilde{\Pi}_{l_k,s_k}^{k,N}f_k\right):N\in\D_{l, s}\right) \right\|_{\ell^p(\Z^k)} 
	\lesssim 2^{-c (l+s)},
	\end{gathered}
	\end{align}
	since summing \eqref{eq:15} over $s_1, l_1,\ldots, s_k,l_k$ yields the desired bound for \eqref{eq:12}. 
	
	\subsubsection{\textbf{Major arc approximations}} We now take advantage of the fact that each function in \eqref{eq:15} is restricted to major arcs, and approximate the average in \eqref{eq:15} by a certain model operator. We first fix some notation.
	
	For $N\ge 1$ we define the exponential sum 
	\begin{align}
	\label{eq:18}
	m_N(\xi) \coloneqq \frac{1}{|(N/2, N]\cap\Z|}\sum_{n\in (N/2, N]\cap\Z}e(\xi\cdot \mathcal P(n)), \qquad \xi\in\T^k,
	\end{align}
	and its continuous counterpart 
	\begin{align}
	\label{eq:6}
	\mathfrak{m}_N(\xi) \coloneqq 2\int_{1/2}^1e(\xi\cdot \mathcal P(Nt))dt, \qquad \xi\in\R^k.
	\end{align}
	Furthermore, given $a\in\Z^k$ and $q\in\Z_+$ such that $\gcd(a, q)=1$, we consider the corresponding complete exponential sum
	\begin{align}
	\label{eq:10}
	G\left(\frac{a}{q}\right) \coloneqq \frac{1}{q}\sum_{n \in [q]} e\left(\frac{a}{q}\cdot \mathcal P(n)\right).
	\end{align}
	We also define the set of rational fractions 
	\begin{align}
	\label{eq:19}
	\Sigma_{l_1,\ldots, l_k} \coloneqq \left\{\frac{a}{q}\in (\Q / \Z)^k: \gcd(a, q)=1  \text{ and } \frac{a}{q}\in \Sigma_{l_1}\times\cdots\times \Sigma_{l_k}\right\},
	\end{align}
	which will be a natural domain for $G$ in our applications. 
	
	\begin{lemma}
		\label{lemma:1}
		Let $\mathcal P \colon \R\to\R^k$ be a polynomial mapping as in \eqref{eq:42} whose components are polynomials with integer coefficients and satisfy \eqref{eq:40}. Then there exists a constant $C_{\mathcal P} \in \R_+$ such that for every $N \geq 1$, $M_{1},\ldots, M_{k} \in \R_+$, and $l_1,\ldots, l_k\in\N$ the following holds. For every $\xi\in\T^k$ and $\theta\in \Sigma_{l_1,\ldots, l_k}$ such that $|\xi_i-\theta_i|\le M_i^{-1}$ for all $i\in[k]$, one has
		\begin{align}
		\label{eq:32}
		|m_N(\xi)-G(\theta)\mathfrak m_N(\xi-\theta)|\le C_{\mathcal P}2^{kl} (\max\{M_i^{-1}N^{d_i-1}:i\in[k]\}+N^{-1})
		\end{align}
		with $m_N, \mathfrak m_N, G, \Sigma_{l_1,\ldots, l_k}$ defined in \eqref{eq:18}--\eqref{eq:19}.
	\end{lemma}
	
	\begin{proof}
		It suffices to prove \eqref{eq:32} for $N\ge 2^{kl}$, because the opposite case $N<2^{kl}$ is trivial. Fix $\theta=\frac{a}{q}\in\Sigma_{l_1,\ldots, l_k}$ and observe that
		\begin{align*}
		m_N(\xi)=\frac{1}{q}\sum_{r \in [q]} e(\theta\cdot \mathcal P(r))\cdot q\sum_{n\in \Z}e((\xi-\theta)\cdot \mathcal P(qn+r))\chi_N(qn+r),
		\end{align*}
		where $\chi_N(n) \coloneqq |(N/2, N]\cap\Z|^{-1}\ind{(N/2, N]}(n)$. We note that
		\[
		\Big|\sum_{n\in \Z}e((\xi-\theta)\cdot \mathcal P(qn+r))\chi_N(qn+r)
		-\int_{\mathbb R}e((\xi-\theta)\cdot \mathcal P(qt+r))\chi_N(q \lfloor t \rfloor +r)dt\Big|
		\]
		is bounded by $C_{\mathcal P} \max\{M_i^{-1} N^{d_i}:i\in[k]\}$ by the mean value theorem. Moreover, we have $\int_{\mathbb R} |\chi_N(q\lfloor t \rfloor +r) - \chi_N(qt+r)| dt \leq 10N^{-1}$. Using this and the restriction $q \leq 2^{kl}$, we obtain \eqref{eq:32}.
	\end{proof}
	
	Next, given a finite subset $\Theta \subset (\Q/\Z)^k$ and functions $m \colon \T^k\to \C$ and $S  \colon (\Q/\Z)^k\to \C$, we will be working with multilinear operators of the form
	\begin{align}
	\label{eq:26}
	B^{\Theta}[S; m](g_1,\ldots, g_k)(x) \coloneqq 
	\sum_{\theta\in\Theta}S(\theta)\sum_{y\in\Z^k}K_{\tau_\theta m}(y)\prod_{i \in [k]} g_i(x-(y\cdot e_i)e_i)
	\end{align}
	for $x\in\Z^k$, where $\tau_\theta m(\xi) \coloneqq m(\xi-\theta)$ and   
	\begin{align}
	\label{eq:27}
	K_m(y) \coloneqq \int_{\T^k}m(\xi)e(-\xi\cdot y)d\xi,
	\end{align}
	and, as before, $\{e_i :i\in [k]\}$ is the standard basis in $\R^k$.
	
	Using \eqref{eq:26} and \eqref{eq:27}, we prove the following technical but very useful lemma, which will often be applied in tandem with
	Lemma~\ref{lemma:1}.
	
	\begin{lemma}
		\label{lemma:2}
		Let $N\ge1$,  $d_1, \dots, d_k \in \Z_+$, and $ M_{1},\ldots,M_{k} \in [3, \infty)$, and let ${\bm m}_N \colon \T^k\to \mathbb C$ be a smooth function supported on $\prod_{i\in[k]}[\pm M_i^{-1}]$ such that
		\begin{align}
		\label{eq:183}
		\left\|\partial_1^{\beta_1}\cdots \partial_k^{\beta_k}{\bm m}_N\right\|_{L^{\infty}(\T^k)}\lesssim_{\beta_1, \ldots, \beta_k}
		N^{\beta_1d_1+\cdots+\beta_kd_k}
		\end{align}
		for all $(\beta_1,\ldots, \beta_k)\in \N^k$. Then, for every $\gamma\in(0, 1)$ and $\beta\in \Z_+$, we have 
		\begin{align}
		\label{eq:184}
		|K_{{\bm m}_N}(y)|\le C_{\beta, \gamma} M_1^{-1}\cdots M_k^{-1} \|{\bm m}_N\|_{L^{\infty}(\T^k)}^{1-\gamma}
		\prod_{i \in [k]} \langle y_i /N^{d_i}  \rangle^{-\gamma \beta}
		\end{align}
		with some constant $C_{\beta, \gamma} \in \R_+$. In particular, if $\gamma\beta>1$, then \eqref{eq:184} implies 
		\begin{align}
		\label{eq:185}
		\|K_{{\bm m}_N}\|_{\ell^1(\Z^k)} \lesssim_{\beta, \gamma}
		(N^{d_1}M_1^{-1})\cdots (N^{d_k}M_k^{-1})\|{\bm m}_N\|_{L^{\infty}(\T^k)}^{1-\gamma}.
		\end{align}
	\end{lemma}
	
	\begin{proof}
		Using \eqref{eq:183}, the support condition for ${\bm m}_N$, and integration by parts (sufficiently many times, say $\beta\in \N$, variable by variable) in the integral in \eqref{eq:27} with $m={\bm m}_N$, we obtain \eqref{eq:184} with $\gamma=1$. A convex combination of this bound and the trivial bound $|K_{{\bm m}_N}(y)|\lesssim_{\beta} M_1^{-1}\cdots M_k^{-1} \|{\bm m}_N\|_{L^{\infty}(\T^k)}$ yields \eqref{eq:184} for $\gamma \in (0,1)$, and \eqref{eq:185} follows whenever $\gamma \beta > 1$.
	\end{proof}
	
	Now we can formulate our first approximation result. 
	
	\begin{proposition}
		\label{prop:1}
		For every $N\ge1$ and $l_1, s_1,\ldots, l_k, s_k\in\N$ such that $l_{(N)}\ge \max\{l, s\}$ and for all $f_1\in \ell^{p_1}(\Z^k),\ldots, f_k\in \ell^{p_k}(\Z^k)$ with $1<p_1,\ldots, p_k<\infty$ such that $\frac{1}{p_1} + \dots + \frac{1}{p_k}=\frac{1}{p}\le 1$ the following estimate holds
		\begin{align}
		\label{eq:28}
		\|E_{l,s}^N (f_1,\ldots, f_k)\|_{\ell^p(\Z^k)}
		\lesssim N^{-\frac{9}{10}}
		\prod_{i \in [k]} \|f_i\|_{\ell^{p_i}(\Z^k)}
		\end{align}
		with $E_{l,s}^N (f_1,\ldots, f_k)$ denoting the error term
		\[
		\tilde A_N\left(\tilde{\Pi}_{l_1,s_1}^{1,N}f_1,\ldots, \tilde{\Pi}_{l_k,s_k}^{k,N}f_k\right)
		-B^{\Sigma_{l_1,\ldots, l_k}}[G; \mathfrak m_Nw_N^{s_1,\ldots, s_k}](f_1,\ldots, f_k),
		\]
		where $w_N^{s_1,\ldots, s_k} \coloneqq \bigotimes_{i \in [k]} \eta_{N}^{i, s_i}$ and $\eta_{N}^{i, s_i}$ is defined in \eqref{eq:34}. 
	\end{proposition}
	
	\begin{proof}
		By homogeneity, we can assume that $\| f_i \|_{\ell^{p_i}(\Z_k)} = 1$ for all $i \in [k]$. 
		
		\paragraph{\bf Step~1} First, we rewrite $\tilde A_N(\tilde{\Pi}_{l_1,s_1}^{1,N}f_1,\ldots, \tilde{\Pi}_{l_k,s_k}^{k,N}f_k)(x)$ as
		\begin{align*}
		\sum_{\theta\in\Sigma_{l_1,\ldots, l_k}} \sum_{y\in\Z^k}K_{m_N\tau_\theta w_N^{s_1,\ldots, s_k}}(y)\prod_{i \in [k]} f_i(x-(y\cdot e_i)e_i).
		\end{align*}
		By the triangle inequality and Minkowski's integral inequality followed by H{\"o}lder's inequality, the left-hand side of \eqref{eq:28} is controlled by
		\begin{align*}
		\sum_{\theta\in\Sigma_{l_1,\ldots, l_k}}
		\left\|K_{m_N\tau_\theta w_N^{s_1,\ldots, s_k}} -G(\theta)K_{\tau_\theta (\mathfrak m_Nw_N^{s_1,\ldots, s_k})}\right\|_{\ell^1(\Z^k)}.
		\end{align*}
		\paragraph{\bf Step~2}
		Fix $\theta = \frac{a}{q} \in\Sigma_{l_1,\ldots, l_k}$ and expand the above norm as follows
		\begin{align*}
		\sum_{y\in\Z^k}
		\left| \int_{\T^k}(m_N(\xi)-G(\theta)\mathfrak m_N(\xi-\theta))w_N^{s_1,\ldots, s_k}(\xi-\theta)e(-\xi\cdot y)d\xi \right|.
		\end{align*}
		Using Lemma~\ref{lemma:1} with $M_i\simeq N^{d_i}2^{-d_is_i}$, together with the support condition for $w_N^{s_1,\ldots, s_k}$ and the restriction $l_{(N)}\ge \max\{l, s\}$, we obtain
		\begin{align}
		\label{eq:29}
		\left\|{m_N\tau_\theta w_N^{s_1,\ldots, s_k}} -G(\theta){\tau_\theta (\mathfrak m_Nw_N^{s_1,\ldots, s_k})}\right\|_{L^{\infty}(\T^k)}\lesssim 2^{2k{d_k}l_{(N)}} N^{-1}.
		\end{align}
		
		\paragraph{\bf Step~3} Now by \eqref{eq:29} and Lemma~\ref{lemma:2}, for any fixed $\gamma\in(0, 1)$, we have
		\begin{align*}
		\left\|K_{m_N\tau_\theta w_N^{s_1,\ldots, s_k}} -G(\theta)K_{\tau_\theta (\mathfrak m_Nw_N^{s_1,\ldots, s_k})}\right\|_{\ell^1(\Z^k)}\lesssim
		2^{3k{d_k}l_{(N)}} N^{-(1-\gamma)}.
		\end{align*}
		Using \eqref{eq:373} we obtain that $\#\Sigma_{l_1,\ldots, l_k}\lesssim 2^{2kl}$, which, combined with the last estimate, yields \eqref{eq:28} upon choosing the parameter $\gamma$ sufficiently small. This completes the proof of Proposition~\ref{prop:1}.
	\end{proof}
	
	Now, in light of Proposition~\ref{prop:1}, instead of \eqref{eq:15} it suffices to prove 
	\begin{align}
	\label{eq:37}
	\left\|{\bf V}^r\left(B^{\Sigma_{l_1,\ldots, l_k}}[G; \mathfrak m_Nw_N^{s_1,\ldots, s_k}](f_1,\ldots, f_k):N\in\D_{l, s}\right) \right\|_{\ell^p(\Z^k)} 
	\lesssim 2^{-c (l+s)}
	\end{align}
	for all normalized functions $f_i \in \ell^{p_i}(\Z^k)$ and all $s_1, l_1,\ldots, s_k,l_k\in\N$. 
	
	\subsubsection{\textbf{Model operators}}
	One difficulty in proving \eqref{eq:37} is the need to obtain the decay in both parameters $l$ and $s$ from \eqref{eq:38}. For this purpose, we distinguish two cases, which will be handled separately later on:
	\begin{itemize}
		\item[(i)] the \textit{high frequency} case $s > Cl$;
		\item[(ii)] the \textit{low frequency} case $Cl \geq s$;
	\end{itemize}
	with some large $C\ge1$ to be specified later. For each fixed $t \in [\frac12, 1]$, using the bumps $\eta_{N}^{i, s_i}$ from \eqref{eq:34}, we define
	\begin{align}
	\label{eq:41}
	\eta_{N, t}^{i, s_i}(\xi_i) \coloneqq e(\xi_i P_i(Nt))\eta_{N}^{i, s_i}(\xi_i), \qquad i\in[k].
	\end{align}
	Using \eqref{eq:41}, we obtain the following identity 
	\begin{align*}
	\mathfrak m_N(\xi)w_N^{s_1,\ldots, s_k}(\xi)=
	2\int_{1/2}^1 \Big( \bigotimes_{i\in[k]}\eta_{N, t}^{i, s_i}\Big) (\xi)dt.
	\end{align*}
	Define
	\begin{align}
	\label{eq:111}
	u\coloneqq 
	\begin{cases}
	100  k (s+1) & \text{ if } s > Cl,\\
	100  k(l+1) & \text{ if } Cl \ge s.
	\end{cases}
	\end{align}
	Note that for any $N\in\D_{l, s}$ we have $N\ge \max\{2^{\max\{l, s\}/\alpha}, C_0\}$, which implies 
	\begin{align}
	\label{eq:117}
	N\ge 2^{10^5  kd_kp_0\max\{l, s\}}C_0^{1/2}\ge 2^{100p_0d_ku},
	\end{align}
	since $C_0 \ge 2^{5/\alpha} \geq 2^{10^6  kd_kp_0}$. By \eqref{eq:117}, we rewrite $B^{\Sigma_{l_1,\ldots, l_k}}[G; \mathfrak m_Nw_N^{s_1,\ldots, s_k}](f_1,\ldots, f_k)(x)$ as
	\begin{align}
	\label{eq:51}
	2 \int_{1/2}^1 B^{\Sigma_{l_1,\ldots, l_k}}[G; \eta_{u}^*]\left(T_{1, \Z^k}^{\Sigma_{l_1}}[\eta_{N, t}^{1, s_1}]f_1,\ldots, T_{k, \Z^k}^{\Sigma_{l_k}}[\eta_{N, t}^{k, s_k}]f_k\right)(x)dt,
	\end{align}
	where $\eta_{u}^* \coloneqq \bigotimes_{i\in[k]}\eta_{\le -p_0 d_i u}$. Introducing new functions 
	\begin{align}
	\label{eq:50}
	F_{N,t}^{i, l_i, s_i} \coloneqq 
	F_{N, t}^{i, l_i, s_i}(f_i) \coloneqq T_{i, \Z^k}^{\Sigma_{l_i}}[\eta_{N, t}^{i, s_i}]f_i, \qquad i\in[k],
	\end{align}
	we will show that the multilinear operator
	\[
	B^{\Sigma_{l_1,\ldots, l_k}}[G; \eta_{u}^*](F_{N, t}^{1, l_1, s_1},\ldots, F_{N, t}^{k, l_k, s_k})\]
	from \eqref{eq:51} can be replaced by the following model multilinear operator
	\[
	\tilde{A}_{2^u}\left(\Pi_{l_1,\le -p_0d_1u}^1(F_{N, t}^{1, l_1, s_1}),\ldots, \Pi_{l_k,\le -p_0d_ku}^k(F_{N, t}^{k, l_k, s_k})\right),
	\]
	with a satisfactory error term.
	\begin{proposition}
		\label{prop:2}
		Fix $l_1,s_1, \ldots, l_k, s_k\in\N$ and $u\in\Z_+$ as in \eqref{eq:111}. Then for every $N\ge 1$ and for all $f_{N}^1\in \ell^{p_1}(\Z^k),\ldots, f_{N}^k\in \ell^{p_k}(\Z^k)$ with $1<p_1,\ldots, p_k<\infty$ satisfying $\frac{1}{p_1} + \cdots + \frac{1}{p_k} =\frac{1}{p}\le 1$  the following estimate holds
		\begin{align*}
		\big\|{\bf V}^r(\tilde E_{l,s}^u(f_{N}^1,\ldots, f_{N}^k): N\in\D_{l, s})\big\|_{\ell^p(\Z^k)}
		\lesssim 2^{-\frac{9u}{10}}\prod_{i\in[k]} \|{\bf V}^r(f_{N}^i:N\in\D_{l, s})\|_{\ell^{p_i}(\Z^k)}
		\end{align*}
		with $\tilde E_{l,s}^u(f_{N}^1,\ldots, f_{N}^k)$ denoting the error term
		\[
		\tilde{A}_{2^u}(\Pi_{l_1,\le -p_0d_1u}^1(f_{N}^1),\ldots, \Pi_{l_k,\le -p_0d_ku}^k(f_{N}^k))
		-B^{\Sigma_{l_1,\ldots, l_k}}[G; \eta_{u}^*](f_{N}^1,\ldots, f_{N}^k).
		\]
	\end{proposition}
	
	\begin{proof} 
		The proof is much the same as the proof of Proposition~\ref{prop:1}.
	\end{proof}
	
	\subsubsection{\textbf{Basic $r$-variational estimates}}
	
	We now gather certain simple $r$-variational estimates that will be useful for our arguments. The key tools will be the Ionescu--Wainger and $r$-variational Ionescu--Wainger theorems; see Theorem~\ref{thm:IW} and Theorem~\ref{thm:IWvar}, respectively.
	
	The following result is a simple exercise that uses the properties of the Fourier transform.\begin{lemma}
		\label{lem:3}
		Let $\phi \colon \R\to \R$ be a compactly supported function with continuous derivative.
		Then there exists a constant $C_{\phi}\in\R_+$ depending on $\phi$ such that for all $n \in \N$ we have
		\begin{align*}
		\|\partial^{n}(\mathcal F_{\R}^{-1}\phi)\|_{L^1(\R)}\le C_{\phi}^{n+1}.
		\end{align*}
	\end{lemma}
	
	\begin{lemma}
		\label{lemma:3}
		Let $1 < p_1, \dots, p_k < \infty$ and $r\in (2, \infty]$. Then we have
		\begin{align}
		\label{eq:60}
		\begin{split}
		\big\|{\bf V}^r(\Pi_{l_i,\le -p_0d_i u}^i(F_{N, t}^{i, l_i, s_i}(f_i): N\in\D_{l, s})\big\|_{\ell^{p_i}(\Z^k)}
		& \lesssim_{p_i, r}C_{p_i}(l, s_i),\\
		\big\|{\bf V}^r(F_{N, t}^{i, l_i, s_i}(f_i): N\in\D_{l, s})\big\|_{\ell^{p_i}(\Z^k)}
		& \lesssim_{p_i, r} 
		C_{p_i}(l, s_i),
		\end{split}
		\end{align}
		for all $i \in [k]$ and all choices of $l_1, s_1, \dots, l_k, s_k \in\N$, uniformly in $t \in [\frac{1}{2},1]$, with normalized $f_i \in \ell^{p_i}(\Z^k)$ and $F_{N, t}^{i, l_i, s_i}(f_i)$ as in \eqref{eq:50}, where
		\[
		C_{p_i}(l, s_i) \coloneqq (s_i + 1) (2^{O_{p_i}(\frac{l \log \log l }{\log l})}\ind{l\ge10}+\ind{l<10}).
		\]
	\end{lemma}
	
	\begin{proof}
		By the definition of $u$ in \eqref{eq:111}, one has $6(l_i+1) \max\{p_i, p_i'\}\le p_0d_iu$ for each $i\in[k]$ so that \eqref{mult-hyp} is satisfied. Consequently, \eqref{eq:116} implies
		\begin{align}
		\label{eq:85}
		\| \Pi_{l_i,\le -p_0d_i u}^i f \|_{\ell^{p_i}(\Z^k)} \lesssim 2^{O_{p_i} (\frac{l \log \log l}{\log l})}\ind{l\ge10}+\ind{l<10}.
		\end{align}
		Therefore, by \eqref{eq:85} it suffices to show the second inequality in \eqref{eq:60}, since
		\[
		\Pi_{l_i,\le -p_0d_iu}^i(F_{N, t}^{i, l_i, s_i}(f_i))=
		F_{N, t}^{i, l_i, s_i}(\Pi_{l_i,\le -p_0d_iu}^if_i).
		\]
		
		If $s_i>0$, then by \eqref{eq:34} and \eqref{eq:41} we dominate the $r$-variational norm by the square function (in fact, the shifted square function; see \cite[Theorem~B.1, p.~1097]{KMT}) and the conclusion in \eqref{eq:60} follows by invoking Theorem~\ref{thm:IW}.
		
		If $s_i=0$, then we replace $\eta_{N, t}^{i, s_i}(\xi_i) = e(\xi_i P_i(Nt))\eta_{N}^{i, s_i}(\xi_i)$ with $\eta_{N}^{i, s_i}(\xi_i)$ and the error term is controlled by the standard Littlewood--Paley arguments (see, for example, \cite{DF} or \cite{MSZ2}) by Lemma~\ref{lem:3}, whereas the $r$-variational estimates for the latter multiplier follow from the L{\'e}pingle's inequality \cite{Lep} (see also \cite{MSZ1}) combined with Theorem~\ref{thm:IWvar} and Theorem~\ref{thm:IW}.
	\end{proof}
	
	\subsubsection{\textbf {Major arcs estimates: high frequency case $s > Cl$}} \label{sec:5}
	We can assume, without loss of generality, that $(s_1,\ldots,s_k)\neq(0,\ldots, 0)$; otherwise we are in the low frequency case. In the high frequency case, since $s > Cl$, it only suffices to obtain \eqref{eq:37} with a decay of the form $2^{-c s}$. Our aim will be to establish the following  maximal theorem.
	
	\begin{theorem}
		\label{thm:2}
		Fix $k\in\Z_+$ and let the exponents $2 \leq p_1,\ldots, p_k<\infty$ such that $\frac{1}{p_1}+\cdots+\frac{1}{p_k}=\frac{1}{p}\le 1$ be given. Then there exists a small constant $c\in(0, 1)$ such that for all normalized $f_1\in \ell^{p_1}(\Z^k),\ldots, f_k\in \ell^{p_k}(\Z^k)$ we have
		\begin{align}
		\label{eq:112}
		\left\|\sup_{N\in\D_{l, s}}\left|B^{\Sigma_{l_1,\ldots, l_k}}[G; \mathfrak m_Nw_N^{s_1,\ldots, s_k}](f_1,\ldots, f_k)\right|\right\|_{\ell^p(\Z^k)} 
		\lesssim 2^{-c s}
		\end{align}
		for all choices of $s_1, l_1,\ldots, s_k,l_k\in\N$ such that $s > Cl$ for some large $C \in \R_+$. 
	\end{theorem}
	
	In contrast to the  bilinear variant of Theorem~\ref{thm:2} from \cite{KMT}, where an integration by parts argument can be employed, the general multilinear case, due to its combinatorial nature, is more delicate. The new robust tool used here is the Sobolev smoothing inequality \eqref{eq:43} from Theorem~\ref{weyl} for $\K = \R$. Once Theorem~\ref{thm:2} is proved, it can be readily used to obtain \eqref{eq:37}. Indeed, we claim that for each fixed $\rho\in (2, \infty)$ and $1 < p_1,\ldots, p_k<\infty$ one has
	\begin{align}
	\label{eq:13}
	\left\|{\bf V}^{\rho}\left(B^{\Sigma_{l_1,\ldots, l_k}}\left[G; \mathfrak m_Nw_N^{s_1,\ldots, s_k}\right](f_1,\ldots, f_k):N\in\D_{l, s}\right) \right\|_{\ell^p(\Z^k)}
	\lesssim 2^{o(s)} 
	\end{align}
	with any normalized $f_1\in \ell^{p_1}(\Z^k),\ldots, f_k\in \ell^{p_k}(\Z^k)$. To prove this claim, we note that $u=100k(s+1)$ in the high frequency case, and using \eqref{eq:51} and Proposition~\ref{prop:2}, we reduce inequality \eqref{eq:13} to establishing
	\begin{align*}
	\left\|{\bf V}^{\rho}\left(\tilde{A}_{2^u}\left(F_{N, t, u}^{1, l_1, s_1},\ldots, F_{N, t, u}^{k, l_k, s_k}\right):N\in\D_{l, s}\right) \right\|_{\ell^p(\Z^k)}
	\lesssim 2^{o(s)},
	\end{align*}
	where
	\begin{align}
	\label{eq:35}
	F_{N, t, u}^{i, l_i, s_i} \coloneqq\Pi_{l_i,\le -p_0d_iu}^i(F_{N, t}^{i, l_i, s_i}(f_i)),
	\qquad i\in[k].
	\end{align}
	The latter inequality follows from \eqref{var-alg}, H{\"o}lder's inequality, and Lemma~\ref{lemma:3}. By interpolating \eqref{eq:13} and \eqref{eq:112}, or rather their unnormalized versions, and renormalizing the resulting bound, we obtain the desired claim in \eqref{eq:37}.
	
	The key ingredient in proving \eqref{eq:112} will be Theorem~\ref{weyl} for $\K = \R$. In order to apply Theorem~\ref{weyl}, we have to pass to the continuous setting. To this end, for $p \in [1,\infty)$, we define an extension operator $E \colon \ell^p(\Z^k)\to L^p(\R^k)$ by
	\begin{align}
	\label{eq:158}
	Ef(x)\coloneqq \sum_{y\in \Z^k}\ind{\mathcal C^k}(x-y)f(y), \qquad x\in \R^k,
	\end{align}
	for each $f\in \ell^p(\Z^k)$, where $\mathcal C^k\coloneqq [-1/2, 1/2)^k$. Observe that $Ef\equiv f$ on $\Z^k$ and $\|Ef\|_{L^{p}(\R^k)}=\|f\|_{\ell^{p}(\Z^k)}$. Thus, in particular, 
	$Ef\in L^{p}(\R^k)$. Moreover,
	$Ef(x+t)=Ef(x)$ holds for all $x\in \Z^k$ and $t\in \mathcal C^k$.
	
	\begin{proof}[Proof of Theorem~\ref{thm:2}]
		The proof is complicated and for the sake of clarity it will be divided into a few simpler blocks. We claim that it suffices to find a small constant $c'\in(0, 1)$ such that for every $\theta\in \Sigma_{l_1,\ldots, l_k}$ we have
		\begin{align}
		\label{eq:161}\left\|\sup_{N\in\D_{l, s}}\left|I_N^{\theta}(f_1,\ldots, f_k)\right|\right\|_{\ell^p(\Z^k)} 
		\lesssim 2^{-c' s},
		\end{align}
		where
		\begin{align*}
		I_N^{\theta}(f_1,\ldots, f_k)(x)\coloneqq \sum_{y\in \Z^k}  K_{\tau_{\theta}(\mathfrak m_Nw_N^{s_1,\ldots, s_k})}(y) \prod_{i \in [k]} f_i(x-(y\cdot e_i)e_i), \qquad x\in \Z^k.
		\end{align*}
		Indeed, since $s>Cl$ and $C\in\R_+$ may be as large as we need, by the triangle inequality and the estimate $\#\Sigma_{l_1,\ldots, l_k}\lesssim 2^{2kl}$ we obtain
		\begin{align*}
		\left\|\sup_{N\in\D_{l, s}}\left|B^{\Sigma_{l_1,\ldots, l_k}}[G; \mathfrak m_Nw_N^{s_1,\ldots, s_k}](f_1,\ldots, f_k)\right|\right\|_{\ell^p(\Z^k)} 
		\lesssim 2^{2kl}2^{-c' s}\lesssim 2^{-cs}.
		\end{align*}
		
		\paragraph{\bf Step~1}
		It suffices to prove \eqref{eq:161}. Fix $\theta\in \Sigma_{l_1,\ldots, l_k}$ and observe that
		\begin{align*}
		I_N^{\theta}(f_1,\ldots, f_k)(x)=e(-k\theta\cdot x)I_N^{0}(M_{\theta}f_1,\ldots, M_{\theta}f_k)(x),
		\end{align*}
		where $M_{\theta}$ is the modulation operator $M_{\theta}f(x)\coloneqq e(\theta\cdot x)f(x)$.
		
		In view of this identity, the proof of \eqref{eq:161} is reduced to showing
		\begin{align}
		\label{eq:162}
		\left\|\sup_{N\in\D_{l, s}}\left|I_N^{0}(g_1,\ldots, g_k)\right|\right\|_{\ell^p(\Z^k)} 
		\lesssim 2^{-c' s}
		\end{align}
		for any normalized $g_i\in \ell^{p_i}(\Z^k)$. We abbreviate $K_{\mathfrak m_Nw_N^{s_1,\ldots, s_k}}$ to $K_N$ and recall the convenient notation $x_{(y_i)}\coloneqq (x_1,\ldots, x_{i-1}, y_i, x_{i+1}, \ldots, x_k)$ for all $x, y\in\R^k$. Using the extension operator $E$ from \eqref{eq:158}, we can rewrite $I_N^{0}(g_1,\ldots, g_k)$ as $I_N^{0}(Eg_1,\ldots, Eg_k)$. Now, the $p$-th power of the norm in \eqref{eq:162} takes the form
		\begin{align}
		\label{eq:167}
		\sum_{x\in\Z^k}\int_{\mathcal C^k} \sup_{N\in\D_{l, s}}
		\Bigg|\sum_{y\in \Z^k} \int_{\mathcal C^k}  K_N(x-y) \prod_{i \in [k]} Eg_i((x+t)_{((y+t')_i)})dt'\Bigg|^pdt,
		\end{align}
		since $Eg_i(x_{(y_i)})=Eg_i((x+t)_{((y+t')_i)})$ for all $t, t'\in\mathcal C^k$ and $x, y\in\Z^k$. 
		
		As in Proposition~\ref{prop:1}, using Lemma~\ref{lemma:2} we can show the bound
		\begin{align}
		\label{eq:163}
		|K_N (x+z)-K_N(x)|\lesssim N^{-1/2}N^{-D}\prod_{i\in[k]}\langle x_i/N^{d_i}\rangle^{-2}
		\end{align}
		for all $x\in\Z^k$ and $z\in 2\mathcal C^k$. Therefore, denoting 
		\[
		\widetilde{K}_N(x, y, t, t') \coloneqq K_N((x+t)-(y+t'))-K_N(x-y),
		\]
		and appealing to \eqref{eq:163}, we readily obtain
		\begin{align*}
		\sum_{x\in\Z^k}\int_{\mathcal C^k} \sup_{N\in\D_{l, s}}
		\Bigg|\sum_{y\in \Z^k} \int_{\mathcal C^k}  \widetilde{K}_N(x, y, t, t') \prod_{i \in [k]} Eg_i((x+t)_{((y+t')_i)})dt'\Bigg|^pdt\lesssim 2^{-pc's}.
		\end{align*}
		Indeed, the above bound holds by Minkowski's integral inequality followed by H{\"o}lder's inequality, as in Proposition~\ref{prop:1}. Next, define
		\begin{align*}
		\tilde I_N^{0}(G_1,\ldots, G_k)(x)\coloneqq \int_{\R^k}  K_N(y) \prod_{i \in [k]} G_i(x-(y\cdot e_i)e_i)dy, \qquad x\in\R^k.
		\end{align*}
		After a change of variables
		$\|\sup_{N\in\D_{l, s}} |\tilde I_N^{0}(Eg_1,\ldots, Eg_k)\|_{L^p(\R^k)}^p$ is equal to
		\begin{align*}
		\sum_{x\in\Z^k}\int_{\mathcal C^k} \sup_{N\in\D_{l, s}}
		\Bigg|\sum_{y\in \Z^k} \int_{\mathcal C^k}  K_N((x+t)-(y+t')) \prod_{i \in [k]} Eg_i((x+t)_{((y+t')_i)})dt'\Bigg|^pdt.
		\end{align*}
		Using \eqref{eq:167} and gathering everything together, we conclude that
		\begin{align*}
		\left\|\sup_{N\in\D_{l, s}}\left|I_N^{0}(g_1,\ldots, g_k)\right|\right\|_{\ell^p(\Z^k)}\lesssim 2^{-c's}+
		\left\|\sup_{N\in\D_{l, s}}\left|\tilde I_N^{0}(Eg_1,\ldots, Eg_k)\right|\right\|_{L^p(\R^k)}.
		\end{align*}
		At this stage, Theorem~\ref{weyl} can be applied, and it remains to establish
		\begin{align}
		\label{eq:169}
		\left\|\sup_{N\in\D_{l, s}}\left|\tilde I_N^{0}(G_1,\ldots, G_k)\right|\right\|_{L^p(\R^k)}\lesssim 2^{-c's}
		\end{align}
		for any normalized $G_i\in L^{p_i}(\R^k)$, since $\|Eg_i\|_{L^{p_i}(\R^k)}=\|g_i\|_{\ell^{p_i}(\Z^k)}=1$.
		
		\paragraph{\bf Step~2}
		Define $V \coloneqq \{i\in[k]: s_i= 0\}$ and $v \coloneqq \# V < k$, and write 
		\[
		V=\{i_1,\ldots, i_v\}
		\quad \text{and} \quad 
		V^c=\{i_{v+1},\ldots, i_k\}.
		\]
		We interpret the polynomial mapping $\mathcal P$ from \eqref{eq:42} as a vector in $\R^k$ writing
		\[
		\mathcal P(t)=\mathcal P_V(t)+\mathcal P_{V^c}(t),
		\qquad t\in\R,
		\]
		where $\mathcal P_V(t) \coloneqq \sum_{i \in V} P_i(t)e_i\in \R^{k}$. For each $i_m\in V$, we can further write
		\[
		e(\xi_{i_m}P_{i_m}(Nt))
		=\sum_{n_m\in\N}\frac{(2\pi {\bm i} N^{d_{i_m}}\xi_{i_m})^{n_m}}{n_m!}  \frac{P_{i_m}(Nt)^{n_m}}{N^{d_{i_m}n_m}},
		\]
		using the Taylor expansion of the function $x \mapsto e(x)$. Since
		\[
		K_N(y)=
		\int_{\R^k}\mathfrak m_N(\xi)w_N^{s_1,\ldots, s_k}(\xi)e(-\xi \cdot y)d\xi,
		\]
		we write $K_N(y)$ as
		\begin{gather*}
		\sum_{(n_1,\ldots, n_v)\in\N^v} \int_{\R^k}\bigg(\prod_{m\in [v]} \frac{\tilde{\eta}_{N, n_m}^{i_m, 0}(\xi\cdot e_{i_m})}{n_m!} \bigg){\mathfrak m}_{N, V^c}^{n_1,\ldots, n_v}(\xi)w_{N, V^c}^{s_1,\ldots, s_k}(\xi) e(-\xi\cdot y) d\xi,
		\end{gather*}
		where for $\zeta \in \R$, $\xi\in\R^k$, and ${i_m}\in V$ we define
		\begin{align*}
		\tilde{\eta}_{N, n_m}^{i_m, 0}(\zeta) & \coloneqq (2\pi {\bm i} N^{d_{i_m}} \zeta)^{n_m}\eta_N^{i_m, 0}(\zeta), \\
		{\mathfrak m}_{N, V^c}^{n_1,\ldots, n_v}(\xi) & \coloneqq 2\int_{1/2}^1e(\xi\cdot \mathcal P_{V^c}(Nt)) B_N^t dt, \\
		w_{N, V^c}^{s_1,\ldots, s_k}(\xi) & \coloneqq \prod_{i\in V^c}\eta_{N}^{i, s_i}(\xi\cdot e_{i}),
		\end{align*}
		with $B_N^t= B_{N}^t(V;n_1,\dots, n_v)$ given by the formula
		\[
		B_N^t \coloneqq \prod_{m\in[v]} \frac{P_{i_m}(Nt)^{n_m}}{N^{d_{i_m}n_m}}.
		\]
		
		These Taylor expansions reveal that
		\[
		\tilde I_N^{0}(G_1,\ldots, G_k)(x)  = \int_{\R^k} K_N (y) \prod_{i \in [k]} G_i(x-(y\cdot e_i)e_i) dy,
		\]
		which we can further rewrite as
		\begin{align}
		\label{eq:120}
		\sum_{(n_1,\ldots, n_v)\in\N^v}
		\bigg(\prod_{m\in[v]} 
		\frac{\tilde G_{N,n_m}^{i_m, 0}(x)}{n_m!} \bigg) A_{N, V^c}^{n_1,\ldots, n_v}(G_{N}^{i_{v+1}, s_{i_{v+1}}},\ldots, G_{N}^{i_{k}, s_{i_{k}}})(x),
		\end{align}
		where
		\begin{align*}
		\tilde G_{N,n}^{i, 0} \coloneqq T_{{i}, \R^k}[\tilde \eta_{N,n}^{{i},0}]G_{i}
		\quad \text{and} \quad  
		G_{N}^{i, s_i} \coloneqq T_{{i}, \R^k}[\eta_{N}^{{i}, s_{i}}]G_{i},
		\end{align*}
		and
		\begin{align*}
		A_{N, V^c}^{n_1,\ldots, n_v}(G_{N}^{i_{v+1}, s_{i_{v+1}}},\ldots, G_{N}^{i_{k}, s_{i_{k}}})(x)\coloneqq 
		2\int_{1/2}^1
		\Big(\prod_{i\in V^c} G_{N}^{i, s_i}(x-P_{i}(Nt)e_i)\Big) B_N^tdt.
		\end{align*}
		
		\paragraph{\bf Step~3}
		Integrating the right-hand side above by parts, we obtain
		\begin{align}
		\label{eq:113}
		A_{N, V^c}^{n_1,\ldots, n_v}(G_{N}^{i_{v+1}, s_{i_{v+1}}},\ldots, G_{N}^{i_{k}, s_{i_{k}}})(x)
		=A_N^1(x) B_N^1-\int_{1/2}^1 A_N^t(x) \frac{\partial B_N^t}{\partial t} dt
		\end{align}
		with $A_{N}^t = A_{N}^t(V^c;G_{N}^{i_{v+1}, s_{i_{v+1}}},\ldots, G_{N}^{i_{k}, s_{i_{k}}})$ given by the formula
		\begin{align*}
		A_{N}^t(x) \coloneqq 2\int_{1/2}^t \prod_{i\in V^c}G_{N}^{i, s_i}(x-P_{i}(Ny)e_i) dy.
		\end{align*}
		Let $1\le q_V, q_{V^c} < \infty$ be defined by $\frac{1}{q_V} = \sum_{i\in V}\frac{1}{p_i}$ and  $\frac{1}{q_{V^c}} = \sum_{i\in V^c}\frac{1}{p_i}$ (if $V = \emptyset$, then we define only $q_{V^c}$). Now by the Sobolev smoothing inequality from Theorem~\ref{weyl} we obtain
		\begin{align}
		\label{eq:118}
		\left\|A_{N; \R^k}^t(G_{N}^{i_{v+1}, s_{i_{v+1}}},\ldots, G_{N}^{i_{k}, s_{i_{k}}})\right\|_{L^{q_{V^c}}(\R^k)}
		\lesssim 2^{-c's}
		\prod_{i\in V^c}\|G_{N}^{i, s_{i}}\|_{L^{p_i}(\R^k)}
		\end{align}
		uniformly in $t \in [1/2,1]$; here, it is crucial that $V^c$ is nonempty when $s > Cl$. By \eqref{eq:113} and \eqref{eq:118}, we obtain
		\begin{align*}
		\begin{gathered}
		\left\|A_{N, V^c}^{n_1,\ldots, n_v}(G_{N}^{i_{v+1}, s_{i_{v+1}}},\ldots, G_{N}^{i_{k}, s_{i_{k}}})\right\|_{L^{q_{V^c}}(\R^k)}
		\lesssim C_{n_1,\dots,n_v}(s) 
		\prod_{i\in V^c}\|G_{N}^{i, s_{i}}\|_{L^{p_i}(\R^k)}
		\end{gathered}
		\end{align*}
		with $C_{n_1,\dots,n_v}(s) \coloneqq C_{i_1}^{n_1}\cdots C_{i_v}^{n_v} 2^{-c's}$ for some large constants $C_{i_m} \in \R_+$.
		
		\paragraph{\bf Step~4} Now, we can bound $\|\sup_{N\in\D_{l, s}}|\tilde I_N^{0}(G_1,\ldots, G_k)|\|_{L^p(\R^k)}$ by
		\begin{align}
		\label{eq:121}
		\begin{split}
		& \sum_{(n_1,\ldots, n_v)\in\N^v} \frac{1}{n_1!\cdots n_v!}
		\Bigg\|\sup_{N\in\D_{l, s}}
		\Bigg|\prod_{m\in[v]}\tilde G_{N,n_m}^{i_m,0}\Bigg|\Bigg\|_{L^{q_V}(\R^k)} \\
		& \quad \times\Bigg\|\sup_{N\in\D_{l, s}}\Bigg|A_{N, V^c}^{n_1,\ldots, n_v}(G_{N}^{i_{v+1}, s_{i_{v+1}}},\ldots, G_{N}^{i_{k}, s_{i_{k}}})\Bigg|\Bigg\|_{L^{q_{V^c}}(\R^k)}
		\end{split}
		\end{align}
		using \eqref{eq:120} and H{\"o}lder's inequality; here, the first line in \eqref{eq:121} does not appear when $V = \emptyset$. Using Lemma~\ref{lem:3} and the properties of $\eta$, we obtain 
		\begin{align} \label{eq:123}
		\left\|\sup_{N\in\D_{l, s}}\left|\tilde G_{N,n_m}^{i_m, 0}\right|\right\|_{L^{p_{i_m}}(\R^k)}\lesssim
		C_{i_m}^{n_m}
		\end{align}
		for every $m\in [v]$. By \eqref{eq:123} and H{\"o}lder's inequality the first norm in \eqref{eq:121} is controlled by $C_{i_1}^{n_1}\cdots C_{i_v}^{n_v}$. By the final estimate from Step~3 and H{\"o}lder's inequality, the second norm in \eqref{eq:121} does not exceed 
		\begin{align} \label{eq:124}
		\begin{split}
		& \Big(\sum_{N\in\D_{l, s}}\big\|A_{N, V^c}^{n_1,\ldots, n_v}(G_{N}^{i_{v+1}, s_{i_{v+1}}},\ldots, G_{N}^{i_{k}, s_{i_{k}}})\big\|_{L^{q_{V^c}}(\R^k)}^{q_{V^c}}\Big)^{1/q_{V^c}}\\
		& \quad \lesssim C_{n_1,\dots,n_v}(s) \Big(\sum_{N\in\D_{l, s}}\prod_{i\in V^c}\| G_{N}^{i, s_{i}} \|_{L^{p_i}(\R^k)}^{q_{V^c}}\Big)^{1/q_{V^c}}\\
		& \quad \leq C_{n_1,\dots,n_v}(s) \prod_{i\in V^c}
		\Big\|\Big(\sum_{N\in\D_{l, s}}|G_{N}^{i,s_{i}}|^{p_i}\Big)^{1/p_i}
		\Big\|_{L^{p_i}(\R^k)} 
		\lesssim C_{n_1,\dots,n_v}(s),
		\end{split}
		\end{align}
		where, in the last line, we have used a standard square function argument \cite{DF} (see also \cite{MSZ2}), which is applicable because $p_i\ge2$. Finally, using \eqref{eq:121}--\eqref{eq:124} we deduce \eqref{eq:169}, since the factors $(C_{i_1}^{n_1}\cdots C_{i_v}^{n_v})^2(n_1!\cdots n_v!)^{-1}$ are summable over $(n_1,\ldots, n_v)\in\N^v$. This completes the proof of Theorem~\ref{thm:2}.
	\end{proof}
	
	\subsubsection{\textbf{Major arcs estimates: low frequency case $Cl \geq s$}}
	\label{sec:6}
	By \eqref{eq:111} in the low frequency case we have
	\begin{align}
	\label{eq:22}
	u \coloneqq 100k(l+1).
	\end{align}
	In this case, since $Cl \geq s$, the proof of \eqref{eq:37} follows from \eqref{eq:51} if for each fixed $r \in (2, \infty)$ we show the estimate
	\begin{align*}
	\left\|{\bf V}^r\left(B^{\Sigma_{l_1,\ldots, l_k}}[G; \eta_{u}^*](F_{N, t}^{1, l_1, s_1},\ldots, F_{N, t}^{k, l_k, s_k}) :N\in\D_{l, s}\right)\right\|_{\ell^p(\Z^k)}
	\lesssim (s+1)^{k} 2^{-c l}
	\end{align*}
	for all normalized $f_i\in \ell^{p_i}(\Z^k)$ and all choices of $s_1, l_1,\ldots, s_k,l_k\in\N$, uniformly in $t \in [\frac12, 1]$, with $F_{N, t}^{i, l_i, s_i} = F_{N, t}^{i, l_i, s_i}(f_i)$ given in \eqref{eq:50}. In view of \eqref{eq:35}, Proposition~\ref{prop:2}, and Lemma~\ref{lemma:3}, the above task reduces to showing
	\begin{align}
	\label{eq:61}
	\left\|{\bf V}^r\left(\tilde{A}_{2^u}(F_{N, t, u}^{1, l_1, s_1},\ldots, F_{N, t,u}^{k, l_k, s_k}):N\in\D_{l, s}\right)\right\|_{\ell^p(\Z^k)}
	\lesssim (s+1)^{k} 2^{-c l},
	\end{align}
	again for all normalized $f_i \in \ell^{p_i}(\Z^k)$ and all choices of  $l_1, s_1,\ldots, l_k, s_k\in\N$, uniformly in $t \in [\frac{1}{2},1]$, with $F_{N, t, u}^{i, l_i, s_i}$ given in \eqref{eq:35}.
	
	The advantage of working with the average $\tilde{A}_{2^u}$ in \eqref{eq:61} is that we can exploit the multilinear Weyl inequality \eqref{eq:43} from Theorem~\ref{weyl} to produce an exponential decay in $l$ as long as the underlying functions have Fourier transforms vanishing on suitable major arcs, like in \eqref{eq:61}. Observe that for each fixed $\rho \in (2,\infty)$, by H{\"o}lder's inequality and Lemma~\ref{lemma:3}, for all normalized $f_i \in \ell^{p_i}(\Z^k)$ and all choices of $l_1, s_1,\ldots, l_k, s_k\in\N$ one has
	\begin{align}
	\label{eq:66}
	\begin{split}
	& \sup_{1/2\le t\le 1}\left\|
	{\bf V}^\rho\left(\tilde{A}_{2^u}(F_{N, t, u}^{1, l_1, s_1},\ldots, F_{N, t,u}^{k, l_k, s_k}):N\in\D_{l, s}\right)\right\|_{\ell^p(\Z^k)}\\
	& \quad \lesssim \prod_{i \in [k]} \sup_{1/2\le t\le 1} \left\|{\bf V}^\rho\left(F_{N, t, u}^{i, l_i, s_i}: N\in\D_{l, s}\right)\right\|_{\ell^{p_i}(\Z^k)}\\
	& \quad \lesssim (s+1)^{k} \left(2^{O_{p_1, \dots, p_k} (\frac{l\log \log l}{\log l})}\ind{l\ge 10}+\ind{l<10}\right).
	\end{split}
	\end{align}
	
	Our aim will be to establish the following maximal theorem.
	
	\begin{theorem}
		\label{thm:1}
		Fix $k\in\Z_+$ and let the exponents $1<p_1,\ldots, p_k<\infty$ such that $\frac{1}{p_1}+\cdots+\frac{1}{p_k}=\frac{1}{p}\le 1$ be given. Then there exists a small constant $\delta\in(0, 1)$ such that for all normalized $f_1\in \ell^{p_1}(\Z^k),\ldots, f_k\in \ell^{p_k}(\Z^k)$ we have
		\begin{align*}
		\sup_{1/2\le t\le 1} \left\|\sup_{N\in\D_{l, s}}\left|\tilde{A}_{2^u}(F_{N, t, u}^{1, l_1, s_1},\ldots, F_{N, t,u}^{k, l_k, s_k})
		\right|\right\|_{\ell^p(\Z^k)}
		\lesssim (s+1)^{k} 2^{-\delta l}
		\end{align*}
		for all choices of $l_1, s_1,\ldots, k_k, s_k\in\N$, with $F_{N, t, u}^{i, l_i, s_i}$ given in \eqref{eq:35}.
	\end{theorem}
	Fixing $r\in (2, \infty)$ and $\rho \in (2,r)$, interpolating the unnormalized ${\bf V}^{\rho}$ bound from \eqref{eq:66} with the unnormalized maximal function bound from Theorem~\ref{thm:1}, and renormalizing, we readily deduce \eqref{eq:61}, as desired. 
	
	\subsubsection{\textbf{Multiparamater norm interchanging inequality}}
	Here we formulate a multiparameter variant of the norm interchanging inequality from \cite[Lemma~9.5, p.~1079]{KMT}. The norm interchanging inequality from \cite{KMT} asserts that for any measure space $(X, \mathcal B(X), \mu)$, any $K\in\mathbb Z_+$ and $1\le R< r\le \infty$, and for any sequence $(f_k)_{k\in[K]}\subseteq L^r(X)$ we have
	\begin{align}
	\label{eq:83}
	\|V^r(f_k: k\in[K])\|_{L^r(X)}\lesssim_{r, R} \|(f_k)_{k\in[K]}\|_{V^R([K]; L^r(X))},
	\end{align}
	where
	\[
	\|(f_k)_{k\in[K]}\|_{V^R([K]; L^r(X))} \coloneqq \sup_{J\in\mathbb Z_+}\sup_{k_0<\dots<k_J}\Big(\sum_{j=0}^{J-1}\|f_{k_{j+1}}-f_{k_j}\|_{L^r(X)}^R\Big)^{1/R}
	\]
	with the inner supremum taken over all increasing sequences from $[K]$ of length $J+1$. Inequality \eqref{eq:83} states that we can insert the $L^r(X)$ norm under the $r$-variation norm, at the cost of worsening the $r$-variation parameter to $R<r$.
	
	We begin with fixing some notation and terminology. Let $\mathbb J\subseteq\mathbb R_+$ be an at most countable set such that $\min\mathbb J\in \mathbb J$. Let ${\bf I}(\mathbb J)$ denote the set of all finite increasing sequences from $\mathbb J$. The elements of ${\bf I}(\mathbb J)$ will be denoted by bold letters ${\bf N}_{<J}=(N_j)_{j\in\N_{<J}}$ for some $J\in\Z_+$.
	
	Fix $k\in\mathbb Z_+$ and let $W_1,\ldots, W_k$ be linear normed spaces. Let $d\in\mathbb Z_+$ and $R_1, R_2\in[1, \infty]$. Suppose that $B \colon W_1\times\cdots\times W_k\to \ell^{R_1}(\mathbb Z^d)$ is a multisublinear operator. By this we mean that for every $f_1, g_1\in W_1,\ldots, f_k, g_k\in W_k$ and any $j\in[k]$ we have the pointwise inequality
	\begin{align*}
	& \big||B(f_1,\ldots,f_{j-1}, f_j, f_{j+1},\ldots,  f_k)|
	- |B(f_1,\ldots,f_{j-1}, g_j, f_{j+1},\ldots,  f_k)|\big| \\ 
	& \quad \le |B(f_1,\ldots,f_{j-1}, f_j-g_j, f_{j+1},\ldots,  f_k)|.
	\end{align*}
	
	For any $(J_1,\ldots, J_k)\in \mathbb Z_+^k$ and $\sigma=(\sigma_1,\ldots, \sigma_k)\in\{0, 1\}^k$ define
	\begin{align*}
	\mathbb N_{<J_i}^{\sigma_i} \coloneqq
	\begin{cases}
	\mathbb N_{<J_i}&\text{ if } \sigma_i=1,\\
	\{-1\}&\text{ if } \sigma_i=0.
	\end{cases}
	\end{align*}
	Moreover, for every $h_i=(h_{i, N})_{N\in \mathbb J}\subseteq W_i$ and $\mathbf N_{<J_i}^i=(N_{j}^i)_{j\in\mathbb N_{<J_i}}\in\mathbf I(\mathbb J)$ for $i\in[k]$ define $\Lambda_{{\mathbf N}_{<J_1}^1,\ldots, {\mathbf N}_{<J_k}^k}^{\sigma, R_1, R_2}(h_1,\ldots, h_k)$ by
	\begin{align*}
	\Big(\sum_{j_1\in \mathbb N_{<J_1}^{\sigma_1}, \ldots, j_k\in \mathbb N_{<J_k}^{\sigma_k}}
	\big\|B(\Delta h_{1, N_{j_1}^1}^{\sigma_1},\ldots, \Delta h_{k, N_{j_k}^k}^{\sigma_k})\big\|_{\ell^{R_1}(\mathbb Z^d)}^{R_2}\Big)^{1/{R_2}},
	\end{align*}
	where we set $N^i_{-1} \coloneqq \min\mathbb J\in \mathbb J$ and
	\begin{align*}
	\Delta h_{i, N_{j_i}^i}^{\sigma_i} \coloneqq
	\begin{cases}
	h_{i, N_{j_i}^i}-h_{i, N_{j_i-1}^i} &\text{ if } \sigma_i=1,\\
	h_{i, N^i_{-1}}&\text{ if } \sigma_i=0.
	\end{cases}
	\end{align*}
	\begin{proposition} 
		\label{nit}
		Let $\mathbb J\subset \mathbb R_+$ be at most countable with $N_{-1} \coloneqq \min\mathbb J\in \mathbb J$. Fix $k\in\mathbb Z_+$ and linear normed spaces $W_1,\ldots, W_k$. Let $d\in\mathbb Z_+$ and $1\le R<r\le\infty$. Suppose that $B \colon W_1\times\cdots\times W_k\to \ell^R(\mathbb Z^d)$ is a multisublinear operator. Then for every $h_i=(h_{i, N})_{N\in \mathbb J}\subseteq W_i$ for $i\in[k]$, we have
		\begin{align}
		\label{eq:84}
		\begin{split}
		& \big\|\sup_{N_1,\ldots, N_k\in\mathbb J}|B(h_{1, N_1},\ldots, h_{k, N_k})|\big\|_{\ell^r(\mathbb Z^d)}\\
		& \quad \lesssim_{r, R}\sum_{\sigma\in\{0, 1\}^k} \ \sup_{{\mathbf N}_{<J_1}^1,\ldots, {\mathbf N}_{<J_k}^k\in{\bf I}(\mathbb J)}
		\Lambda_{{\mathbf N}_{<J_1}^1,\ldots, {\mathbf N}_{<J_k}^k}^{\sigma, R, R}(h_1,\ldots, h_k).
		\end{split}
		\end{align}
	\end{proposition}
	\begin{proof}
		The proof follows from simple iterative applications of \eqref{eq:83}.
	\end{proof}
	
	The advantage of using Proposition~\ref{nit} is that we do not distinguish between small and large scales, which was an essential argument in \cite{KMT}. This is another important simplification of our method. Proposition~\ref{nit} will be the key tool in our approach and will be applied with $d=k$,  $B=\tilde A_{2^u}$, and $W_i=\ell^{R_i}(\mathbb Z^k)$ for $i\in[k]$, where $\frac{1}{R}=\frac{1}{R_1}+\cdots+\frac{1}{R_k}$.
	
	\subsubsection{\textbf{Uncertainty principle}}
	We will need a smooth partition of unity $(\phi_I: I\in\mathcal I_h)$ adapted to the partition $\mathcal I_h \coloneqq \{[lh, (l+1)h): l\in\Z\}$ of $\R$ into disjoint intervals $I$ of equal length $h \in \R_+$. In fact, we can take
	\begin{align}
	\label{eq:68}
	\phi_{I_l}(x) \coloneqq \phi(h^{-1}(x-c(I_l)))
	\end{align}
	for each $I_l \coloneqq [lh, (l+1)h)\in\mathcal I_h$, where $c(I_l) \coloneqq (l+\frac{1}{2})h$ is the center of $I_l$, and $\phi \colon \R\to[0, 1]$ is a suitable smooth bump function. Using \eqref{eq:68} we define
	\begin{align}
	\label{eq:70}
	\phi_I^{i}(x) \coloneqq \phi_I(x\cdot e_i), \qquad x\in\Z^k,
	\end{align}
	for every $i\in[k]$. We will need two technical lemmas. Here, the symbol $\delta_0$ stands for the Dirac delta distribution at the origin.
	
	\begin{lemma}
		\label{lemma:5}
		Fix $i\in[k]$, $\varepsilon\in(0, 1)$, and an interval $I$ of length comparable to $2^{d_i u}$. Then, for $M\in\N$ and $x\in\Z^k$, one has
		\begin{align}
		\label{eq:73}
		|\phi_{I}^i*_i(\delta_0-\mathcal F^{-1}_{\R}\eta_{\le -{(d_i-\varepsilon)u}})(x)|\lesssim_{\varepsilon, M} 2^{-Mu}
		\end{align}
		and, for $c(I)$ being the center of $I$, one also has
		\begin{align}
		\label{eq:74}
		|\phi_{I}^i*_i\mathcal F^{-1}_{\R}\eta_{\le -{(d_i-\varepsilon)u}}(x)|\lesssim_M\langle 2^{-d_i u}|x\cdot e_i-c(I)|\rangle^{-M}.
		\end{align}
		In particular, combining \eqref{eq:73} and \eqref{eq:74} one obtains
		\begin{align}
		\label{eq:75}
		|\phi_{I}^i*_i(\delta_0-\mathcal F^{-1}_{\R}\eta_{\le -{(d_i-\varepsilon)u}})(x)|\lesssim_{\varepsilon, M} 2^{-Mu} \langle 2^{-d_i u}|x\cdot e_i-c(I)|\rangle^{-M}.
		\end{align}
	\end{lemma}
	
	\begin{proof}
		Regarding \eqref{eq:73}, by the Poisson summation formula, we have 
		\begin{align*}
		|\mathcal F_{\Z}\phi_{I}(\xi)|
		= \Big| \sum_{n \in \Z} \mathcal F_{\R}^{-1} \phi_{I}(\xi+n) \Big|
		\lesssim_{N} \frac{2^{d_i u}}{(2^{d_i u}\|\xi\|)^{N}}
		\end{align*}
		for every $N\in\N$, with $\|\xi\| \coloneqq \min_{n\in\Z}|\xi+n|$. Observe that
		\begin{align*}
		\sup_{x\in\Z^k}|\phi_{I}^i*_i(\delta_0-\mathcal F^{-1}_{\R}\eta_{\le -{(d_i-\varepsilon)u}})(x)|
		\le \int_{\T}|\mathcal F_{\Z}\phi_{I}(\xi)||1-\eta_{\le -{(d_i-\varepsilon)u}}(\xi)|d\xi.
		\end{align*}
		Since $2^{d_i u}|\xi|\gtrsim 2^{\varepsilon u}$ whenever $1-\eta_{\le -{(d_i-\varepsilon)u}}(\xi) \neq 0$, applying the bound for $|\mathcal F_{\Z}\phi_{I}(\xi)|$ with some fixed $N\simeq (M+d_i)/\varepsilon$, we can control the integral above by $2^{-Mu}$, as desired. Regarding \eqref{eq:74}, we can control its left-hand side by
		\begin{align*}
		\sum_{y\in\Z} \left\langle 2^{-d_i u}|x\cdot e_i-y-c(I)|\right\rangle^{-M}2^{-(d_i-\varepsilon)u}\left\langle 2^{-(d_i-\varepsilon)u }|y|\right\rangle^{-2M}
		\end{align*}
		and use the fact that $|x\cdot e_i-c(I)|\le |x\cdot e_i-y-c(I)| + |y|$ for each $y \in \Z$.
	\end{proof}
	
	\begin{lemma}
		\label{lemma:4}
		Fix $i\in[k]$, $\varepsilon\in(0, 1)$, $f\in\ell^2(\Z^k)$, and an interval $I$. Then
		\begin{align}
		\label{eq:72}
		(\phi_{I}^i*_i\mathcal F^{-1}_{\R}\eta_{\le -{(d_i-\varepsilon)u}})
		\Pi_{l_i,\le -p_0d_iu}^i(f) \in \ell^2(\Z^k)
		\end{align}
		has its $i$-th Fourier transform supported on the major arcs $\mathcal M^i_{l_i, \le -(d_i-\varepsilon) u+1 }$. 
	\end{lemma}
	\begin{proof}
		Our claim readily follows, since the integral
		\begin{align*}
		\int_{\T}\sum_{\theta\in\Sigma_{l_i}}\eta_{\le -p_0d_iu}(\xi-\zeta-\theta)\mathcal F_{i, \Z^k}f(x_{(\xi-\zeta)}) \mathcal F_{\Z}\phi_{I}(\zeta) \eta_{\le -{(d_i-\varepsilon)u}}(\zeta)d\zeta
		\end{align*}
		is the $i$-th Fourier transform of the function from \eqref{eq:72}, evaluated at $x_{(\xi)}$.
	\end{proof}
	
	\subsubsection{\textbf {All together: proof of Theorem~\ref{thm:1}}}
	We now focus our attention on proving Theorem~\ref{thm:1}. To simplify notation, we write
	\begin{align}
	\label{eq:76}
	F_{N}^{i}(f_i)\coloneqq F_{N, t, u}^{i, l_i, s_i}(f_i)= \Pi_{l_i,\le -p_0d_iu}^i(F_{N, t}^{i, l_i, s_i}(f_i)), \qquad i\in[k],
	\end{align}
	following \eqref{eq:35}. By Lemma~\ref{lemma:3}, we have
	\begin{align}
	\label{eq:FN}
	\| \sup_{N\in\D_{l, s}}|F_{N}^{i}(f_i)| \|_{\ell^{p_i}(\Z^k)} \lesssim
	(s+1) 2^{o_{p_i}(l)}
	\end{align}
	for all $i \in [k]$ and all normalized $f_i \in \ell^{p_i}(\Z^k)$.
	
	\begin{proof}[Proof of Theorem~\ref{thm:1}]
		Using notation from \eqref{eq:76}, we see that it suffices to prove that, with some constant $\delta_p\in(0, 1)$, we have 
		\begin{align}
		\label{eq:134}
		\sup_{1/2\le t\le 1}\|\sup_{N\in\D_{l, s}}|\tilde{A}_{2^u}(F_{N}^{1}(f_1),\ldots, F_{N}^{k}(f_k))|\|_{\ell^p(\Z^k)}\lesssim (s+1)^k 2^{-\delta_p l} 
		\end{align}
		for all normalized $f_1\in\ell^{p_1}(\Z^k), \dots, f_k\in\ell^{p_k}(\Z^k)$. Since our argument is fairly intricate, it will be divided into several steps.
		
		\paragraph{\bf Step~1}
		It suffices to prove that, with some constant $\gamma\in(0, 1)$, we have 
		\begin{align}
		\label{eq:125}
		\sup_{1/2\le t\le 1}\|\sup_{N\in\D_{l, s}}|\tilde{A}_{2^u}(F_{N}^{1}(f_1),\ldots, F_{N}^{k}(f_k))|\|_{\ell^4(\Z^k)}
		\lesssim (s+1)^k 2^{-\gamma l}
		\end{align}
		for all normalized $f_1, \dots, f_k \in\ell^{4k}(\Z^k)$, with $F_{N}^{i}(f_i)$ as in \eqref{eq:76} for all $i\in[k]$. Any other parameter \( p \in (2, \infty) \) in place of \( 4 \) in inequality \eqref{eq:125} would also work. However, using the number \( 4 \) would reduce the number of parameters we need to declare later, making our argument clearer. Note that inequality \eqref{eq:125} implies inequality \eqref{eq:134}. To see this, suppose that \eqref{eq:125} holds. Combining inequality \eqref{q<1} from the proof of Theorem~\ref{weyl} with \eqref{eq:FN}, we obtain
		\begin{align}
		\label{eq:int q<1}
		\left\|\sup_{N\in\D_{l, s}}|\tilde{A}_{2^u}(F_{N}^{1}(f_1),\ldots, F_{N}^{k}(f_k))|\right\|_{\ell^v(\Z^k)}
		\lesssim (s+1)^k 2^{o_{v_1,\ldots, v_k, v}(l)}
		\end{align}
		for all normalized $f_1\in\ell^{v_1}(\Z^k), \dots, f_k\in\ell^{v_k}(\Z^k)$, where $1<v_1,\ldots, v_k<\infty$ satisfy $\frac{1}{v_1}+\cdots+\frac{1}{v_k}=\frac{1}{v}$ with $v<1$ close to $1$. Interpolating the unnormalized bounds from \eqref{eq:int q<1} and \eqref{eq:125}, using the multilinear interpolation result from \cite[Theorem~1]{GrTao}, and renormalizing yields \eqref{eq:134} for all $1<p_1,\ldots, p_k<\infty$ from Theorem~\ref{thm:1} satisfying
		$\frac{1}{p_1}+\cdots+\frac{1}{p_k}=\frac{1}{p}\le 1$.
		
		By a simple density argument, it suffices to prove \eqref{eq:125} for all normalized $f_1,\ldots, f_k\in\ell^{4k}(\Z^k)$ that are also compactly supported. This additional assumption ensures that all functions also belong to $\ell^2(\Z^k)$. Hence, from now on we assume that $f_1,\ldots, f_k\in\ell^{4k}(\Z^k)$ are compactly supported. 
		
		\paragraph{\bf Step~2}
		
		By inequality \eqref{eq:43} from Theorem~\ref{weyl}, there exists a constant $c \coloneqq c_k\in(0, 1)$ such that for all $g_{1},\ldots, g_k\in \ell^{2k}(\mathbb Z^k)$, if $g_i\in \ell^{2}(\mathbb Z^k)$ for some $i\in[k]$ and the $i$-th Fourier transform of $g_i$ vanishes on $\mathcal M^i_{\le w_1, \le - d_i u +w_2}$, then 
		\begin{align}
		\label{eq:135}
		\|\tilde{A}_{2^u}(g_1,\ldots, g_k)\|_{\ell^2(\Z^k)}
		\lesssim  (2^{-c \min\{w_1, w_2 \}}+2^{-cu}) \prod_{i\in[k]}\|g_i\|_{\ell^{2k}(\Z^k)}.
		\end{align}
		
		Thus, our aim is to establish \eqref{eq:125} by reducing the matter to \eqref{eq:135}. For this purpose, we choose $r\in(4, \infty)$ such that
		\begin{align}
		\label{eq:131}
		100kD\bigg(\frac{1}{4}-\frac{1}{r}\bigg)<\frac{c}{1000},
		\end{align}
		where $c \coloneqq c_{k} \in (0,1)$ is the exponent from inequality \eqref{eq:135}. 
		
		\paragraph{\bf Step~3}
		
		For each $i \in[k]$ we now fix a partition of unity $(\phi_{I^i}:I^i\in\mathcal I^i)$ as in \eqref{eq:68}, where $\mathcal I^i \coloneqq \mathcal I_{2^{d_i u}}$
		is the partition of $\R$ into disjoint intervals $I^i$ of equal length
		$2^{d_i u}$. Next, define $\mathcal I \coloneqq \mathcal I^1\times\cdots\times\mathcal I^k$ and observe that for every $I=I^1\times\cdots\times I^k\in\mathcal I$  we have
		\begin{align}
		\begin{gathered}
		\label{eq:39}
		\tilde{A}_{2^u}(\phi_{I^1}^1F_{N}^{1}(f_1),\ldots,\phi_{I^k}^kF_{N}^{k}(f_k))=
		\tilde{A}_{2^u}(\phi_{I^1}^1F_{N}^{1}(f_1^I),\ldots,\phi_{I^k}^kF_{N}^{k}(f_k^I))\ind{ \tilde{I}},
		\end{gathered}
		\end{align}
		with $\phi_{I^i}^i$ defined in \eqref{eq:70}, $\tilde{I} \coloneqq \tilde{I}^1\times\cdots\times \tilde{I}^k\subseteq \R^k$, and
		\begin{align}
		\label{eq:25}
		f_i^I \coloneqq \big(\ind{\tilde{I}^1}\otimes\cdots\otimes\ind{\tilde{I}^{i-1}}\otimes \ind{\R}\otimes\ind{\tilde{I}^{i+1}}\otimes\cdots\otimes\ind{\tilde{I}^{k}}\big)f_i,\qquad i\in[k],
		\end{align}
		where $\tilde{I}^i$ is an interval centered at $c(I^i)$ of length $O(2^{d_i u})$ such that, whenever $\phi_{I^i}(x_i-P_i(n)) \neq 0$ for some $x \in \R^k$ and $n \in [2^u]$, the arguments $x_i - P_i(n)$ and $x_i$ both belong to $\tilde{I}^i$. Indeed, with this choice we have $\phi_{I^i}\ind{\tilde{I}^i}=\phi_{I^i}$ and we can also localize the operator to $\tilde{I}$. Finally, note that for each $i\in [k]$, the family of intervals $\{\tilde{I}^i : I^i\in\mathcal I^i\}$ has bounded overlap.
		
		For each $i\in[k]$, define
		\begin{align}
		\label{eq:136}
		G_{N, I}^i \coloneqq G_{N, I}^i(f_i) \coloneqq \phi_{I^i}^iF_{N}^{i}(f_i^I)
		\end{align}
		with $f^I_i$ as in \eqref{eq:25}. Moreover, for each $\omega  = (\omega_1, \dots, \omega_k) \in\{0, 1\}^k$, define
		\begin{align*}
		G_{N, I}^{i, \omega_i} \coloneqq
		G_{N, I}^{i, \omega_i}(f_i) \coloneqq 
		\begin{cases}
		(\phi_{I^i}^i*_i\mathcal F^{-1}_{\R}\eta_{\le -{(d_i-\varepsilon )u}})F_{N}^{i}(f_i^I) & \text{ if } \omega_i=1,\\
		(\phi_{I^i}^i*_i(\delta_0-\mathcal F^{-1}_{\R}\eta_{\le -{(d_i-\varepsilon )u}}))F_{N}^{i}(f_i^I)& \text{ if }  \omega_i=0.
		\end{cases}
		\end{align*}
		Note that
		\[
		G_{N, I}^{i}=G_{N, I}^{i, 0}+G_{N, I}^{i, 1}, \qquad i \in [k].
		\]
		Using \eqref{eq:39} and \eqref{eq:136}, and the bounded overlap of $\{\tilde{I} : I\in\mathcal I\}$, we may write
		\begin{align}
		\label{eq:157}
		\begin{split}
		&\|\sup_{N\in\D_{l, s}}|\tilde{A}_{2^u}(F_{N}^{1}(f_1),\ldots, F_{N}^{k}(f_k))|\|_{\ell^4(\Z^k)}^4\\
		&\quad\lesssim\sum_{\omega\in\{0, 1\}^k}
		\sum_{I\in\mathcal I}\|\sup_{N\in\D_{l, s}}|\tilde{A}_{2^u}(G_{N, I}^{1, \omega_1},\ldots, G_{N, I}^{k, \omega_k})|\|_{\ell^4(\tilde{I}\cap\Z^k)}^4.
		\end{split}
		\end{align}
		
		\paragraph{\bf Step~4}
		We need a few properties of the functions $G_{N, I}^{i,\omega_i}$. We fix $i\in[k]$ and $p_i\in(1, \infty)$. Since $G_{N, I}^{i, 1}=G_{N, I}^{i}-G_{N, I}^{i, 0}$, we have  
		\begin{align*}
		\|\sup_{N\in\D_{l, s}}|G_{N, I}^{i,1}|\|_{\ell^{p_i}(\Z^k)}^{p_i}\lesssim
		\|\sup_{N\in\D_{l, s}}|
		G_{N, I}^{i}|\|_{\ell^{p_i}(\Z^k)}^{p_i}+
		\|\sup_{N\in\D_{l, s}}|G_{N, I}^{i, 0}|\|_{\ell^{p_i}(\Z^k)}^{p_i}.
		\end{align*}
		Observe that $ \sup_{N\in\D_{l, s}}| G_{N,I}^{i}| = \phi_{I^i}^i \sup_{N\in\D_{l, s}}|
		F_{N}^{i}(f_i^I)|$ by \eqref{eq:136}. Moreover, $\{\tilde{I} : I\in\mathcal I\}$ has bounded overlap. Thus, using Lemma~\ref{lemma:3} we obtain
		\begin{align}
		\label{eq:138}
		\|\sup_{N\in\D_{l, s}}|
		G_{N, I}^{i}|\|_{\ell^{p_i}(\Z^k)}^{p_i}\lesssim
		C'_{p_i}(l, s_i) \|f_i\|_{\ell^{p_i}(\Z^k)}^{p_i} 
		\end{align}
		with a constant $C'_{p_i}(l, s_i) \in \R_+$ of the form $ (s_i + 1)^{p_i} (2^{O_{p_i}(\frac{l \log \log l }{\log l})}\ind{l\ge10}+\ind{l<10})$. 
		
		Next, taking $M =1$ in \eqref{eq:75} from Lemma~\ref{lemma:5}, we have
		\begin{align*}
		\sup_{N\in\D_{l, s}}|G_{N, I}^{i, 0}(x)|^{p_i}
		\lesssim_{\varepsilon}
		2^{-p_i u} w_i (x\cdot e_i-c(I^i))\sup_{N\in\D_{l, s}}|
		F_{N}^{i}(f_i^I)(x)|^{p_i}, \qquad x\in\Z^k,
		\end{align*}
		where $w_i(y) \coloneqq \langle 2^{-d_i u}|y|\rangle^{-p_i}$ for $y\in\Z$. Since $w_i (x\cdot e_i-c(I^i))$ are summable over $I^i \in \mathcal I^i$, uniformly in $x\in \mathbb Z^k$, we see that
		\begin{align*}
		\sum_{I^i\in\mathcal I^i}\|\sup_{N\in\D_{l, s}}|G_{N, I}^{i, 0}(f_i)|\|_{\ell^{p_i}(\Z^k)}^{p_i}
		\lesssim2^{-p_i u}
		\|\sup_{N\in\D_{l, s}}|F_{N}^{i}(f_i^I)| \|
		_{\ell^{p_i}(\Z^k)}^{p_i}
		\end{align*}
		and, since $\{\tilde{I} : I\in\mathcal I\}$ has bounded overlap, using Lemma~\ref{lemma:3} we obtain
		\begin{align}
		\label{eq:140}
		\sum_{I\in\mathcal I}\|\sup_{N\in\D_{l, s}}|G_{N, I}^{i,0}|\|_{\ell^{p_i}(\Z^k)}^{p_i}
		\lesssim 2^{-p_i u} C'_{p_i}(l, s_i) \|f_i\|_{\ell^{p_i}(\Z^k)}^{p_i}.
		\end{align}
		
		Finally, by Lemma~\ref{lemma:4} with $\varepsilon = \frac{1}{10}$ for each $i\in[k]$ the $i$-th Fourier transform of the function $G_{N, I}^{i, 1}(f_i)$ is supported on the major arcs $\mathcal M^i_{l_i, \le -(d_i-\varepsilon) u +1}$. Thus, it vanishes on the major arcs $\mathcal M^i_{\le l_i-1, \le -d_i u+\varepsilon(l_i-1)}$.
		
		\paragraph{\bf Step~5}
		Fix $\omega \in \{0,1\}^k \setminus \{1\}^k$. By the previous discussion, we have
		\begin{align}
		\label{eq:130}
		\begin{split}
		\sum_{I\in\mathcal I}\|\sup_{N\in\D_{l, s}}|\tilde{A}_{2^u}(G_{N, I}^{1, \omega_1},\ldots, G_{N, I}^{k, \omega_k})|\|_{\ell^4(\tilde{I}\cap\Z^k)}^4
		\lesssim 2^{-4u}(s+1)^{4k} 2^{o_{k}(l)}
		\end{split}
		\end{align}
		for all normalized $f_1,\ldots, f_k\in\ell^{4k}(\Z^k)$. Indeed, to prove inequality \eqref{eq:130}, we use the boundedness of $\tilde{A}_{2^u}$ and H{\"o}lder's inequality, together with inequalities \eqref{eq:138} and \eqref{eq:140} applied with $p=4$ and $p_i=4k$ for all $i\in[k]$. 
		
		\paragraph{\bf Step~6} It remains to consider $\omega = (1,\dots,1)$. We shall prove the bound
		\begin{align}
		\label{eq:80}
		\sum_{I\in\mathcal I} \|\sup_{N\in\D_{l, s}}|\tilde{A}_{2^u}(G_{N, I}^{1, 1},\ldots, G_{N, I}^{k, 1})|\|_{\ell^4(\tilde{I}\cap\Z^k)}^4 \lesssim
		(s+1)^{4k} 2^{o_{k}(l)} 2^{-\frac{4cl}{40}}
		\end{align}
		for all normalized $f_1,\ldots, f_k\in\ell^{4k}(\Z^k)$. This is the main part of the argument, and once inequality \eqref{eq:80} is verified, the proof follows. Indeed, it suffices to combine inequalities \eqref{eq:157} and \eqref{eq:130} for all $\omega \in \{0,1\}^k \setminus \{1\}^k$ (recall that now $u=100k(l+1)$ by \eqref{eq:22}), and \eqref{eq:80} to deduce \eqref{eq:125} with $\gamma \coloneqq c/50$.
		
		By H{\"o}lder's inequality (with $r>4$ as in condition \eqref{eq:131}), we obtain
		\begin{align}
		\label{eq:155}
		\begin{split}
		&\sum_{I\in\mathcal I}\|\sup_{N\in\D_{l, s}}|\tilde{A}_{2^u}(G_{N, I}^{1, 1},\ldots, G_{N, I}^{k, 1})|\|_{\ell^4(\tilde{I}\cap\Z^k)}^4\\
		&\quad \lesssim 2^{4uD(\frac{1}{4}-\frac{1}{r})}\sum_{I\in\mathcal I}\|\sup_{N\in\D_{l, s}}|\tilde{A}_{2^u}(G_{N, I}^{1, 1},\ldots, G_{N, I}^{k, 1})|\|_{\ell^r(\Z^k)}^4.
		\end{split}
		\end{align}
		In view of condition \eqref{eq:131}, the bound \eqref{eq:80} follows from \eqref{eq:155} and
		\begin{align}
		\label{eq:64}
		\begin{split}
		\sum_{I\in\mathcal I} \|\sup_{N\in\D_{l, s}}|\tilde{A}_{2^u}(G_{N, I}^{1, 1},\ldots, G_{N, I}^{k, 1})|\|_{\ell^r(\Z^k)}^4
		\lesssim
		(s+1)^{4k} 2^{o_{k}(l)} 2^{-\frac{4cl}{30}}.
		\end{split}
		\end{align}
		
		\paragraph{\bf Step~7} From now on we use the notation introduced above Proposition~\ref{nit}. By the norm interchanging inequality from Proposition~\ref{nit} with $R=4<r$ applied to the $\ell^r(\Z^k)$ norm in \eqref{eq:64}, followed by H{\"o}lder's inequality, we have
		\begin{align}
		\label{eq:90}
		\begin{split}
		& \sum_{I\in\mathcal I}\|\sup_{N\in\D_{l, s}} |\tilde{A}_{2^u}(G_{N, I}^{1, 1},\ldots, G_{N, I}^{k, 1})|\|_{\ell^r(\Z^k)}^4\\
		& \quad \lesssim_{r}
		\sum_{\sigma\in\{0, 1\}^k}\sum_{I\in\mathcal I}\Big(\sup_{{\mathbf N}_{<J_1}^1,\ldots, {\mathbf N}_{<J_k}^k\in{\bf I}(\mathbb D_{l, s})}
		\Lambda_{{\mathbf N}_{<J_1}^1,\ldots, {\mathbf N}_{<J_k}^k}^{I,\sigma, 4, 4}(f_1,\ldots, f_k)\Big)^4,
		\end{split}
		\end{align}
		where for $\sigma=(\sigma_1,\ldots, \sigma_k)\in\{0, 1\}^k$, we define $\Lambda_{{\mathbf N}_{<J_1}^1,\ldots, {\mathbf N}_{<J_k}^k}^{I,\sigma, 4, 4}(f_1,\ldots, f_k)$ by
		\begin{align*}
		\Big(\sum_{(j_1,\ldots, j_k)\in\mathbb X_{\sigma}}\|\tilde{A}_{2^u}(\Delta G_{N_{j_1}^1, I}^{1, 1, \sigma_1},\ldots, \Delta G_{N_{j_k}^k, I}^{k, 1, \sigma_k})\|_{\ell^4(\Z^k)}^4\Big)^{1/4},
		\end{align*}
		and for each $i\in[k]$, setting $N^i_{-1} \coloneqq \min \D_{l, s}$, we have 
		\begin{align*}
		\mathbb X_{\sigma} \coloneqq \prod_{i\in [k]}\N_{<J_i}^{\sigma_i}
		\quad \text{with} \quad
		\N_{<J_i}^{\sigma_i} \coloneqq 
		\begin{cases}
		\N_{<J_i}  & \text{ if } \sigma_i=1,\\
		\{-1\} & \text{ if } \sigma_i=0,
		\end{cases}
		\end{align*}
		and
		\begin{align*}
		\Delta G_{N^i_{j_i}, I}^{i,1, \sigma_i} \coloneqq 
		\Delta G_{N^i_{j_i}, I}^{i,1, \sigma_i}(f_i) \coloneqq 
		\begin{cases}
		G_{N^i_{j_{i}}, I}^{i, 1}(f_i)-G_{N^i_{j_i-1}, I}^{i, 1}(f_i) & \text{ if } \sigma_i=1,\\
		G_{N^i_{-1}, I}^{i, 1}(f_i) & \text{ if } \sigma_i=0.
		\end{cases}
		\end{align*}
		In view of the final remark in Step~4, for each $i\in[k]$, the $i$-th Fourier transform of $\Delta G_{N^i_{j_i}, I}^{i,1, \sigma_i}$ is supported on the major arcs $\mathcal M^i_{l_i, \le -(d_i-\varepsilon) u +1}$ with $\varepsilon=\frac{1}{10}$. Thus, it vanishes on $\mathcal M^i_{\le l_i-1, \le -d_i u+\varepsilon(l_i-1)}$ and we may write
		\begin{align*}
		\Delta G_{N^i_{j_i}, I}^{i,1, \sigma_i}=\Delta G_{N^i_{j_i}, I}^{i,1, \sigma_i}-\Pi^i_{\le l_i-1, \le -d_i u+\varepsilon(l_i-1)}\Delta G_{N^i_{j_i}, I}^{i,1, \sigma_i}.
		\end{align*}
		
		\paragraph{\bf Step~8} Fix $\sigma\in\{0, 1\}^k$ and ${\mathbf N}_{<J_1}^1,\ldots, {\mathbf N}_{<J_k}^k\in{\bf I}(\mathbb D_{l, s})$, and let $\mathfrak S(\sigma) \coloneqq \{i\in[k]: \sigma_i=1\}$. We may assume, without loss of generality, that $\mathfrak S(\sigma)=[m]$ for some $m\in\N_{\le k}$. In order to simplify notation, for each $i\in[k]$, we define
		\begin{align*}
		H_{i}(f) \coloneqq f-\Pi^i_{\le l_i-1, \le -d_i u+\varepsilon(l_i-1)}f, 
		\end{align*}
		and, assuming momentarily that $I \in \mathcal I$ is fixed, we set
		\begin{align*}
		g_{j_i}^{i} \coloneqq g_{j_i}^{i}(f_i) \coloneqq\Delta G_{N^i_{j_i}, I}^{i,1, 1}(f_i) \quad \text{and} \quad
		g_{i} \coloneqq g_{i}(f_i) \coloneqq G_{N^i_{-1}, I}^{i,1}(f_i).
		\end{align*}
		Given $p\in[1, \infty]$, define $\Theta^{\sigma, p}_{{\mathbf N}_{<J_1}^1,\ldots, {\mathbf N}_{<J_k}^k}({\bm h}_1,\ldots, {\bm h}_m, h_{m+1},\ldots, h_k)$ by setting
		\begin{align*}
		\Big(\sum_{(j_1,\ldots, j_k)\in \mathbb X_{\sigma}}\|\tilde{A}_{2^u}(H_{1}(h_{j_1}^{1}),\ldots,H_{m}(h_{j_m}^{m}), H_{m+1}(h_{m+1}),\ldots, H_{k}(h_k))\|_{\ell^p(\Z^k)}^p\Big)^{1/p}
		\end{align*} 
		for all sequences of functions ${\bm h}_1=(h_{j_1}^{1})_{j_1\in\mathbb N_{<J_1}},\ldots, {\bm h}_m=(h_{j_m}^{m})_{j_m\in\mathbb N_{<J_m}}$ and all functions $h_{m+1}, \dots, h_k$ for which this formula makes sense. Note that
		\begin{align*}
		\Theta^{(0,\dots,0), p}_{{\mathbf N}_{<J_1}^1,\ldots, {\mathbf N}_{<J_k}^k}(h_{1},\ldots, h_k) =
		\|\tilde{A}_{2^u}(H_{1}(h_{1}),\ldots, H_{k}(h_k))\|_{\ell^p(\Z^k)}.
		\end{align*}
		We also observe that $\Lambda_{{\mathbf N}_{<J_1}^1,\ldots, {\mathbf N}_{<J_k}^k}^{I,\sigma, 4, 4}(f_1,\ldots, f_k) = \Theta^{\sigma, 4}_{{\mathbf N}_{<J_1}^1,\ldots, {\mathbf N}_{<J_k}^k}({\bm g}_1,\ldots, {\bm g}_m, g_{m+1},\ldots, g_k)$
		for ${\bm g}_1=(g_{j_1}^{1})_{j_1\in\mathbb N_{<J_1}},\ldots, {\bm g}_m=(g_{j_m}^{m})_{j_m\in\mathbb N_{<J_m}}$ and $g_{m+1}, \dots, g_k$ as before. 
		
		Given $\sigma\in\{0, 1\}^k$ and ${\mathbf N}_{<J_1}^1,\ldots, {\mathbf N}_{<J_k}^k\in{\bf I}(\mathbb D_{l, s})$, we will prove that
		\begin{align}
		\label{eq:81}
		\begin{split}
		& \Theta^{\sigma, 4}_{{\mathbf N}_{<J_1}^1,\ldots, {\mathbf N}_{<J_k}^k} ({\bm h}_1,\ldots, {\bm h}_m, h_{m+1},\ldots, h_k)\\
		& \quad \lesssim 2^{o_k(l)}2^{-\frac{cl}{30}}
		\Big(\prod_{i\in\mathfrak S(\sigma)}\big\|\big(\sum_{j_i\in\mathbb N}|h_{j_i}^{i}|^{4}\big)^{1/4}\big\|_{\ell^{4k}(\mathbb Z^k)}\Big)
		\Big( \prod_{i\in\mathfrak S(\sigma)^c}\|h_i\|_{\ell^{4k}(\mathbb Z^k)} \Big).
		\end{split}
		\end{align}
		for every ${\bm h}_i=(h_{j_i}^{i})_{j_i\in\mathbb N}\in\ell^{4k}\big(\Z^k; \ell^{4}\big(\mathbb N)\big)$ for $i\in[m]$ and every $h_i\in\ell^{4k}(\Z^k)$ for $i\in[k]\setminus[m]$. If \eqref{eq:81} holds, then we can control $\Lambda_{{\mathbf N}_{<J_1}^1,\ldots, {\mathbf N}_{<J_k}^k}^{I,\sigma, 4, 4}(f_1,\ldots, f_k)$ by 
		\begin{align*}
		2^{o_k(l)}
		2^{-\frac{cl}{30}} \Big(\prod_{i\in\mathfrak S(\sigma)}\big\|\big(\sum_{j_i\in\mathbb N_{<J_i}}|\Delta G_{N^i_{j_i}, I}^{i,1, 1}|^{4}\big)^{1/4}\big\|_{\ell^{4k}(\mathbb Z^k)}\Big)
		\Big( \prod_{i\in\mathfrak S(\sigma)^c}\|G_{N^i_{-1}, I}^{i, 1}\|_{\ell^{4k}(\mathbb Z^k)} \Big)
		\end{align*}
		so that $\sup_{{\mathbf N}_{<J_1}^1,\ldots, {\mathbf N}_{<J_k}^k\in{\bf I}(\mathbb D_{l, s})}\Lambda_{{\mathbf N}_{<J_1}^1,\ldots, {\mathbf N}_{<J_k}^k}^{I,\sigma, 4, 4}(f_1,\ldots, f_k)$ is controlled by
		\begin{align*}
		\begin{split}
		2^{o_k(l)}2^{-\frac{cl}{30}} \Big(\prod_{i\in\mathfrak S(\sigma)}\|V^{4}(G_{N, I}^{i,1}: N\in\mathbb D_{l, s})\|_{\ell^{4k}(\mathbb Z^k)}\Big) \Big(
		\prod_{i\in\mathfrak S(\sigma)^c}\|G_{N^i_{-1}, I}^{i, 1}\|_{\ell^{4k}(\mathbb Z^k)} \Big).
		\end{split}
		\end{align*}
		Inserting this bound into inequality \eqref{eq:90} and applying Hölder's inequality, we sum over \(I \in \mathcal{I}\). Since \(\mathcal{I}\) has bounded overlap, we can use Lemma~\ref{lemma:5} and proceed as in Step~4. Additionally, invoking Lemma~\ref{lemma:3}, we conclude that inequality \eqref{eq:64} holds, as desired. 
		
		\paragraph{\bf Step~9}
		It remains to establish \eqref{eq:81}. For this purpose, we first prove that
		\begin{align}
		\label{eq:89}
		\begin{split}
		& \Theta^{\sigma, 2}_{{\mathbf N}_{<J_1}^1,\ldots, {\mathbf N}_{<J_k}^k} ({\bm h}_1,\ldots, {\bm h}_m, h_{m+1},\ldots, h_k)\\
		&\quad \lesssim 2^{o_k(l)}2^{-\frac{cl}{10}}
		\Big(\prod_{i\in\mathfrak S(\sigma)}\big\|\big(\sum_{j_i\in\mathbb N}|h_{j_i}^{i}|^{2}\big)^{1/2}\big\|_{\ell^{2k}(\mathbb Z^k)}\Big)
		\Big( \prod_{i\in\mathfrak S(\sigma)^c}\|h_i\|_{\ell^{2k}(\mathbb Z^k)} \Big)
		\end{split}
		\end{align}
		for every ${\bm h}_i=(h_{j_i}^{i})_{j_i\in\mathbb N}\in\ell^{2k}\big(\Z^k; \ell^{2}\big(\mathbb N)\big)$ for $i\in[m]$ and every $h_i\in\ell^{2k}(\Z^k)$ for $i\in[k]\setminus[m]$. To this end, we apply Khinchine's inequality (see, for instance, \cite[Section~C.5, p.~589]{Grafakos1}), which allows us to linearize the $\ell^2(\bar{\mathbb X}_{\sigma})$ square function that appears in the definition of $\Theta^{\sigma, 2}_{{\mathbf N}_{<J_1}^1,\ldots, {\mathbf N}_{<J_k}^k} ({\bm h}_1,\ldots, {\bm h}_m, h_{m+1},\ldots, h_k)$. Inserting the summation over $\bar{\mathbb X}_{\sigma}$ under the operator $\tilde A_{2^u}$ and applying Weyl's inequality \eqref{eq:135} with $w_1=l-1$ and $w_2=\varepsilon(l-1)$ for $\varepsilon=\frac{1}{10}$, we obtain \eqref{eq:89}.
		
		Similarly, we prove that
		\begin{align}
		\label{eq:115}
		\begin{split}
		& \Theta^{\sigma, 8}_{{\mathbf N}_{<J_1}^1,\ldots, {\mathbf N}_{<J_k}^k} ({\bm h}_1,\ldots, {\bm h}_m, h_{m+1},\ldots, h_k)\\
		& \quad \lesssim 2^{o_k(l)}
		\Big(\prod_{i\in\mathfrak S(\sigma)}\big\|\big(\sum_{j_i\in\mathbb N}|h_{j_i}^{i}|^{8}\big)^{1/8}\big\|_{\ell^{8k}(\mathbb Z^k)}\Big)
		\Big( \prod_{i\in\mathfrak S(\sigma)^c}\|h_i\|_{\ell^{8k}(\mathbb Z^k)} \Big)
		\end{split}
		\end{align}
		for every ${\bm h}_i=(h_{j_i}^{i})_{j_i\in\mathbb N}\in\ell^{8k}\big(\Z^k; \ell^{8}\big(\mathbb N)\big)$ for $i\in[m]$ and every $h_i\in\ell^{8k}(\Z^k)$ for $i\in[k]\setminus[m]$. This inequality readily follows by applying Minkowski's integral inequality and the boundedness of the operator $\tilde{A}_{2^u}$. 
		
		Then, noting that $\frac{1}{4}=\frac{\theta}{2}+\frac{1-\theta}{8}$ and $\frac{1}{4k}=\frac{\theta}{2k}+\frac{1-\theta}{8k}$ for $\theta=\frac{1}{3}$, we can interpolate \eqref{eq:89} and \eqref{eq:115}, using a multilinear vector-valued interpolation similar to a Riesz--Thorin argument, to derive inequality \eqref{eq:81}. Specifically, we can adapt the proof of the multilinear complex interpolation theorem \cite[Theorem~7.2.9, p.~514]{Grafakos2} to the multilinear vector-valued setting, following the approach in \cite[Exercise~5.5.2, p.~339]{Grafakos1}.
		
		Now, the proof of Theorem~\ref{thm:1} follows.
	\end{proof}
	
	\begin{remark}
		\label{rem:5} 
		Finally, we note that if one wants to use the method from \cite{KMT}, then the argument in Step~6 would fail badly. In \cite{KMT} the parameter $u$ is exponential in $l$, while in our case it is linear; see \eqref{eq:22}. This is because in \cite{KMT} the bilinear Weyl inequality was proved with the decay $\delta^c+ \langle \log N \rangle^{-c\widetilde C}$, which includes a logarithmic term that requires exponential $u=2^{O(l)}$. However, thanks to the improvement in the Ionescu--Wainger theorem, our Weyl's inequality holds with $\delta^c+N^{-c}$, allowing us to have a linear dependence between $l$ and $u$. This is crucial in deriving inequality \eqref{eq:80}, which follows from choosing a parameter $r$ in condition \eqref{eq:131}. This choice renders the contribution of $u$ in \eqref{eq:155} negligible compared to the size of $c$ arising in \eqref{eq:135}. If $u$ were exponential in $l$, our approach in Step~6 would not be feasible at all.
	\end{remark}

\end{document}